\RequirePackage{fix-cm}
\documentclass{svjour3}                     % onecolumn (standard format)
\smartqed  % flush right qed marks, e.g. at end of proof
\usepackage{graphicx}

\newtheorem{thm}{\bf Theorem} 
\newtheorem{cor} {\bf Corollary}

\newcommand{\AT}[1]{{\color{red} Alternate Theorem:\\  #1 }}

\newtheorem{defn} {\bf Definition}

 \newcommand{\ignore}[1]{}
 \usepackage{tikz}
 \usepackage{ifthen}
 \usepackage{xcolor}
 \usepackage{extarrows}

\newcommand{\gp}[1]{\gamma_{_{p_{#1}}}}

\newcommand{\eop}{{\hfill $\blacksquare$} }

\newcommand{\ForTR}{\ifthenelse{1<0}} 
\newcommand{\ATR}{\ifthenelse{1<0}}

\newcommand{\x}{{\bar {\bf x}}}
\newcommand{\X}{{\bar {\bf X}}}

\renewcommand{\u}{{\bar {\bf u}}} 
\usepackage{tikz}
\usetikzlibrary{matrix}
\usetikzlibrary{positioning}
\usetikzlibrary{arrows}
\usepackage{xcolor}
\usepackage{amsmath, amsfonts, amssymb,  enumerate}
\usepackage[final]{pdfpages}

\usepackage{mathrsfs}
\usepackage{multirow,array}
\usepackage{caption}
\usepackage{subcaption}
\usepackage{hyperref}
\newcommand{\f}{{\bf  f}}

\newcommand{\indc}[1]{ \mathbf{1}_{#1} }

\renewcommand{\u}{{\bar {\bf u}}} 
\newcommand{\re}{r_{e_2}}

\newcommand{\UV}{V}

\usepackage{changes}
\usepackage{algorithm,algpseudocode}
\newcounter{algsubstate}
\renewcommand{\thealgsubstate}{\alph{algsubstate}}
\newenvironment{algsubstates}
  {\setcounter{algsubstate}{0}%
   \renewcommand{\State}{%
     \stepcounter{algsubstate}%
     \Statex {\footnotesize\thealgsubstate:}\space}}
  {}

\begin{document}

\title{Random fixed points, systemic risk and resilience of heterogeneous  financial network %\footnote{A preliminary version of this work was presented at the CDC conference, 2018 \cite{Systemicrisk}.}    %\thanks{Grants or other notes
}


\author{Indrajit Saha and Veeraruna Kavitha %etc.
 }

\institute{Indrajit Saha \at
              IIT Bombay, India \\
              \email{indrojit@iitb.ac.in}           %  \\
           \and
           Veeraruna Kavitha \at
           IIT Bombay, India \\
           \email{vkavitha@iitb.ac.in}  
}

\date{Received: date / Accepted: date}

\maketitle
\begin{abstract} 
We consider a large random network, in which  the performance of a node depends upon that of its neighbours  and some  external random influence factors.  This results in  random vector valued  fixed-point (FP) equations in large dimensional spaces, and our aim is to study their   almost-sure solutions. An underlying directed random graph  defines the connections between various components of the FP equations.  Existence of an edge between  nodes  $i,j$ implies the $i$-th FP equation  depends on the $j$-th component. We consider a special case where any component of the FP equation depends upon an appropriate aggregate of that of the random `neighbour'  components. We obtain   finite dimensional limit FP equations in a much   smaller dimensional space, whose solutions aid to approximate the solution of FP equations for   almost all realizations,  as the number of nodes increases.  We use Maximum theorem for non-compact sets to prove this convergence.
 
We apply the results to study systemic risk  in an example financial network with large number of  heterogeneous entities. We utilized the simplified limit system to analyse  trends of default probability (probability that an entity fails to clear its liabilities) and expected surplus (expected-revenue after clearing  liabilities) with varying degrees of interconnections between two diverse groups.  We illustrated the accuracy of  the approximation using exhaustive Monte-Carlo simulations. 

Our approach can be  utilized for a   variety of financial networks (and others); the developed methodology   provides approximate small-dimensional solutions to  large-dimensional  FP equations that represent the  clearing vectors in case of financial networks.

\ignore{
We consider a large random network, in which   the performance of a node depends upon that of its neighbours  and some  external random influence factors.  This results in   random vector valued  fixed point (FP) equations in large dimensional spaces, and our aim is to study their   almost sure solutions. 
An underlying directed random graph  defines the connections between various components of the FP equations.  Existence of an edge between  nodes  $i,j$ implies the $i$-th FP 
equation  depends on the $j$-th component. 
We consider a special case where any component of the FP equation depends upon an appropriate aggregate of that of the random `neighbour'  components. We obtain   finite dimensional limit FP equations in a much   smaller dimensional space, whose solutions aid to approximate the solution of the random FP equations for   almost all realizations,  as the number of nodes increases.  We use Maximum theorem for non-compact sets to prove this convergence.
Our techniques are  different from the traditional mean-field  methods,  which  deal with
 stochastic FP equations in the space of distributions   to describe  the   stationary distributions. In contrast  our focus is on almost sure FP solutions.  

We apply the results to study systemic risk  
in an example financial network with large number of  heterogeneous entities.
We utilized the simplified limit system to analyse  trends of default probability (the probability that an entity fails to clear its liabilities) and expected surplus (expected revenue after clearing the liabilities) with varying degrees of interconnections between two diverse groups.  
We illustrated the accuracy of  the approximation using exhaustive Monte-Carlo simulations. 

Our approach is  general and can be applied to a   variety of financial networks (and others); the developed methodology   provides approximate small dimensional solutions to  large-dimensional  FP equations that represent the  clearing vectors in case of financial networks.

%We consider a complex network with many nodes; these nodes are interconnected and shared a resource. The nodes' performance depends on the other connected nodes' aggregated performance and appears as vector-valued random fixed point equations in the large dimensional space. This paper asymptotically solves this random fixed point equation under suitable graph structure assumption. We proved that the solution converges to a limit in an almost sure sense. This asymptotic simplification helps dimension collapse in the limit.

%One can apply this result in a financial network with heterogeneous financial entities to study systemic risk. We used the simplified limit system and analyse important trends of default probability and expected surplus. We also tested the accuracy of the proposed approximation with the practical number of financial entities by performing Monte Carlo simulations; Simulation suggests that the limit system's performance and the actual system well matched with the proposed approximation. This paper also provides a simple algorithm to compute the clearing vector of the banks. 
}
\end{abstract}
\section{Introduction}
Random fixed points (FPs) are generalization of classical deterministic  fixed points,  and arise  when one considers systems with uncertainty. One can think of two types of  fixed points under uncertainty.   There is considerable literature that  considers  stochastic fixed point equations on the space of probability distributions  (e.g., \cite{Urns,WeightedBranch}). These equations typically arise as a limit of some iterative schemes, or  as asymptotic (stationary) distribution of stochastic systems.  Alternatively, one might be interested in sample wise (almost sure)  fixed points as in \cite{Measure_FP,Measure_Approx};  for each realization of the random quantities describing the system, there is  one deterministic fixed point equation. These kind of equations can arise  when the performance/status of an  agent  depends upon that of  other agents. For example, a financial network with any given liability graph   is affected by individual/common random economic shocks received by the agents (see Section \ref{sec_finance} for more details).  The amount cleared (full/fraction of liability) by an agent  depends upon: a) the  random shocks  it receives; and  b) the liabilities cleared by the other agents. Our  focus in this paper is on the second type of equations, defined in almost sure sense.  Current literature primarily considers the existence  of measurable fixed point solutions, given the existence of realization-wise  fixed point solutions (e.g., \cite{Measure_FP,Measure_Approx}). In \cite{Measure_Approx} (and reference therein) authors consider the idea of random proximity points.

To the best of our knowledge, there are no (common) techniques that provide `good' solutions to (even some special type of) these equations. We consider  special type of fixed point equations,  and provide a procedure to compute the approximate almost sure solutions; 
 here the performance/status of an agent is influenced only by the aggregate  performance/status of its neighbours. 

 We consider a random graph where nodes denote agents and the edges denote interaction between the agents. For  example in a financial setting, nodes may be banks and edges may denote liability structure between banks.  A set of fixed point equations (one per realization of the random quantities, e.g., economic shocks)  describe certain  performance vectors of  the agents. The performance of each agent is influenced by aggregate of the performance  of  its  neighbours, with the aggregate defined using the random edges.  For example  the clearing vectors in the financial setting. 

The key  idea is to study these fixed points, asymptotically as the number of agents increase.    Towards this, we first analyze the aggregate influence factors, with an aim to reduce the  dimensionality of the problem. But due to random connections,  
the aggregate influence  factors can also depend upon the nodes. However   the  aggregates might converge towards the same limit almost surely (e.g., as in  law of large numbers). Considering such scenarios,    the random fixed points are shown to   converge to  that of  a limit system, under certain conditions. The performance of  the agents   in the limit system, depends upon finitely many  `aggregate'  limits.  For some examples,    closed-form expressions are derived for   approximate almost sure solutions.

The mean-field theory    primarily deals with a system of large number of agents, wherein the state/behaviour of an individual agent  is influenced by its own (previous) state and the mean (aggregate)  field  seen by it (e.g., \cite{Mean_WLAN} and reference therein). The mean-field is largely described in terms of  occupation (empirical) measures representing the fraction of agents in different states. 
 The theory shows the convergence of the mean state trajectories as well as the stationary (time limit) distributions of the  original system towards that of a limit deterministic system. 

  The stationary distribution   can be described by  fixed point  equations in the space of distributions (e.g., \cite{Mean_WLAN}).  As opposed to that, 
we consider  a set of fixed point  equations, which are defined in almost sure sense.
%and which depend upon the realization-wise  `mean' performance.   

 We consider fixed point equations  with possibly multiple solutions, and, show that any chosen sequence of the fixed points  converge  almost surely to a fixed point of the limit system (along a sub-sequence).   Towards this, we construct an appropriate parameterized optimization problem and apply the relatively recent result (\cite{Feinberg}) on Maximum Theorem for non-compact sets to show almost sure convergence of the aggregate random fixed points; the main idea is to construct appropriate topological spaces (e.g., Tychonoff's topology) and an appropriate objective function.  Under some additional (mild) conditions, we show the uniqueness of the fixed points; we further derive limit solutions using that of  a significantly low-dimensional system. The results are derived for the case with two diverse groups (homogeneous within the group) of agents, for which one has to solve three-dimensional equations; one can easily extend the results to any finite number of groups.

\subsection*{Application to financial networks} 
 We   apply our   results to study  systemic risk  related aspects in a  large financial network. 
The institutions borrow/lend money from/to other institutions, and will have to clear their obligations at a later time.  These systems are subjected to economic shocks, 
because of which some entities default (do not clear their obligations).  Because of inter-dependencies, this can lead to further defaults and the cascade of these reactions can lead to  the
(partial/full)  {\it collapse} of the system.

After the financial crisis of $2007$-$2008$, there is a surge of activity towards studying systemic risk (e.g.,\cite{acemoglu2015systemic,allen2000financial,eisenberg2001systemic}).  The focus in these papers has been on several aspects including, measures to capture systemic risk,   influence of network structure on systemic risk, phase transitions etc. 
Some   papers  (e.g., \cite{allen2000financial,blume2011networks,eisenberg2001systemic,Systemicrisk,Freixas,Haldane}) consider network-based approach, while  \cite{carmona2013mean,Garnier} considers mean-field analysis based approach.
%to study various aspects related to   measuring and managing systemic risk. The former set of papers discuss the influence of network structure on systemic risk and  
%, for example, a) network-based approach \cite{allen2000financial,blume2011networks,eisenberg2001systemic,Systemicrisk,Freixas,Haldane} b) mean-field analysis \cite{carmona2013mean}  etc. 
%
%our work also focuses on similar aspects. 
%
Further these   papers primarily discuss homogeneous systems, although  heterogeneity 
 is a crucial feature of real world networks.    As already mentioned, 
 the clearing vectors are represented by FP equations and  one must analyze the same to study the more realistic heterogeneous networks;   our asymptotic solution can be of significant relevance in this context.  

%The systemic risk is re-highlighted after the financial crisis of $2007$-$2008$. There is a vast literature on systemic risk;   some   papers  (e.g., \cite{allen2000financial,blume2011networks,eisenberg2001systemic,Systemicrisk,Freixas,Haldane}) consider network-based approach, while  \cite{carmona2013mean} considers mean-field analysis based approach to study various aspects related to   measuring and managing systemic risk. The former set of papers discuss the influence of network structure on systemic risk and  
%our work also focuses on similar aspects.  

The seminal work in this line of research is provided by  \cite{allen2000financial}, which   shows that incomplete financial networks are less resilient and more vulnerable to contagion than complete networks (all nodes are interconnected as in complete graph).  A similar kind of conclusions  are derived in   \cite{Freixas}, in  the context  of liquidity shocks.
Another piece of pioneering work is \cite{eisenberg2001systemic}, wherein, the authors show that the clearing payment vector is unique under mild conditions. The paper also provides a fictitious default algorithm to compute the clearing vector. In  recent years, the authors of  \cite{acemoglu2015systemic} extended the work of  \cite{eisenberg2001systemic}, to accommodate the external shocks; they also  showed   
the stability of complete graph (when the magnitude of the negative shock is below a specific range) and vulnerability of the ring graphs  among all regular class of networks.

The previous papers consider time-static models, while \cite{acemoglu2015systemic} also considers three time-period model; at time $t=0$  the portfolio is chosen, partial returns  and liability repayment is at $t=1$ and the final returns (in case of no default)  are at  $t=2$. 

Majority of the papers discussed above consider  deterministic  networks. Real world networks are seldom deterministic, it is more appropriate to model them using random quantities.   Authors  in \cite{Glasserman} consider random networks and derive  a network independent bound on the probability of financial contagion. 
The authors in \cite{Amini} also consider random networks, and derive analysis under the assumption that the recovery rate is negligible for the  defaulted nodes. 

\subsection*{Our results related to financial networks}
We consider random networks with diverse groups  (homogeneous agents within each group), two-time period  return model and with random economic shocks. Further the defaulted banks  pay-back their liabilities to the best extent possible. Under  certain growth condition on the number of neighbours (in each group) we derived a very general technique to obtain approximate  closed form expressions (easily computable) for clearing vectors.  Our methods can handle a  large variety of  networks and the  approximate clearing vectors can be used in computing a variety of performance measures, e.g.,  default probability, expected surplus. For example, in this paper we consider a network with two sets of users,  the first group takes measured risk and the second group is aggressive while choosing their portfolios at time $t=0$. We identified a regime of parameters (interest rates, parameters of economic shocks, percentage of taxes etc.) in which both the groups benefit by small amount of inter-connections between the two groups; for the rest of the regimes, only one of the groups benefits.

One can use our clearing vector based results to study various other aspects. In \cite{Saha} 
we used 
these results to study the convergence of replicator dynamics in a financial network 
  where the agents alter their choices   between  risk-free or risky  portfolios (based on their experiences and observations). We showed that all the agents eventually revert either to risky or risk-free portfolios, unless the agents choose  their strategies  based on large number of observations. In the former case the dynamics converges to pure evolutionary stable strategy (ESS), while the latter converges to a mixed ESS.

Some initial results of this flavour are available in our conference paper
  \cite{Systemicrisk}.   However, the current paper is a sufficient generalization; we consider a more complex network/graph   with a larger variety of entities to define the FP equations and also prove the results using alternate assumptions on graph structure. In addition, the current paper includes all the relevant proofs. We also analyze a more complex financial network.  
  Further using   exhaustive Monte-Carlo simulations, we illustrate  good accuracy of approximation 
  even for moderate-sized networks. 
  To summarise,  our analysis helps identify important patterns in a complex structure, since the structure (often) simplifies when large number of constituents are involved.
  %
  %We performed Monte Carlo simulations with the practical number of financial entities, suggesting  the  approximation is well-matched with the Monte Carlo simulations. We also provided an algorithm to compute the clearing vector of the financial entities.
  %We illustrated the accuracy of  the approximation using exhaustive Monte-Carlo simulations. 

  \ignore{
  with a more complex structure of the financial entities and derived the limiting result of the equilibrium wealth of the entities.
 
We consider one stylized example of a heterogeneous financial network with one big bank and two groups of small banks. {\it Our key contribution is that we develop a methodology to arrive at simplified asymptotic representation to large bank networks.}
This allows the easy resolution of many practical what-if scenarios. For instance, we showed few interesting  properties in the limiting financial network: a phase transition behaviour of the financial network below a certain threshold of the function of the connectivity parameters the whole network is stable as soon as it crosses the threshold the network jumps to default regime; a non-monotone trend of the expected surplus of the financial entities;  provided the regime of the connectivity where both the groups expected surplus improves with inter-lending which highlights the connectivity improves the systemic risk performance. However, in the large shock regime the expected surplus degrades, which suggest interconnectivity plays an adverse effect to the network; proved that group $1$ banks are more robust to the economic shock than that of group $2$ banks.
}

%The proposed methodology can be similarly used to provide insights into many other practical scenarios. {\it To summarise,  our analysis helps identify important patterns in a complex structure, since the structure simplifies when large number of constituents are involved.}

%{\color{red}
 %We consider one  stylized example of  \underline{heterogeneous financial } network, that of one big bank and numerous small banks. 
%{\it Our key contribution is that we develop a methodology to arrive at simplified asymptotic representation to large bank networks.}
%This allows easy resolution of many practical what-if  scenarios. For instance, in a simple framework we observe that
%having a big bank in an economy well connected to the small banks can stabilize the small banks as well as the big bank.  
%We consider a conditional analysis, conditioned on the shocks of the big bank. 
%There exists a range of the parameters describing connectivity between the big bank and small banks for which the system behaves the best: the fraction of defaults is minimum and the conditional surplus is the best.  
%The proposed methodology can be similarly used to provide insights into many other practical scenarios. We  hope to analyse these in future. {\it To summarise,  our analysis helps identify important patterns in a complex structure, since the structure simplifies when large number of constituents are involved.}
%}
\ignore{
\subsection{Related Literature}
The systemic risk is re-highlighted after the financial crisis of $2007$-$2008$. There is a vast literature on systemic risk, which focuses on measuring and mitigating the systemic risk, for example, a) network-based approach \cite{allen2000financial,eisenberg2001systemic,Freixas,Haldane} b) mean-field analysis \cite{carmona2013mean} c) empirically-based work. Our work focuses on a network-based system; usually, in tradition portfolio optimization analysis, financial entities do not have a feedback effect on the other financial entities. In contrast, the network approach is one method to capture systemic risk, which captures the feedback effect.

The seminal work for the contribution of the systemic risk is provided by the author \cite{allen2000financial} and shows how symmetric financial network leads to contagion. The paper's finding is that incomplete financial networks are less resilient and more vulnerable to contagion than complete networks.  A similar kind of conclusions is addressed in the paper \cite{Freixas}.
Another piece of pioneering work \cite{eisenberg2001systemic} wherein the author develops a static model of contagion of the financial network where each agent holds one another liabilities. This is the first work that shows a clearing payment vector is unique under mild conditions to the best of our knowledge. The paper also provides a fictitious default algorithm to compute the clearing vector. In the recent year, the author in  \cite{acemoglu2015systemic} extended the work of  \cite{eisenberg2001systemic}, to accommodate the external shock. This paper \cite{acemoglu2015systemic} studies that as long as the negative shock is beyond a specific range, the complete graph is most stable while ring graph is least stable among all regular class of networks. But as soon as the negative shock crosses some threshold, the network becomes more vulnerable to the default. 

The author in \cite{Glasserman}  addressed similar problem as in \cite{acemoglu2015systemic}. The main focus of the work based on without detailed knowledge of network structure as opposed to the existing literature, the author considers three pieces of information about each node: net worth,  outside leverage, and its financial connectivity. Using these pieces of information, established a network independent bound on the probability of financial contagion.

In \cite{Gai}  also addressed a network model of inter-bank lending with unsecured claims. By performing the numerical simulations, the paper shows how greater complexity and concentration in the financial network may amplify the system's fragility.

In recent years \cite{Systemicrisk},  proposes a random graph approach to model the systemic risk in a heterogeneous financial network. It shows that the finite banking network well approximates with an appropriate limiting financial system which can capture many real-world scenarios. It also shows that the network is resilient as long as the connectivity below some connectivity parameter range. When the operational cost and the common shock are large enough, the network is more prone to the systemic risk event.

Another piece of work originates from the result of \cite{Systemicrisk}. The paper \cite{Saha} shows how replicator dynamics useful in a financial problem where agents have to choose between two strategies, risk-free or risky. In contrast, in the classical social behaviour models, the population chooses between aggressive and passive strategies.
The author studies the emerging strategies when different types of replicator dynamics capture inter-agent interactions. The author proved that the equilibrium strategies converge almost surely to that of an attractor of the ordinary differential equation (ODE).

 The authors in \cite{Elliott} show how discontinuous changes in asset values causes the failure of financial entities and how the network plays a role in it. In contrast to the work of \cite{acemoglu2015systemic}, this paper shows how the consequences of a given moderate shock depend on diversification and integration. The results show that intermediate levels of diversification and integration can be the most problematic.
 
  In \cite{Anand}, the author addressed the problem that inter-bank exposure is unobserved in a financial network. This paper suggests a methodology to estimate the inter-bank network by information-theoretic arguments. It shows that for a stress testing purpose, the minimum density-based approach overestimates the contagion while maximum entropy approach underestimates the contagion.
 
  %The problem of partial information available in inter-bank discussed 
  
In \cite{Mastromatteo}  address the issue, the maximum entropy method leads to underestimating the contagion. Also proposed an efficient message-passing algorithm for reconstructing the network by using the information of partially unknown credit networks, to estimate their robustness.

The author in \cite{Amini}  proposed a cascade of default model through random graph approach. Under the constant recovery rate assumption, shows the analytical expression for the asymptotic fraction of defaults in terms of network parameter.}
%\ignore{
%{\color{blue} What  is systemic risk with respect to graphical model? Why measuring systemic risk is important? What is the connection with the Fixed point? Mean field analysis,
%Modelling description in brief way.}
%\\
%\section*{Some important links}
% \url{https://epubs.siam.org/doi/pdf/10.1137/S0040585X97T988599}
 
 %\url{https://arxiv.org/pdf/1811.00141.pdf}
 
 %\url{https://www.econ.pitt.edu/sites/default/files/erol.selman.pdf}
 
 %\url{http://eprints.lse.ac.uk/66042/1/__lse.ac.uk_storage_LIBRARY_Secondary_libfile_shared_repository_Content_LSE%20SRC%20Discussions%20papers_2016_dp-53.pdf}

%Discuss the literature where the paper models the systemic risk , through asymptotic sparse regime.
%Mention some of the literature on the stochastic block model (SBM ) which is applied to financial network. This graphical model is like a stochastic block model where we want to study how does inter-connectivity among the groups make the stability of the network or rather which kind of link is benefit to the financial network. Does  inter-bank network structure is alone stable or intra-bank connection make robust to the financial network.
%\begin{enumerate}
%\item Daron Acemoglu
%\item Eisenberg
%\item Gai Kapadia
%\item Allen and Gale
%\item Bloom
%\item Hamini et al
%\item Rene Carmona : Mean field game and systemic risk
%\item Rochet
%\item Literature from optimization perspective
%\item Literature on Stochastic fixed point
%\item Literature on Mean field 
%\item Stochastic block model and systemic risk
%\item  CDC work
%\item Netgcoop paper
%\end{enumerate}

%
%}
\noindent{\bf Organization of the  paper:}
The rest of the paper is organized as follows: Sections \ref{Graphicalmodel}  and  \ref{sec_alternate}  provide random fixed point  almost sure results for two  different structures of the network. Section \ref{Section_examples}  provides various other graphical models.
Section \ref{sec_finance}  describes the large  financial system, while, Section \ref{sec_asymptotic} provides  its asymptotic analysis. We have Monte-Carlo simulations in Section \ref{sec_MCsimulation} and Section \ref{sec_conclusion}  concludes the paper. All the proofs are provided in the Appendices.
\section{ Graphical model and fixed point equations}
\label{Graphicalmodel}
We consider $(n+1)$  nodes in a random network indexed by the set $N$
% = \lbrace 1,2, \cdots	n,b \rbrace $  
whose directed   edges,   have random weights    $\{W_{j,i}\}$, representing  influence factors. The  node $b\in N$  is  a big node and   has significant influence on the network. The  remaining $n$ nodes are small nodes and are  classified into two groups $\mathcal{G}_1$  and $\mathcal{G}_2$ respectively. The size of the group $\mathcal{G}_1$ is $n_1 = n\gamma$ and  that  of   $\mathcal{G}_2$ group is  $n_2 = n(1-\gamma)$, where $0 <\gamma < 1$ is a positive fraction\footnote{For any given fraction $\gamma$,  $n$ is  chosen such that $n_1$ and $n_2$ are integers. Further note that ${\cal G}_m$ depends upon $n$, but $n$ is avoided in notation for simplicity and that at limit system (discussed in later parts) ${\cal G}_m$  would have countably infinite elements.  }.
Let
$
{\cal G}_1 = \{1, 2, \cdots, n_1\} \mbox{ and }  {\cal G}_2 = \{1,  \cdots, n_2\}.$
Any small node in $N$ is represented by   pair  $(m, i)$ with $i \le n_m$ and $m \in \{1, 2\}$.

The probability of an edge connecting  two small nodes belonging to the same group (say  $\mathcal{G}_m$)  is  $p_m$, while, that of an edge connecting two nodes belonging to   different groups is  $p_{c_1}$ or $p_{c_2}$. 
 All the edge forming events are independent of the others (or need to satisfy  {\bf B.2(C)} defined later, if some connections are correlated) and let $\{I_{j,i}\}$, be the corresponding indicators.
To summarize,  for any $i \in {\cal G}_m$ and $j \in {\cal G}_{m'}$:
\begin{eqnarray}
\label{Eqn_connections}
P(I_{j,i} =1 ) = p_{mm^{'}} = \left \{  
\begin{array}{llll}
 p_m  &\mbox{ if }  m= m^{'}  \mbox{, }   \mbox{ for any }  m \in \{ 1, 2 \}  \mbox{ and }  \\
 p_{c_1}  &  \mbox{ if } ~~  j\in \mathcal{G}_1 ~~ and ~~ i\in    \mathcal{G}_2 ,\\
  p_{c_2}  &  \mbox{ if } ~~  j\in \mathcal{G}_2~~ and ~~ i\in \mathcal{G}_1 .
\end{array}
\right  .
\end{eqnarray}

 From any small node $j\in \mathcal{G}_m$, there is a dedicated fraction $\eta^{sb}_j$ (of weight) towards the b-node while the  remaining  ($1-\eta^{sb}_j$) fraction is shared   by  other connected small nodes. %Similarly, for $j \in \mathcal{G}_2$, $(1-\zeta^{sb})$ fraction towards the other  small nodes.  
This fraction, for example, can represent an investment to a  particular stock of a big player or to a government security or to a nationalized bank (more details are in Section \ref{sec_finance}, that discusses finance based application).
In all,  weights from a small node $j\in \mathcal{G}_m $ %and  $j\in \mathcal{G}_2 
    are  the respective  fractions  as below\footnote{Note that
 $\displaystyle \sum_{i \in  \mathcal{G}_1 \cup  \mathcal{G}_2} W_{j,i} + W_{j, b} = 1$ for all $j \in  \mathcal{G}_1 \cup  \mathcal{G}_2.$}:
\begin{eqnarray}
\label{Eqn_Weights 1}
W_{j, b} &=&  \eta_j^{sb} \mbox{ (to b-node),}  \ \ 
W_{j,i} =  \frac{I_{j,i} (1-\eta_j^{sb}) } {\displaystyle \sum_{i'  \in \mathcal{G}_1 \cup \mathcal{G}_2} I_{j, i'}  } ,  \mbox{ (to another small node $i$)},  \mbox{ with} \nonumber\\
p_m^{sb} &:= & E[\eta_j^{sb}]  \mbox{ for any } j \in {\cal G}_m. 
\end{eqnarray}
In the above, $\{\eta_j^{sb} \}_{j \in {\cal G}_m}$,
 are i.i.d. (independent and  identically distributed) random variables with values between $0$ and $1$ for any $m$ and are independent of all other random variables. The weights from $b$-node are  given by $\{\eta_j^{bs}\}$,
 %
% the fractions, 
%$$
%\hspace{19mm}
%W_{b, j}  = \frac{\eta_j^{bs} }{\displaystyle \sum_i  \eta_i^{bs} } \mbox{ with } E[\eta_j^{bs}] = p_m^{bs}  \mbox{ (when $j\in \mathcal{G}_m $}),
%$$
where $\{\eta_j^{bs} \}_{j \in {\cal G}_m}$, are   {\it bounded i.i.d.  random variables   for any $m$} and are independent of others. We consider an alternate form of interconnections in Section~\ref{sec_alternate}.

We are interested in the performance of the nodes, which depends upon the weighted aggregate of the performance of other nodes with weights as given by $\{W_{j, i}\}_{j, i \in N}$, $\{\eta^{bs}_j\}$ and \ $\{W_{j,b} \}$.
As mentioned,  the weights of the performance measures may be stochastically different for the two groups. % because of heterogeneous connectivity. % structure in the graph.
We consider the following fixed-point (FP) equation  (in $\mathcal{R}^{n+1}$)  constructed using  functions $(f^{1},f^{2}, f^b)$, which in turn depend upon weighted averages $\{{\bar X}^m_i \}_i$ and ${\bar X}^b$ (constructed using weights $\{W_{j,i} \}$ and $\{W_{j,b}\}$),  and 
whose solution  ($i$-th component) represents   important performance measure of the nodes (node-$i$), as below:
\begin{eqnarray}
\label{Eqn_FixedeqGen 1}
X_i^{m}  &=& f^{m}(G^m_i,   {\bar X}^m_i , \eta_i^{bs} X^b) \mbox{ {\normalsize  for each }} i \in \mathcal{G}_m,  \ \mbox{and},  \\  
%\label{Eqn_FixedeqGen 2}
%X_i^{2}  &=& f^{2}(G^2_i,   {\bar X}^2_i , \zeta_i^{bs} X^b) \mbox{ {\normalsize  for each }} i \in \mathcal{G}_2,   \\  
\label{Eqn_Fixedeq2Gen 3}
X^b  &=&    f^b (   {\bar X}^b ) \ \ \ \mbox{ {\normalsize with aggregates}}    \\
{\bar X}^m_i  & := &  \displaystyle \sum_{j \in \mathcal{G}_1} X_j^1 W_{j, i} + \displaystyle \sum_{j \in \mathcal{G}_2} X_j^2  {W}_{j, i} \mbox{ {\normalsize  for each }} i \in \mathcal{G}_m, 
\label{Eqn_aggragate_forsn}\   \     \\
%\tilde{X}^s_i  & := &  \sum_{j \in \mathcal{G}_1} X_j^s W_{j, i} + \sum_{j \in \mathcal{G}_2} X_j^s \tilde{W}_{j, i}\mbox{ {\normalsize  for each }} i \in \mathcal{G}_2, \   \  \nonumber  
{\bar X}^b & := & \frac{1}{n}  \sum_{j 	\in \mathcal{G}_1} {X}^1_j W_{j, b} + \frac{1}{n}  \sum_{j \in \mathcal{G}_2} {X}^2_j W_{j, b} . \label{Eqn_Fixedeq2Gen 4}
\end{eqnarray}
In the above, $\{G^m_i\}_i$  is an i.i.d.  sequence for any fixed $m$ and further is independent of the sequence  corresponding to other $m$ and other random variables; the performance of the big node $X^b$ is defined  per small node (performance divided by $n$) and influences that of the small nodes via terms $\{\eta_i^{bs} X^b \}$.  
\noindent For any $n = n_1+n_2$ define mapping 
\begin{eqnarray}
\f := (f^b, \underbrace{f^{ 1},   f^{ 1} , \cdots ,  f^{ 1}}_{ n_1 \mbox{ times } }, \underbrace{f^{ 2}, f^{ 2}\cdots f^{2}}_{n_2 \mbox{ times }}), \mbox{ with } \label{Eqn_func1}  
{\bf x}:=  ( x_1^{1}, x_2^{1}, \cdots  x_{n_1}^{1}, x_{1}^{2}, \cdots x_{n_2}^{2}),
\end{eqnarray}
 component wise as below ($m \in \{1,2\}$):    
 \begin{eqnarray}
f_{1}( x_b, {\bf x} )  &:=& f^b ( {\bar x}_b),  \mbox{ with, }   {\bar x}_b  := \frac{1}{n}  \displaystyle \sum_{j 	\in \mathcal{G}_1} {x}^1_j W_{j, b} + \frac{1}{n}  \sum_{j \in \mathcal{G}_2} {x}^2_j  {W}_{j, b} \mbox{ and,} 
\label{Eqn_func2}\\
f^m_i ( x_b ,{\bf x}) &:=&   f^m(G_i^m,  {\bar x}_i^m, \eta_i^{bs}x_b ),\ \   {\bar x}^m_i  :=  \displaystyle \sum_{j \in \mathcal {G}_1} x^1_j     W_{j,i} + \displaystyle \sum_{j \in \mathcal{G}_2} x^2_j     { W}_{j,i} 
\ \forall  i \in {\mathcal{G}_m} . \hspace{5mm}
\label{Eqn_func3}
%\\
%f_i ( {\bf x}, x_b ) :=   f^s(\tilde{G}_i,  {\tilde {x}}_i, \zeta_i^{bs}x_b ) , \   {\bar x}_i  :=  \sum_{j \in G_1} x_j     W_{j,i} + \sum_{j \in \mathcal{G}_2} x_j    \tilde{ W}_{j,i} 
%\ \forall  i \in {\mathcal{G}_2},
%\vspace{-1mm}
\end{eqnarray}
It is clear that the  above mapping 
  represents the fixed point equations  corresponding to the random  operator (\ref{Eqn_FixedeqGen 1})-(\ref{Eqn_Fixedeq2Gen 4}), sample path wise (i.e., for each realization of the random variables, $\{G_j^m\}_{j, m},   \{\eta_j^{bs}\}_j, \{\eta_j^{sb}\}_j, \{I_{j,i}\}_{j,i}$). 

%\noindent
We assume the following:\\ 
{\bf B.1} The functions $f^{1}(.),f^{2}(.), f^b(.)$ are non-negative, continuous and are bounded by a constant $y< \infty$, % i.e., for e.g.,
   $$
0 \le f^{1}(g_1,x,x_b) ,  f^{2}(g_2,x,x_b) , f^b (x_b) \le y \mbox{ for all }  g_1,g_2, x, x_b. 
$$

This is a typical assumption required for existence of fixed points; we also require that the functions are bounded. This is a reasonable assumption as  many applications satisfy this, including our financial network.

Under the above assumption,  we have a measurable fixed point solution: 
\begin{lemma}
\label{CLemma_exist 1}
For any $n$ consider mapping $\f$ defined as in \eqref{Eqn_func1}-\eqref{Eqn_func3}. %:= (f^b, f^{1},   \cdots f^{1},f^{2}\cdots ,f^{2}  )$,
%with ${\bf x}:=  ( x_1^{1}, x_2^{1}, \cdots  x_{n_1}^{1}, x_{1}^{2}, \cdots x_{n_2}^{2}$),    component wise:   
%
%
%\begin{eqnarray*}
%\begin{array}{llllll}
%f_{1}( {\bf x}, x_b ) & :=& f^b ( {\bar x}_b),  &  {\bar x}_b &  :=&  \frac{1}{n} \sum_{j 	\in \mathcal{G}_1} {x}^1_j W_{j, b} + \frac{1}{n}  \sum_{j \in \mathcal{G}_2} {x}^2_j  {W}_{j, b},   \\ \\
%f_{1+i} ( {\bf x}, x_b )  & :=  &  f^m (G_i^m,  {\bar x}_i^m, \eta_i^{bs}x_b ) ,   \  \    &{\bar x}^m_i  & :=&  \sum_{j \in \mathcal {G}_1} x^1_j     W_{j,i} + \sum_{j \in \mathcal{G}_2} x^2_j     { W}_{j,i}  \\ \\
%&& && & \hspace{20mm}
%\mbox{ for all }   i \in {\cal G}_m.
%\end{array}
%\end{eqnarray*}
 Then we have (almost sure) random fixed point $(X_b^*, {\bf X}^{*})$ for each $n$ (see \cite{Measure_FP}).   
 \end{lemma}
 {\bf Proof:}
 Each component of this function 
 is a mapping from $[0, y ]^{n+1}  \to [0, y ]$   
  for almost all $\{G^{m}_i\}$,  $\{W_{j,i}\}$. Thus the function  $\f$  from  $[0,y]^{n+1} \to [0, y ]^{n+1}$.
 Further by continuity of $\f$, using the well known Brouwer’s fixed point Theorem, we have a deterministic fixed point   for  all realizations  of $\{G^{m}_i\}$  , $\{W_{j,i}\}$, $ \{W_{j,b}\}$ and $\{\eta_j^{bs}\}_j$ under {\bf B.1}. Then the overall measurability result follows by \cite[Theorem 8]{Measure_FP}. \hspace{-2mm}. \hspace{-2mm} \eop

{\bf Assumptions on the graph structure:}  We require that the number of  nodes  influencing any given  node,  grows asymptotically linearly with $n$ for almost all   sample paths. Towards this, first define the following set:

\vspace{-2mm}
{\small
\begin{eqnarray}
{\cal E} &: =&  \left\lbrace  \omega : \lim_{n \to \infty } \displaystyle \sum_{j \in \mathcal{G}_1}    \left|  
    \frac{1 }{  \displaystyle \sum_{i \in   \mathcal{G}_1 \cup  \mathcal{G}_2 } I_{j, i}}  - \frac{1}{ n\gp{1} }    \right |  =0,  ~ 
\lim_{n \to \infty } \displaystyle \sum_{j \in \mathcal{G}_2}    \left|    \frac{1 }{ \displaystyle \sum_{i \in   \mathcal{G}_1 \cup  \mathcal{G}_2 } I_{j, i}}  - \frac{1}{ n\gp{2} }    \right | = 0  
\right\rbrace, \nonumber \\
&&
\hspace{-6mm}
\mbox{\normalsize with, } 
\gp{1} := \gamma p_1+ (1-\gamma)p_{c_1}  \mbox{\normalsize  and  } \gp{2} :=  \gamma p_{c_2}+ (1-\gamma) p_2.
\hspace{2mm}
\label{Eqn_calE}
\end{eqnarray}}
%
% We require the following assumption:\\
% \noindent {\bf B.2} Assume $\gp{1}  > 0$ and $\gp{2} >0 $.  Also consider only graphs for which, $P({\cal E}) = 1$. \\
% The initial results are with this assumption, latter (in sub-section \ref{sec_assumptions}) we have  results under more general conditions (with $P({\cal E}) < 1$). We also provide an equivalent  assumption on the  growth pattern of the  graphs in the same sub-section. 

We require the following assumption, which has two parts. In part {\bf(C)} we consider that $\{I_{j,i}\}$ need not be independent; however they  remain independent of the other quantities like $\{G_j^m\}$, $\{\eta_j^{sb} \}_{j \in {\cal G}_m}$ etc:\\
\noindent {\bf B.2} 
 Assume $\gp{1}  > 0$ and $\gp{2} >0 $.  Also consider only graphs for which, $P({\cal E}) = 1$. 
\\
{\bf {\bf B}.2(C), extra assumption for correlated $\{I_{j,i}\}$}: When $\{I_{j, i}\}_j$ are not i.i.d. we additionally require: 
$$
\frac{1}{n \gamma_m} \sum_{j \in {\cal G}_m} I_{j, i} \to p_{m', m} \mbox{  a.s. } \mbox{ for any } m, m' \mbox{ and } i \in \mathcal{G}_{ m'}.
$$
The initial results are with this assumption, latter (in sub-section \ref{sec_assumptions}) we have  results under more general conditions (with $P({\cal E}) < 1$). We also provide an equivalent  assumption on the  growth pattern of the  graphs in the same sub-section.

Our assumptions on graph structure are quite general. 
We firstly require that the  marginal probabilities related to random connections, are the same within a group, i.e., $P(I_{j,i}=1)= p_{mm'}$  as in  equation \eqref{Eqn_connections} (for each $i,j$). Further the joint probabilities are supposed to satisfy {\bf B.2}, i.e., mainly $P({\cal E}) = 1$. 
Furthermore, when the connections $\{I_{j,i}\}$ are not independent,  the results are still true under the most natural assumption {\bf B.2(C).} Regular graphs constructed in sub-section \ref{sec_MCsimulation}  are some example graphs that satisfy our assumptions. Our results can also be extended partially to Erd\H{o}s-R\'enyi graphs, if required via  Theorem \ref{Thm_onsetD}.
The most restrictive assumption is {\bf B.2} and one can avoid such an assumption by considering a different  structure on weights as discussed in sub-section \ref{Section_examples} (see e.g., equation \eqref{Eqn_reWeights 1}). With this our results   can cover many more graphs.
\subsection{Aggregate fixed points:}
We rewrite the  fixed point equations in terms of weighted averages and  first analyze the aggregate system.
Towards this, define the following random variables, that depend upon real constants $({\bar x}^m_i, x_b)$: 
\begin{eqnarray}\label{Eqn_xi 1}
\xi^m_i ({\bar x}^m_i, x_b) :=   f^{ m} (	G_i^m,   {\bar x} _i^m,   \eta_{i}^{bs} x_b)  \mbox{ for any }  i \in {\cal G}_m,
\end{eqnarray} and assume for each value of $m \in  \{1,2\}$ (see \eqref{Eqn_connections}, \eqref{Eqn_Weights 1}):   \\

\noindent  {\bf B.3}    $|\xi^{m}_i (x, x_b)  - \xi^{m}_i (u, u_b) | \le \sigma ( |x- u| + \varsigma|x_b- u_b|)$  with   $\sigma \le 1$ and 
$0 <\varsigma \le 1$. 
Basically we require the following: 
\begin{eqnarray*}
|f^b({\bar x}_b) - f^b ({\bar u}_b) | &\le &  |{\bar x}_b-{\bar u}_b | \mbox{ for all } {\bar x}_b, {\bar u}_b, \mbox{ and, } \\
|f^m (g, x, \eta x_b) - f^m(g, u, \eta u_b)| 
&\le &\sigma ( |x- u| + \varsigma|x_b- u_b|) \mbox{ for all } x, u, x_b, u_b, g,  \eta. 
\end{eqnarray*} 
{\bf B.3} is  a typical contraction mapping type of  assumption 
that ensures the existence (and uniqueness) of fixed points. Observe this assumption   does not imply strict contraction mapping (as $\sigma \le 1$ and not $\sigma <1$), but is nonetheless sufficient. 
\\
\noindent{\bf B.4}  Assume  $\varrho \le 1$, where 

\vspace{-4mm}
 {\small
 \begin{eqnarray*} \hspace{4mm}
\varrho &:=&  \max \bigg \lbrace \frac{\gamma p_{c_1} (1-p_1^{sb})  }{\gp{1}} + \frac{(1-\gamma)p_2  (1-p_2^{sb}) }{\gp{2} } ,  \ \ \  \ \  
   \frac{\gamma p_1   (1-p_1^{sb}) }{\gp{1}} + \frac{(1-\gamma)p_{c_2}  (1-p_2^{sb})  }{\gp{2} } \bigg \rbrace  \\ && + \left ( \gamma p_1^{sb} + (1-\gamma) p_2^{sb} \right ).
 \end{eqnarray*}}
 Observe that {\bf B.4} is readily satisfied, for symmetric conditions, for example, when
 $\gamma =0.5$, $p_1^{sb} = p_2^{sb}$, $p_1 =p_2$ and $p_{c_1} = p_{c_2}$.

 Let $\x^m := ({\bar x}^m_1, {\bar x}^m_2 \cdots )  \mbox{ for } m  \in \{ 1,2\}  \mbox{ and }  {\x} := (\x^1 , \x^2)$.
 Consider  the  following  operators   on  Banach  space\footnote{
\label{footnote_sinf}
 Here $s^\infty$ is the space (subset) of bounded sequences  equipped  with $l_\infty$ norm  $ |\x |_\infty :=  \sup_{i} |x_i|$, 
$$
s^\infty := \{  \x =  (x_1, x_2, \cdots ) : x_i \in [0, y ]  \mbox{ for all }  i    \}. 
$$ We also consider different  other norms (and/or topologies) on $s^{\infty}$ for various parts of the   proofs in the appendices and the same is mentioned at the relevant parts} $  [0, y] \times s^\infty \times s^\infty$,   one for each $ n= n_1+ n_2$: 
\begin{eqnarray}
\label{Eqn_bar_fixedpoint_randGen_1}
{\bar {\bf f}}^n ({\bar x}_b, {\x} ) & = & ({\bar f}^n_b, \ \ \ {\bar f}_1^{n,1}, {\bar f}_2^{n,1}, \cdots, {\bar f}_{n_1}^{n,1}, \cdots,  \ \ \ {\bar f}_{1}^{n,2}, \cdots,  {\bar f}_{n_2}^{n,2}, \cdots),\end{eqnarray}
where for any $(n,m)$ we have, 
\begin{eqnarray} 
\label{Eqn_key }
\label{Eqn_bar_fixedpoint_randGen_1_more}
{\bar f}^{n, m}_i ( {\bar x}_b, \x ) &:=&  
\left \{  \begin{array}{lll}
 \displaystyle \sum_{j \in \mathcal{G}_1}   \xi^1_{j} ({\bar x}^1_j, x_b)  W_{j, i} + \displaystyle \sum_{j \in \mathcal{G}_2}   \xi_{j}^2  ({\bar  x}^2_j, x_b)   {W}_{j, i} &\mbox{ if }  i \in \mathcal{G}_m ,    \\ 
 0  &\mbox{\normalsize else,}  \mbox{ \normalsize and, }  \\
\end{array} \right.   \\ 
 {\bar f}^n_b ({\bar x}_b,{\x})& := & 
 \frac{1}{n}
\displaystyle \sum_{j \in \mathcal{G}_1}   \xi^{1}_{j} ({\bar x}^{1}_j, x_b) W_{j, b} +  \frac{1}{n}
\displaystyle \sum_{j \in \mathcal{G}_2}   \xi^{2}_{j} ({\bar x}^{2}_j, x_b) W_{j, b} \ \mbox{\normalsize with } \  x_b :=  f^b (  {\bar x}_b) . 
\nonumber 
\end{eqnarray}

\ignore{
\begin{eqnarray}\nonumber
\label{Eqn_key }
{\bar f}^{n, m}_i ( {\bar x}_b, \x ) &:=&  
\left \{  \begin{array}{lll}
 \displaystyle \sum_{j \in \mathcal{G}_1}   \xi^1_{j} ({\bar x}^1_j, x_b)  W_{j, i} + \displaystyle \sum_{j \in \mathcal{G}_2}   \xi_{j}^2  ({\bar  x}^2_j, x_b)   {W}_{j, i} &\mbox{ if }  i \in \mathcal{G}_m ,    \\ \\
 0  &\mbox{\normalsize else,}   \hspace{5mm} \mbox{\normalsize and, }
\end{array} \right .  \\
{\bar f}^n_b ({\bar x}_b,{\x})& := & 
 \frac{1}{n}
\displaystyle \sum_{j \in \mathcal{G}_1}   \xi^{1}_{j} ({\bar x}^{1}_j, x_b) W_{j, b} +  \frac{1}{n}
\displaystyle \sum_{j \in \mathcal{G}_2}   \xi^{2}_{j} ({\bar x}^{2}_j, x_b) W_{j, b} \ \mbox{\normalsize with } \  x_b     :=  f^b (  {\bar x}_b).
\nonumber
\end{eqnarray}}
\noindent  Thus we require the fixed point of the  operator:
$$
{\bar {\bf f}}^n   \mbox{ where }   {\bar{\bf f}}^n  :  {[0, y]}   \times  {s}^\infty   \times  {s}^\infty  \to   {[0, y]}   \times  {s}^\infty  \times  {s}^\infty ,
$$which provides the   aggregate vectors, $({\bar X}_b, \{{\bar X}_i^m \}_{i \le n_m})$ given in \eqref{Eqn_aggragate_forsn}-\eqref{Eqn_Fixedeq2Gen 4}.
{\it Observe here that, for uniformity we have infinite dimensional mappings even for finite $n$, where the extra components  are set to zero functions (i.e., ${\bar f}^{n, m}_i  \equiv 0$ for all $i > n_m$)}. 

Recall that the weights sum up to one, i.e., $ \sum_{i} W_{j, i} +W_{j,b} = 1$  for all $j$.  Thus
the {\it idea is to derive a kind of mean-field analysis where their expected values will approximate the aggregates.}
%
%
%It is clear that the  fixed points   of  above operators  equal  the   aggregate vectors, $({\bar X}_b, \{{\bar X}_i^m \}_{i \le n})$ given in \eqref{Eqn_Fixedeq2Gen 3 1}-\eqref{Eqn_Fixedeq2_Gen 3}.
Towards this,  as  a first step,  we   analyze the point-wise limits of the above operator.
\begin{lemma}{\bf [Constant sequences]}
\label{lem: LLN_fixed point}
%{\color{red} (i) For any $\x^m \in s^\infty \times s^\infty$ and ${\bar x}_b$, the sequence   ${\bar {\bf f}}^n ({\bar x}_b, \x^m)$ (with  ${\bar {\bf f}}^n$ defined in \eqref{Eqn_bar_fixedpoint_randGen_1}) 
%converges along a sub sequence ( This is immediate because of Bolzano Weierstrass theorem)}.  \\ 
%(ii)
Assume  {\bf B.1}-{\bf B.2}. Consider any constant sequence ${\x}=(\x^1 , \x^2)$, i.e., sequence with   $\x^m  = ( {\bar x}^m,   {\bar x}^m,   \cdots )$ for some ${\bar x}^m < \infty$,   for each $m$. The functions ${\bar {\bf f}}^n ({\bar x}_b, {\x} )$  defined in \eqref{Eqn_bar_fixedpoint_randGen_1}-\eqref{Eqn_bar_fixedpoint_randGen_1_more}  converge component-wise and the limits
%
%the limit superiors in the definition of   ${\bar {\bf f}}^\infty$    given by   \eqref{Eqn_bar_fixedpoint_limtGen 1}  become limits  (under {\bf B.2})  and 
equal almost surely (with $x_b = f^b ({\bar x}_b)$ as in \eqref{Eqn_Fixedeq2Gen 3}, and see \eqref{Eqn_connections}):
%\begin{multline*}
%{\bar f}_i^{\infty,1}  ({\bar x}_b,\x)  =      E_{G^{1}_i, \eta_i^{bs}} \left [ \xi^{1}_i ({\bar x}^1, x_b)   \right ] \frac{ \gamma p_1}{\gp{1} }(1-p_1^{sb})  \\
%+  E_{G^{2}_i, \eta_i^{bs}} \left [ \xi^{2}_i ({\bar x}^2, x_b)   \right ] \frac{(1-\gamma)  p_{c_2}}{ \gp{2}}  (1-p_2^{sb})\   \mbox{for} \  i \in \mathcal{G}_1,
%\end{multline*}
%
\begin{eqnarray}
{\bar f}_i^{\infty,m}  ({\bar x}_b,\x)  &=&      E_{G^{1}_i, \eta_i^{bs}} \left [ \xi^{1}_i ({\bar x}^1, x_b)   \right ] \frac{ \gamma p_{1m}}{\gp{1} }(1-p_1^{sb}) \nonumber \\
& &+  E_{G^{2}_i, \eta_i^{bs}} \left [ \xi^{2}_i ({\bar x}^2, x_b)   \right ] \frac{(1-\gamma)  p_{2m}}{ \gp{2}}  (1-p_2^{sb})\   \mbox{for} \  i \in \mathcal{G}_m, \ \forall m, \nonumber\\
{\bar f}_b^\infty  ({\bar x}_b,\x) &= &  \gamma E_{ G^{1}_i, \eta_i^{bs}}  [\xi^{1}_i ({\bar x}^1, x_b)  ]  p_1^{sb}  + (1-\gamma) E_{ {G}^{2}_i, \eta_i^{bs}}  [\xi^{2}_i ({\bar x}^2, x_b)  ] p_2^{sb}  . \label{Eqn_barf_limit 1}
 \end{eqnarray}
\ignore{
\begin{multline*}
{\bar f}_i^{\infty,2}  ({\bar x}_b,\x)  =      E_{G^{1}_i, \eta_i^{bs}} \left [ \xi^{1}_i ({\bar x}^1, x_b)   \right ] \frac{\gamma p_{c_1}}{ \gp{1}}(1-p_1^{sb})  + \\
 E_{G^{2}_i, \eta_i^{bs}} \left [ \xi^{2}_i ({\bar x}^2, x_b)   \right ] \frac{(1-\gamma ) p_2}{ \gp{2} }  (1-p_2^{sb}) \   \mbox{for} \  i \in \mathcal{G}_2,
\end{multline*} }
\ignore{
\begin{multline}
{\bar f}_b^\infty  ({\bar x}_b,\x) =   \gamma E_{ G^{1}_i, \eta_i^{bs}}  [\xi^{1}_i ({\bar x}^1, x_b)  ]  p_1^{sb}  + (1-\gamma) E_{ {G}^{2}_i, \eta_i^{bs}}  [\xi^{2}_i ({\bar x}^2, x_b)  ] p_2^{sb}  . \label{Eqn_barf_limit 1}
\end{multline}}
%\end{eqnarray}}}
\end{lemma}
{\bf Proof:} is available in Appendix A.
\eop
\\
In the above, $E_{X,Y}$ represents the  expectation with respect to  $X,Y$.
For constant sequences
%(as in part (ii) of the lemma)
the random variables $\{\xi^{m}_i (\bar {x}^m, x_b)\}_{i \in {\cal G}_m}$ are i.i.d., and hence the first equation of  \eqref{Eqn_barf_limit 1} has   same  right hand side  value for all $i \in \mathcal{G}_m$.
We now 
define
the  following `limit' operator, which 
in view of the above lemma equals a limit for constant sequences:

\vspace{-3mm}
 \begin{eqnarray}
\label{Eqn_bar_fixedpoint_limtGen 1}
\bar {\bf f}^\infty ({\bar x}_b,\x)  &=&  ({\bar f}^\infty_b, \ \ \  {\bar f}^{\infty,1}_1,  \cdots , \ \ \ {\bar f}^{\infty,2}_1, \cdots) \mbox{ with} \\
{\bar f}_b^{\infty} ({\bar x}_b,\x)&: =& \displaystyle \lim \sup_{n}{\bar f}^n_b ({\bar x}_b,\x) \  \mbox{and}, \nonumber \\
 {\bar f}_i^{\infty,m} ( {\bar x}_b,\x ) &:=
&\displaystyle \lim \sup_{n} {\bar f}^{n,m}_i ({\bar x}_b,\x) \mbox{ for all } i \mbox{ and } m.\nonumber
\end{eqnarray}
The idea is to show that the aggregate fixed points of the original system converge towards the fixed point of this `limit' system/operator (\ref{Eqn_bar_fixedpoint_limtGen 1}) (more specifically the fixed point of the three dimensional system \eqref{Eqn_barf_limit 1}), as $n \to \infty$. %
%
%By Theorem \ref{Thm_MainGen 1},  given below, 
%one such sequence would be the almost sure limit of the solutions of the aggregate fixed point equations (\ref{Eqn_bar_fixedpoint_randGen_1}).  
%Thus one will have to solve a three-dimensional fixed point equation corresponding to the above function (\ref{Eqn_barf_limit 1}) in the limit. 
%And then random fixed points  (\ref{Eqn_FixedeqGen 1})-(\ref{Eqn_Fixedeq2Gen 3}) are asymptotically independent depending upon the other nodes only via the aggregate fixed point, as given by the Theorem below.
We require another assumption:

\noindent{\bf B.5} The limit system ${\bar f}^\infty$ given by \eqref{Eqn_bar_fixedpoint_limtGen 1} has a fixed point among constant sequences.

The assumption  {\bf B.5} demands that  the limit system has a fixed point among constant sequences. This kind of an assumption can be restrictive, but provides required convergence results in the most general settings. Further it might be readily satisfied by some future applications; thus this assumption gives more flexibility for future applications.
Furthermore, we prove this assumption  (along with others)  in   Theorem \ref{Lemma_limit_system_uniquness} and the assumptions of the latter theorem are readily satisfied by the financial network based case studies of Section \ref{sec_finance}.

We now prove one of the the main results of this paper:
\begin{thm}
\label{Thm_MainGen 1}
Assume {\bf B.1}-{\bf B.5}.
 %Also assume either   $0 < p_1^{sb}, p_2^{sb}   < 1$      or   $\sigma < 1$ in {\bf B.3}.
The aggregates  of the random system \eqref{Eqn_FixedeqGen 1}-\eqref{Eqn_Fixedeq2Gen 4},  which are   FPs of    \eqref{Eqn_xi 1}-\eqref{Eqn_bar_fixedpoint_randGen_1}, denoted by $({\bar X}_b, {\X}) (n) = ({\bar X}_b, \{{\bar X}_i^m\}_{i,m}) (n) $, converge   as $n \to \infty$ along a sub-sequence. That is, there exists $k_n \to \infty$ such that:
\begin{eqnarray}
 {\bar X}_i^m (k_n) \to {\bar x}_i^{m\infty*} \   \forall i,m ,  \mbox{ and }  {\bar X}_b (k_n)  \to {\bar x}_b^{\infty*} \mbox{ almost surely (a.s.),}
 \label{Eqn_conv_aggr}
\end{eqnarray}
 where $({\bar x}_b^{\infty*}, \x^{\infty*})$ with $\x^{\infty*} := ({\bar x}_1^{1\infty*}, {\bar x}_2^{1\infty *} , \cdots {\bar x}_1^{2\infty*}, {\bar x}_2^{2\infty*}, \cdots )$ is an FP of 
  the limit system given by  \eqref{Eqn_bar_fixedpoint_limtGen 1}.
Further    (any sequence of) FPs of the  original system \eqref{Eqn_FixedeqGen 1}-\eqref{Eqn_Fixedeq2Gen 3} converge almost surely (along the sub-sequence of \eqref{Eqn_conv_aggr}, i.e., as $k_n \to \infty$): %to the limit :
\begin{eqnarray}
\label{Eqn_Thm_1_orginaleq_convg}
X^b (k_n)  & \to &  f^b (  {\bar x}_b^{\infty*}) \mbox{ and } \hspace{3mm}\label{Eqn_Act_Fixed_with_m}   \
X^m_i (k_n)  \ \to \   f^m(G_i^m,  {\bar x}_i^{m\infty*}, \eta_i^{bs} X_b^{\infty*})\ \forall i,m   .
\end{eqnarray}
\end{thm}
{\bf Proof:} available in Appendix B. 
\eop
\\
By the above theorem, under minimal conditions on the fixed points of the limit system \eqref{Eqn_bar_fixedpoint_limtGen 1},  the fixed points of the finite $n$-system can be studied using that of the limit system, the latter is an approximation and the approximation would be better for larger $n$.  %
We now consider an additional assumption under which   one will have unique fixed points, and in addition, the FP is a constant sequence for limit system:
\ignore{
{\color{red}
\newcommand{\yb}{{\underline y}} 
\newcommand{\yu}{{\bar y}}
\newcommand{\ybb}{{\underline {\tilde y}}} 
\newcommand{\yuu}{{\bar {\tilde y}}}
{\bf Unique fixed points proof, trying out:}
Consider any finite $n$-system and consider the following norm for this proof:
\begin{eqnarray}
|| (\x, {\bar x}_b) ||_1 := \frac{1}{n}
\sum_m \sum_{j \in {\cal G}_m}\left ( |x^m_j|+ \varsigma  |{\bar x}_b) | \right )
\end{eqnarray}
where,  $0 < \varsigma < 1$.
Assume that the function is upper and lower bounded 
$$0 < \yb \le  f \le \yu.$$
Note that for any $ ( {\bar x}_b, \x ) $, 
\begin{eqnarray*}
||{\bar {\bf f}}^n ( {\bar x}_b, \x )  ||_1 &=& \frac{1}{n}\sum_m 
\sum_{i \in {\cal G}_m} | \left ( |{\bar f}^{n, m}_i ( {\bar x}_b, \x )|+ \varsigma  |  {\bar f}^{n}_b ( {\bar x}_b, \x ) | \right ) \\
&=&  \frac{1}{n}\sum_m 
\sum_{i \in {\cal G}_m}  \left (  \sum_{j \in \mathcal{G}_1}   \xi^1_{j} ({\bar x}^1_j, x_b)  W_{j, i} + \displaystyle \sum_{j \in \mathcal{G}_2}   \xi_{j}^2  ({\bar  x}^2_j, x_b)   {W}_{j, i}  \right )
\\
&&
+ \varsigma \frac{1}{n}  \sum_{j \in \mathcal{G}_1}   \xi^{1}_{j} ({\bar x}^{1}_j, x_b) W_{j, b} + \varsigma  \frac{1}{n} 
\displaystyle \sum_{j \in \mathcal{G}_2}   \xi^{2}_{j} ({\bar x}^{2}_j, x_b) W_{j, b} \\
&=&  \frac{1}{n}    \sum_{j \in \mathcal{G}_1}   \xi^1_{j} ({\bar x}^1_j, x_b) \sum_m 
\sum_{i \in {\cal G}_m}  W_{j, i} +   \frac{1}{n}  \sum_{j \in \mathcal{G}_2}   \xi_{j}^2  ({\bar  x}^2_j, x_b) \sum_m 
\sum_{i \in {\cal G}_m}   {W}_{j, i}   
\\
&&
+\varsigma \frac{1}{n} \sum_{j \in \mathcal{G}_1}   \xi^{1}_{j} ({\bar x}^{1}_j, x_b) W_{j, b} +\varsigma   \frac{1}{n}  
\displaystyle \sum_{j \in \mathcal{G}_2}   \xi^{2}_{j} ({\bar x}^{2}_j, x_b) W_{j, b} \\
&=&   \frac{1}{n}     \sum_{j \in \mathcal{G}_1}   \xi^1_{j} ({\bar x}^1_j, x_b)  \left ( \sum_m 
\sum_{i \in {\cal G}_m}  W_{j, i} +\varsigma W_{j,b} \right ) +   \frac{1}{n}  \sum_{j \in \mathcal{G}_2}   \xi_{j}^2  ({\bar  x}^2_j, x_b) \left (  \sum_m 
\sum_{i \in {\cal G}_m}   {W}_{j, i}   +\varsigma W_{j,b} \right ) 
 \\
&=&   \left (\frac{1}{n}     \sum_{j \in \mathcal{G}_1}   \xi^1_{j} ({\bar x}^1_j, x_b) \left ( 1- \eta_j^{sb} + \varsigma \eta^{sb}_j  \right )
+ \frac{1}{n}  \sum_{j \in \mathcal{G}_2}   \xi_{j}^2  ({\bar  x}^2_j, x_b)  \left ( 1- \eta_j^{sb} + \varsigma \eta^{sb}_j  \right ) \right ) 
\end{eqnarray*}
Thus clearly, $ \yb \le ||{\bar {\bf f}}^n (\x, {\bar x}_b) ||_1
\le \yu $ and hence
\begin{eqnarray*}
{\bar {\bf f}}^n : {\cal D}_n \to {\cal D}_n, \mbox{ where the domain }  {\cal D}_n := \left \{ ( {\bar x}_b, \x )  : 
  \ybb \le || ( {\bar x}_b, \x )  ||_1  \le \yuu  \right \} \mbox{ with } \\
  \ybb := \yb \frac{1}{n} \sum_m \sum_{j \in {\cal G}_m} \left ( 1- \eta_j^{sb} + \varsigma \eta^{sb}_j  \right ) \mbox{ and }
   \yuu := \yu \frac{1}{n} \sum_m \sum_{j \in {\cal G}_m} \left ( 1- \eta_j^{sb} + \varsigma \eta^{sb}_j  \right )
\end{eqnarray*}
Working in exactly similar lines we have:

{\small 
\begin{eqnarray*}
||{\bar {\bf f}}^n ( {\bar x}_b, \x ) - {\bar {\bf f}}^n ( {\bar u}_b, \u )  ||_1  & \le & \frac{1}{n}     \sum_{j \in \mathcal{G}_1}   |  \xi_{j}^1  ({\bar  x}^1_j, x_b) - 
 \xi_{j}^1  ({\bar  u}^1_j, u_b) | + \frac{1}{n}  \sum_{j \in \mathcal{G}_2} |  \xi_{j}^2  ({\bar  x}^2_j, x_b) - 
 \xi_{j}^2  ({\bar  u}^2_j, u_b) |  \\
 &\le & \sigma \frac{1}{n}  \sum_m \sum_{j \in {\cal G}_m} \left ( | {\bar x}_j^m - {\bar u}_j^m| +\varsigma | {\bar x}_b -  {\bar u}_b| \right )    \left ( 1- \eta_j^{sb} + \varsigma   \eta^{sb}_j  \right )   % \\ 
% &=&  \sigma  || ( {\bar x}_b, \x )   - ( {\bar u}_b, \u )  ||_1
\end{eqnarray*}}
It is easy to see that 

{\small \begin{eqnarray*}
\left | \frac{1}{n}  \sum_m \sum_{j \in {\cal G}_m} \left ( | {\bar x}_j^m - {\bar u}_j^m| +\varsigma | {\bar x}_b -  {\bar u}_b| \right )    \left ( 1- \eta_j^{sb} + \varsigma   \eta^{sb}_j  \right )  - \frac{1}{n}  \sum_m \sum_{j \in {\cal G}_m} \left ( | {\bar x}_j^m - {\bar u}_j^m| +\varsigma | {\bar x}_b -  {\bar u}_b| \right )    \left ( 1- p_m^{sb} + \varsigma   p^{sb}_m  \right )  \right |  \hspace{-100mm}  \\
&\le & (1+\varsigma) \bar{y} \frac{1}{n}  \sum_m \sum_{j \in {\cal G}_m}    \left | \left ( 1- \eta_j^{sb} + \varsigma   \eta^{sb}_j  \right )  -      \left ( 1- p_m^{sb} + \varsigma   p^{sb}_m  \right ) \right |  
\end{eqnarray*}}

For any $(i, m)$ we have:
\begin{eqnarray*}
  |{\bar f}^{n, m}_i ( {\bar x}_b, \x )| - |{\bar f}^{n, m}_i ( {\bar u}_b, \u )|  &\le & 
\sum_m  \sum_{j \in {\cal G}_m}   |  \xi_{j}^1  ({\bar  x}^1_j, x_b) - 
 \xi_{j}^1  ({\bar  u}^1_j, u_b) |  W_{j,i}  \\
 &\le & \sigma \max_{m, i \in {\cal G}_m} \left ( |{\bar x}_i^m - {\bar u}_i^m | + | {\bar x}_b -  {\bar u}_b| \right )
 \sum_m  \sum_{j \in {\cal G}_m} \frac{ \left (  1-\eta_j^{sb} \right ) }{ \sum_{i' } I_{j,i'}}
\end{eqnarray*}
Now we know by the Lemma \ref{Master lemma}  that 
$$
 \sum_m  \sum_{j \in {\cal G}_m} \frac{ \left (  1-\eta_j^{sb} \right ) }{ \sum_{i' } I_{j,i'}} \to  \frac{
 \gamma (1-p_1^{sb} ) } {\gp{1}} + \frac{(1-\gamma) (1-p_2^{sb} ) } {\gp{2}}
$$
}
}

%, i.e., it would also be a fixed point of the simpler (3 dimensional) system given by \eqref{Eqn_barf_limit 1}. We will also have converge along the  entire sequence as given below: 
%\begin{thm}
 %\label{Lemma_limit_system_uniquness}
%Assume {\bf B.1} to {\bf B.4}. Also assume,  either  $\varrho <1$ or  $\sigma < 1$ in {\bf B.3}. Then we have unique fixed point   of the limit system \eqref{Eqn_bar_fixedpoint_limtGen 1}, which is a constant sequence, i.e., we have 
%${\bar x}_i^{m\infty*} = {\bar x}^{m\infty*} $ for all $i \in {\cal G}_m$ in equation %\eqref{Eqn_conv_aggr}. This limit is the fixed point of the three dimensional system given by \eqref{Eqn_barf_limit 1}.  Further there exists an ${\bar N} < \infty$, such that  equations 
 %\eqref{Eqn_xi 1}-\eqref{Eqn_bar_fixedpoint_randGen_1}  (whose FPs provide the aggregates) have unique fixed point for all $n \ge {\bar N}$ and then the convergence  in \eqref{Eqn_conv_aggr} and \eqref{Eqn_Act_Fixed_with_m} is along the original sequence, i.e., as $n \to \infty$. 

%\end{thm}
%{\bf Proof:} available in Appendix B.
%\eop

\begin{thm}{\bf [Unique fixed points]}
\label{Lemma_limit_system_uniquness}
Assume {\bf B.1}-{\bf B.4} and also assume  $\eta_j^{sb} \ge  \underline{\eta} \ge 0$ a.s., with
$\sigma \left ( 1- \underline{\eta} + \varsigma   \underline{\eta} \right ) < 1$.   Then we have unique fixed point of the finite $n$-system \eqref{Eqn_bar_fixedpoint_randGen_1_more} for each $n$. We also  have unique fixed point for the limit system \eqref{Eqn_bar_fixedpoint_limtGen 1}, which is a constant sequence, i.e., we have 
${\bar x}_i^{m\infty*} = {\bar x}^{m\infty*} $ for all $i \in {\cal G}_m$ and for each $m$ in equation \eqref{Eqn_conv_aggr}. This limit is the unique  fixed point of the three dimensional system given by \eqref{Eqn_barf_limit 1}.
\end{thm}
{\bf Proof:} available in Appendix B.
\eop\\
The above theorem immediately implies the following corollary: {\it one can solve three dimensional system \eqref{Eqn_barf_limit 1} and derive the fixed points for large dimensional system given by \eqref{Eqn_FixedeqGen 1}-\eqref{Eqn_Fixedeq2Gen 3} almost surely.}

\begin{cor}{\bf [Three dimensional approximation]}
\label{Cor_Three dimensional approximation}
 Assume the conditions of Theorem  \ref{Lemma_limit_system_uniquness}.  Then convergence in the equations \eqref{Eqn_conv_aggr}-\eqref{Eqn_Thm_1_orginaleq_convg} is along the  original sequence, i.e., as $n\to \infty$. Further,  $({\bar x}_b^{\infty*}, \x^{\infty*})$
 is a constant sequence and is the unique fixed point of the three dimensional system given by  \eqref{Eqn_barf_limit 1}. 
 \end{cor}
 {\bf Proof:} available in Appendix B. \eop \\
\noindent \textbf{Remarks:} We have several remarks regarding the above results.  \\
\noindent $\bullet$ Observe  that  the  aggregate  fixed  points  converge almost surely  to the same limit; further the aggregates at limit are also constant across the agents of the same group as given by Theorem \ref{Lemma_limit_system_uniquness} and Corollary \ref{Cor_Three dimensional approximation}.\\
$\bullet$
The fixed points of the finite $n$ system converge to that of the limit system.  From \eqref{Eqn_Thm_1_orginaleq_convg} and Corollary \ref{Cor_Three dimensional approximation}, the fixed points are asymptotically independent and depend upon the other nodes only via an almost sure constant (representing the aggregate), which is the same for all $i$ in a group. %, which is also the same for all nodes from the same group.

\noindent $\bullet$ Under the more general assumptions of the Theorem~\ref{Thm_MainGen 1}  \underline{the aggregate fixed}  \underline{points need not be unique}, for initial $n$. However, any sequence of fixed points (one for each $n$) converges towards that of the limit system (when it has a unique fixed point).  If the limit system has many fixed points then every such sub-sequence converge to one among these fixed points.
Thus the three dimensional fixed points of  \eqref{Eqn_barf_limit 1}  (if any) are useful even under  general conditions. 

\noindent
$\bullet$ The graphical  model of the current paper  is a significantly generalised version of our  previous model considered in \cite{Systemicrisk}; and it  reduces to  the model considered in \cite{Systemicrisk},  when $p_1= p_2=p_{c_1}=p_{c_2}= p_{ss}$ and $\gamma \in \lbrace 0,1 \rbrace$. Also, the current graphical model is heterogeneous in many more aspects,  e.g.,  the interconnection probabilities, the connections to b-node etc. 
\\
$\bullet$
From \eqref{Eqn_barf_limit 1},
the limiting fixed point is dependent on the interconnection probability ($p_{c_m}$) as well as the group-wise connectivity parameters ($\{p_{m}\}$). While in
our initial model of 
\cite{Systemicrisk}, with only one group and a big node,     the limiting fixed point   is independent of the exact value of the connectivity parameter (referred to as $p_{ss}$ in \cite{Systemicrisk});  it only requires that $p_{ss} > 0 $ (together with other  appropriate  assumptions on graph structure). Thus it establishes that when more groups are coupled, the limits depend upon inter-group as well as intra-group connectivity parameters. 
%, basically   it only requires that every node can potentially influence every other node directly or indirectly.

%{\color{red}
%\noindent $\bullet$ When one has only existence of fixed point for the limit aggregate system (as in {\bf B.5}), we could prove  required convergence only along a sub-sequence. However under the extra conditions of  Lemma \ref{Lemma_limit_system_uniquness}, we established the existence of unique fixed point for limit system as well as for finite systems with all big enough  $n$; this also   gave us convergence along $n \to \infty $.  Further with this extra condition, one can derive the fixed points by just solving the three-dimensional fixed point equations given by \eqref{Eqn_barf_limit 1}.
%}
%Further one can easily generalize this result to a case without uniqueness for limit system also, using the results available for Maximum Theorem  (\cite{Berge,Equilibria}). Here, the fixed points  of the finite systems converge along a sub-sequence to one of the fixed points of the limit system (see \cite{Berge}).  

\subsection{Assumption {\bf B.2}}
\label{sec_assumptions}

In the previous sub-sections we consider graphs that satisfy assumption {\bf B.2}. In this sub-section we  generalize the assumption. We first show that uniform convergence is equivalent to assumption {\bf B.2}. We begin with a definition, followed by the result.
\begin{defn}
\label{main_def}
Any  property is said to hold a.s.  on a set $A$ if there exist a set $B$ such that  $B \subset A $  and $P(B)=P(A)$. 
 \end{defn}
\begin{lemma}
\label{Lemma_uniformconvergence}
Define $ A^{n}_j :=  \sum_{i \in \mathcal{G}_1\cup \mathcal{G}_2} I_{j,i}$ for any  $j$. Also,  define the following set:

\vspace{-4mm}
{\small \begin{eqnarray*}
{\cal D} :=  \Big \lbrace  w :   \frac{A^{n}_j}{n} \stackrel{ n \to \infty}{\to}  \gp{1} ~  \mbox{uniformly in  j} \in \mathcal{G}_1 ~  \mbox{and } \frac{A^{n}_k}{n} \stackrel{ n \to \infty}{\to}
\gp{2} ~\mbox {uniformly in k} \in\mathcal{G}_2 \Big \rbrace . \end{eqnarray*}}Then (a.s.) convergence on set ${\cal D}$ is equivalent to the (a.s.) convergence on the   set ${\cal E}$ defined in assumption {\bf B.2}. 
\end{lemma}
\textbf{Proof:}  available in Appendix C.
\eop

It is possible that the graphs may not satisfy uniform convergence with probability one as defined in set  ${\cal D}$. However  we will show  that   the results of Theorem \ref{Thm_MainGen 1} are true almost surely on set ${\cal D}$ or ${\cal  E}$ as given below:
\begin{thm}\label{Thm_onsetD}
Consider a scenario satisfying the assumptions of Theorem \ref{Thm_MainGen 1} and Theorem \ref{Lemma_limit_system_uniquness}, except for assumption {\bf B.2}. Then the respective  conclusions of the   theorems  hold almost surely on the set    ${\cal E}$. 
\end{thm}    
\textbf{Proof:}  available in Appendix C.
\eop
\ignore{
 {\color{red}
 \subsection{An example of multiple fixed points in pre-limit with unique fixed point in  the limit}
Consider a ring network with $n$  nodes and set $\gamma =1$. The first node is liable to the second node, and the second node is liable to third and so on. The $n^{th}$ node is liable to   node $1$. Consider the sample paths in which all the nodes are connected only to the small nodes even when $p^{sb}_1 > 0$. Basically, in these sample paths, the first $n$ nodes are only connected from other small nodes. For any value of $n$, the set of these sample paths has a positive probability.

Let  us consider the  fixed point equation for the given ring network structure, after realization of the random shocks  and for the sample paths as mentioned above. The fixed point equation for any bank $i \in \lbrace 1,2, \cdots , n-1 \rbrace$ is:
\begin{eqnarray*}
 X_i = \min \bigg \lbrace C_i + X_{i-1}, \bar{y}_1 \bigg \rbrace ~ \mbox{a.s.} 
\end{eqnarray*}
 In the above equation the $C_i$'s are obtained by the realization of the shocks and  after paying the tax. But  for the $n^{th}$ nodes  fixed point equation is:
 \begin{eqnarray*}
X_n = \min \bigg \lbrace C_n + X_{1}, \bar{y}_1 \bigg \rbrace ~\mbox{a.s.}  
\end{eqnarray*}
We are computing the fixed point solution under large shock regime i.e.,  if the shock value is  large enough and the nodes are unable to pay the liability amount i.e.,  all nodes are paying:  $X_i = C_i +  X_{i-1}$. Adding all the fixed point solutions we have,  $\sum_{i=1}^{n} X_i = \sum_{i=1}^{n}( C_i + X_{i-1})$ which leads to  $\sum_{i=1}^{n} C_i = 0$. Thus if we choose $C_n =   -\sum_{i=1}^{n-1} C_i$,  one can obtain multiple solutions for ${X_i}$, $i{\le n}$, for any finite value  of $n$. The above sample path  can also occur with positive probability for a sufficiently large value of $n$. However this probability reduces as $n$ increases, beyond certain $n$ some of the small  nodes   lend to big node. Now in the limit we will have a  unique solution for almost all paths as $0 < p^{sb}_1 < 1$ . \\}}

Before proceeding further, we consider   a simple example to illustrate the idea of random fixed point equations and the relevant coupling arguments that are crucial in  comparing the networks of different sizes   sample-path-wise. 
\subsection*{An Example of Random fixed point equations}
\label{Section_examplesoffp}
To illustrate the idea of random fixed points (FP), we consider an example finance network  first with $3$ nodes, and then with a fourth node attached to the existing connections. Basically we consider a realization of $\{I_{i,j}\}_{i, j \le 3}$, $\{\eta_i^{sb}\}_{i \le 3}$  and a realization of the shocks $\{V_i\}_{i \le 3}$ and write the FP  for the corresponding clearing vector (more details are in Section $5$). We then consider the fourth node, extend the realization of the connections $\{I_{i,j}\}_{i = 4 \mbox{ or } j =4}$, $\eta^{sb}_4$ and the shocks $V_4$.

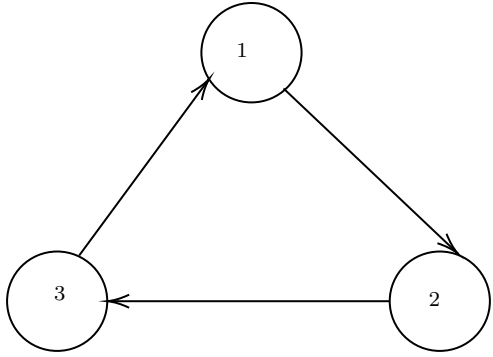
\begin{figure}[htbp!]
\begin{center}
\tikzset{every picture/.style={line width=0.75pt}} %set default line width to 0.75pt        

\begin{tikzpicture}[x=0.75pt,y=0.75pt,yscale=-1,xscale=1]
%uncomment if require: \path (0,390); %set diagram left start at 0, and has height of 390

%Shape: Circle [id:dp5541205930686315] 
\draw   (290,63) .. controls (290,49.19) and (301.19,38) .. (315,38) .. controls (328.81,38) and (340,49.19) .. (340,63) .. controls (340,76.81) and (328.81,88) .. (315,88) .. controls (301.19,88) and (290,76.81) .. (290,63) -- cycle ;
%Shape: Circle [id:dp6801867297656192] 
\draw   (193,188) .. controls (193,174.19) and (204.19,163) .. (218,163) .. controls (231.81,163) and (243,174.19) .. (243,188) .. controls (243,201.81) and (231.81,213) .. (218,213) .. controls (204.19,213) and (193,201.81) .. (193,188) -- cycle ;
%Shape: Circle [id:dp37910788093146985] 
\draw   (384,188) .. controls (384,174.19) and (395.19,163) .. (409,163) .. controls (422.81,163) and (434,174.19) .. (434,188) .. controls (434,201.81) and (422.81,213) .. (409,213) .. controls (395.19,213) and (384,201.81) .. (384,188) -- cycle ;
%Straight Lines [id:da2735383333098629] 
\draw    (331,81) -- (416.55,162.62) ;
\draw [shift={(418,164)}, rotate = 223.65] [color={rgb, 255:red, 0; green, 0; blue, 0 }  ][line width=0.75]    (10.93,-3.29) .. controls (6.95,-1.4) and (3.31,-0.3) .. (0,0) .. controls (3.31,0.3) and (6.95,1.4) .. (10.93,3.29)   ;
%Straight Lines [id:da5325541390997666] 
\draw    (229,165) -- (292.82,77.62) ;
\draw [shift={(294,76)}, rotate = 126.14] [color={rgb, 255:red, 0; green, 0; blue, 0 }  ][line width=0.75]    (10.93,-3.29) .. controls (6.95,-1.4) and (3.31,-0.3) .. (0,0) .. controls (3.31,0.3) and (6.95,1.4) .. (10.93,3.29)   ;
%Straight Lines [id:da6317581713258629] 
\draw    (384,188) -- (245,188) ;
\draw [shift={(243,188)}, rotate = 360] [color={rgb, 255:red, 0; green, 0; blue, 0 }  ][line width=0.75]    (10.93,-3.29) .. controls (6.95,-1.4) and (3.31,-0.3) .. (0,0) .. controls (3.31,0.3) and (6.95,1.4) .. (10.93,3.29)   ;

% Text Node
\draw (306,57) node [anchor=north west][inner sep=0.75pt]   [align=left] {1};
% Text Node
\draw (402,182) node [anchor=north west][inner sep=0.75pt]   [align=left] {2};
% Text Node
\draw (215,179) node [anchor=north west][inner sep=0.75pt]   [align=left] {3};

\end{tikzpicture}
\end{center}
\caption{Ring Graph
\label{Fig_ringgraph}}
\end{figure}

\begin{example}[Multiple fixed points with $n=3$ nodes]
\label{Example_multiple_fixed_point}
In the financial network, agent $i$ borrows from agent $j$ if $I_{i,j}=1$. It can also borrow $y \eta_i^{sb} $ from a big bank (BB), for some $y >0$;   the amount borrowed by $i$ from each of its lenders equals, $y (1-\eta^{sb}_i)/ \sum_{j} I_{i,j}$. Further each agent lends to others in a similar way. Further more, each agent invests the remaining amount in outside risky ventures; thus in all, agent $i$ invests the following amount  in risky ventures (where $k_0$ is the initial wealth)
\begin{equation*}
\Omega_i = \bigg( k_0 + y - \sum_{j} I_{j,i} \frac{y}{\sum_{i'} I_{j,i'}} \bigg).
\end{equation*}
The agent  receives a shock $V_i$ which is either upward $V_i = u$  or downward $V_i = d$, in other words the returns equal $K_i = \Omega_i (1+V_i)$. Each agent has to clear its liability using these returns as well as the claims from the other agents.  

We now discuss a realization of a network with three agents which resulted in a ring liability graph as in Figure \ref{Fig_ringgraph}:  
\begin{eqnarray}
\label{Eqn_example_1}
I_{1,2} = I_{2,3} = I_{3,1} = 1, \ V_1 = d, V_2 = d \mbox{ and,   } V_3 = u,  \eta_1^{sb} = \eta_2^{sb} = \eta_3^{sb} = 0.
\end{eqnarray}
The rest of the $I_{i,j} = 0$ with $i \in \{1,2,3\}$ and $j \in \{1,2,3\}$.

% We consider the ring liability graph as in Figure \ref{Fig_ringgraph} and the corresponding  liability entries are $L_{12}= L_{23} = L_{31} = y$. The other liability entries are set to be zero.

The fixed point (clearing vector in case of financial network) for the three nodes are governed  by the following:
  \begin{eqnarray}
      X_1 &=& \min\left \{  (K_1 + X_3)^+, y \right \}, X_2 =  \min\left \{ (K_2 + X_1)^+, y \right \}, \nonumber \\
      X_3 &=&  \min\left \{ (K_3 + X_2)^+, y \right \}.
  \end{eqnarray}
 
 %  At the time zero the initial wealth of each agent is $k_0$ and every agent borrow and lend equal amount. 
%  
%  The agent receives return from the outside investment  either with upward return $(1+u)$ or with the downward return $1+d$ with $u>d$.  
 
%  The random returns are $a= -k_0(1+d-\kappa)$ and $b= k_0(1+u-\kappa)$ with respective probabilities $w$ and $1-w$, where $\kappa> 0$ be the scaled deposit factor; we also let $b = 2a$. 
Note that $K_i$ is the random return of the agent $i$. Say two of them get shocks, i.e., $K_1 = K_2 = -a$ and $ K_3 = b$ with $a= -k_0(1+d)$ and $b= k_0(1+u)$,  as in  equation \eqref{Eqn_example_1}; we also let $b = 2a$. Then there are multiple fixed point (FP) solutions $(X_1, X_2, X_3) = (a+\epsilon, \epsilon, b+ \epsilon)= (a+\epsilon, \epsilon, 2a+ \epsilon)$ for any $ 0 \le \epsilon < y-b$.
\end{example}   

 %{\color{red}
 % \newcommand{\C}{C}
%   With three nodes in ring graph and consider the liability  with  $L_{12}= L_{23} = L_{31} = y$ 

\begin{example}[Unique fixed point with $n=4$ nodes]
\label{Example_unique_fixed_point}
%   When fourth node is added, the node say borrows from node 1 and lends to node 3.  

We now extend the previous example to include a fourth node, where  the previous quantities are still applicable, with   additional  new details as below. The new connections are as below along with old connections as in \eqref{Eqn_example_1}:
$$
I_{4,1}=1, \eta^{sb}_4 = 1/2,  V_4 = u \mbox{ and }
$$the rest are zero. 

This results in the  liability graph as in  Figure \ref{Fig_4nodesgraph}, with four nodes. The fixed point equations corresponding to  Figure \ref{Fig_4nodesgraph} are as below:

%and with liability structure be $L_{12}= L_{23} = y$ and  $L_{31} = L_{34} =  y/2 = L_{41}$ while $\eta_4^{sb} = y/2$, where $L_{ij}$ represents liability of node $i$ towards node $j$.

\begin{figure}[htbp!]
\begin{center}
\tikzset{every picture/.style={line width=0.75pt}} %set default line width to 0.75pt        

\begin{tikzpicture}[x=0.75pt,y=0.75pt,yscale=-1,xscale=1]

%Shape: Circle [id:dp143786857155225] 
\draw   (144,98) .. controls (144,84.19) and (155.19,73) .. (169,73) .. controls (182.81,73) and (194,84.19) .. (194,98) .. controls (194,111.81) and (182.81,123) .. (169,123) .. controls (155.19,123) and (144,111.81) .. (144,98) -- cycle ;
%Shape: Circle [id:dp3254120518147159] 
\draw   (147,200) .. controls (147,186.19) and (158.19,175) .. (172,175) .. controls (185.81,175) and (197,186.19) .. (197,200) .. controls (197,213.81) and (185.81,225) .. (172,225) .. controls (158.19,225) and (147,213.81) .. (147,200) -- cycle ;
%Shape: Circle [id:dp922123623810787] 
\draw   (316,203) .. controls (316,189.19) and (327.19,178) .. (341,178) .. controls (354.81,178) and (366,189.19) .. (366,203) .. controls (366,216.81) and (354.81,228) .. (341,228) .. controls (327.19,228) and (316,216.81) .. (316,203) -- cycle ;
%Straight Lines [id:da668330917268781] 
\draw    (197,200) -- (315.44,104.26) ;
\draw [shift={(317,103)}, rotate = 141.05] [color={rgb, 255:red, 0; green, 0; blue, 0 }  ][line width=0.75]    (10.93,-3.29) .. controls (6.95,-1.4) and (3.31,-0.3) .. (0,0) .. controls (3.31,0.3) and (6.95,1.4) .. (10.93,3.29)   ;
%Straight Lines [id:da6792164942725476] 
\draw    (194,98) -- (313.5,95.05) ;
\draw [shift={(315.5,95)}, rotate = 178.59] [color={rgb, 255:red, 0; green, 0; blue, 0 }  ][line width=0.75]    (10.93,-3.29) .. controls (6.95,-1.4) and (3.31,-0.3) .. (0,0) .. controls (3.31,0.3) and (6.95,1.4) .. (10.93,3.29)   ;
%Straight Lines [id:da1894297853739797] 
\draw    (340.5,120) -- (340.98,176) ;
\draw [shift={(341,178)}, rotate = 269.51] [color={rgb, 255:red, 0; green, 0; blue, 0 }  ][line width=0.75]    (10.93,-3.29) .. controls (6.95,-1.4) and (3.31,-0.3) .. (0,0) .. controls (3.31,0.3) and (6.95,1.4) .. (10.93,3.29)   ;
%Shape: Circle [id:dp7260871117771865] 
\draw   (315.5,95) .. controls (315.5,81.19) and (326.69,70) .. (340.5,70) .. controls (354.31,70) and (365.5,81.19) .. (365.5,95) .. controls (365.5,108.81) and (354.31,120) .. (340.5,120) .. controls (326.69,120) and (315.5,108.81) .. (315.5,95) -- cycle ;
%Straight Lines [id:da2715781355981217] 
\draw    (172,175) -- (169.12,125) ;
\draw [shift={(169,123)}, rotate = 86.7] [color={rgb, 255:red, 0; green, 0; blue, 0 }  ][line width=0.75]    (10.93,-3.29) .. controls (6.95,-1.4) and (3.31,-0.3) .. (0,0) .. controls (3.31,0.3) and (6.95,1.4) .. (10.93,3.29)   ;
%Straight Lines [id:da20082441483363778] 
\draw    (316,203) -- (199.5,204.97) ;
\draw [shift={(197.5,205)}, rotate = 359.03] [color={rgb, 255:red, 0; green, 0; blue, 0 }  ][line width=0.75]    (10.93,-3.29) .. controls (6.95,-1.4) and (3.31,-0.3) .. (0,0) .. controls (3.31,0.3) and (6.95,1.4) .. (10.93,3.29)   ;

% Text Node
\draw (330,87) node [anchor=north west][inner sep=0.75pt]   [align=left] {1};
% Text Node
\draw (333,194) node [anchor=north west][inner sep=0.75pt]   [align=left] {2};
% Text Node
\draw (166,195) node [anchor=north west][inner sep=0.75pt]   [align=left] {3};
% Text Node
\draw (162,93) node [anchor=north west][inner sep=0.75pt]   [align=left] {4};
\end{tikzpicture}
\end{center}
\caption{ Graph with 4 nodes
\label{Fig_4nodesgraph}}
\end{figure}
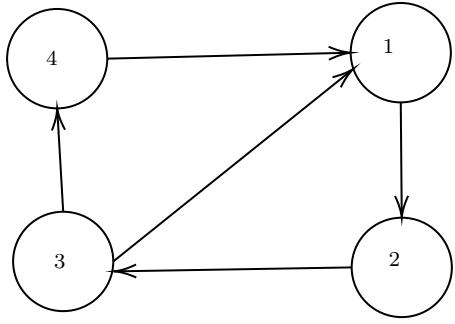
      
\begin{eqnarray}
\label{Eqn_example_with4nodes}
X_1 &=& \min\left \{  (K_1 + (X_3 + X_4)/2)^+,y \right \}, \nonumber \\
X_2 &=&  \min\left \{ (K_2 + X_1)^+, y \right \}, \nonumber \\
X_3 &= & \min\left \{ (K_3 + X_2)^+,y \right \}, \nonumber \\
X_4 &=& \min\left \{ ( K_4 + X_3/2)^+,y \right \}.
 \end{eqnarray}
%   where ${\tilde C} = (K-v)\frac{ w+y/2}{w}.$ and $C = K-v $.

% Once  again consider the stressed scenario i.e.,  when three of them receives the shock and unable to clear the full liability. Adding all the above we have  the following:
% \begin{eqnarray}
% \sum_{i=1}^{4} X_i &=&  \sum_{i=1}^{4} C_i + \sum_{i=1}^{3} X_i +\frac{X_4}{2} \nonumber \\
% \implies X_4 &=& 2\sum_{i=1}^{4} C_i.
% \end{eqnarray}
% From \eqref{Eqn_example_with4nodes} we can compute the fixed point for the other entities.  
% \begin{eqnarray}
% X_3&=& 2(X_4-C_4), %= 4\sum_{i=1}^{4} C_i -2C_4,
% \nonumber\\
% X_1&=& C_1+ \frac{X_3+X_4}{2} %= C_1-C_4+3 \sum_{i=1}^{4} C_i
% , \nonumber\\
% X_2&=& C_2+X_1.
% %=  C_2 +C_1-C_4+3 \sum_{i=1}^{4} C_i.
% \end{eqnarray}
Note that $K_i$ are the  random returns as in the Example \ref{Example_multiple_fixed_point}. Say  node $1$ and $2$ receives the shock and the realization of the returns are  $K_1= K_2=-a, K_3=b, K_4=b+ \frac{by}{2k_0}$, with $a< b< y$   then the unique  fixed point solution becomes:
 $$
\bigg(X_1, X_2, X_3, X_4\bigg)=  \bigg (y-a, y-2a, y,y\bigg).
$$
In the above example the defaulted nodes are $\{1,2\}$ (more details in  Section \ref{sec_finance}). 

Thus by adding an extra node we have unique fixed point. A close observation  at the two networks indicates that this is possible because the fourth  (new)  node is connected to BB; this resulted in a contraction mapping. By law of large numbers, as $n$ increases, we will have networks with sufficient connections to BB, the resultant of which will be a contraction mapping (and then the existence of unique fixed point). 
\end{example}

\section{Another Graphical Model}
\label{sec_alternate}
\label{Alternate_graphical_model}
In the previous section, we discussed a random graphical model with   large number of interacting nodes. In this  section, we extend the methodology to an alternate model. One can extend the results to many such variants in a similar way. 

Both   the graphical models,  are important on their own. The first model is relevant when  say the resources are shared  equally across all the connected neighbours of the entire network (see \eqref{Eqn_Weights 1}, where the weights are divided by $ \sum_{i'  \in \mathcal{G}_1 \cup \mathcal{G}_2} I_{j, i'} $). In the first  model, the resources are shared  equally across the entire network. In contrast,  in this section,  the resources are shared equally only within the group  after allocating dedicated fractions  ($\{\lambda_m\}$) to each group.
We again have  convergence of the random fixed point equations in almost sure sense (like in Theorems \ref{Thm_MainGen 1}-\ref{Thm_onsetD} and Corollary \ref{Cor_Three dimensional approximation}), as the number of nodes increases to infinity.  The graphical model of this section is used to study     the financial network-based application in the next section. 

As before, we  have two groups and a big $b$-node.  
We  begin  by providing the details of   the random weights between various nodes    in the following:
\begin{eqnarray}
\label{Alt_Eqn_Weights 1}
W_{j, b} &=&  \eta_j^{sb}  \lambda_m \mbox{ (towards $b$-node),} \mbox{ with } p_m^{sb} \ :=  \ E[\eta_j^{sb}]  \mbox{ for any } j \in {\cal G}_m, \mbox{ and, }  \ \\ 
W_{j,i} &=& \left \{ \begin{array}{lll}   \frac{I_{j,i} (1-\eta_j^{sb}) \lambda_m } {\displaystyle \sum_{i'  \in \mathcal{G}_m \ } I_{j, i'}  } ,  & \mbox{ if } i \in {\cal G}_m \\ \\
  \frac{I_{j,i}  (1-\lambda_m) } {\displaystyle \sum_{i'  \notin \mathcal{G}_m \ } I_{j, i'}  } 1_{p_{c_m} > 0} ,  & \mbox{ if } i \notin {\cal G}_m \\
  \end{array} \right .  
    \   \   \ 
\mbox{ (towards small node $i,$)} 
\end{eqnarray}
where $\lambda_m$ (with $0 \le \lambda_m \le 1$) is a non-negative fraction. Consider any small node $j\in \mathcal{G}_m$, $\lambda_m$ fraction is dedicated towards the nodes of its own group, while, the remaining $(1-\lambda_m)$  fraction is towards the other group.  Further $\lambda_m \eta_j^{sb}$   fraction   is  towards the big node $b$, while the remaining fraction is equally shared within the group ${\cal G}_m$ (among the interested members, interests represented by $\{I_{j,i}\}$ flags).   All the edge formulation events are independent of one other (or   satisfy an assumption like {\bf B.2(C)}) and the corresponding probabilities are as in the previous section. Observe 
the weights  sum up to one, i.e., $\sum_{i \in  \mathcal{G}_m} W_{j,i} + \sum_{i \notin  \mathcal{G}_m} W_{j,i} + W_{j, b} = 1$ for all $j .$  Further recall  that 
the weights  from $b$-node are  given by  $\{\eta_j^{bs}\}.$
%: 
%$$
%\hspace{19mm}
%W_{b, j}  = \frac{\eta_j^{bs} }{\displaystyle \sum_i  \eta_i^{bs} } \mbox{ with } E[\eta_j^{bs}] = p_m^{bs}  \mbox{ (when $j\in \mathcal{G}_m $}).
%$$
As before we are interested in the performance of the nodes which depends  upon the performance of the other nodes  via a set of fixed point equations (which are of different structure to that considered in previous section):
\begin{eqnarray}
\label{Alt_Eqn_FixedeqGen 1}
X_i^{m}  &=& f^{m}(G^m_i,   {\bar X}^{m1}_i ,   {\bar X}^{m2}_i ,  \eta_i^{bs} X^b) \mbox{ {\normalsize  for each }} i \in \mathcal{G}_m \mbox{ and any } m,  \ \mbox{and},  \\  
\label{Alt_Eqn_Fixedeq2Gen 3}
X^b  &=&    f^b (   {\bar X}^b ), \ \ \ \mbox{ {\normalsize with aggregates,}}    \\
{\bar X}^{m1}_i  & := &  \displaystyle \sum_{j \in \mathcal{G}_1} X_j^1 W_{j, i} , \  \ \ \
{\bar X}^{m2}_i \  :=  \  \displaystyle \sum_{j \in \mathcal{G}_2} X_j^2  {W}_{j, i} \mbox{ {\normalsize  for each }} i \in \mathcal{G}_m, 
\label{Alt_Eqn_aggregatesn} \   \     \\
%\tilde{X}^s_i  & := &  \sum_{j \in \mathcal{G}_1} X_j^s W_{j, i} + \sum_{j \in \mathcal{G}_2} X_j^s \tilde{W}_{j, i}\mbox{ {\normalsize  for each }} i \in \mathcal{G}_2, \   \  \nonumber  
{\bar X}^b & := & \frac{1}{n}  \sum_{j 	\in \mathcal{G}_1} {X}^1_j W_{j, b} + \frac{1}{n}  \sum_{j \in \mathcal{G}_2} {X}^2_j W_{j, b} . \label{Alt_Eqn_Fixedeq2Gen 4}
\end{eqnarray}

% Note: In financial set up we have $\lambda_1 = 1$, $p_1^{sb}=0$,$p_2^{sb}> 0$ and $\sigma=1$. Therefore, the Lemma \ref{Lemma_limit_system_uniquness} condition does not  applicable.

We are interested in solving these random fixed point equations asymptotically, and   we begin with aggregate
fixed points.   {\it We would like to mention here  that the proofs for this section follow  exactly in a similar way as in the previous section and we are only stating the differences in the relevant expressions and assumptions.} % here. 

 \subsection{Aggregate fixed points}
 We rewrite the  fixed point equations in terms of weighted averages and first analyze the aggregate system   as before.
Define the following random variables, that depend upon real constants $({\bar x}^{m1}_i, {\bar x}^{m2}_i, {\bar x}_b), \mbox{ (with } x_b:= f^b (  {\bar x}_b))$ : 
\begin{eqnarray}\hspace{15mm} \label{Alt_Eqn_xi 1}
\xi^m_i ({\bar x}^{m1}_i, {\bar x}^{m2}_i, x_b) :=   f^{ m} (	G_i^m,   {\bar x} _i^{m1},  {\bar x} _i^{m2},   \eta_{i}^{bs} x_b)  \mbox{ for any }  i \in {\cal G}_m.
\end{eqnarray} 
 Let ${\x}^{m} := ({\bar x}^{m1}_1, {\bar x}^{m2}_1, {\bar x}^{m1}_2, {\bar x}^{m2}_2,  \cdots )$,
%$ \bar{x}^{m2} := ({\bar x}^{m2}_1, {\bar x}^{m2}_2 \cdots )$
 for  $m  \in \{ 1,2\}$  and  let
 ${\x} := ({\x}^{1} , {\x}^{2})$.
 Consider  the  following  operators  on infinite  sequence space
 $s^\infty$, one for each $ n= n_1+ n_2$:

 \vspace{-2mm}
{\small
\begin{eqnarray}
\label{Alt_Eqn_bar_fixedpoint_randGen_1}
{\bar {\bf f}}^n ({\bar x}_b,{\x} )  =  ({\bar f}^n_b, \ \  {\bar f}_1^{n,11},
{\bar f}_1^{n,12},
%{\bar f}_2^{n,11}, {\bar f}_2^{n,12} %
\cdots, 
{\bar f}_{n_1}^{n,11}, {\bar f}_{n_1}^{n,12}, \cdots,  \ \ {\bar f}_{1}^{n,21}, {\bar f}_{1}^{n,22}, \cdots,  {\bar f}_{n_2}^{n,21}, {\bar f}_{n_2}^{n,22}, \cdots),     \hspace{4mm}
\end{eqnarray}}where, for any $ m$,

\vspace{-1mm}
{\small 
\begin{eqnarray}
\label{Eqn_altG1_aggregate}
{\bar f}^{n, m 1}_i ({\bar x}_b,\x ) &=& \left \{  \begin{array}{lll} \sum_{j \in\mathcal{G}_1}   \xi^1_{j} ({\bar x}^{11}_j, {\bar x}^{12}_j, x_b)  W_{j, i}, &\mbox{ \normalsize if }     i \in \mathcal{G}_m  \\
 0  &\mbox{\normalsize else,}    \hspace{20mm} 
\end{array} \right .   \\ \nonumber
\\
{\bar f}^{n, m 2}_i ({\bar x}_b,\x ) & = 
& \left \{  \begin{array}{lll} \sum_{j \in\mathcal{G}_2}   \xi^2_{j} ({\bar x}^{21}_j, {\bar x}^{22}_j, x_b)  W_{j, i} &\mbox{ \normalsize if }   i \in \mathcal{G}_m  \\
 0  &\mbox{\normalsize else, \hspace{10mm}  and, }   \hspace{5mm}
 \end{array} \right. 
 \label{Eqn_alt_aggregate}  \\
 {\bar f}^n_b ({\bar x}_b,{\x})& := & 
 \frac{1}{n}
\displaystyle \sum_{j \in \mathcal{G}_1}   \xi^{1}_{j} ({\bar x}^{11}_j,{\bar x}^{12}_j, x_b) W_{j, b} +  \frac{1}{n}
\displaystyle \sum_{j \in \mathcal{G}_2}   \xi^{2}_{j} ({\bar x}^{21}_j, {\bar x}^{22}_j, x_b) W_{j, b} . 
 \label{Alt_Eqn_key} 
 \end{eqnarray}}
\ignore{
{\color{red}
We get from (27) with $m = 1$
$$
\sup_{j \in {\cal G}_1} \left ( |{\bar x}^{m1}_j - {\bar u}^{m1} | + |{\bar x}^{m2}_j - {\bar u}^{m2} |
\right ) $$
And from (28) we get with same $m$
$$
\sup_{j \in {\cal G}_2} \left ( |{\bar x}^{m1}_j - {\bar u}^{m1} | + |{\bar x}^{m2}_j - {\bar u}^{m2} |
\right ) $$
By adding the appropriate terms we get for $i \in {\cal G}_{m'}$
\begin{eqnarray*}
|{\bar f}^{n, m' 2}_i ({\bar x}_b,\x ) - {\bar f}^{n, m' 2}_i ({\bar u}_b,\u ) |
+ |{\bar f}^{n, m' 1}_i ({\bar x}_b,\x ) - {\bar f}^{n, m' 1}_i ({\bar u}_b,\u )| 
+ | {\bar f}^n_b ({\bar x}_b,{\x}) - {\bar f}^n_b ({\bar u}_b,{\u}) |
\hspace{-80mm} \\
&\le& \sup_{j \in {\cal G}_1, m = 1} \left ( |{\bar x}^{m1}_j - {\bar u}^{m1} | + |{\bar x}^{m2}_j - {\bar u}^{m2} |
\right )  \sum_{j \in\mathcal{G}_1}    W_{j, i} 
+
\sup_{j \in {\cal G}_2, m = 2} \left ( |{\bar x}^{m1}_j - {\bar u}^{m1} | + |{\bar x}^{m2}_j - {\bar u}^{m2} |
\right )
\sum_{j \in\mathcal{G}_2}    W_{j, i}  \\
&\le & 
\max_{m \in \{1, 2\} }
\sup_{j \in {\cal G}_m } \left ( |{\bar x}^{m1}_j - {\bar u}^{m1} | + |{\bar x}^{m2}_j - {\bar u}^{m2} | 
\right )  
\left ( \sum_{j \in\mathcal{G}_1}    W_{j, i} + \sum_{j \in\mathcal{G}_2}    W_{j, i} \right )
\end{eqnarray*}

\begin{eqnarray*}
| {\bar f}^n_b ({\bar x}_b,{\x}) - {\bar f}^n_b ({\bar u}_b,{\u}) |
&\le & 
\max_{m \in \{1, 2\} }
\sup_{j \in {\cal G}_m } \left ( |{\bar x}^{m1}_j - {\bar u}^{m1} | + |{\bar x}^{m2}_j - {\bar u}^{m2} |  +|\bar{x}_b -\bar{u}_b |\right )  
\frac{1}{n}\left ( \sum_{j \in\mathcal{G}_1}    W_{j, b} + \sum_{j \in\mathcal{G}_2}    W_{j, b} \right )
\end{eqnarray*}

}

{\color{red} For alternate system --
 \begin{eqnarray*} 
     || ({\bar x}_b ,\x)- ({\bar u}_b, \u)  ||_\infty 
      &=&  \max_{m,m' \in \lbrace 1, 2\rbrace} \sup_{i\in \mathcal{G}_m } \left (|\bar{x}^{m,m'}_i - \bar{u}^{m,m'}_i| + | {\bar x}_b -   {\bar u}_b |\right) \\
       &=&  \max_{m \in \lbrace 1, 2\rbrace} \sup_{i\in \mathcal{G}_m } \left (|\bar{x}^{m1}_i - \bar{u}^{m1}_i|+ |\bar{x}^{m2}_i - \bar{u}^{m2}_i| + | {\bar x}_b -   {\bar u}_b |\right)
\end{eqnarray*}

}
}
As before, we define the limit operators using limits for any $m, m'$ and $i \in {\cal G}_m$ 
\begin{eqnarray}
\label{Eqn_lim_sup_alt}
{\bar f}^{\infty, m m'}_i := \limsup_{n \to \infty} {\bar f}^{n, m m'}_i
\mbox{ and }  {\bar f}^\infty_b ({\bar x}_b,{\x}) :=  \limsup_{n \to \infty} {\bar f}^n_b ({\bar x}_b,{\x}).
 \end{eqnarray}
Now with this  description  we are ready to define the  aggregate convergence of the random fixed points. Towards this, exactly as in Lemma \ref{lem: LLN_fixed point}, under constant sequences,  the above limit operator is given by (with $\x = ({\x}^{1} , {\x}^{2})$):
 
 \vspace{-2mm}
 {\small\begin{eqnarray}
 \label{Eqn_aggregatefunc_altgraph}
 {\bar f}_i^{\infty,11}  ({\bar x}_b,\x )  &=&      E_{G^{1}_i, \eta_i^{bs}} \left [ \xi^{1}_i ({\bar x}^{11},{\bar x}^{12}, x_b)   \right ]  \lambda_1 (1-p_1^{sb}), \label{Alt_Eqn_barf_limit 1} \\
 {\bar f}_i^{\infty,12}  ({\bar x}_b,\x )  &=&  E_{G^{2}_i, \eta_i^{bs}} \left [ \xi^{2}_i ({\bar x}^{21},{\bar x}^{22}, x_b)   \right ] \frac{1-\gamma}{\gamma}(1-\lambda_2) 1_{p_{c_2}> 0}, \nonumber \\
 {\bar f}_i^{\infty,21}  ({\bar x}_b,\x)  &=&      E_{G^{1}_i, \eta_i^{bs}} \left [ \xi^{1}_i ({\bar x}^{11},{\bar x}^{12}, x_b)   \right ] (1-\lambda_1)\frac{\gamma}{(1-\gamma)}1_{p_{c_1}> 0}, \nonumber  \\
  {\bar f}_i^{\infty,22}  ({\bar x}_b,\x)  &=&  E_{G^{2}_i, \eta_i^{bs}} \left [ \xi^{2}_i ({\bar x}^{21},{\bar x}^{22}, x_b)   \right ] \lambda_2(1-p_2^{sb}) ,\nonumber  \\
 {\bar f}_b^\infty  ( {\bar x}_b,\x) &=&   \gamma E_{ G^{1}_i, \eta_i^{bs}}  [\xi^{1}_i ({\bar x}^{11},  {\bar x}^{12},x_b)  ]\lambda_1  p_1^{sb}  + (1-\gamma) E_{ {G}^{2}_i, \eta_i^{bs}}  [\xi^{2}_i ({\bar x}^{21}, {\bar x}^{22}, x_b)  ]\lambda_2 p_2^{sb}. \hspace{4mm} \label{Alt_Eqn_barf_limit e}
\end{eqnarray}}
Recall that one such constant sequence will form the (sequence of) aggregate fixed points. 
% Once we have convergence of the aggregate random fixed points,  which in turn would provide us   the convergence of the original  random  fixed points equations (see equation \eqref{Alt_Eqn_FixedeqGen 1} -\eqref{Alt_Eqn_Fixedeq2Gen 4}) as below. 
 We are interested in solving these random fixed point equations asymptotically, and towards that, we make the following modified assumptions.  We would like to mention again that the proofs for this section follow  exactly in a similar way and we are only stating the  modifications here; we begin with the assumptions: 
\\
{\bf B.2$'$}.
The assumption \textbf{B.2} is modified to use the following definition: 

\vspace{-4mm}
{\small\begin{equation}
{\cal E} : =   \cap_m \left\lbrace \omega: \lim_{n \to \infty } \displaystyle  \hspace{-1mm}  \sum_{j \in \mathcal{G}_m}    \left|  
    \frac{1 }{ \displaystyle  \sum_{i \in   \mathcal{G}_m  } I_{j, i}}  - \frac{1}{ np_m }    \right |  =0 , ~   
    %\right .  \hspace{-20mm} \nonumber \\ &
%    \left .
 \displaystyle \lim_{n\to \infty} {\displaystyle \hspace{-1mm} \sum_{j \in \mathcal{G}_m} }   \left|    \frac{1 }{ \displaystyle \sum_{i \notin   \mathcal{G}_m } I_{j, i}}  - \frac{1}{ np_{c_m} }    \right | 1_{p_{c_m}> 0} = 0 
\right\rbrace.  \hspace{4mm}
\label{Alt_Eqn_calE}
\end{equation}}
We will again require equivalent of {\bf B.2(C)} when $\{I_{j,i}\}$ are not i.i.d.\\
\noindent {\bf B.3$'$:} For each  $m \in  \{1,2\}$ and  $i \in \mathcal{G}_m$ assume (for some $\sigma \le 1$, $\varsigma \le 1$): 
\begin{eqnarray*}
|\xi^{m}_i (x^{m1},x^{m2}, x_b)  - \xi^{m}_i (u^{m1}, u^{m2}, u_b) |
\hspace{-13mm}
\\ &\le & \sigma ( |x^{m1}- u^{m1}| + |x^{m2}- u^{m2}|+ \varsigma|x_b- u_b|).
\end{eqnarray*}  
Basically we require
for all  $x^1, x^2, u^1, u^2,  x_b, u_b, g,  \eta$,  ${\bar x}_b$ and ${\bar u}_b$:
\begin{eqnarray*}|f^b({\bar x}_b) - f^b ({\bar u}_b) | &\le &  |{\bar x}_b-{\bar u}_b | \\
f^m (g, x^1, x^2, \eta x_b) - f^m(g, u^1,u^2, \eta u_b)| 
&\le &\sigma ( |x^1- u^1| +|x^2- u^2| + \varsigma|x_b- u_b|) .
\end{eqnarray*}
 
{\bf B.4$'$:} The assumption \textbf{B.4} modified as below, we now assume:
{ 
 \begin{eqnarray}
 \label{assumption_alt_modified_contraction}
 \varrho : &=&\max \bigg\lbrace  \lambda_1(1- p_1^{sb})+ \frac{1-\gamma}{\gamma}(1-\lambda_2) 1_{p_{c_2}> 0}, \frac{\gamma}{1-\gamma}(1-\lambda_1)1_{p_{c_1}> 0}  + \lambda_2 (1-p_2^{sb})\bigg \rbrace \nonumber\\
 && +  \left ( \gamma \lambda_1 p_1^{sb} + (1-\gamma) \lambda_2 p_2^{sb} \right )  \leq  1.
 \end{eqnarray}}
It is easy to observe that $\varrho =1$  for symmetric conditions, i.e.,  $\gamma =0.5$, $\lambda_1 = \lambda_2$, $p_1^{sb} = p_2^{sb}$, $p_{c_1} =p_{c_2}$.
 
The  assumptions \textbf{B.1} and \textbf{B.5}    remain unaltered,  
 except that the quantities are redefined; for example,  limit function $\bar {\bf f}^\infty $ is now given by \eqref{Eqn_lim_sup_alt}.  With these modified assumptions we have:

\begin{thm}
\label{Thm_MainGen 2}
Assume {\bf B.1}, {\bf B.2$'$}-{\bf B.4$'$} and {\bf B.5}.
The aggregates  of the random system (see  \eqref{Alt_Eqn_FixedeqGen 1}-\eqref{Alt_Eqn_Fixedeq2Gen 4}),  which are   FPs of   \eqref{Alt_Eqn_bar_fixedpoint_randGen_1}-\eqref{Alt_Eqn_key}   denoted by 
$$({\bar X}_b,{\X}) (n) := ({\bar X}_b, \{{\bar X}_i^{m1}\}_{i,m}, \{{\bar X}_i^{m2}\}_{i,m}) (n) $$
converge   as $n \to \infty$, along a sub-sequence. That is, there exists $k_n \to \infty$ such that:
\begin{eqnarray}
\label{Alt_Eqn_aggregate_groupwise}
{\bar X}_i^{m1}(k_n)  \to {\bar x}^{\infty,m1}, \
{\bar X}_i^{m2}(k_n)  \to  {\bar x}^{\infty,m2} \
\mbox{ and }  {\bar X}_b(k_n)  \to {\bar x}_b^{\infty} \mbox{   (a.s.), } \ 
\end{eqnarray}
 where $( {\bar x}_b^{\infty},\x^{\infty} )$
 with   $\x^{\infty} := ({\bar x}^{\infty,11},
 {\bar x}^{\infty,12},  {\bar x}^{\infty,11},
 {\bar x}^{\infty,12} , \cdots, \  \   {\bar x}^{\infty,21},
 {\bar x}^{\infty,22}, \cdots )$ is an  FP of
  the limit system given by  \eqref{Eqn_lim_sup_alt}.
Further    (any sequence of) FPs of the  original system \eqref{Alt_Eqn_FixedeqGen 1}-\eqref{Alt_Eqn_Fixedeq2Gen 3} converge along the sub-sequence in almost sure sense: %to the limit :
\begin{eqnarray}
\label{Alt_Eqn_main_groupwise}
X^b (k_n)  &\to&  f^b (  {\bar x}_b^{\infty}) \mbox{ as $n \to \infty$    and } \hspace{3mm}\label{Eqn_Act_Fixed 1} \\
X^m_i (k_n)  &\to& f^m(G_i^m,  {\bar x}^{\infty,m1},  {\bar x}^{\infty,m2},\eta_i^{bs} X_b^{\infty})\ \forall \mbox{i,m}.\nonumber
\end{eqnarray}
\end{thm}
\textbf{Proof:}  available in Appendix D.
\eop

Now we  state the theorem related to the uniqueness of the fixed point analogous of the  Theorem \ref{Lemma_limit_system_uniquness} as below:
\begin{thm}{\bf [Unique Fixed points]}
\label{Thm_modified_forThm_2_Alt_system}
Assume {\bf B.1}, {\bf B.2$'$}-{\bf B.4$'$} and also assume  $\eta_j^{sb} \ge  \underline{\eta} \ge 0$ a.s.,  with $\sigma \left ( 1- \underline{\eta} + \varsigma   \underline{\eta} \right ) < 1$. Then we have unique fixed point of the finite $n$-system \eqref{Alt_Eqn_bar_fixedpoint_randGen_1} for each $n$. We also have unique fixed point  for the limit system \eqref{Eqn_lim_sup_alt}, which is  a constant sequence, i.e., we have 
${\bar x}_i^{\infty,m1} = {\bar x}^{\infty,m1} $, ${\bar x}_i^{\infty,m2} = {\bar x}^{\infty,m2} $  for all $i \in {\cal G}_m$ and for each $m\in\{1,2\}$ in equation \eqref{Alt_Eqn_aggregate_groupwise}. This limit is the unique  fixed point of the five dimensional system given by \eqref{Alt_Eqn_barf_limit 1}-\eqref{Alt_Eqn_barf_limit e}.
\end{thm}
\textbf{Proof:}  available in Appendix D.
\eop

\begin{cor}{\bf [Three dimensional approximation]}
\label{Cor_Three dimensional approximation_Alt_system}
Assume the conditions of Theorem \ref{Thm_modified_forThm_2_Alt_system}.  Then convergence in  \eqref{Alt_Eqn_aggregate_groupwise}-\eqref{Alt_Eqn_main_groupwise} is along the  original sequence, i.e., as $n\to \infty$. Further,  $({\bar x}_b^{\infty}, \x^{\infty})$
 is a constant sequence and is given by: 
\begin{eqnarray}
\label{Eqn_aggrragate_al_system}
{\bar x}^{\infty,12}  ( {\bar x}_b,\x)  &=& {\bar x}^{\infty,22}  ({\bar x}_b,\x ) \mu_1, \mbox{ with, } \mu_1 := \frac{1- \gamma}{\gamma} \frac{1-\lambda_2}{\lambda_2} \frac{1}{(1- p^{sb}_2)}  1_{p_{c_2}> 0},\nonumber
\\
{\bar x}^{\infty,21}  ({\bar x}_b,\x) &=&  {\bar x}^{\infty,11}  ({\bar x}_b,\x) \mu_2, \mbox{ with, } \mu_2 := \frac{\gamma}{1-\gamma} \frac{1-\lambda_1}{\lambda_1} \frac{1}{(1- p^{sb}_1)}1_{p_{c_1}> 0}, \hspace{3mm}
\end{eqnarray}
where   $({\bar x}^b, {\bar x}^{\infty, 11}, {\bar x}^{ \infty, 22 } $)  is the unique fixed point of the   three dimensional system:

\vspace{-4mm}
{\small
\begin{eqnarray}
 \label{Alt_eqn_limit_system}
 {\bar f}_i^{\infty,1}  ({\bar x}_b, {\bar x}^{1}, {\bar x}^{2})  &=&      E_{G^{1}_i, \eta_i^{bs}} \left [ \xi^{1}_i ({\bar x}^{1}, \mu_1 {\bar x}^{2}, x_b)   \right ]  \lambda_1 (1-p_1^{sb}), 
 \  x_b:= f^b (  {\bar x}_b))
 % \mbox{\normalsize with, } \x = ( {\bar x}^{1}, {\bar x}^{2}) 
 \nonumber, 
 \\
 {\bar f}_i^{\infty,2}  ( {\bar x}_b,{\bar x}^{1}, {\bar x}^{2})  &=&  E_{G^{2}_i, \eta_i^{bs}} \left [ \xi^{2}_i (\mu_2 {\bar x}^{1},{\bar x}^{2}, x_b)   \right ] \lambda_2(1-p_2^{sb}) \nonumber, \\ 
 {\bar f}_b^\infty  ( {\bar x}_b, {\bar x}^{1}, {\bar x}^{2}) &=&   \gamma E_{ G^{1}_i, \eta_i^{bs}}  [\xi^{1}_i ({\bar x}^1, x_b)  ] \lambda_1 p_1^{sb}  + (1-\gamma) E_{ {G}^{2}_i, \eta_i^{bs}}  [\xi^{2}_i ({\bar x}^2, x_b)  ] \lambda_2 p_2^{sb}.
\end{eqnarray}}
\end{cor}
\textbf{Proof:}  available in Appendix D.
\eop

\ignore{
\begin{eqnarray*}
{\bar x}^{\infty,12}  ( {\bar x}_b,\x)  &=& {\bar x}^{\infty,22}  ({\bar x}_b,\x ) \mu_1, \mbox{ with, } \mu_1 := \frac{1- \gamma}{\gamma} \frac{1-\lambda_2}{\lambda_2} \frac{1}{(1- p^{sb}_2)}  1_{p_{c_2}> 0},
\\
{\bar x}^{\infty,21}  ({\bar x}_b,\x) &=&  {\bar x}^{\infty,11}  ({\bar x}_b,\x) \mu_2, \mbox{ with, } \mu_2 := \frac{\gamma}{1-\gamma} \frac{1-\lambda_1}{\lambda_1} \frac{1}{(1- p^{sb}_1)}1_{p_{c_1}> 0}.
\end{eqnarray*}
In all,  first  one has to solve the   three dimensional fixed point equation as given below ($\x = ({\bar x}^{1}, {\bar x}^{2}$)): 
{\small
\begin{eqnarray}
 \label{Alt_eqn_limit_system}
 {\bar f}_i^{\infty,1}  ({\bar x}_b,\x)  &=&      E_{G^{1}_i, \eta_i^{bs}} \left [ \xi^{1}_i ({\bar x}^{1}, \mu_1 {\bar x}^{2}, x_b)   \right ]  \lambda_1 (1-p_1^{sb}) \nonumber, 
 \\
 {\bar f}_i^{\infty,2}  ( {\bar x}_b,\x)  &=&  E_{G^{2}_i, \eta_i^{bs}} \left [ \xi^{2}_i (\mu_2 {\bar x}^{1},{\bar x}^{2}, x_b)   \right ] \lambda_2(1-p_2^{sb}) \nonumber, \\ 
 {\bar f}_b^\infty  ( {\bar x}_b,\x) &=&   \gamma E_{ G^{1}_i, \eta_i^{bs}}  [\xi^{1}_i ({\bar x}^1, x_b)  ] \lambda_1 p_1^{sb}  + (1-\gamma) E_{ {G}^{2}_i, \eta_i^{bs}}  [\xi^{2}_i ({\bar x}^2, x_b)  ] \lambda_2 p_2^{sb},
\end{eqnarray}}and then the limits of fixed points of the original system are given by \eqref{Eqn_Act_Fixed 1}.
 \eop \\
 }
 
 \ignore{
 {\color{red}
\begin{cor}
Assume {\bf B.1}, {\bf B.2$'$} to {\bf B.4$'$} and also assume,  either  $\varrho <1$ or  $\sigma < 1$ in {\bf B.3$'$}. Then we have unique fixed point   of the limit system \eqref{Eqn_lim_sup_alt}, which is a constant sequence, i.e., we have 
${\bar x}_i^{\infty,m1} = {\bar x}^{\infty,m1} $ , ${\bar x}_i^{\infty,m2} = {\bar x}^{\infty,m2} $  for all $i \in {\cal G}_m$ in equation \eqref{Alt_Eqn_aggregate_groupwise}. This limit is the fixed point of the three dimensional system given by \eqref{Alt_eqn_limit_system}.
\end{cor}
\textbf{Proof:} To proof the contraction mapping for the limit system we use the following norm:
\begin{eqnarray*} 
     || ({\bar x}_b ,\x)- ({\bar u}_b, \u)  ||_\infty 
       &=&  \max_{m \in \lbrace 1, 2\rbrace} \sup_{i\in \mathcal{G}_m } \left (|\bar{x}^{m1}_i - \bar{u}^{m1}_i|+ |\bar{x}^{m2}_i - \bar{u}^{m2}_i| + | {\bar x}_b -   {\bar u}_b |\right).
\end{eqnarray*}

First observe that, adding the terms  from the equations \eqref{Alt_Eqn_key} - \eqref{Eqn_alt_aggregate} we get for $i \in {\cal G}_{m'}$

\vspace{-2mm}
{\scriptsize
\begin{eqnarray*}
|{\bar f}^{\infty, m' 1}_i ({\bar x}_b,\x ) - {\bar f}^{\infty, m' 1}_i ({\bar u}_b,\u ) |
+ |{\bar f}^{\infty, m' 2}_i ({\bar x}_b,\x ) - {\bar f}^{\infty, m' 2}_i ({\bar u}_b,\u )|  + | {\bar f}^\infty_b ({\bar x}_b,{\x}) - {\bar f}^\infty_b ({\bar u}_b,{\u}) | \hspace{-120mm} & & \\
& =& |\limsup_n({\bar {\bf f}_i}^{n,m'2} ( {\bar x}_b, \x ) - {\bar {\bf f}}_i^{n,m'2} ( {\bar u}_b, \u ) ) | + |\limsup_n({\bar {\bf f}_i}^{n,m'2} ( {\bar x}_b, \x ) - {\bar {\bf f}}_i^{n,m'2} ( {\bar u}_b, \u ) ) | \\
& & +|\limsup_n ( {\bar f}^n_b ({\bar x}_b,{\x}) - {\bar f}^n_b ({\bar u}_b,{\u})| \\
&\le& \sigma \sup_{j \in {\cal G}_1, m = 1} \left ( |{\bar x}^{m1}_j - {\bar u}^{m1} | + |{\bar x}^{m2}_j - {\bar u}^{m2} | +|\bar{x}_b-\bar{u}_b|
\right )  \limsup _n \sum_{j \in\mathcal{G}_1}    W_{j, i} 
\\
& & + \sigma \sup_{j \in {\cal G}_2, m = 2} \left ( |{\bar x}^{m1}_j - {\bar u}^{m1} | + |{\bar x}^{m2}_j - {\bar u}^{m2} | + |\bar{x}_b-\bar{u}_b|
\right )
\limsup_n \sum_{j \in\mathcal{G}_2}    W_{j, i}  \\
& & + \sigma \max_{m \in \{1,2\}}\sup_{j \in {\cal G}_m } \left ( |{\bar x}^{m1}_j - {\bar u}^{m1} | + |{\bar x}^{m2}_j - {\bar u}^{m2} |  +|\bar{x}_b -\bar{u}_b |\right )  
\limsup_n\frac{1}{n}\left ( \sum_{j \in\mathcal{G}_1}    W_{j, b} + \sum_{j \in\mathcal{G}_2}    W_{j, b} \right )  \\
&\le & 
 \sigma || ({\bar x}_b ,\x)- ({\bar u}_b, \u)  ||_\infty \lim_n \left ( \sum_{j \in\mathcal{G}_1}    W_{j, i} + \sum_{j \in\mathcal{G}_2}    W_{j, i}  +  \frac{1}{n}\sum_{j \in\mathcal{G}_1}    W_{j, b} + \frac{1}{n} \sum_{j \in\mathcal{G}_2}    W_{j, b}\right ) \\
&\le & 
\sigma\varrho|| ({\bar x}_b ,\x)- ({\bar u}_b, \u)  ||_\infty.  
\end{eqnarray*}}
The first inequality holds due to assumption  {\bf B.3$'$} and   the last inequality
is by Lemma  \ref{Master lemma} (with $M_j = I_{j,i} $). Therefore it is a contraction mapping under the given hypothesis. Hence the proof.

\eop
}

\noindent \textbf{Remarks:}  {\color{red}{}Using similar arguments, one can show that  Theorem \ref{Lemma_limit_system_uniquness}, Lemma \ref{Lemma_uniformconvergence} and Theorem \ref{Thm_onsetD} are true even for this alternate system. Also, observe that the remarks of section \ref{Graphicalmodel} holds for this Theorem as well.  In the alternate graphical model to prove the equivalent of Theorem \ref{Lemma_limit_system_uniquness} we use following norm:
\begin{eqnarray*} 
     || ({\bar x}_b ,\x)- ({\bar u}_b, \u)  ||_\infty 
       &=&  \max_{m \in \lbrace 1, 2\rbrace} \sup_{i\in \mathcal{G}_m } \left (|\bar{x}^{m1}_i - \bar{u}^{m1}_i|+ |\bar{x}^{m2}_i - \bar{u}^{m2}_i| + | {\bar x}_b -   {\bar u}_b |\right).
\end{eqnarray*} 
The rest of the details are same as in the proof of Theorem \ref{Lemma_limit_system_uniquness}}.
}
One can also prove equivalent conditions for the graph structure {\bf B.2$'$} exactly as in sub-section \ref{sec_assumptions}. 
\section{Various other graphical models}
\label{Section_examples}
In this section we consider some more variants of the graphical model, obtained after modifications of the models  discussed in  Sections \ref{Graphicalmodel} and \ref{Alternate_graphical_model}. All the previous results will be valid after some minor modifications to the proof.

 \subsection{Less Randomized Weights}
\label{subsection_withlesserandomizedweights}
Previously we assumed that the random weights satisfy  $\sum_{i} W_{j, i} +W_{j,b} = 1$ for all $j$.
Towards this the weights were normalized with sum of all the involved random quantities (see \ref{Eqn_Weights 1}). 
Now we consider a generalization in which such an equality is true only in limit. Further we don't require such a normalization. 
In all, 
we again consider a random graph with two sets of nodes and one big node as before, but now with the following modification to the connectivity details given in  \eqref{Eqn_Weights 1}
for $i\in {\cal G}_m$ with $m \in \{ 1, 2\}$ as follows:
\begin{eqnarray}
\label{Eqn_reWeights 1}
W_{j, b} &=&  \eta_j^{sb} \mbox{ (to b-node),}  \ \ 
W_{j,i} =  \frac{I_{j,i} (1-\eta_j^{sb}) } { n\gp{m}} ,  \mbox{ (to another small node $i$)},  \mbox{ with} \nonumber\\
p_m^{sb} &:= & E[\eta_j^{sb}]  \mbox{ for any } j \in {\cal G}_m.  
\end{eqnarray}
In the above $I_{j,i}$ are as before, i.e., as in \eqref{Eqn_connections} and so are the remaining details, i.e.,
  $\gp{1} := \gamma p_1+ (1-\gamma)p_{c_1}  \mbox{  and  } \gp{2} :=  \gamma p_{c_2}+ (1-\gamma) p_2.$
With   weights as in \eqref{Eqn_reWeights 1} the assumption  {\bf B.2} is readily satisfied because the denominators in \eqref{Eqn_calE}, $\sum_{i \in   \mathcal{G}_1 \cup  \mathcal{G}_2 } I_{j, i}$ is now replaced by $n \gp{m}$. Thus again under {\bf B.1},  {\bf B.3},  {\bf B.4} and  {\bf B.5}, the  Theorem \ref{Thm_MainGen 1}, Theorem \ref{Lemma_limit_system_uniquness} and Corollary \ref{Cor_Three dimensional approximation} are applicable. 
We will require some minor changes in the proof given in Appendix: for example  the term $\sum_{i\in \mathcal{G}_1 \cup \mathcal{G}_2} W_{j,i} $ is no more upper bounded by 1; towards achieving upper bound ($c'$)   in last inequality of  \eqref{Eqn_Lip_cont}, one can upper bound (sample-path wise) $\sum_{i\in \mathcal{G}_1 \cup \mathcal{G}_2} W_{j,i} $  for all $n > N_w$ as in \eqref{eqn_inequality_joint1}- \eqref{eqn_inequality_forcontinuity}.

In a similar way the weights in \eqref{Alt_Eqn_Weights 1} can be modified to the following and results of Section \ref{sec_alternate},  Theorem \ref{Thm_MainGen 2}, Theorem \ref{Thm_modified_forThm_2_Alt_system} and Corollary \ref{Cor_Three dimensional approximation_Alt_system}  are again applicable under {\bf B.1}, {\bf B.3$'$}, {\bf B.4$'$} and {\bf B.5}:
\begin{eqnarray}
\label{Alt_modEqn_Weights 1}
W_{j, b} &=&  \eta_j^{sb}  \lambda_m \mbox{ (towards $b$-node),} \mbox{ with } p_m^{sb} \ :=  \ E[\eta_j^{sb}]  \mbox{ for any } j \in {\cal G}_m, \mbox{ and, }  \ \\ 
W_{j,i} &=& \left \{ \begin{array}{lll}   \frac{I_{j,i} (1-\eta_j^{sb}) \lambda_m } { n\gamma_m p_m} ,  & \mbox{ if } i \in {\cal G}_m \\ \\
  \frac{I_{j,i}  (1-\lambda_m) } { n(1-\gamma_m)p_{c_m} } 1_{p_{c_m} > 0} ,  & \mbox{ if } i \notin {\cal G}_m  \nonumber \\
  \end{array} \right .  
    \   \   \ 
\mbox{ (towards small node $i$)}.
\end{eqnarray}

\subsection{Single Group with big node}
\label{RandomizedWeightswithSingleGroup}
When one requires performance related to a single group, one can deduce the required results by letting $p_{c_1} = p_{c_2} = 0$ and letting both the groups have the parameters of the single group, i.e., consider \eqref{Eqn_connections}-\eqref{Eqn_aggragate_forsn} with:
\begin{eqnarray}
\label{Eqn_SingleGroup}
p_1 = p_2 = p_{ss}, \  \gamma = 0.5, \   p_{c_1} = p_{c_2 } = 0, \ \   p_1^{sb} = p_2^{sb} = p^{sb} 
\end{eqnarray}
and where $\{G_i^m\}$, $\{\eta_j^{sb}\}$ and $\{\eta_j^{bs}\}$   are distributed alike for both the groups ($m = 1 $ or 2). 
In this case $\rho$ of {\bf B.4} is given by:
\begin{eqnarray}
\label{Eqn_SingleGroup1}
\rho = (1-p^{sb}) + p^{sb} = 1
\end{eqnarray}
and hence assumption {\bf B.4} is readily satisfied. Observe here that results of both Sections \ref{Graphicalmodel} and \ref{Alternate_graphical_model} coincide for single group. 
Further, one can again consider that the denominators of the weight factors are constant values $n \gamma_{p_m}$ as in \eqref{Eqn_reWeights 1}  of previous sub-section. 
The system in \cite{Systemicrisk} requires results of single group with weights as in \eqref{Eqn_Weights 1}, while  the model considered in  \cite{anof_ess} requires single group results with constant denominators as in \eqref{Eqn_reWeights 1}.

\subsection{Aggregates of some functions of other components}
\label{Subsec_Aggregatesofsomefunctions}
We now consider  another variant where
the component functions depend upon random aggregates as in previous sections, but now depend upon given functions of the other components. 
In this regard only equations \eqref{Alt_Eqn_aggregatesn}, \eqref{Alt_Eqn_FixedeqGen 1}-\eqref{Alt_Eqn_Fixedeq2Gen 4} and \eqref{Eqn_altG1_aggregate}-\eqref{Alt_Eqn_key} of Section \ref{sec_alternate} and  equations  \eqref{Eqn_aggregatefunc_altgraph}-\eqref{Alt_Eqn_barf_limit e}  of Section \ref{sec_alternate} will change. The results are true even for this model with one additional assumption as explained below. We provide the precise details for the model of Section \ref{sec_alternate}, the same can be  done for the other model.  

The fixed point equations \eqref{Alt_Eqn_FixedeqGen 1}-\eqref{Alt_Eqn_Fixedeq2Gen 4} modify to the following depending upon the given functions $h_m^{m'} (\cdot)$ as below (observe only third equation is different):
\begin{eqnarray}
\label{Alt_Eqn_FixedeqGen 1_opt}
X_i^{m}  &=& f^{m}(G^m_i,   {\bar X}^{m1}_i ,   {\bar X}^{m2}_i ,  \eta_i^{bs} X^b) \mbox{ {\normalsize  for each }} i \in \mathcal{G}_m \mbox{ and any } m,  \ \mbox{and},  \\  
\label{Alt_Eqn_Fixedeq2Gen 3_opt}
X^b  &=&    f^b (   {\bar X}^b ), \ \ \ \mbox{ {\normalsize with aggregates,}}    \\
{\bar X}^{m1}_i  & :=&   \displaystyle \sum_{j \in \mathcal{G}_1} h^m_1(X_j^1) W_{j, i} , \  \ \ \
{\bar X}^{m2}_i \  :=  \  \displaystyle \sum_{j \in \mathcal{G}_2} h^m_2(X_j^2)  {W}_{j, i} \mbox{ {\normalsize  for each }} i \in \mathcal{G}_m
\label{Alt_Eqn_aggregatesn_opt} \   \     \\
%\tilde{X}^s_i  & := &  \sum_{j \in \mathcal{G}_1} X_j^s W_{j, i} + \sum_{j \in \mathcal{G}_2} X_j^s \tilde{W}_{j, i}\mbox{ {\normalsize  for each }} i \in \mathcal{G}_2, \   \  \nonumber  
{\bar X}^b & := & \frac{1}{n}  \sum_{j 	\in \mathcal{G}_1} {X}^1_j W_{j, b} + \frac{1}{n}  \sum_{j \in \mathcal{G}_2} {X}^2_j W_{j, b} . \label{Alt_Eqn_Fixedeq2Gen 4_opt}
\end{eqnarray}
We will require that $ h^{m'}_m(.)$ are  Lipschitz continuous functions.  Further we require the following additional assumption for Theorem \ref{Thm_MainGen 2} counter part:\\
{\bf B.6$'$} 
The functions $h^{m'}_m(.)$ for any $m, m' 
\in \{1, 2\}$ are  Lipschitz continuous with  Lipschitz co-efficient 1.

We will now have the following fixed point equations for aggregate vectors (the remaining   quantities as in section \ref{sec_alternate}):

\vspace{-1mm}
{\small 
\begin{eqnarray}
\label{Eqn_alt_aggregateg1_opt}
{\bar f}^{n, m 1}_i ({\bar x}_b,\x ) &=& \left \{  \begin{array}{lll} \sum_{j \in\mathcal{G}_1}  h_1^m \left ( \xi^1_{j} ({\bar x}^{11}_j, {\bar x}^{12}_j, x_b) \right )  W_{j, i}, &\mbox{ \normalsize if }     i \in \mathcal{G}_m  \\
 0  &\mbox{\normalsize else,}    \hspace{20mm} 
\end{array} \right .   \\ \nonumber
\\
{\bar f}^{n, m 2}_i ({\bar x}_b,\x ) & = 
& \left \{  \begin{array}{lll} \sum_{j \in\mathcal{G}_2}  h_2^m \left ( \xi^2_{j} ({\bar x}^{21}_j, {\bar x}^{22}_j, x_b)  \right ) W_{j, i} &\mbox{ \normalsize if }   i \in \mathcal{G}_m  \\
 0  &\mbox{\normalsize else, \hspace{10mm}  and, }   \hspace{5mm}
 \end{array} \right. 
 \label{Eqn_alt_aggregate_opt}  \\
 {\bar f}^n_b ({\bar x}_b,{\x})& := & 
 \frac{1}{n}
\displaystyle \sum_{j \in \mathcal{G}_1}   \xi^{1}_{j} ({\bar x}^{11}_j,{\bar x}^{12}_j, x_b) W_{j, b} +  \frac{1}{n}
\displaystyle \sum_{j \in \mathcal{G}_2}   \xi^{2}_{j} ({\bar x}^{21}_j, {\bar x}^{22}_j, x_b) W_{j, b} . 
 \label{Alt_Eqn_key_opt} 
 \end{eqnarray}}
 The rest of the details are exactly the same, after modifying Lemma \ref{lem: LLN_fixed point}  with the following five dimensional fixed point equations:
\vspace{-2mm}
 {\small\begin{eqnarray}
 {\bar x}_i^{\infty,11}  ({\bar x}_b,\x )  &=&      E_{G^{1}_i, \eta_i^{bs}} \left [ h^1_1(\xi^{1}_i ({\bar x}^{11},{\bar x}^{12}, x_b))   \right ]  \lambda_1 (1-p_1^{sb}), \label{Alt_Eqn_barf_limit 1_opt} \\
 {\bar x}_i^{\infty,12}  ({\bar x}_b,\x )  &=&  E_{G^{2}_i, \eta_i^{bs}} \left [h^1_2(\xi^{2}_i ({\bar x}^{21},{\bar x}^{22}, x_b) ) \right ] \frac{1-\gamma}{\gamma}(1-\lambda_2) 1_{p_{c_2}> 0}, \nonumber \\
 {\bar x}_i^{\infty,21}  ({\bar x}_b,\x)  &=&      E_{G^{1}_i, \eta_i^{bs}} \left [h^2_1(\xi^{1}_i ({\bar x}^{11},{\bar x}^{12}, x_b) )  \right ] (1-\lambda_1)\frac{\gamma}{(1-\gamma)}1_{p_{c_1}> 0}, \nonumber  \\
  {\bar x}_i^{\infty,22}  ({\bar x}_b,\x)  &=&  E_{G^{2}_i, \eta_i^{bs}} \left [h^2_2( \xi^{2}_i ({\bar x}^{21},{\bar x}^{22}, x_b))\right ] \lambda_2(1-p_2^{sb}) ,\nonumber  \\
 {\bar x}_b^\infty  ( {\bar x}_b,\x) &=&   \gamma E_{ G^{1}_i, \eta_i^{bs}}  [\xi^{1}_i ({\bar x}^{11},  {\bar x}^{12},x_b)  ]\lambda_1  p_1^{sb}  + (1-\gamma) E_{ {G}^{2}_i, \eta_i^{bs}}  [\xi^{2}_i ({\bar x}^{21}, {\bar x}^{22}, x_b)  ]\lambda_2 p_2^{sb}. \hspace{4mm} \label{Alt_Eqn_barf_limit e_opt}
\end{eqnarray}}
 The proofs will also go through in a similar way and we will have the  results, i.e., Corollary \ref{Cor_Three dimensional approximation_Alt_system}, Theorems \ref{Thm_MainGen 2}-\ref{Thm_modified_forThm_2_Alt_system} (i.e., almost sure convergence given in \ref{Alt_Eqn_main_groupwise} is true), now using the limit fixed point  equations provided in the above equations, when additionally {\bf B.6$'$} is assumed.
 
Further one can have lesser random weight and single group modifications (the previous models of this section) for this case also. 

 \subsection{Single group with variations}
 \label{subsectionSingle group with variations}
 As in previous sub-sections, one can easily extend this analysis to the case when  the random fixed point equations are in the following form:
 \begin{eqnarray}
 \label{Eqn_fixedpoint_forsinglefuncton}
 X_i &=& f(G_i,   {\bar X}_i , \eta_i^{bs} X^b) \mbox{ {\normalsize  for each }} i \in \mathcal{G},  \ \mbox{ and } X^b =  f_b ({\bar X}^b)
 \mbox{ with, } \\
 {\bar X}_i &:=& \sum_{j \in {\cal L}_i^1} h_1 (X_j) \frac{I_{j,i} (1-\eta_j^{sb}) }{n p (1-\alpha)} + \sum_{j \in {\cal L}_i^2} h_2 (X_j) \frac{I_{j,i}(1-\eta_j^{sb})}{n p  \alpha}, \\
 {\bar X}^b & := & \frac{1}{n}  \sum_{j  } {X}_j\eta_j^{sb}.
 \end{eqnarray}
 where ${\cal L}_i^1$ and ${\cal L}_i^2$ are random subsets of size approximately  $n p \alpha$ and $n p (1-\alpha)$.
 We will require that $ h^{}_1(.),h^{}_2(.)$ are  Lipschitz continuous functions.
 One can alternatively replace the denominators with $\sum_{i'} I_{j, i'} $ also (would require {\bf B.2} type of assumption). 
 The limiting fixed point  in this case would be given by:

 \vspace{-4mm}
 {\small\begin{eqnarray}
 \label{Eqn_single_group_aggrgatefuncdependentfixedpoint}
{\bar x}^\infty  &=&     \bigg( E_{G_i, \eta_i^{bs}} \left [ h_1(\xi_i ({\bar x}^\infty, x_b))   \right ]   +  E_{G_i, \eta_i^{bs}} \left [ h_2(\xi_i ({\bar x}^\infty, x_b))   \right ] \bigg) (1-p^{sb}), 
\mbox{\normalsize  where } \nonumber \\
h_1(\xi_i ({\bar x}^\infty, x_b)) &=&h_1( f (	G_i,   {\bar x} _i, \eta_{i}^{bs} x_b ) ),  \ 
h_2(\xi_i ({\bar x}^\infty, x_b)) = h_2( f (	G_i,   {\bar x} _i, \eta_{i}^{bs} x_b ) ), \nonumber \mbox{\normalsize and, }\\
 x_b &:=&  f^b (\bar{x}_b), \  E[\eta^{sb}_j] = p^{sb}.
\end{eqnarray}}
The proofs will also go through in a similar way and again Corollary \ref{Cor_Three dimensional approximation}, Theorems \ref{Thm_MainGen 1}-\ref{Lemma_limit_system_uniquness} are true for this model; we will only require the following change in the   step \eqref{Eqn_Lemma13} of Lemma \ref{Master lemma},

\vspace{-4mm}
{\scriptsize
\begin{eqnarray*}
 \bigg| \sum_{j \in \mathcal{L}_i^1} M_j \frac{(1- \eta_j^{sb})}{n p (1-\alpha)}- E[M_j]    \frac{1-p^{sb}}{  p }   \bigg|  =  \bigg| \sum_{j \in {\cal G}} \indc{j \in {{\cal L}_i^1}} M_j \frac{(1- \eta_j^{sb})}{n p (1-\alpha)}- E[M_j] (1-\alpha) \frac{1-p^{sb}}{  p (1-\alpha) }   \bigg| 
\end{eqnarray*}}
which converges to zero  because $\left \{ \indc{j \in {{\cal L}_i^1}}\right \}_j$  are either i.i.d. for each $i$, or should satisfy an assumption like {\bf B.2(C)}.

\section{Financial Network }
\label{sec_finance}
 In the previous section, we described a graphical model with a large number of nodes.  As the number of nodes increases, one can approximate the system by a simplified limit system described by Theorems \ref{Thm_MainGen 1}-\ref{Thm_MainGen 2} and their corollaries. In this section we apply Theorem \ref{Thm_MainGen 2}, more importantly Corollary \ref{Cor_Three dimensional approximation_Alt_system}, to study systemic risk aspects in a large complex \textbf{financial} network.

In our previous work \cite{Systemicrisk}, we considered an example of a  large heterogeneous financial network with one big bank and a large number of small (identical) entities, to study the systemic risk.
In this paper, we consider further heterogeneity, with two large groups of small entities and a big bank. The entities within a group are identical but are different from those of the other group. 

%This model is analyzed using a similar approach as in the current paper.  
The network in \cite{Systemicrisk} consists of $n$ small banks and one big bank.  The small banks,  borrow some money from the big bank, also borrow from their neighbouring small banks at time $t=0$ and invest the total borrowed money along with their initial wealth into risky investments. Basically, the banks select a portfolio at time $t=0$.  They would get returns at two instances of time (time slot $1$ and $2$), depending on their portfolio and subject to the economic shocks. They attempt to clear their liabilities at time slot $t=1$, some of them may default because of the economic shocks. This can result in further defaults,  and,  these effects can percolate throughout the network. 
Using an asymptotic analysis, which uses similar flavour like that in the current paper, we showed some interesting conclusions:  a) when the banks borrow more from big bank and neighbours and invest more in risky assets, (as anticipated) the probability of defaults increases;  
b) however the expected surplus of the network increases with the increase in the investment towards risky assets; and c)   more interestingly,  the increase is possible only till a certain threshold on the investment; after this threshold the expected surplus reduces.  Thus we observed interesting non-monotone trends in expected surplus (this quantity is influenced by percolation of shocks) as the amount of investment towards risky investment is varied. 
  
  In the example of \cite{Systemicrisk}, the small banks are homogeneous. But this may not be true in many scenarios. Here we consider a network with two groups of homogeneous (within the group) entities (each as in \cite{Systemicrisk}), but the two groups have different characteristics.  For example, one group might be aggressive and might consider more risky portfolios, while the other group could be cautious. 
  We are now interested in the influence of one group on the other when some interconnections ($p_{c_m}$, $\lambda_m$ of  previous sections can represent  interconnection parameters) between the two groups are formed. 
  
  In this paper  we focus on the influence of interconnection parameters, however  one can study many other relevant and important aspects using this approach. For example, one can study the two time-period model of \cite{Systemicrisk} to further study the expected surplus as a function of inter and intra connection parameters.   One can  study the effect of the entire network on big bank, by studying its performance.  
  As already mentioned, one can study the evolutionary trends of aggressive and recessive  behaviours using replicator dynamics and random fixed point theorems of this paper as in \cite{Saha}, etc. 
  
The network of \cite{Systemicrisk} can also be analysed  using Theorem  \ref {Thm_MainGen 1} of this paper when $p_1 = p_2 = p_{c_1}=p _{c_2}$.  To analyze the  more complicated network of the current paper, we had to  extend the results  of  \cite{Systemicrisk} to Theorems \ref{Thm_MainGen 1} and \ref{Thm_MainGen 2} (and the corresponding corollaries).

We will begin by providing precise modelling  details of the system. 

%\paragraph{\textbf{Financial network description:}} 

In this paper, we  consider a heterogeneous financial network with a large number of  entities. %{\color{red}We  refer these entities as `banks', henceforth for  convenience}.
In particular, we consider a stylized example of a financial network  with $n$ small banks and one big bank (BB): a)  one group of banks are willing to take more risk, while, the other group prefers less risky portfolios, b) financial linkages of the banks are different across the groups. %Depending on the willingness of risk taken; 
Thus the network is classified into two groups, namely, $\mathcal{G}_1$ and $\mathcal{G}_2$. Let $\gamma$ and $(1-\gamma)$ be the  respective fractions of  banks in  groups $\mathcal{G}_1$ and  $\mathcal{G}_2$. %while in $\mathcal{G}_2$ contains the remaining (1-$\gamma$) fraction.
The banks are sustained in the economy for two time periods, namely, $t= 0,1$. In the initial period, (i.e., at $t=0$) banks form links by lending and/or borrowing, are also investing in  risky assets (outside the network). In the next period, i.e., at $t=1$, they  have to clear their liabilities.

\subsection{Connectivity Details} As already mentioned, we consider two groups of banks. These banks are interconnected by the credit instruments borrowed from each other or by direct cash lending. Any bank from the group $\mathcal{G}_m$ (with $m$ =$1$ or $2$)  provides loan to any  other bank in  its group with probability $p_m$, independent of other banks. Also any member from $\mathcal{G}_1$ can lend to any bank in $\mathcal{G}_2$ with probability $p_c$. The members of $\mathcal{G}_2$ prefer  risky investments, do not lend to    $\mathcal{G}_1$ (rather prefer to invest in risky assets). Further the $\mathcal{G}_2$ banks borrow (more) funds from the BB.
%for larger-scale risky investment.
To summarize, we have the following  connectivity between various entities ($j\in \mathcal{G}_m$ and $i\in \mathcal{G}_{m'}$):
\begin{eqnarray}
\label{Eqn_connectivity}
P(I_{j,i} =1 ) = p_{mm^{'}} = \left \{  
\begin{array}{llll}
 p_m  &\mbox{ if }  m= m^{'}  \mbox{, }   \mbox{ for any }  m = 1, 2   \\
 p_c  &  \mbox{ if } j \in \mathcal{G}_2~ \mbox{and} \ i\in\mathcal{G}_1  \\
  0  &  \mbox{ else, } 
\end{array}
\right  .
\end{eqnarray}
where $I_{j,i}$ is the   indicator that bank $j$ is liable to bank $i$.

\subsection{Initial Investments and Liabilities} 
We assume that each small bank has initial wealth\footnote{Most of these quantities can be changed to i.i.d.  random variables, but they are kept constants  to keep the discussions simple.}  $k_0>0$ and while that of a BB is $nk_b$, where $n$ is the number  of  small banks. %At the time $t=0$, all the values of the initial wealth are known, although our theory can analyze even when the initial capital is random (see  Section {\ref{Graphicalmodel}} ). At   time $t=0$,  banks form links to one another and also borrow funds. 
Each bank  chooses a portfolio at   $t= 0$:  by investing the borrowed amount (borrowed from other members of the network) and  the initial wealth, in  outside  risky investments, and also in lending to the other entities, as explained below.
%As already mentioned we consider a two-period model $t=0,1$, where investments are made in the initial period ($t=0$) and  returns are realized in the next period ($t=1$).

\newcommand{\gc}[1]{\mbox{\color{green} $#1$}}
%\paragraph{\textbf{Dynamics of 
{\bf $\mathcal{G}_1$ banks:}
At the time $t=0$, $\mathcal{G}_1$ banks borrow funds and lend to within the group. Also, some portion of the available wealth is lent to the banks of  $\mathcal{G}_2$, and the remaining is invested in risky assets. 
Consider a typical bank in $j \in \mathcal{G}_1$, say it  borrows a total amount of $y_1(1-\eta_j^{sb})$ from all its $\mathcal{G}_1$ lenders, in particular  from $i \in \mathcal{G}_1$ it derives an amount:
\begin{equation}
\label{eqn_lendingtowardsG_1}
I_{j, i}
\frac{y_1 (1-\eta_j^{sb})}{\sum_{i'\in \mathcal{G}_1 } I_{j, i'} }.
\end{equation}

The $\mathcal{G}_1$ banks also  borrow  an amount, $y_1 \eta_j^{sb}$, from  BB.
Thus total liability of $j \in \mathcal{G}_1$ at $t = 0$, towards all  the banks in $\mathcal{G}_1$ equals:
\begin{equation}
\label{eqn_lenadingforintergroup}
\sum_{ i \in \mathcal{G}_1 }  I_{j, i}  \frac{y_1(1-\eta_j^{sb})}{\sum_{i'\in \mathcal{G}_1 } I_{j, i'} } = y_1(1-\eta_j^{sb}).
\end{equation}
This borrowed amount has to be repaid at $t=1$,  with  interest rate $r_1 >0$. Thus  the total liability of any  $\mathcal{G}_1$ bank at time period $t=1$ equals, $\bar{y}_1 = y_1(1+r_1)$.

A typical agent    $j \in \mathcal{G}_2$ borrows a total of  $y_c$ from all its $\mathcal{G}_1$ lenders,  in particular from $i \in \mathcal{G}_1$, it   borrows an amount:
\begin{equation}
I_{j, i}
\frac{y_c}{\sum_{i'\in \mathcal{G}_1} I_{j, i'} }.
\label{Eqn_yc}
\end{equation}
The rate of  interest of this liability is  $r_2$. We   assume $r_1 < r_2$; the agents of ${\mathcal G}_1$ agree to lend to   ${\mathcal G}_2$, only at a higher interest rate, as they are aware  of the risky nature of the latter group.   % intra  group liability is more costly than inter-bank liability i.e., $r_2 > r_1$.
%Thus $G_1$ banks have claims at $t=1$  is from the within $G_1$ as well as form the $G_2$ banks.

In all, the money lent by a typical agent $i \in \mathcal{G}_1 $, at $t=0$, towards other banks of the network approximately equals (for large $n$ and with assumption\footnote{\label{foot_approx_exact}One can make this convergence rigorous    exactly as in the proofs of the previous section and the expressions would be exact at limit (a.s.). } {\bf B.2$'$}):

\vspace{-4mm}
{\small \begin{eqnarray}
\label{eqn_approximatelendingG_1}
\hspace{-2mm} \sum_{j  \in \mathcal{G}_1 } \hspace{-1mm}
\frac{I_{j, i} y_1 (1-\eta_j^{sb})}{\sum_{i'\in \mathcal{G}_1} I_{j, i'} } + \hspace{-.2mm} \sum_{j\in \mathcal{G}_2 } 
\hspace{-.1mm}\frac{ I_{j, i}y_c}{\sum_{i'\in \mathcal{G}_1} I_{j, i'} } \approx
%y_1 (1-p_1^{sb}) + \frac{(1-\gamma) y_c  p_c}{\gamma p_c}= 
y_1 (1-p_1^{sb}) + \frac{(1-\gamma)y_c  }{\gamma  }   . \hspace{1.5mm}
\end{eqnarray}}
\paragraph {\textbf{ $\mathcal{G}_2$ banks:}} The $\mathcal{G}_2$ banks borrow and lend from within the group, as well as, borrow from  $\mathcal{G}_1$ banks. In addition,   the $\mathcal{G}_2$ banks   borrow more form the BB ($p_2^{sb} > p_1^{sb}$) for a bigger risky investment. Say $y_2$ is the total amount that a typical $\mathcal{G}_2$ bank borrows at $t=0$ (from BB and from $\mathcal{G}_2$): the   agent   $j \in \mathcal{G}_2$ borrows $y_2\eta^{sb}_j$ from BB,  and, 
 \begin{equation}
 \label{eqn_lendingG_2toG_1}
 I_{j, j''}\frac{ y_2 (1-\eta^{sb}_j) }{
 \sum_{j'\in \mathcal{G}_2} I_{j, j'}
 },
 \end{equation}
 from each   $j''\in \mathcal{G}_2$, that is interested in giving loan to agent $j$. Here $\{\eta_j^{sb}\}$ are i.i.d. random variables with expected value $p_2^{sb}$. Thus the total loan (from BB and $\mathcal{G}_2$) taken by agent $j$ at $t= 0$ is given by:
 \begin{equation}
 \label{eqn_overalllendingG_2}
 y_2 \eta^{sb}_j + \sum_{j'' \in \mathcal{G}_2} I_{j, j''}\frac{ y_2 (1-\eta^{sb}_j) }{
 \sum_{j'\in \mathcal{G}_2} I_{j, j'} } = y_2.
 \end{equation}
In all, from \eqref{Eqn_yc},  the  total liability of the $\mathcal{G}_2$ banks at time period $t=1$ becomes $\bar{y}_2:= (y_2+y_c)(1+r_2)$. Similarly, the money lent by a typical agent $i \in \mathcal{G}_2 $, at $t=0$, towards other banks of the network approximately equals (for large $n$ and with assumption  {\bf B.2$'$}):\vspace{-4mm}
\begin{equation}
\label{eqn_approxG_2withinlending}
      \sum_{j\in \mathcal{G}_2 }  I_{j, i}
\frac{y_2 (1-\eta^{sb}_j) }{\sum_{i'\in \mathcal{G}_2} I_{j, i'} } \approx y_2 (1-p^{sb}_2)   .
\end{equation}

\paragraph {\textbf{ Big bank:}} The BB only provides loans to the  small banks, and    has zero liability. 

\paragraph {\textbf{ Risky investments:}}
As already mentioned, banks select a portfolio  along with liability connections. The banks invest the remaining money (after lending and borrowing) in risky investments (at $t=0$). Thus node $i \in \mathcal{G}_1$ invests (see   (\ref{eqn_lendingtowardsG_1}) - (\ref{eqn_approximatelendingG_1})),  with $(x)^+ := \max\{0, x\}$:

\vspace{-2mm}
{\small
\begin{eqnarray}
\Omega^1_{i}  = \left (k_0 + y_1 \eta^{sb}_i + \sum_{ j \in \mathcal{G}_1 }    \frac{I_{ i,j} y_1 (1-\eta_i^{sb})}{\sum_{j'\in \mathcal{G}_1 } I_{i, j' } } -   \left (   \sum_{j  \in \mathcal{G}_1 }
\frac{ I_{j, i} y_1 (1-\eta_j^{sb})}{\sum_{i'\in \mathcal{G}_1} I_{j, i'} } + \sum_{j\in \mathcal{G}_2 } 
\frac{ I_{j, i} y_c}{\sum_{i'\in \mathcal{G}_1} I_{j, i'} } \right ) \right )^+\hspace{-2mm}. \hspace{4mm} \label{Eqn_Omega1} \end{eqnarray}
}
For  large $n$, the risky (out-side) investment of   any   bank from $\mathcal{G}_1$ approximately  equals (see  (\ref{eqn_lendingtowardsG_1}) - (\ref{eqn_approximatelendingG_1})), which is exact at limit as in footnote \ref{foot_approx_exact}: 
\begin{eqnarray}
\Omega_{1}  
&\approx&   \bigg (k_0 +  y_1 -   y_1 (1-p_1^{sb}) - \frac{ (1-\gamma) y_c} {\gamma}\bigg )^+= \bigg(k_0+ y_1p_1^{sb} -  \frac{ (1-\gamma) y_c} {\gamma} \bigg)^+\hspace{-2mm}. \hspace{4mm}
\end{eqnarray}

  In a similar way, the  risky investment by  bank   $j \in \mathcal{G}_2$  equals 
  (see \eqref{eqn_lendingG_2toG_1}-\eqref{eqn_approxG_2withinlending}):
  
\vspace{-4mm}
{\small\begin{eqnarray}
\label{Eqn_Omegaj2}
\Omega_j^2  =   k_0 + y_2 \eta^{sb}_j + \hspace{-2mm}\sum_{j'' \in \mathcal{G}_2}\frac{ I_{j, j''}  y_2 (1-\eta^{sb}_j) }{
 \sum_{j'\in \mathcal{G}_2} I_{j, j'} }    + \hspace{-2mm} \sum_{i \in \mathcal{G}_1 }
\frac{ I_{j, i} y_c}{\sum_{i'\in \mathcal{G}_1} I_{j, i'} } - \hspace{-2mm} \sum_{i\in \mathcal{G}_2 }  
\frac{I_{i, j} y_2 (1-\eta^{sb}_i) }{\sum_{i'\in \mathcal{G}_2} I_{i, i'} }, \hspace{3mm}
\end{eqnarray}}
which  approximately equals (the same for any $j \in \mathcal{G}_2$ at limit):
\begin{eqnarray}
\Omega_2 \approx  k_0 + y_2    + y_c - y_2 (1-p^{sb}_{2})
 = k_0 + y_2 p^{sb}_{2} + y_c.
 \label{eqn_omega_2}
\end{eqnarray}
\ignore{
{\color{red}
Similarly, the risky investment   of the BB  equals:

 \begin{eqnarray}
 \Omega_b \approx (k_b - y_2p_2^{sb})^+.
 \label{eqn_omega_b}
 \end{eqnarray}}
 }
 \subsection {Economic shocks at $t=1$}
 The banks receive returns from their risky investments at time period, $t=1$. These returns can have shocks. We assume   binomial distribution to  model the shocks,  as is majorly considered in literature (see  e.g.,  \cite{acemoglu2015systemic,Saha,Gai,Goldstein}).  The  (risky) asset prices at time period $t=1$, can have upward movement with rate $u$ and this happens  with   probability $1-w$,  while, the price can have downward movement (rate $d$)  with probability $w$. By   {\it standard no-arbitrage principle it is reasonable to assume that  $d < r_1 < r_2 < u$}. Thus the (random) returns of   the risky investments at $t=1$ equal (for $m = 1,2 $): 

 \begin{equation}
 \label{eqn_assetpriceG_1}
 \mbox{ $K^m_j  = \Omega_j^m (1+\UV^m_j-d_c) $, where $\UV^m_j  = \begin{cases}  u  &w.p. (1-w) \\ d &else \end{cases}$ }  
\end{equation}
  where {\it $d_c$ is the common shock} (for example created by COVID-19 pandemic) which can affect all the banks. The shocks $\{V_j^m \}$ are i.i.d across all banks, irrespective of $m$.
\subsection{Returns and Clearing vector} 
At time $t=1$ all the entities receive returns from their risky (outside) investments. Using these returns the banks attempt to clear the liabilities,  created during the time period  $t=0$, further using the returns from the other banks. The final payments made by the banks, after clearing the liabilities to the maximum extent possible,  are called the clearing vector  (e.g., \cite{acemoglu2015systemic,eisenberg2001systemic,Systemicrisk}).   However,  the risky investments   are subjected to economic shocks (see \eqref{eqn_assetpriceG_1}), which could significantly reduce the returns of some (or all) banks. This in turn can potentially reduce the clearing capacity of the connected banks, and this goes on. {\it Systemic risk precisely studies this aspect, basically  micro-level  (entity-level) shocks could trigger  cascade of defaults, which can eventually  lead to the collapse  of the entire system.} 
%We are assuming the loss of the banks is proportional to the investments.
 Let $X^{m}_i$ denote the clearing value of the $i$-th bank of group $\mathcal{G}_m$, which indicates the maximum possible amount (out of the liability), cleared by $i$-th bank. 
 The clearing vector  $X = (X_i^m)_{i,m}$ is obtained by the standard bankruptcy rule, i.e., under the assumption of limited liability and pro-rata basis repayment of the debts in case of default  (e.g.,  \cite{acemoglu2015systemic,eisenberg2001systemic,Systemicrisk}); here the amounts returned are proportional to their liability ratios; the  bank  $j\in \mathcal{G}_2$ pays\footnote{We drop the group notation $m$, when there is no ambiguity, to keep notations simple.} back $X_j W_{j,i}$ towards bank $i$, where $W_{j,i}$,  the liability fraction borrowed during the initial period, equals (see \eqref{eqn_lendingtowardsG_1}-\eqref{eqn_approxG_2withinlending}):
 \begin{equation}
 \label{eqn_liabilityG_1}
   W_{j,i}= 
    \begin{cases}
      \frac{I_{j,i}}{\sum_{j'' \in \mathcal{G}_1} I_{j, j''} }\frac{y_c }{y_2 + y_c} , & \text{for~ all}\ i \in \mathcal{G}_1 , \ j \in \mathcal{G}_2 \\
       \frac{I_{j,i} (1-\eta^{sb}_{j})}{\sum_{j'' \in \mathcal{G}_2} I_{j, j''} } \frac{y_2 }{y_2 + y_c} ,
       & \text{for~ all}\ i,j \in \mathcal{G}_2.
    \end{cases}
  \end{equation}
  Similarly the liability  fractions   for the entities of  $\mathcal{G}_1$  are given by:
  \begin{equation}
  \label{Eqn_liabilityforG1G2}
    W_{j,i}=
     \frac{I_{j,i} (1-\eta_j^{sb})}{\sum_{j'' \in \mathcal{G}_1} I_{j, j''} },  
    \text{ for~ all}\ i,j  \in \mathcal{G}_1.
\end{equation}
 Thus the maximum possible amount cleared by any agent  $i \in \mathcal{G}_1$ is  given by the following (random fixed point) equation,
 \begin{equation}
 \label{eqn_clearingvectorG_1}
X_i = \min \bigg \lbrace   \bigg(K^{1}_i +  \displaystyle\sum_{j\in \mathcal{G}_1} X_j W_{j,i}+ \displaystyle\sum_{j\in \mathcal{G}_2} X_j W_{j,i}-v_1 \bigg)^+,\bar{y}_1 \bigg \rbrace ,  \mbox{ where, }
\end{equation}
\begin{enumerate}[(a)]
    \item the first term denotes the return from  risky investment,  note that $\{K^{1}_i\}_{i \in {\cal G}_1}$ are i.i.d. random variables distributed according to $K_i^1$ defined in equations \eqref{Eqn_Omega1}-\eqref{eqn_assetpriceG_1};
    \item  the claims from the other  banks are given by the second and third term, ${\sum_{j\in \mathcal{G}_1}} X_j W_{j,i}+  {\sum_{j\in \mathcal{G}_2}} X_j W_{j,i}$;    
    \item the fourth term, $v_1$, is the taxes/security deposits/senior debt;  and
    \item the banks repay at maximum ${\bar y}_1$, their total liability.
\end{enumerate}  
 The banks first have to clear the taxes, the remaining  money can then be distributed to its creditors according to pro-rata basis (as in  \cite{acemoglu2015systemic,Systemicrisk} and see  \eqref{eqn_liabilityG_1}, \eqref{Eqn_liabilityforG1G2}).
%We assume the taxes  are proportional to the upward return, i.e.,   $v_1 =\kappa \Omega_1(1+u-d_c)$ for some appropriate $0< \kappa < \infty$.
%Similarly,  the clearing vector for   $\mathcal{G}_2$  group of banks is given by: 
\ignore{
\begin{equation}
\label{eqn_assetpriceG_2}
    K_{2}=
    \begin{cases}
     \Omega_2(1+u-d_c), & \text{wp}\ 1-w  \\
      \Omega_2(1+d-d_c), & \text{otherwise}.
    \end{cases}
  \end{equation}
  }
  %
  %
%The liability fraction between the entities of group $\mathcal{G}_2$ are given as:
%\begin{equation}
    %\tilde{L}_{j'j}=
    %\begin{cases}
      %I_{j',j} \frac{1- \eta^{sb}_{j'}}{\sum_{j'' \in \mathcal{G}_2} I_{j', j''} }\frac{y_2}{y_2 + y_c}, & \text{for~ all}\ i \in \mathcal{G}_2  \\
       %0,  & \text{for~ all}\ i \in %\mathcal{G}_1.
    %\end{cases}
  %\end{equation}
In a similar way,  the clearing vector for a typical entity from  $\mathcal{G}_2$ equals,
\begin{equation}
\label{eqn_clearingvectorG_2}
X_i = \min \bigg \lbrace \bigg( K^{2}_i  +  \ \displaystyle\sum_{j\in \mathcal{G}_2} X_jW_{j,i} -v_2 \bigg)^+,\bar{y}_2 \bigg \rbrace,  \mbox{ for any } ~i \in \mathcal{G}_2.
\end{equation}
%where,  $v_2 = \kappa \Omega_2(1+u-d_c)$.\\
\ignore{
{\color{red}
Let $W_{jb}$ be the liability fraction of  $\mathcal{G}_2$ banks towards BB and it is given by:
\begin{equation}
    W_{jb}:= \eta^{sb}_j  \frac{y_2  }{y_2+y_c} ~\mbox{for} ~ j \in \mathcal{G}_2.
\end{equation}
Also the clearing money  for the big bank is given by:
\begin{equation}
\label{eqn_clearingvectorBB}
X_b =  \bigg ( K_{b}  +  \ \displaystyle\sum_{j\in \mathcal{G}_2} X_j \frac{W_{jb}}{n(1-\gamma)}-v_b \bigg )^+.
\end{equation}}
}
%where,  $v_b = \kappa \Omega_b (1+u-d_c).$
\subsection{ Systemic risk performance measures} 
We consider three important performance measures related to systemic risk.\\
\begin{enumerate}
\item {\textbf{Probability of default:}}
 We say a bank  defaults  when it is unable to settle the liability amount at period $t=1$. The probability of such an event is an important aspect for the network and let:% $P^{m,n}_{D,i}$ be this probability:  
\begin{eqnarray}
\label{eqn_prelimitPD}
 P^{m,n}_{D,i} : = P(X_i < \bar{y}_m) \ \mbox{with} \ m=1,2,  \  i\in \mathcal{G}_m \ \ignore{{\color{red}\mbox{and}, \ P^{b}_D: =  1_{X_b < 0 }}}.
 \end{eqnarray}
%{\color{red} We are also interested in the fraction of small banks that defaulted in a given group (say group $m$), i.e., ${\mathbb F}^{m,n} := \frac{\sum_{i \in {\cal G}_m}  1_{\{X_i < \bar{y}_m\} } }{n_m}$.
% }

\item {\textbf{Expected Surplus:}}  The surplus of  any bank is the total income of the bank
(small banks), after clearing the liabilities and taxes. 
Let $E_i^{n}[S_m]$ be the expected surplus of a typical agent of group $\mathcal{G}_m$ with  $m= 1,2$ and it is formally defined as follows: 
\begin{eqnarray}
\label{eqn_prelimitsurplussb}
E_i^{n}[S_1]: &= &E\bigg (  K^{1}_i  +  \displaystyle\sum_{j\in \mathcal{G}_1} X_j W_{j,i}+ \displaystyle\sum_{j\in \mathcal{G}_2} X_j W_{j,i}-v_1- \bar{y}_1 \bigg )^+ \ \mbox{for} \ i\in \mathcal{G}_1, \mbox{ and, }\nonumber \\
E_i^{n}[S_2]:& =& E \bigg ( K^{2}_i  +  \ \displaystyle\sum_{j\in \mathcal{G}_2} X_j W_{j,i} -v_2- \bar{y}_2 \bigg )^+ \ \mbox{for}  \ i\in \mathcal{G}_2,
\end{eqnarray}
\ignore{
{\color{red}
while that of  the BB (per small bank)  is defined as:
\begin{equation}
\label{eqn_prelimitsurplusbb}
E^n[S_b]:=  E\bigg ( K_{b}  +  \ \displaystyle\sum_{j\in \mathcal{G}_2} X_j \frac{W_{jb}}{n(1-\gamma)}-v_b \bigg )^+.
\end{equation}}
}
\item \textbf{Returns with upward movement:} 
 The banks within the network have heterogeneous belief towards the asset price returns from the risky investment. The $\mathcal{G}_2$ banks believe that the asset price will go up at period $t=1$ with high probability. Hence these banks would be interested in best possible returns from their investments. In this regard, we define  a third performance measure as the best possible surplus (one achieved with upward movement of risky asset) as below (see  \eqref{eqn_assetpriceG_1} and (\ref{eqn_prelimitsurplussb})):
 \begin{eqnarray}
\label{eqn_prelimitupwardreturnsb}
%S_{1,u}: &= & \bigg (  K_{1,u}  +  \displaystyle\sum_{j\in \mathcal{G}_1} X_j L_{ji}+ \displaystyle\sum_{j\in \mathcal{G}_2} X_j L_{ji}-v_1- \bar{y}_1 \bigg )^+ , \nonumber \\
\hat{S}^{n}_{2,u,i}:& =&  \bigg ( K^{2}_{u,i}  +  \ \displaystyle\sum_{j\in \mathcal{G}_2} X_j W_{j,i} -v_2- \bar{y}_2 \bigg )^+,\  K^2_{u,i} := \Omega^2_i (1+u-d_c). \hspace{2mm}
\end{eqnarray}
{\it We refer this as SaU, the Surplus at Upward movement.} 
 As opposed to this,   the $\mathcal{G}_1$ banks are interested  only in the expected surplus.  
 \end{enumerate}
% Also discuss the expected surplus and return of banks with the upward movement Appropriate notations..
\section{Asymptotic approximation of the banking network}\label{sec_asymptotic}
 The finite banking network   is complicated to analyze;  the most complex aspect  being the derivation of the clearing vector. Further, usually, the number of entities in such a network is sufficiently large. Thus, we obtain  asymptotic (as $n\to \infty$) analysis; 
we derive the approximate closed-form expression for the clearing vector using Theorem \ref{Thm_MainGen 2} and Corollary \ref{Cor_Three dimensional approximation_Alt_system}.
This becomes instrumental in deriving the systemic risk performance measures (discussed above).  In Section \ref{sec_MCsimulation} using exhaustive Monte-Carlo simulations, 
we   demonstrate the accuracy of this approximation  even for moderate values of $n$ (the number of banks).

The clearing vector equations (\ref{eqn_clearingvectorG_1})-(\ref{eqn_clearingvectorG_2}) can be viewed as  random fixed point equations, which depend upon the realizations of the economic shocks $\{K^{1}_i \}_{i \in {\cal G}_1}$, $\{K^{2}_i\}_{i \in {\cal G}_2}$ to the   network. This financial system is exactly like the graphical model discussed in Section \ref{Alternate_graphical_model} with the following mapping   details (see  \eqref{Alt_Eqn_FixedeqGen 1} - \eqref{Alt_Eqn_Fixedeq2Gen 4} and \eqref{Eqn_connectivity}):
\begin{eqnarray}
 G^m_i &=& K^{m}_i, \  \eta_i^{bs} =0  \  (a.s.), \ \lambda_1 = 1 ,\  \ \lambda_2 = \frac{y_2}{y_2+y_c} , \  p_{c_1} = 0 , \ p_{c_2} = p_c > 0  , \nonumber \\ 
\mu_1 & = & \frac{1-\gamma}{\gamma}\frac{y_c}{y_2}\frac{1}{(1-p^{sb}_2)}, \  \mu_2 = 0, \  p^{sb}_1 > 0 \mbox{ and }  p^{sb}_2 >0. \label{Eqn_mu1mu2}
\end{eqnarray} 
Further observe from \eqref{eqn_clearingvectorG_1}-\eqref{eqn_clearingvectorG_2} and equation \eqref{Alt_Eqn_xi 1} that:
$$
\xi_i^m (x^{m1}, x^{m2}, x_b) =  \left \{ 
\begin{array}{ll}
\min \bigg \lbrace   \bigg(K^{1}_i +  x^{11} + x^{12} -v_1 \bigg)^+,\bar{y}_1 \bigg \rbrace     &  \mbox{ if }  m = 1, i \in {\cal G}_1, \\
   \min \bigg \lbrace \bigg( K^{2}_i  +  \ x^{22} -v_2 \bigg)^+,\bar{y}_2 \bigg \rbrace   &  \mbox{ else. }
\end{array}
\right .
$$
Thus assumption {\bf B.3$'$} is satisfied with $\sigma = 1$ and   any $0 \le \varsigma < 1$ (does not depend upon $x_b$). 
We assume {\bf B.2$'$}, {\bf B.4$'$} and that $\eta_j^{sb} \ge {\bar \eta} >0$ a.s. for all $j$. It is easy to verify that  assumption {\bf B.1} is satisfied (see  \eqref{eqn_clearingvectorG_1}-\eqref{eqn_clearingvectorG_2}). % 
Thus by   Corollary \ref{Cor_Three dimensional approximation_Alt_system},  the aggregate clearing vector   converges almost surely (see   \eqref{Alt_Eqn_aggregate_groupwise}-\eqref{Eqn_aggrragate_al_system} and $\mu_1$  as in \eqref{Eqn_mu1mu2}): 

\vspace{-4mm}
{\small 
\begin{eqnarray*} 
\displaystyle\sum_{j\in \mathcal{G}_1} X_j W_{j,i}+ \displaystyle\sum_{j\in \mathcal{G}_2} X_j W_{j,i}
\to \ {\bar x}_1^\infty + \mu_1 {\bar x}_2^\infty  \mbox{ for } i \in {\cal G}_1 
 \mbox{, \normalsize and } 
 \displaystyle\sum_{j\in \mathcal{G}_2} X_j W_{j,i} \ \to \ {\bar x}_2^\infty   \mbox{   for } i \in {\cal G}_2, 
\end{eqnarray*}}
  where  ${\bar x}_1^\infty $, ${\bar x}_2^\infty $  satisfy the following deterministic fixed point equations:
  
 \vspace{-4mm}
{\small
\begin{eqnarray}
 \label{eqn_avgclearingvectorG1}
{\bar x}_1^\infty &=& E \bigg[  \min\bigg \lbrace  \bigg ( K^{1}_i +  {\bar x}_1^\infty +  {\bar x}_2^\infty\frac{1- \gamma}{\gamma} \frac{y_c}{y_2}\frac{1}{(1-p^{sb}_2)} - v_1 \bigg)^+  , \ {\bar y}_1 \bigg \rbrace  \bigg] (1-p_1^{sb}),\\
\label{eqn_avgclearingvectorG2}
  {\bar x}_2^\infty &=& E \left [ \min \bigg  \lbrace \bigg ( K^{2}_i +  {\bar x}_2^\infty - v_2 \bigg)^+  , \ {\bar y}_{2}  \bigg \rbrace \right ] (1-p^{sb}_{2}) \frac{y_2}{y_2 + y_c}.\ 
\end{eqnarray} }
Observe that  the aggregate clearing vector   converges (in   almost sure sense) to a constant value, which is the same for all the banks in the same group. By the same corollary,  the clearing vector  converges almost surely to: 
\begin{eqnarray}
X_i &\to& \min \bigg \lbrace \bigg ( K^{2}_i +  {\bar x}_2^\infty   - v_2 \bigg)^+  , \ {\bar y}_{2}   \bigg \rbrace, \ \mbox{   for all } i \in \mathcal{G}_2, \mbox{ and, } \label{Eqn_Finance_Conv}
 \\ 
X_i &\to &  \min \bigg \lbrace   \bigg ( K^{1}_i +  {\bar x}_1^\infty +  {\bar x}_2^\infty \frac{1-\gamma}{\gamma}\frac{y_c}{y_2}\frac{1}{(1-p^{sb}_2)}- v_1 \bigg)^+  , \ {\bar y}_1  \bigg \rbrace, \  \mbox{   for all } i \in \mathcal{G}_1.
\nonumber
\end{eqnarray}
  
  %where, ${\bar x}_1^\infty$ satisfies the following fixed point equation:
 %\begin{eqnarray}
 %\label{eqn_avgclearingvectorG1}
%{\bar x}_1^\infty = E \bigg[  \min\bigg \lbrace  \bigg ( K_{1,i} +  {\bar x}_1^\infty +  {\bar x}_2^\infty \frac{y_c }{y_2 + y_c} 
%\frac{1-\gamma }{\gamma}- v_1 \bigg)^+  , \ {\bar y}_1 \bigg %\rbrace  \bigg]  \ \mbox{a.s.}
%\end{eqnarray}

%\vspace{-2mm}
%{\small
%\begin{eqnarray}
 %\label{eqn_avgclearingvectorG1}
%{\bar x}_1^\infty = E \bigg[  \min\bigg \lbrace  \bigg ( K^{1}_i +  {\bar x}_1^\infty +  {\bar x}_2^\infty\frac{1- \gamma}{\gamma} \frac{y_c}{y_2}\frac{1}{(1-p^{sb}_2)} - v_1 \bigg)^+  , \ {\bar y}_1 \bigg \rbrace  \bigg] (1-p_1^{sb}).
 %\ \mbox{a.s.}
%\end{eqnarray} }

\ignore{
{\color{red}
And the asymptotic return of the big bank is given by:

\begin{equation*}
X_b \to \bigg ( K_{b}  +  {\bar x}_2^\infty p^{sb}_{2} \frac{y_2}{y_2+y_c}-v_b\bigg )^+\ \mbox{a.s.}
\end{equation*}
}
}
\noindent \textbf{Asymptotic default probability:} The probability of default of any bank from group $\mathcal{G}_m$ converges  by bounded convergence theorem and by \eqref{Eqn_Finance_Conv} as $n \to \infty$: %{\color{red} and indicator that the BB defaults} 
\begin{eqnarray}
\label{eqn_limitPD}
P^{m,n}_{D, i} \to 
 P^{m}_D: =P(X_i({\bar x}_1^\infty , {\bar x}_2^\infty ) 
 < \bar{y}_m) \ \mbox{ for any }  \  i\in \mathcal{G}_m, \mbox{with} \ m\in \{1,2\}.\
 %\mbox{and}, \ {\color{red}P^{b}_D: = 1_{X_b(\bar{x}^{\infty}_2) < 0}}. 
 \end{eqnarray}
 Observe here that $\{ X_i({\bar x}_1^\infty , {\bar x}_2^\infty)\}_i$ are identical for all $i$ from the same group and hence the right hand side is the same for any $i$  of the same group.
%{\color{red} Further
  %the fraction of small banks that defaulted, ${\mathbb F}^{m,n}$,  also converges,
%by bounded convergence Theorem and Theorem \ref{Thm_MainGen 2},   
%$$
%{\mathbb F}^{m,n} = 
%\frac{\sum_{j \in {\cal G}_m}  1_{\{X_j < \bar{y}_m\} } }{n_m} \to  P^{m}_D \mbox{ a.s.}
%$$}
 
\noindent \textbf{Asymptotic expected surplus:} 
 By again using   \eqref{Eqn_Finance_Conv} and bounded convergence theorem the expected surplus of any bank of each group    is obtained as below:
 \begin{eqnarray}
\label{eqn_limitsurplussb}
E_i^{n}[S_1] \to 
E[S_1]: &= &E\bigg (  K^{1}_i  +   {\bar x}_1^\infty +  {\bar x}_2^\infty \frac{1-\gamma}{\gamma}\frac{y_c}{y_2} \frac{1}{(1-p^{sb}_2)}-v_1- \bar{y}_1 \bigg )^+ \ \mbox{for} \ i\in \mathcal{G}_1, \nonumber \\
E_i^{n}[S_2] \to E[S_2]:& =& E \bigg ( K^{2}_i  +  {\bar x}_2^\infty   -v_2- \bar{y}_2 \bigg )^+ \ \mbox{for}  \ i\in \mathcal{G}_2.
\end{eqnarray}
\textbf{Asymptotic SaU:} In a similar way, by \eqref{Eqn_Finance_Conv} the asymptotic surplus with upward movement (SaU) is  obtained as follows (almost surely):
\begin{equation}
\label{Eqn_surplus_atu_G_2}
\hat{S}^{n}_{2,u,i} \to
\hat{S}_{2,u}    := \bigg ( K^{2}_u  +   \bar{x}^{\infty}_2  -v_2- \bar{y}_2 \bigg )^+.
\end{equation}
Thus we have a simplified limit system, conditioned on the common shock ($d_c$), and one can  compute the performance measures of systemic risk. 
\ignore{{\color{blue}
The following  are the governing equations of the limit system ($p_{c_1} = 0$, $p_{c_2} \ne 0$):
\begin{eqnarray*}
X_i &\to& \min \bigg \lbrace \bigg ( K_{2,i} +  {\bar x}_2^\infty \frac{y_2}{y_2 + y_c}   - v_2 \bigg)^+  , \ {\bar y}_{2}   \bigg \rbrace \ \mbox{a.s., for all } i \in \mathcal{G}_2, \\
X_i &\to &   \min \bigg \lbrace   \bigg ( K_{1,i} +  {\bar x}_1^\infty  - v_1 \bigg)^+  , \ {\bar y}_1  \bigg \rbrace \ \mbox{a.s., for all } i \in \mathcal{G}_1,
\end{eqnarray*}
The aggregate clearing vector satisfies the following fixed point equations:

\begin{eqnarray*}
{\bar x}_2^\infty &=& E \left [ \min \bigg  \lbrace \bigg ( K_{2,i} +  {\bar x}_2^\infty  - v_2 \bigg)^+  , \ {\bar y}_{2}  \bigg \rbrace \right ](1-p^{sb}_{2})  \frac{(1-\gamma ) p_2}{ \gp{2} } \\
{\bar x}_1^\infty&=& E \left [ \min \bigg  \lbrace \bigg ( K_{1,i} +  {\bar x}_1^\infty  - v_1 \bigg)^+  , \ {\bar y}_{1}  \bigg \rbrace \right ] + E \left [ \min \bigg  \lbrace \bigg ( K_{2,i} +  {\bar x}_2^\infty  - v_2 \bigg)^+  , \ {\bar y}_{2}  \bigg \rbrace \right ](1-p^{sb}_{2})\frac{(1-\gamma ) p_2}{ \gp{2}}.
\end{eqnarray*}
Also note that $\gp{2} = \gamma p_c  +(1-\gamma)p_2$, $\gp{1} = \gamma p_1$
}}
\ignore{
\newpage
Consider a network with $n$ small banks and one big bank (one can talk without as well). The banks are partitioned into two groups based on the willingness of the different investments in the market. Let say $\gamma > 0$  fraction of banks in group 1 and rest in group 2.
The economy sustained for two time period $t=0,1$. At the initial stage, investments are made, and at the next stage, returns are realized.
{
\color{red}
\begin{equation*}
    \tilde{K}^1=
    \begin{cases}
      k_0- y_{c}, & \text{if lent to group 2 banks}  \\
      k_0, & \text{otherwise}.
    \end{cases}
  \end{equation*}
  Thus $\tilde{K}^1$ is  invested in the outside security for more significant return, and the following equation gives its probability distribution:
  
  \begin{equation*}
    K^{1}=
    \begin{cases}
     (\tilde{K}^1(1+u)-Z_c)^+, & \text{wp}\ 1-w  \\
        (\tilde{K}^1(1+d)-Z_c)^+, & \text{otherwise}.
    \end{cases}
  \end{equation*}
  Thus the clearing vector for this banks are given as follows:
  \begin{equation}
X_i = \min \bigg \lbrace   K^{1}  +  \displaystyle\sum_{j\in G_1} X_j w_{ji}+ \displaystyle\sum_{j\in G_2} X_j w_{ji}-v^1,Y_1 \bigg \rbrace ,  \forall i \in G_1.
\end{equation}
Where, $v^1$ be the senior debt that each bank has to clear first and $Y_1= y_0(1+r_1)$ be the liability amount that bank has to make at the end of the contract period, and the random weights are given by:
\begin{eqnarray}
w_{j, b} &=&  \eta_j^{sb}  \ \ 
w_{j,i} =  \frac{I_{j,i} (1-\eta_j^{sb}) } {\displaystyle \sum_{i' \le n} I_{j, i'}  } , 
\end{eqnarray}

\begin{eqnarray}
P(I_{i,j} =1 ) = p_{mm^{'}} = \left \{  
\begin{array}{llll}
 p_m  &\mbox{ if }  m= m^{'}  \mbox{, }   \mbox{ for any }  m = 1, 2   \\
 p_c  &  \mbox{ if group 1 lends to group 2 } \\
  0  &  \mbox{ else. } 
\end{array}
\right  .
\end{eqnarray}
Dynamics: Now consider the dynamics of the group 2 bank, the group 2 banks have the same initial wealth as group 1. Group 2 banks borrow funds  $y_b$ from the big bank, group 2  or group 1 to make a  bigger investments return. The total accumulated  wealth of the group 2 banks invested in the outside security   and its return are as follows:

\begin{equation*}
    \tilde{K}^{2}=
    \begin{cases}
      k_0+y_b/(1+r_b) +y_{c}, & \text{if borrow from group 1 banks}  \\
      k_0+y_b/(1+r_b), & \text{otherwise}.
    \end{cases}
  \end{equation*}

\begin{equation*}
    K^{2}=
    \begin{cases}
     (\tilde{K}^{2}(1+u)-Z_c)^+, & \text{wp}\ 1-w  \\
        (\tilde{K}^{2}(1+d)-Z_c)^+, & \text{otherwise}.
    \end{cases}
  \end{equation*}
  Thus the clearing vector for this banks are given as follows:
  \begin{equation}
X_i = \min \bigg \lbrace K^{2}  +  \ \displaystyle\sum_{j\in G_2} X_j w_{ji}-v^2,Y_2 \bigg \rbrace , \forall i \in G_2.
\end{equation}
 where $Y_2= y_b + y_0(1+r_2) +y_{c}(1+r_2)1_{G_1\rightarrow G_2}$.
 
 Dynamics: The dynamics of the big bank can be similarly explained, the big bank total initial wealth is $nk_b$ and out of which  $\alpha$ fraction of wealth is invested towards the loans for the group 2 banks and remaining as the outside investments. Each agents receives $k_b \alpha /(1-\gamma)$ amount of money ( $y_b = k_b \alpha /(1-\gamma)*(1+r_b)$) .
 Thus the fixed point equation  for the big bank is given by as follows:
 
 \begin{equation}
X_b =  \bigg ( K^{b}  +  \ \displaystyle\sum_{j\in G_2} X_j w_{jb}-v^b\bigg )^+.
\end{equation}
with 
\begin{equation*}
    K^{b}=
    \begin{cases}
     ((1-\alpha)k_b(1+u)-Z_c)^+, & \text{wp}\ 1-w  \\
        ((1-\alpha)k_b(1+d)-Z_c)^+, & \text{otherwise}.
    \end{cases}
  \end{equation*}
}
 
 Note: The above quantities are erroneous because the scaling is not considered, below is the precise details of  the  system.
 
 \subsection*{\bf At time $t=0$}

Say a member of $j \in G_2$ borrows totally $y_c$ from all $G_1$ lenders, in particular from $i \in G_1$ it derives:
$$
I_{j, i}
\frac{y_c}{\sum_{i'\in G_1} I_{j, i'} }
$$
Say a member $j \in G_1$ borrows totally $y_1$ from all $G_1$ lenders, in particular  from $i \in G_1$ it derives:
$$
I_{j, i}
\frac{y_1}{\sum_{i'\in G_1 } I_{j, i'} }
$$
Thus total Liability of $j \in G_1$ at $t = 0$ towards all $G_1$ equals:
$$
\sum_{ i \in G_1 }  I_{j, i}  \frac{y_1}{\sum_{i'\in G_1 } I_{j, i'} } = y_1.
$$
This equals the total liability of agent  $j \in G_1$, also equals the total money borrowed at $t= 0.$

The total amount lent by an agent $i \in G_1$ equals:
\begin{eqnarray*}
\sum_{j \in G_1}
 I_{j, i}
\frac{y_1}{\sum_{i'\in G_1 } I_{j, i'} }
+
\sum_{j \in G_2}
 I_{j, i}
\frac{y_c}{\sum_{i'\in G_1} I_{j, i'} } \approx \frac{p_1 \gamma y_1 }{\gamma p_1} +  \frac{p_c (1-\gamma) y_c }{\gamma p_c } =   y_1 + \frac{ (1-\gamma) y_c} {\gamma}
\end{eqnarray*}
So, the money invested by any agent from $G_1$ towards out-side investment equals:
{\small\begin{eqnarray*}
\Omega_1  = {\tilde K}_1 &\approx& k_0 +  y_1 -   y_1 - \frac{ (1-\gamma) y_c} {\gamma}= \bigg(k_0 -  \frac{ (1-\gamma) y_c} {\gamma} \bigg)^+
\end{eqnarray*}}

Total
liability of $j \in G_2$ at $t= 0$ towards all $G_1$ equals
$$
\sum_{ i \in G_1 }  I_{j, i}
\frac{y_c}{\sum_{i'\in G_1} I_{j, i'} } = y_c.
$$

\ignore{
\newpage
Any member from Group 1 (say $i$-agent) lends $(y_1 + y_c)$ totally, so a member of $j \in G_2$ receives (after equal distribution to all small banks)
$$
\frac{(y_1+y_c) I_{j, i }}{\sum_{j'} I_{ j', i}},
$$similar is the case with a member $j \in G_1$. 
Thus the total money lent by an agent $i$ from $G_1$ at $t= 0$ equals:
$$
\sum_{j \in G_2}
\frac{(y_1+y_c) I_{  j,i}}{\sum_{j'} I_{j',i}} +  
\sum_{j \in G_1}
\frac{(y_1+y_c) I_{  j, i}}{\sum_{j'} I_{j', i}} = y_1 + y_c
$$
The liability of agent $i \in G_1$ at $t=0$, 
$$
Y_{1,i} := \sum_{j \in G_1} I_{i,j} \frac{(y_1+y_c)}{\sum_{j'} I_{j',j}} \approx \frac{p_1 \gamma }{ p_1 \gamma + p_c (1-\gamma)}(y_1+y_c) := y_{1,i}
$$

So, the money invested by any agent from $G_1$ towards out-side investment equals:
{\small\begin{eqnarray*}
\Omega_1  = {\tilde K}_1 &\approx& k_0 + \frac{p_1 \gamma }{ p_1 \gamma + p_c (1-\gamma)}(y_1+y_c) - (y_1+y_c) \\
&=
&  k_0 - \frac{p_c (1-\gamma)}{ p_1 \gamma + p_c (1-\gamma)}(y_1+y_c) \\
&
  = & k_0 - \frac{h_c (1-\gamma)}{\gamma},
\end{eqnarray*}}where 
$$
h_c := \frac{p_c \gamma}{ p_1 \gamma + p_c (1-\gamma)}(y_1+y_c)
$$}
 
 An agent   $j \in G_2$ borrows $y_2\eta_j$ from BB and 
 $$
 I_{j, j''}\frac{ y_2 (1-\eta_j) }{
 \sum_{j'\in G_2} I_{j, j'}
 }
 $$from each agent $j''\in G_2$ that is interested in giving loan to agent $j$. Thus the total loan (from BB and $G_2$) taken by agent $j$ at $t= 0$ is given by:
 $$
 y_2 \eta_j + \sum_{j'' \in G_2} I_{j, j''}\frac{ y_2 (1-\eta_j) }{
 \sum_{j'\in G_2} I_{j, j'} } = y_2.
 $$
 Agent $j \in G_2$  also borrows the following   (repeat)
 \ignore{
 $$
 \sum_{i \in G_1}  \frac{I_{i, j}    (y_1+y_c) }{ \sum_{j''} I_{ j'',i}} \approx \frac{p_c \gamma }{p_c (1-\gamma) + p_1 \gamma} (y_1+y_c)  = h_c 
 $$}
 $$
\sum_{ i \in G_1 }  I_{j, i}
\frac{y_c}{\sum_{i'\in G_1} I_{j, i'} } = y_c.
$$from all agents of $G_1$  totally.
 
 The total money lent by $j \in G_2$ is given by:
 $$
 \sum_{j' \in G_2}
  I_{j', j}\frac{ y_2 (1-\eta_{j'}) }{
 \sum_{j''\in G_2} I_{j', j''}
 } \approx y_2 (1-p^{sb}_{2})
 $$
 Thus the investment towards risky regime by $j \in G_2$ is given by:
 $$
 \Omega_2 = {\tilde K}_2 = k_0 + y_2    + y_c - y_2 (1-p^{sb}_{2})
 = k_0 + y_2 p^{sb}_{2} + y_c.
 $$

 \subsection*{\bf At time $t=1$}
 Note that $y_1$, $y_2$ and $y_c$ are time $t= 0$ quantities. 
Liabilities of $j \in G_2$ are given by:
\ignore{
\begin{eqnarray*}
Y_{2,j} &=& \sum_{i \in G_1}  I_{j,i} \frac{y_1+y_c}{ \sum_{j'} I_{ j',i}} (1+r_2) + \sum_{i \in G_2} I_{j,i} \frac{y_2(1-\eta_j)}{\sum_{j'\in G_2 } I_{j, j'}} (1+r_2) +y_2 \eta_j (1+r_b) \\
&\approx & y_{2,j} \mbox{ where } \\
y_{2,j} &:=& \frac{\gamma p_c (y_1+y_c) (1+r_2) }{ p_c (1-\gamma) + p_1  \gamma } + y_2 (1-\eta_j) (1+r_2) +y_2 \eta_j (1+r_b)  \\
&=& h_c (1+r_2) + y_2 (1-\eta_j) (1+r_2) +y_2 \eta_j (1+r_b)
\end{eqnarray*}}
$$
{\bar y}_{2} = y_2 (1+r_2) + y_c (1+r_2).
$$

Thus the clearing vector for agents of  $ G_2$ is given by (for $j \in G_2$):
\begin{eqnarray*}
X_j = \left <  \bigg(K_2 + \sum_{j' \in G_2} X_{j'} I_{j',j} \frac{1- \eta_{j'}}{\sum_{j'' \in G_2} I_{j', j''} } \frac{y_2}{y_2 + y_c}- v_2 \bigg)^+   , \ {\bar y}_{2}  \right >
\end{eqnarray*}
The fixed point approximation  using the Theorem 1, for the same should be
 \begin{eqnarray*}
X_j = \left < \bigg ( K_2 +  {\bar x}_2^\infty (1-p^{sb}_{2}) \frac{y_2}{y_2 + y_c}  - v_2 \bigg)^+  , \ {\bar y}_{2}  \right >
\end{eqnarray*}
 where, ${\bar x}_2^\infty$ satisfies the following fixed point equation:
  \begin{eqnarray*}
  {\bar x}_2^\infty = E \left < \bigg ( K_2 +  {\bar x}_2^\infty (1-p^{sb}_{2}) \frac{y_2}{y_2 + y_c} - v_2 \bigg)^+  , \ {\bar y}_{2}   \right >
  \end{eqnarray*}
  Let ${\bar y}_{1}  = y_1 (1+r_1)$.
The clearing vector for the the $ G_1$ is given by (for $j \in G_1$):
{\small\begin{eqnarray*}
X_j = \left < \bigg( K_1 + \sum_{j' \in G_1} X_{j'} I_{j',j} \frac{1}{\sum_{j'' \in G_1} I_{j', j''} }
+\sum_{j' \in G_2} X_{j'} I_{j',j} \frac{1 }{\sum_{j'' \in G_1} I_{j', j''} } \frac{y_c }{y_2 + y_c} - v_1 \bigg)^+  , \ {\bar y}_1  \right >
\end{eqnarray*}}

The fixed point approximation  using the Theorem 1, for the same should be
 \begin{eqnarray*}
{\bar x}_1^\infty = E \left < \bigg ( K_1 +  {\bar x}_1^\infty +  {\bar x}_2^\infty \frac{y_c }{y_2 + y_c}
\frac{1-\gamma }{\gamma}- v_1 \bigg)^+  , \ {\bar y}_1  \right > 
\end{eqnarray*}

 The clearing money  for the big bank is given by as follows:
 
 \begin{equation}
X_b =  \bigg ( K_{b}  +  \ \displaystyle\sum_{j\in G_2} X_j \frac{\eta_j}{n(1-\gamma)} \frac{y_2}{y_2+y_c}-v_b\bigg )^+.
\end{equation}
Using the approximation we have:

\begin{equation}
X_b^{\infty }= \bigg ( K_{b}  +  {\bar x}_2^\infty p^{sb}_{2} \frac{y_2}{y_2+y_c}-v_b\bigg )^+.
\end{equation}
The invested money towards the outside security  for the big bank is left out money after giving loans to the group 2 banks: $(k_b - y_2p_2^{sb})$. This has again upward and downward movement.\\

The following are the  returns of the  banks  if upward and downward shocks are realized.
\begin{eqnarray*}
k_{d1}&=&\bigg(k_0 -  \frac{ (1-\gamma) y_c} {\gamma} \bigg)^+(1+d - d_c) \\
k_{u1}&=&\bigg(k_0 -  \frac{ (1-\gamma) y_c} {\gamma} \bigg)^+(1+u - d_c) \\
k_{d2}&=&\bigg(k_0 + y_2 p_2^{sb} + y_c \bigg ) (1+d - d_c)\\
k_{u2}&=&\bigg(k_0 + y_2 p_2^{sb} + y_c \bigg ) (1+u - d_c)\\
k_{db}&= &\bigg (k_b - y_2 p_2^{sb} \bigg)(1+d-d_c)\\
k_{ub}&= &\bigg (k_b - y_2 p_2^{sb} \bigg)(1+u-d_c)\\
\end{eqnarray*}

\newpage
\subsection{Performance measures  of systemic risk}
We are interested in probability  of default and the expected surplus

The surplus for the group 1 and banks,   when there is no default for group 1 as well as group 2: 
 
\begin{eqnarray}
E(S_1)&=   E \bigg (K_1 + \bar{y}_2\frac{y_c}{y_2+y_c}\frac{1-\gamma}{\gamma} -v_1 \bigg)^+ \\
 &= E(K_1) - v_1 +y_c(1+r_2)\frac{1-\gamma}{\gamma}
\end{eqnarray} \label{Eqn_surplus 1}
The expected surplus for the group 2 banks are given by :
\begin{eqnarray}
E(S_2)&=   E \bigg (K_2 +  {\bar x}_2^\infty \frac{y_2}{y_2+y_c}(1-p^{sb}_2) -v_2 - \bar{y}_2\bigg )^+ \\
 &= E(K_2) - v_2 -(1+r_2)\bigg(y_c+ y_2 p^{sb}_2 \bigg)
\end{eqnarray} \label{Eqn_surplus 2}
}
\subsection{Analysis of the limit system}
\label{Subsec_analysis_limit_system}
We obtain the performance measures of the financial network  by analyzing the simplified limit system derived in the above. %We analyze the system for a given realization of the random shocks. 
We begin with few  more notations (for any $m$): 
\begin{eqnarray*}
 k_{dm} := \Omega_m(1+d-d_c), \
k_{um} := \Omega_m(1+u-d_c),\mbox{ and, }
\bar{l}_m:= wk_{dm}+(1-w)k_{um}~. 
\end{eqnarray*}
We first  derive  ${\bar x}_2^\infty $ and $P^2_D$, the aggregate clearing vector and the default probability of ${\cal G}_2$, for a given set of  system parameters in the following: 
\begin{lemma}
\label{Lemma_G2 default}
Consider $k_{d2} > v_2$.
There is a unique solution to \eqref{eqn_avgclearingvectorG2} and
the asymptotic aggregate
clearing vector and the  default probability of   $\mathcal{G}_2$   is  given by:

\vspace{-4mm}
{\small
\begin{eqnarray}
\hspace{-4mm}
\mbox{{\small$( {\bar x}_2^\infty, P_D^2)$}}  &=& \left \{ \begin{array}{lllll}
	\bigg(\bar{y}_2 (1-p_2^{sb})\lambda_2) ,  \ 0\bigg)  & \mbox{ if }    \lambda_2  \ge  \frac{\bar{y}_2 +v_2-k_{d2} }{\bar{y}_2 (1-p_2^{sb})} \\
	\bigg( \bigg(\bar{y}_2- \frac{  (\bar{y}_2-k_{d2}+v_2)w -\bar{y}_2  (1 - p_2^{sb})w\lambda_2  }{ 1 -w (1-p_2^{sb})\lambda_2  } \bigg)(1-p_2^{sb})\lambda_2 , \ w  \bigg)   & \mbox{ if } \beta_0 <  \lambda_2
	< \frac{\bar{y}_2+v_2-k_{d2}}{\bar{y}_2 (1-p_2^{sb})} \nonumber
	\\ 
\bigg (\frac{ ( \bar{l}_2 -v_2)^+}{ 1- (1-p_2^{sb})\lambda_2 }(1-p_2^{sb})\lambda_2, \ 1 \bigg) & \mbox {  if }
                     \lambda_2<   \beta_0,  \hspace{6mm} \mbox{\normalsize with, }
\end{array}
\right . \ \label{Eqn_PDSB2} \\ \nonumber 
\beta_0 &:=& \frac{  (\bar{y}_2 +v_2-k_{u2})
} { \big (\bar{y}_2-w(k_{u2}-k_{d2}) \big ) \left  (1-p_2^{sb} \right ) } 
\mbox{ $1_{{\bar y}_2  > w(k_{u2}-k_{d2}) }$}.
\end{eqnarray}}
 %where,{\small\begin{eqnarray*}E[L] & := &wk_{d2}+(1-w)k_{u2}.
% H(y_2,y_c,p_2^{sb})&:= &(1- p^{sb}_{2})\frac{y_2}{y_2+y_c}.
 %\end{eqnarray*}}
 \end{lemma}
 \textbf{Proof:}  available in Appendix E.
 \eop

 \begin{lemma}
\label{Lemma_G2 default with v2> kd}
Consider $v_2> k_{d2}$ and   $\bar{y}_2  > w(k_{u2}-k_{d2})$.
There is a unique solution to \eqref{eqn_avgclearingvectorG2} and $( {\bar x}_2^\infty, P_D^2)$ are    given by:
\vspace{-4mm}

\begin{eqnarray}
\hspace{-4mm}
\mbox{{\small$( {\bar x}_2^\infty, P_D^2)$}}  &=& \left \{ \begin{array}{cllll}
	\bigg(\bar{y}_2(1-w) (1-p_2^{sb})\lambda_2) ,  \ w\bigg)  & \mbox{ if }  \beta_4  < \lambda_2 \leq \beta_1  \\
		\bigg (\frac{(k_{u2}-v_2)^+ (1-w)(1-p^{sb}_2)\lambda_2}{1-(1-p^{sb}_2)\lambda_2(1-w)}, \ 1 \bigg) & \mbox {  if }
                     \lambda_2 < \min  \bigg \lbrace \beta_4,  \beta_3 \bigg\rbrace\\ 
                     \\
	\bigg( \bigg(\frac{(k_{d2}-v_2)w +\bar{y}_2(1-w) }{1-(1-p^{sb}_2)w\lambda_2} \bigg)(1-p^{sb}_2)\lambda_2, \ w  \bigg)   & \mbox{ if } \lambda_2 > \max \bigg \lbrace \beta_2, \beta_1 \bigg \rbrace \nonumber
	\\ 
\bigg (\frac{(\bar{l}_2-v_2)^+(1-p^{sb}_2)\lambda_2}{1-(1-p^{sb}_2)\lambda_2},   \ 1 \bigg) & \mbox {  if }
                   \beta_3 < \lambda_2   <  \beta_2
                       \hspace{6mm} \mbox{\normalsize with } ,
\end{array}
\right . \nonumber 
\end{eqnarray}

\vspace{-4mm}
\begin{align*} 
\beta_ 1 &:=\frac{v_2- k_{d2}}{\bar{y} _2(1-w)(1-p_2^{sb})}, & \beta_2 &:= \frac{  (\bar{y}_2 +v_2 -k_{u2})
} { \big (\bar{y}_2-w(k_{u2}-k_{d2})  \big ) \left  (1-p_2^{sb} \right ) } ,  \\  
\beta_3&:=\frac{v_2-k_{d2}}{(1-w)(k_{u2}-k_{d2})(1-p^{sb}_2)} \mbox{and} ,
& \beta_4 &:= \frac{\bar{y}_2 - k_{u2}+v_2}{\bar{y}_2(1-w)(1-p_2^{sb})}.
\end{align*} 

\end{lemma}
\textbf{Proof:}  available in Appendix E.
\eop

One can derive closed form expressions for the remaining case also, but the expressions could be more complicated; for such cases,  numerically solving the fixed point equations \eqref{eqn_avgclearingvectorG2} for limit system is not complicated and the same is considered in sub-section \ref{sub_section_numerical}  for some numerical examples.    We now analyze the aggregate clearing vector and the default probability of the $\mathcal{G}_1$ banks.
 \begin{lemma}
 \label{Lemma_G1 default}Assume $k_{d1}-v_1 + \mu_1 {\bar x}_2^\infty \ge 0$. 
Given the (unique) asymptotic  aggregate clearing vector ${\bar x}_2^\infty$ for   $\mathcal{G}_2$, the 
fixed point equation \eqref{eqn_avgclearingvectorG1} has a unique solution; 
the aggregate clearing vector and the  default probability of $\mathcal{G}_1$ is  given by
 \begin{eqnarray}
\hspace{-4mm}
\mbox{$( {\bar x}_1^\infty, P_D^1)$}  = \left \{ \begin{array}{cllll}
 	\bigg(\bar{y}_1 (1-p_1^{sb}) ,  \ 0\bigg)  & \mbox{ if } \mu_1 {\bar x}_2^\infty \ge  e_1 \\ \\
	\bigg( \frac{ \bar{y}_1  (1 - w )   + w(k_{d1}-v_1+\mu_1 {\bar x}_2^\infty)}{ 1 -w(1-p_1^{sb})  }(1-p_1^{sb})  , \  w  \bigg)   & \mbox{ if }  e_2 \le \mu_1 {\bar x}_2^\infty
	<  e_1 \nonumber
	\\ 
	
	%(0,1) & \mbox{if} \
	%\beta < v_1-k_{d1}

\bigg ( \frac{(\bar{l}_1 -v_1 + \mu_1 {\bar x}_2^\infty )(1-p_1^{sb})}{p_1^{sb}}  , \ 1 \bigg) & \mbox { if }
                   \mu_1 {\bar x}_2^\infty  \le  e_2,
\end{array}
\right .
\end{eqnarray}
where $\mu_1$ is in \eqref{Eqn_mu1mu2},  $
   e_1 := v_1-k_{d1} +\bar{y}_1p_1^{sb}, $ and, $
e_2 := v_1-\bar{l}_1 +p_1^{sb}  (\bar{y}_1+ w(k_{d1} -k_{u1})).
$
\end{lemma}
\textbf{Proof:}  available in Appendix E.
\eop
\\
Like before, one can derive fixed points \eqref{eqn_avgclearingvectorG1} even for other cases; the expressions can be more complicated, it is rather   easier to solve the fixed point equations of limit system;  this is considered in sub-section \ref{sub_section_numerical}.

%\noindent\textbf{Remarks:} 
%From the above three lemmas, it is easy to derive the   performance measures given in sub-section \ref{Subsec_analysis_limit_system},  for any given set of  system parameters.  One can estimate  the aggregate clearing vector for $\mathcal{G}_1$  using Lemma \ref{Lemma_G1 default}, once that for  $\mathcal{G}_2$ is known, using the same approach.   

%Another important aspect is that one could also consider the performance analysis of the BB to study the effect of the entire network on BB.

%{\color{blue}\noindent{\textbf Note:} In the systemic risk event regime, our model does not satisfies the   assumptions of Corollary \ref{corolary_limit_system_uniquness} i.e.,  $p^{sb}_1 > 0$ or $\sigma < 1$ is not satisfying and hence there is no fixed point solution in that regime. This is evident from the Banach fixed point Theorem; it is not strict contraction. We avoid this regime in the financial model because the failure of all the entities in the group $\mathcal{G}_1$ is not an interesting regime.}

We now consider an interesting sub-case and derive some more analysis related to  the network. 

\subsubsection*{Taxes proportional to investments} From now on, we assume that the taxes are proportional to risky investments, i.e., $v \propto \Omega$, in particular we assume 
$v = \kappa \Omega$ for some $\kappa >0$.
Thus, $v_m = \kappa \Omega_m   $ for any agent from group ${\cal G}_m.$ This is a natural assumption.
%For notational  convenience we define another constant, $\kappa := \kappa' / (1+u-d_c)$, i.e., 
%$v_m = \kappa (1+u-d_c)  \Omega_m.$
For this sub-case we have some more interesting observations,  we begin with some definitions followed by an interesting property: 
\begin{defn}
{\bf Resilient regime}: A group of banks is said to be in resilient regime if none of them  default (pay back their liabilities completely),   irrespective of the economic shocks that they receive.  
The financial system is said to be in  resilient regime, if all its groups are in resilient regime.  
\end{defn}
\begin{defn}
{\bf Systemic risk regime:} A financial system is said to be in   systemic risk regime when  the local shocks   trigger cascade of defaults and   all the agents default, i.e.,  the entire system collapses.
\end{defn}

\begin{lemma}
\label{Lemma_default}{\bf[$\mathcal{G}_1$  is more robust]}
Assume   proportional taxes, i.e., $v =  \kappa \Omega $ and $y_1p_1^{sb} < y_2p_2^{sb}$. Then,    if the  $\mathcal{G}_2$ banks are resilient, so are the $\mathcal{G}_1$ banks. 
\end{lemma}
\textbf{Proof:}  is available in Appendix E.
\eop 

The condition $y_1p_1^{sb} < y_2p_2^{sb}$ further implies that  ${\cal G}_1$  borrows lesser from   BB and hence invests even lesser in risky investments. Under this condition, when ${\cal G}_2$ banks are resilient, so are ${\cal G}_1$ banks.
Further more, $e_2$ given in Lemma \ref{Lemma_G1 default} is  usually a small value (as usually taxes are less than the expected returns from risky investments, i.e.,  $v_1 < {\bar l}_1$ and $p_1^{sb}$ is typically  a small value) and hence $P_D^1 \le w$, however the ${\cal G}_2$ banks can enter into ``Systemic risk regime".  

%\noindent\textbf{Remarks:} The above lemma highlights  that  the  $\mathcal{G}_1$ banks are more robust towards economic shock; the only possible default patterns are following:
%$$
%(P_D^1, P_D^2) \in \{ (0,0),(0,w),(0,1),(w,w),(w,1),(1,1) \}.
%$$

%The condition required for the  second part of the lemma is quite general, observe that the negation of it would imply even the expected increase ($1+{\bar r}_r$) with risky investments is smaller  than the common shock plus the rate of taxes $\kappa$. One can always  work the details for this case also, but we do not consider this as it is not a practically relevant scenario. 

We now discuss the trends of the  performance measures of the two groups as a function of the inter-lending parameter $y_c$.

\begin{lemma}
\label{lemma_surplus1}
Assume   proportional taxes, i.e., $v =  \kappa \Omega $ and   let ${\bar r}_r := u(1-w) +dw$,  
 $\Delta_u := (1+u) - (1+r_2)$ and $\Delta_r := (1+\bar{r}_r) - (1+r_2)$. Under the resilient regime,
%The expected surplus of  $\mathcal{G}_1$  under the resilient regime is :
\begin{enumerate}[a)]
    \item  the expected surplus $E[S_1] \geq E[S_2]$  if  and only if  $ \Delta_r  \leq  d_c+ \kappa$,
    \item the expected surplus ($E[S_1]$) of  $\mathcal{G}_1$  increases with inter lending amount $y_c$ if  and only if  
    $\Delta_r < d_c + \kappa$,   and remains constant when $\Delta_r = d_c + \kappa $; and; 
     \item the surplus at upward movement SaU, $\hat{S}_{2,u}$   of $\mathcal{G}_2$, increases with $y_c$ if  and only if $\Delta_u > d_c +\kappa$, remains unaltered when $\Delta_u = d_c +\kappa$.  
     \ignore{
   \item if further $ {\bar r}_r<  \re  <u $
    then,  SaU
$\hat{S}_{2,u}$   of $\mathcal{G}_2$  and the expected surplus $E[S_1]$ of $\mathcal{G}_1$ are both  increasing with inter lending amount $y_c$ .}
\end{enumerate}
\end{lemma}
\textbf{Proof:}  available in Appendix E.
\eop

\noindent\textbf{Remarks:}  
If the banks operate in  resilient regime, i.e., even the banks with economic shocks manage to clear their liabilities completely, then   the expected surplus of the first group is larger, only when the expected return rate from risky investment ($\Delta_r$) is smaller than the 
`burden factor' $d_c +\kappa$.

The SaU of the $\mathcal{G}_2$ banks and $\mathcal{G}_1$ banks expected surplus   increases with the inter lending parameter $y_c$ in the  regime $\Delta_r< d_c +\kappa < \Delta_u$.
{\it In such scenarios,  it is beneficial for both groups to increase the inter-lending amount $y_c$}.  
%The $\mathcal{G}_2$  banks portfolio comprised of risky (such as derivatives, stocks, etc.) instruments, while as opposed to that the $\mathcal{G}_1$ banks portfolio consist of the less risky instrument. Therefore, SaU for $\mathcal{G}_2$ is more crucial than the expected surplus of $\mathcal{G}_1$ banks.

It would be interesting to study similar aspects in default regime, we consider the same using numerical computations in the next sub-section.

\ignore
{\color{red}
\noindent\textbf{Remarks:}
\label{fig_remark}
The conclusions of  Lemma \ref{lemma_surplus1}  hold  true  for the  default regime   i.e., when    $P^{m}_D > 0$ if the probability of economic shock is sufficiently small  $0 \leq w \le \bar{w}$, for some $\bar{w} < 1$.
First consider the $\mathcal{G}_2$ banks aggregate clearing vector in the default regime. Define, $\delta_0 : =  (1-p_2^{sb})\lambda_2$ and it is immediate that  $0 < \delta_0 < 1$. For all $w \le {\bar w}$, from Lemma \ref{Lemma_G2 default}  we have the following:
\begin{eqnarray*}
& & \displaystyle{\sup _{\delta_0}}|{\bar x}_2^\infty - \bar{y}_2 \delta_0|\\
& < &\displaystyle {\sup _{\delta_0}} \bigg|  \frac{\bar{y}_2(1-w) + w(k_{d2}-v_2)}{1-w\delta_0}- \bar{y}_2 \bigg|  \\
&  = & \displaystyle{\sup_{\delta_0}} ~  \bigg |\frac{\bar{y_2}(\delta_0-1) +w\delta_0(k_{d2}-v_2)}{1-w\delta_0}\bigg| \\
&< & \displaystyle{\sup _{\delta_0}} \bigg |\frac{\bar{y_2}(\delta_0-1) +w\delta_0(k_{d2}-v_2)}{1-w}\bigg|  < \frac{w}{1-w}\eta_0
\end{eqnarray*}
for some $\eta_0 > 0$, which decreases to 0 as  $w \to 0$ (the first and last two inequality is due to $\delta_0 < 1$ and then taking  supremum over  $\delta_0$). 
By continuity of  the above upper bound with respect to $(1-w)$  
and, $\exists$  a ${\bar w} > 0 $ (further small if required) such that  the above arguments are true for all $w \le {\bar w}$.
 We can repeat  the same argument for the aggregate clearing vector of the $\mathcal{G}_1$ banks.
 
 {\color{red}
 The system defaults with $v_2 < k_{d2}$ when the following is satisfied
 $$ k_0 (1+d-d_c-\kappa) < (y_c + y_2 p^{sb}_2) (\re - d) $$}
 }

 \subsection{Numerical observations}
 \label{sub_section_numerical}
As discussed before, by Theorem \ref{Thm_MainGen 2} and Corollary \ref{Cor_Three dimensional approximation_Alt_system},   a large banking network can be well approximated  by an appropriate  limit system almost surely.  For large networks, it is complicated to derive the performance directly, one can rather use the limit system. 
In the next section, we reaffirm the accuracy of this approximation,  using Monte-Carlo simulation-based results; in this sub-section
we obtain some interesting performance trends  using   the limit system (given by \eqref{eqn_avgclearingvectorG1}-\eqref{Eqn_Finance_Conv}).   Our key objective is to analyze the role of the inter lending parameter $y_c$ and its feedback effect in the network;  we numerically solve the limit fixed point equations to study the trends of probability of default and surplus   based measures (for different groups)   with the inter-lending parameter.

\begin{figure}[hbt!]
\begin{subfigure}{\textwidth}
\begin{subfigure}{.5\textwidth}
  \centering
  % include first image
  \includegraphics[width=.6\linewidth]{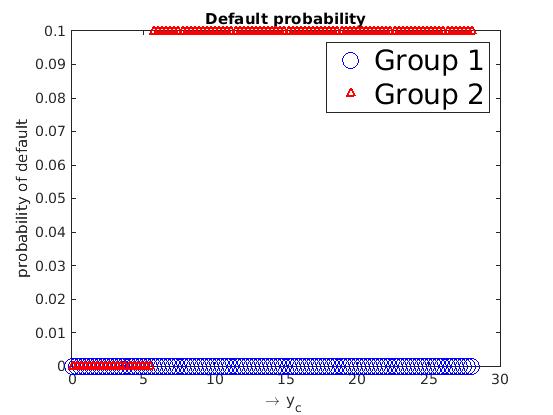}  
 % \caption{ $u =0.5$ , $d= -0.35$, $d_c =0.1$, $\kappa = %0.175$}
 % \label{fig:PD_6}
\end{subfigure}
\begin{subfigure}{.5\textwidth}
  \centering
  % include second image
  \includegraphics[width=.6\linewidth]{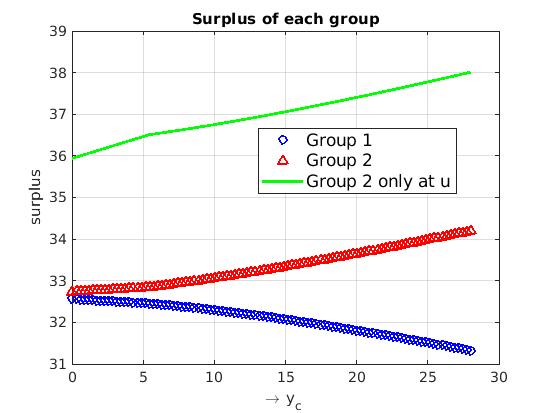}  
 
\end{subfigure}
 \caption{\small $u =.5$ , $d= -.35$, $d_c =.1$, $\kappa = .175$, with, $d_c+\kappa=.275< \Delta_r =.295< \Delta_u=.38$
  \label{fig:ES_6} }
  \end{subfigure}

   \begin{subfigure}{\textwidth}
\begin{subfigure}{.5\textwidth}
  \centering
  % include first image
  \includegraphics[width=.6\linewidth]{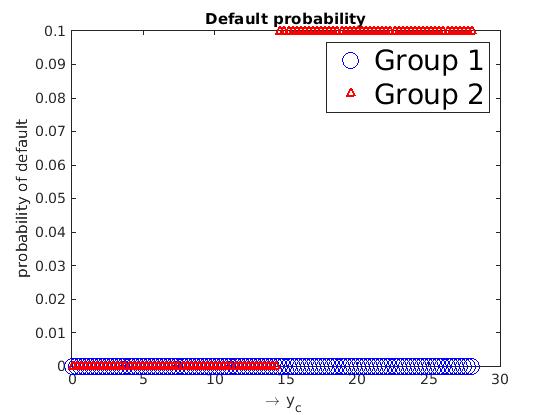}  
  %\caption{ $u =.5$ , $d= -.15$, $d_c =.1$, $\kappa = .245$}
  %\label{fig:PD_8}
\end{subfigure}
\begin{subfigure}{.5\textwidth}
  \centering
  % include second image
  \includegraphics[width=.6\linewidth]{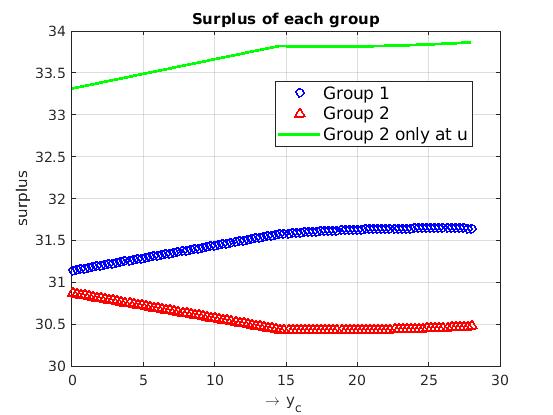}  
 
\end{subfigure}
 \caption{\small $u =.5$ , $d= -.15$, $d_c =.1$, $\kappa = .245$, with, $\Delta_r = .315 <  d_c+\kappa= .345 < \Delta_u =.38$  \label{fig:ES_8}}
 \end{subfigure}

  \begin{subfigure}{\textwidth}
\begin{subfigure}{.5\textwidth}
  \centering
  % include first image
  \includegraphics[width=.6\linewidth]{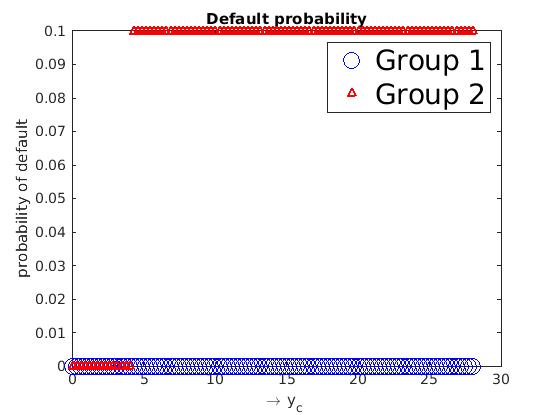}  
  %\caption{ $u =.5$ , $d= -.2$, $d_c =.1$, $\kappa = .35$}
  %\label{fig:PD_7}
\end{subfigure}
\begin{subfigure}{.5\textwidth}
  \centering
  % include second image
  \includegraphics[width=.6\linewidth]{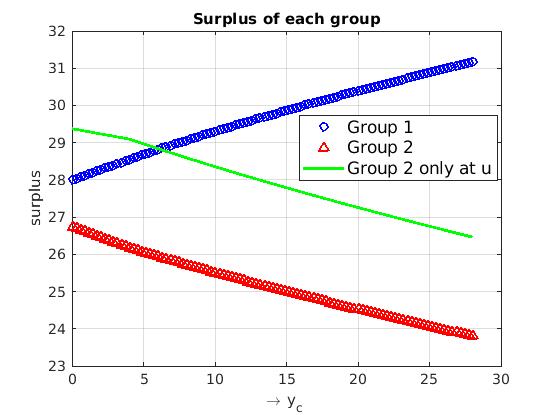}  
  \end{subfigure}
  \caption{ \small $u =.5$ , $d= -.2$, $d_c =.1$, $\kappa = .35$, with, $\Delta_r = .31 < \Delta_u= .38 < d_c+\kappa= .45$ \label{fig:ES_7}}
  \end{subfigure}
  
\begin{subfigure}{\textwidth}
\begin{subfigure}{.5\textwidth}
  \centering
  % include first image
  \includegraphics[width=.6\linewidth]{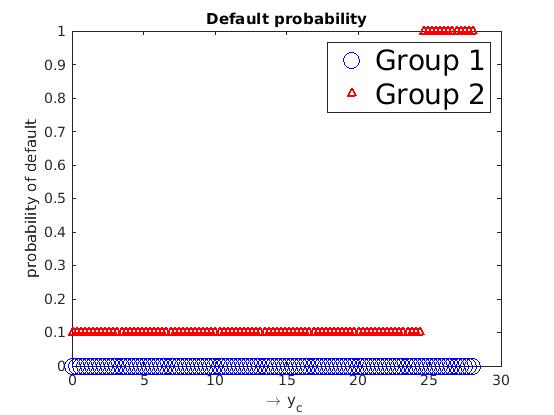}  
  %\caption{ $u =.5$ , $d= -.01$, $d_c =.1$, $\kappa = .77$}
  %\label{fig:PD_9}
\end{subfigure}
\begin{subfigure}{.5\textwidth}
  \centering
  % include second image
  \includegraphics[width=.6\linewidth]{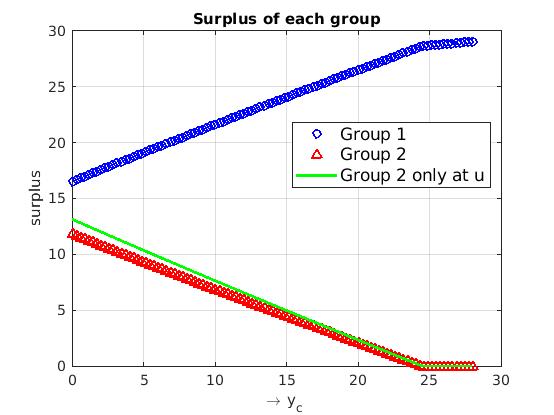}  
 \end{subfigure}
  \caption{\small $u =.5$ , $d= -.01$, $d_c =.1$, $\kappa = .77$, with, $\Delta_r = .329 < \Delta_u = .38 < d_c + \kappa = .87$
  \label{fig:ES_9}}
   \end{subfigure}
\caption{Small shock regime: $w=.1$, $r_1= .1$, $r_2=.12$, $p^{sb}_2 =.2$, $p^{sb}_1 =.01$.}
\label{fig:2}
\vspace{-1mm}
\end{figure}
 
 We have used the following common set of  parameters for our numerical examples:
 $k_0 = 40/(1+u-d_c)$,
$k_b = 50/(1+u-d_c)$, $y_1 = 50/(1+r_1)$, $y_2 = 50/(1+r_2)$, $\bar{y} _1 = y_1(1+r_1)$, $\bar{y}_2 = (y_2+y_c)(1+r_2)$, $\gamma =0.5$ and the rest of the system parameters are given in the captions of the respective figures.
In Figure \ref{fig:2}  (and its sub-figures),  we consider small shock regime ($v_{2} < k_{d2}$), while, Figure \ref{fig:1}  studies  the large shock regime ($v_2> k_{d2}$).

%Based on simulations and theoretical results given by Lemmas \ref{Lemma_default} and \ref{lemma_surplus1}  and  we have the following observations for surplus:

\noindent{\bf Moderate and small shock regime ($v_2 <  k_{d2}$):}
Lemma \ref{lemma_surplus1} characterizes the trends in the performance under resilient regime (which is possible only under small-shock regime); and interestingly the trends continue even in default regime (when $P^{2}_D > 0$) for the case with small-shocks:
\vspace{-2mm}
\begin{figure}[hbt!]
  %\vspace{-1mm}
  \begin{subfigure}{\textwidth}
\begin{subfigure}{.5\textwidth}
  \centering
  % include first image
  \includegraphics[width=.6\linewidth]{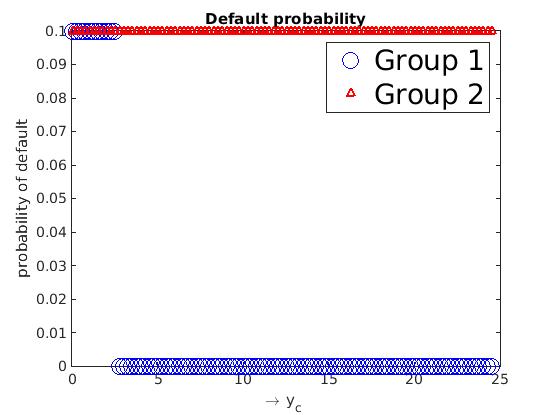}  
  %\caption{ $u =.7$ , $d= -.7$, $d_c =.1$, $\kappa = .32$}
  %\label{fig:PD_3}
\end{subfigure}
\begin{subfigure}{.5\textwidth}
  \centering
  % include second image
  \includegraphics[width=.6\linewidth]{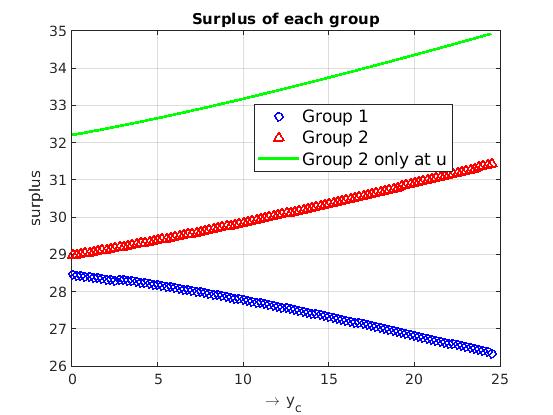}  
  \end{subfigure}
  \caption{ \small $u =.7$ , $d= -.7$, $d_c =.1$, $\kappa = .32$, with, $d_c +\kappa=.42 < \Delta_r=.44 <\Delta_u=.58$
  \label{fig:ES_3}}
   \end{subfigure}

\begin{subfigure}{\textwidth}
\begin{subfigure}{.5\textwidth}
  \centering
  % include first image
  \includegraphics[width=.6\linewidth]{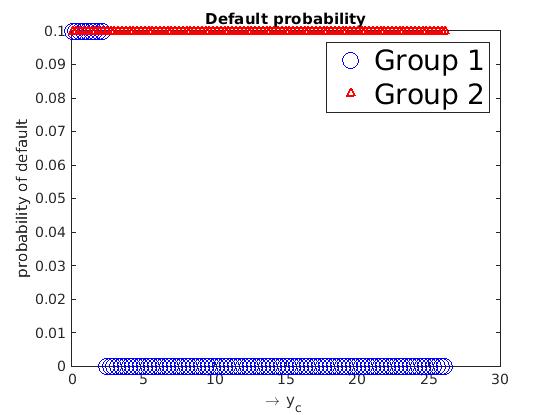}  
  %\caption{ $u =.6$ , $d= -.7$, $d_c =.1$, $\kappa = .3$}
  %\label{fig:PD_4}
\end{subfigure}
\begin{subfigure}{.5\textwidth}
  \centering
  % include second image
  \includegraphics[width=.6\linewidth]{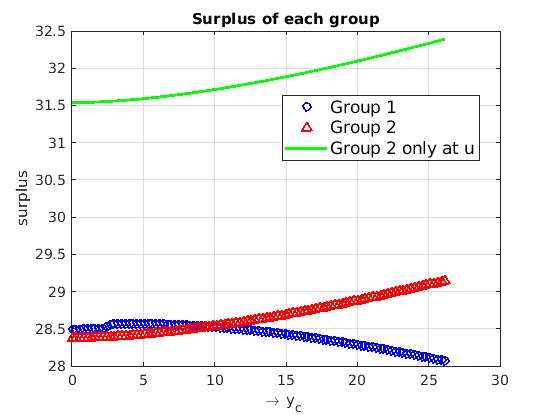}  
 \end{subfigure}
 \caption{\small $u =.6$ , $d= -.7$, $d_c =.1$, $\kappa = .3$, with, $\Delta_r=.35 < d_c +\kappa=.4< \Delta_u = .48$
  \label{fig:ES_4}}
   \end{subfigure}

  \begin{subfigure}{\textwidth}

\begin{subfigure}{.5\textwidth}
  \centering
  % include first image
  \includegraphics[width=.6\linewidth]{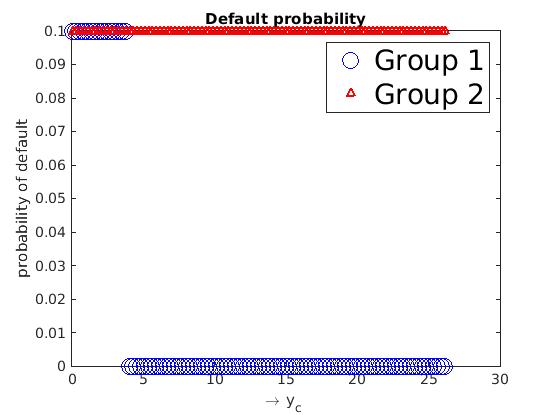}  
  %\caption{ $u =.6$ , $d= -.7$, $d_c =.1$, $\kappa = .375$}
  %\label{fig:PD_2}
\end{subfigure}
\begin{subfigure}{.5\textwidth}
  \centering
  % include second image
  \includegraphics[width=.6\linewidth]{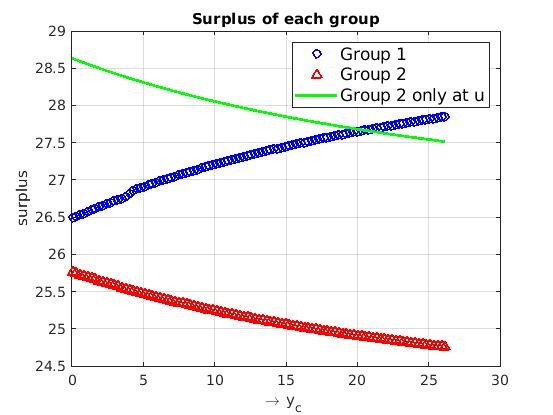}  
 \end{subfigure}
 \caption{\small$u =.6$ , $d= -.7$, $d_c =.1$, $\kappa = .375$, with,  $\Delta_r=.35 < d_c +\kappa=.475< \Delta_u = .48$
  \label{fig:ES_2}}
   \end{subfigure}
 \begin{subfigure}{\textwidth}
\begin{subfigure}{.5\textwidth}
  \centering
  % include first image
  \includegraphics[width=.6\linewidth]{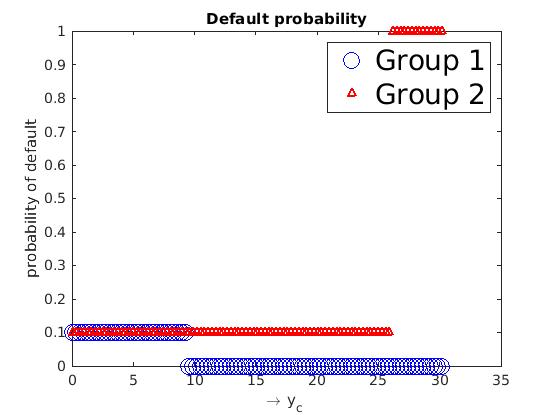}  
  %\caption{ $u =.4$ , $d= -.7$, $d_c =.1$, $\kappa = 0.65$}
  %\label{fig:PD_1}
\end{subfigure}
\begin{subfigure}{.5\textwidth}
  \centering
  % include second image
  \includegraphics[width=.6\linewidth]{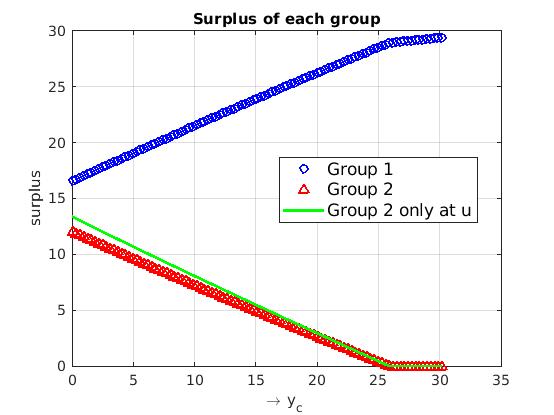}  
 \end{subfigure}
 \caption{\small $u =.4$ , $d= -.7$, $d_c =.1$, $\kappa = 0.65$, with, $\Delta_r=.17 < \Delta _u=.28 < d_c +\kappa=.75$
  \label{fig:ES_1}}
   \end{subfigure}
 
\caption{Large shock regime: $w=.1$, $r_1= .1$, $r_2=.12$, $p^{sb}_2 =.2$, $p^{sb}_1 =.01$.}
\label{fig:1}
\vspace{-4mm}
\end{figure}

%with the following further observations: 
 %  
  %
%
%  The  following regime for v2 < k_d2
%% end here with the regime v_2 < k_d2
%  In this regime the observations (even in default regime) majorly follow the results of the two Lemmas  \ref{Lemma_default} and  \ref{lemma_surplus1}:
\begin{enumerate}[a)]
 \item  When, $\Delta_r$, the difference in the expected rate of return from risky assets and  the rate  $r_2 $ of liabilities of group ${\cal G}_2$
is smaller than the system `burden factor' $d_c + \kappa$, 
%i.e., if $\Delta_r < d_c +\kappa$
then the expected surplus of ${\cal G}_1$ banks improves with increase in $y_c$ (see figures (\ref{fig:ES_8}), (\ref{fig:ES_7}), (\ref{fig:ES_9})). Moreover, this trend   continues in   default regime also.

\item  On the other hand, if $\Delta_u$ (the  difference    between upward rate and $r_2$) is bigger than $d_c + \kappa$,  then SaU of group ${\cal G}_2$  improves with $y_c$ as seen in figures (\ref{fig:ES_6}) and  (\ref{fig:ES_8}). This trend also continues in the default regime.  
\item For the case study of figure (\ref{fig:ES_8}), both the groups improve;  SaU as well as $E[S_1]$ increase with $y_c$, even in default regime. % the trend is opposite in the default regime as for example see in the . Therefore we observed  non-monotone trends in expected surplus in the  default regime. %this quantity is influenced by percolation of shocks.

%But the above is true only in resilience regime, while in default regime one of the  performance measure continue to increase while the other decreases (see in  Figures  (\ref{fig:ES_6})  and (\ref{fig:ES_8})).

\item When $\Delta_r > d_c+\kappa$,   $E[S_1]$   decreases, while, $E[S_2]$ as well as the SaU of ${\cal G}_2$ banks  improves with $y_c$ (figure (\ref{fig:ES_6})). 
While with  $d_c+\kappa > \Delta_u$, only group 1 improves (figures (\ref{fig:ES_7}), (\ref{fig:ES_9})). These trends also continue   (even with   $P_D^2 = 1$). 

%This trend also continues   in default regime, as in the previous cases. 

%\item In figures (\ref{fig:ES_7}), (\ref{fig:ES_9}) with $d_c+\kappa > \Delta_u$, only group 1 improves and the trend continues even when   $P_D^2 = 1$. 

%\item In the first sub-region i.e., when  $\Delta_r < d_c+ \kappa$, both the groups can simultaneously improve their respective utilities by inter lending,  we showed a sub-region of $y_c >0$, in which both groups do better than without inter-lending (see Figure (\ref{fig:ES_8})).
\end{enumerate}

\noindent{\bf Large shock regime ($v_2 > k_{d2}$):} 
This regime is considered in Figure \ref{fig:1}. 
The observations are almost similar to that in the previous case, except for the switch-over points: 
a) when $d_c+\kappa$ is small, only group 2 banks benefit; b) with  larger burden factor, both the groups   benefit; c)  when the burden factor is increased further, only group 1 banks benefit with increase in $y_c$; and d) the switch over points of $d_c+\kappa$ for  above three types of regimes are given in Lemma \ref{lemma_surplus1},   are valid even in default regime with small shocks; however e) the switch-over points can be different with larger shocks (see figure  (\ref{fig:ES_2})).

\ignore{
% finishes here
 \begin{enumerate}
    % \item  The group  $\mathcal{G}_2$ always defaults with non-zero probability  i.e.,  $P^2_{D} \ge w$.  
     \item In the large shock regime with $\Delta_r <  d_c +\kappa  <  \Delta_   u$  the expected surplus $E[S_1]$ of the $\mathcal{G}_1$ banks decreases with $y_c$ while $\mathcal{G}_2$ banks  expected surplus  $E[S_2]$ and SaU increases with  $y_c$. We plotted such scenario  in the  Figure  (\ref{fig:ES_4}). However, when the difference between any two of quantities  $\Delta_u$, $d_c+\kappa$, $\Delta_r$ are negligible then 
      we have `opposite' trend for the expected surplus of both the groups. For instance 
     expected surplus $E[S_1]$ of the $\mathcal{G}_1$ banks increases while $\mathcal{G}_2$ banks  expected surplus $E[S_2]$ and $\hat{S}_{2,u}$ decreases with  $y_c$  see in Figure (\ref{fig:ES_2}).
  \item Also if we consider  when $ \Delta_r < \Delta_u < d_c +\kappa$,  expected surplus of the $\mathcal{G}_1$ banks increases  with $y_c$ while  $\mathcal{G}_2$ banks expected surplus decreases with the inter lending parameter $y_c$. We plotted this in the Figure  (\ref{fig:ES_1}).
  \item When we consider the regime $\Delta_r > d_c+\kappa$ expected surplus $E[S_1]$ decreases with $y_c$ and surplus at the upward return, expected surplus  of the $\mathcal{G}_2$ banks increases with the inter lending parameter $y_c$. We observed this in the Figure (\ref{fig:ES_3}).
  \item The inter lending parameter $y_c$ can't simultaneously help both the groups. It is evident from the  Figures  (\ref{fig:ES_3})- (\ref{fig:ES_1}), and in all those scenarios one of the performance measure increases while other strictly decreases with  $y_c$. Therefore in the large shock regime,  inter-lending parameter plays an adverse effect. Hence it is beneficial for both the groups of banks to become isolated.
 \end{enumerate}
 }

\ignore{
 Now by the Lemmas \ref{Lemma_G2 default} - \ref{Lemma_G1 default} and numerical simulations,  under the range of small shock regime ($v_2 <  k_{d2}$) to the large shock regime  ($v_2 <  k_{d2}$) we have the following observations of the default probability: 
\begin{enumerate}
\item As observed in Lemma \ref{Lemma_default}, the $\mathcal{G}_1$ banks are   resilient, once   $\mathcal{G}_2$ banks are resilient.  Further from all the figures, $\mathcal{G}_2$ has more defaults, i.e.,   $P_D^{1} \le P_D^{2}$.  
\item For large-shocks (when the common shock and operational costs (taxes or deposit) are high),  $\mathcal{G}_2$ banks are  prone to enter `Systemic risk regime'.  But in none of the figures  $\mathcal{G}_1$ has entered  the `Systemic risk regime' ($P_D^1 = 1$). On the contrary, $\mathcal{G}_2$ enjoys better SaU (and some times better expected surplus) than $\mathcal{G}_1$.  However  the surplus based measures increase for initial values of $y_c$ and decrease after certain value and the trend is similar to the one observed in \cite{Systemicrisk}. 
\item This is because of the liability of  $\mathcal{G}_2$ banks increase with an increase of $y_c$. In contrast,   $\mathcal{G}_1$ banks are more robust to economic shocks. Initially, with large shocks, $\mathcal{G}_1$ banks tend to default, but as the inter-lending parameter increases, $\mathcal{G}_1$ banks manage to become resilient. Such occurs because  $\mathcal{G}_2$ banks clear the liability as much as possible, and as a result, it improves the liquidity of  $\mathcal{G}_1$ banks. We plotted such scenario in Figure (\ref{fig:ES_1}).
\item We also have $\mathcal{G}_1$ banks are more `resilient' to the downward movement and the common shock.  At the same time,  $\mathcal{G}_2$ banks are more sensitive towards common shock (as for example, see Figure (\ref{fig:ES_9})).
\end{enumerate}
 }
\section{ Monte Carlo Simulations}
\label{sec_MCsimulation}
In the previous sections, we derived asymptotic performance and  systemic-risk analysis of a financial network using the fixed-point convergence 
theorems of Section \ref{Graphicalmodel} and  \ref{Alternate_graphical_model}.
This section reinforces the approximation demonstrated by Theorem \ref{Thm_MainGen 2}, using exhaustive Monte-Carlo (MC) simulations. Alongside, we discuss the rate of convergence, which in turn discusses the accuracy of the convergence result for the smaller (and practical) number of entities. 

We consider an example system with $\gamma = 0$, i.e., with one group in this section.   For each run of the simulation, we first generate a realization of the random graph by generating independent binary random variables (with probability $p_2$) $\{I_{j,i}\}$  for all $j,i$ and another set of independent binary random variables (with probability $p_2^{sb}$) $\{\eta_j^{sb}\}$ for all $j$, also independent of the former set. Thus we have a realization of the financial network along with the portfolios of each entity.   We then generate the realization of economic shocks by generating the two-valued random variables  with upward movement $u$  (with probability $(1-w)$) or downward movement $d$.  

For each random sample generated as above, we compute $\{W_{j,i}\}$, $\{\Omega_j^2 \}$ (as in \eqref{Eqn_Omegaj2}) and $\{K^{2}_i\}$ and solve the fixed point equations given in  \eqref{eqn_clearingvectorG_2}, (with $y_c = 0$) to obtain the clearing vector  $\{X_i\}_i$.
The corresponding fixed point equation modifies to the following with $\gamma = 0$ (with ${\bf X} := (X_1, \cdots, X_n)$):
\begin{eqnarray}
\label{eqn_clearingvectorG_2_modified}
f_i ({\bf X}) &:= & \min \bigg \lbrace \bigg( K^{2}_i  +  \ \displaystyle\sum_{j\in \mathcal{G}_2} X_j W_{j,i} -v_2 \bigg)^+,\bar{y}_2 \bigg \rbrace,  \mbox{ for any } ~i.
%, \mbox{ with } .
 \mbox{
where,  }
 \\ \nonumber 
    K^{2}_i&=&
    \begin{cases}
    k_{u2} & \text{with probability }\ 1-w \\
      k_{d2}, & \text{otherwise,}  \mbox{ with }    
    k_{u2} =
       \Omega_i^2(1+u) \mbox{ and }  k_{d2} = \Omega_i^2 (1+d). 
    \end{cases}
  \end{eqnarray}In particular we are considering the scenarios with  $P^2_D = w$ and $v_2 > k_{d2}$ for this case-study; under these conditions by Lemma \ref{Lemma_G2 default with v2> kd} (case 3) we have:
  \begin{eqnarray}
  \label{Eqn_sim_barxth}
  {\bar x}_{th} &=&  \frac{ {\bar y}_2 (1-w) +( k_{d2}-v_2) w} {1- w (1-p_2^{sb})} (1-p_2^{sb}),\\
  E[S]_{th} &=& ( k_{u2} - v_2 + {\bar x}_{th}-{\bar y}_2 ) (1-w).
  \label{Eqn_sim_Es}
  \end{eqnarray}
 Throughout the simulations, we use the following set of common parameters, and any additional changes  of the  parameters are mentioned in the respective table itself:
$u=0.2$, $d=-0.6$, $r_2= 0.12$,
   $\kappa= 0.56$, $\Omega_2= 12.5$, $\bar{y}_2 =35$, $v_2=7$, $w= 0.2$, $p_2^{sb}=0.001$, $y_c=0$, $d_c=0$.
 \ignore{ 
  {\color{red}
  We set the following to get one to one correspondence:
   $$
   u = \bigg(\frac{k}\Omega_2 -1\bigg)^+ ,  \ 
   d = \bigg(\frac{k-\epsilon}{\Omega_2} -1\bigg) =\bigg( u- \frac{\epsilon}{\Omega_2}\bigg)\  \mbox{ and } \
   k= \bigg(k_0 +y_2p_2^{sb} +c^{'}\bigg)
   $$
   where, $c^{'}$  is a any positive constant such that $u\geq r_2 > d$.
   Also set $y_2 =\frac{c}{1+r_2}$ where $c$ is an appropriate constant. Recall with $y_c =0$,  $\Omega_2$ reduces to $\Omega_2 = (k_0 +y_2 p_2^{sb})$.
   
   In the table 1 experiment we have: $u=0.2$, $d=-0.68$, $r_2= 0.12$,
   $\kappa= 0.56$, $\Omega_2= 12.5$, $y_2= 31.25$, $\bar{y}_2 =35$,$v_2=7$, $w= 0.2$. $p=0.05$, $p_2^{sb}=0.001$, $y_c= d_c=0$.
   
    In the table 2 experiment we have: $u=0.2$, $d=-0.6$, $r_2= 0.12$,
   $\kappa= 0.56$, $\Omega_2= 12.5$, $y_2= 31.25$, $\bar{y}_2 =35$,$v_2=7$, $w= 0.2$. $p=0.05$, $p_2^{sb}=0.001$, $y_c= d_c=0$

}}

We use an iterative algorithm to minimize $\sum_{i\le n} \left ( X_i - f_i({\bf X})\right )^2$ (observe FP is the minimizer)  to obtain the fixed point, i.e., the clearing vector of any sample path. The same is provided in  Algorithm \ref{algo_Fixed_point_algo}. %Note that for simulation purpose we have taken the parameter choice of  $p^{sb}_2 =0.001$  to have the minimal influence of the big bank in the network. However one can take any values between $(0,1)$. 

\begin{algorithm}

      \caption{Fixed-point algorithm to compute  the clearing vector \label{algo_Fixed_point_algo}}
      \begin{algorithmic}[1]
       \State Inputs: $n$, $k$, $K_i^2$, $\{W_{j,i}\}$, $\bar{y}_2$, $v_2$, $\delta$, $T$.
        
            \State Initialize $  X^0_i  = {\bar y}_2$ for all $i\le n$
            
           \State Iteration:  $t= 1,2,\cdots T$

        \begin{algsubstates}
            \State  update: $X_i^{t+1} = X_i^{t} - \epsilon_t(f_i({\bf X}^t)-X_i^{t})$
            
            \State if $\sum_{i=1}^{n}|X_i^{s+1} - X_i^{s}| < n\delta $, for all $s = t-k, t-k+1, \cdots, t$
            
            \State  algorithm  converged and end

        \end{algsubstates}
       \State end
        \end{algorithmic}
    \end{algorithm} 
    
    \begin{table}[hbt!]
\centering
\begin{tabular}{|l|l|l|l|l|l|l|}
\hline
$n$ & $\bar{x}_{th}$ &$\widehat{\bar{x}}$  & Error($\%$) & 
 $ E[S]_{th}$& $ \widehat{E}[S]$  &  Error($\%$)  \\  \hline
600  & 34.5  & 34.2414 & 0.7496 & 6      & 5.9036  & 1.6061 \\ \hline
700  & 34.5  & 34.2760  & 0.6493 & 6      & 5.9551  & 0.7484 \\ \hline
800  & 34.5  & 34.3653 & 0.3904 & 6      & 5.9717  & 0.4722 \\ \hline
900  & 34.5  & 34.4262 & 0.2139 & 6      & 6.0245  & 0.4078 \\ \hline
1000 & 34.5  & 34.4395 & 0.1754 & 6      & 6.0267  & 0.4452 \\ \hline
\end{tabular}
\caption{Sample path wise estimates for ER-graphs.}
\label{Table_1 withbarx}
\end{table}

\ignore{
We demonstrate the simulation model with simplified modelling assumption, i.e., $\gamma=1$ and with a large number of homogeneous entities (say $n$). %and without the big bank. 
%By the homogeneity we refer all the banks are making same decisions during investment. 
We are considering a  simple two time period framework $t=0,1$.
%{\color{red} Recall at the initial period, borrow money from each other and invest in the outside project. In the next period, the bank pays back the loans to the other institutions. If any bank fails to cover the liability, it is declared  a defaulted bank. We are assuming simple bankruptcy rules to pay the liability.
%}
Let $X_i$  denote clearing vector which represents maximum possible clearance of the bank $i$ towards the network.
\begin{equation}
\label{Eqn_clearing_vector_simulation}
X_i = \min \bigg \lbrace \bigg( (K_1 -Z_i)^+ +  \displaystyle\sum_{j=1}^{n} X_j \frac{I_{ji}}{\displaystyle\sum _{i' \leq n} I_{ji'}}-v_1 \bigg)^+,\bar{y}_1 \bigg \rbrace
\end{equation}
where, $Z_i$ be the (random) idiosyncratic shock. We implemented the following algorithm in MATLAB :
\begin{algorithm}
      \caption{Fixed-point Algorithm }
      \begin{algorithmic}[1]
        \State Input: $n$, $\epsilon$, $w$, $p_{1}$, $K_1$, $\bar{y}_1$, $v_1$, $\delta$, $T$
        \State Generate a random sample of $(K_1 -Z_i)^+$, with probability $w$, where $Z_i$ is binary valued
        \State Generate  $I_{ij}$ with probability $p_{1}$
        
            \State Estimate  $|\sum _{i' \leq n} I_{ji'}|$  $  \forall j= 1,2, \cdots n$
            
            \State Initialize $\lbrace X^0_i \rbrace$
            
            \State Iteration:  $t= 1,2,\cdots  T$

        \begin{algsubstates}
            \State  update: $X_i^{t+1} = X_i^{t} - \epsilon_t(f_i(X_i^t)-X_i^{t})$
            
            \State if $\sum_{i=1}^{n}|X_i^{t+1} - X_i^{t}| < \delta $
            
            \State  fixed point converged

        \end{algsubstates}
       \State end
        \end{algorithmic}
    \end{algorithm} \label{Fixed point algo}
 \subsection{Simulation  framework}
We built a simulation model to check the accuracy of the approximation of the convergence of the fixed-point solution. In our simulation framework we consider binary distribution for the shock values $Z_i$  as below:
\begin{equation}
\label{eqn_simulation_shockdist}
    Z_i=
    \begin{cases}
      \epsilon, & \text{wp}\ w  \\
      0, & \text{otherwise}.
    \end{cases}
  \end{equation}

\noindent The connection between banks  taken as  binary random variable  i.e.,
\begin{equation}
\label{eqn_connection_distribution}
    I_{ji}=
    \begin{cases}
     1, & \text{wp}\ p_{1}\\
     0, & \text{otherwise}.
    \end{cases}
  \end{equation}
  where $I_{ji}$ indicates that bank $j$ and bank $i$ are connected. The anticipated return, i.e. $K_1 = k$, $\forall i \in \mathcal{G}_1$ and liability amount $\bar{y}_1$ and tax $v_1$ are kept as a constant.

  First, we generate  random liability matrix and the asset return distribution (see \eqref{eqn_simulation_shockdist}- \eqref{eqn_connection_distribution}). By using  stochastic gradient algorithm, we calculate the  fixed point equations for all the banks %In the above, we describe the algorithm for computing the random fixed point equations. 
  and solve them using Algorithm \ref{Fixed point algo} described above.}
  
 Once we ensure the convergence of the estimates in the algorithm   (when the difference of step 3(b) of Algorithm \ref{algo_Fixed_point_algo} is below $\delta=0.0001$ for $k=100$ consecutive steps), we compute the performance measures related to the systemic risk.   We tabulate these estimates (represented using $\widehat{\cdot}$), along with aggregate fixed point 
 in Table \ref{Table_1 withbarx}. We also tabulate the theoretical 
  $\bar{x}_{th}$ (computed using \eqref{Eqn_sim_barxth})
 and the theoretical expected surplus $E[S]_{th}$ (computed using \eqref{Eqn_sim_Es}) in the same table. We further included an index by name ``Error(\%)" that compares the two sets by computing the normalized error as below, for example for expected surplus:
 $$
 \frac{|E[S]_{th} -\widehat{E}[S] |}{E[S]_{th}}\times 100\%
$$ 
Our theoretical results well-match  the Monte-Carlo estimates (sample path-wise), even with a few hundred banks in the network. Observe here that the table is for one sample path of the graph model.

  %%%%%%%%%%%%
  
{\bf Random Graphs:} In the previous example we considered only Erd\H{o}s-R\'enyi (ER)  graphs. 
 For these  graphs  the generated sample paths  are highly irregular, i.e., the variance  
 in the  number of connections (e.g., $\sum_i {I}_{j,i}$ is the number of lenders and $\sum_i {I}_{i,j}$ is the number of borrowers for entity $j$) across different  entities 
 of the network is high.
Thus, 
we include another set of regular graphs.   We   generated the (correlated) regular graphs by discarding the samples if  the   number  of  connections  deviated  significantly  from  the  true  average.  Such  a  controlled  generation  of  the  random  graphs  leads  to  more  regular graphs with lesser variations in   the     number  of  lenders of various entities of  each  sample  path.  For example, with  $np_2= 500 \times 0.05 = 25$,  we allowed $\pm 2$ variations in the number of lenders.  Observe here that the variations in the number of borrowers is still significantly high.

% New table with respect to the theory and surplus
\begin{table}[hbt!]
\vspace{-2mm}
\centering
\begin{tabular}{r|r|r|r|r|r|r|r|r|}
\cline{2-9}
\multicolumn{1}{c|}{} & \multicolumn{2}{c|}{\textbf{Regular}} & \multicolumn{2}{c|}{\textbf{ER}} & \multicolumn{2}{c|}{\textbf{Regular}} & \multicolumn{2}{c|}{\textbf{ER}} \\ \cline{2-9} 
\multicolumn{1}{l|}{} & \multicolumn{4}{c|}{$p_2$=0.03} & \multicolumn{4}{c|}{$p_2$=0.05} \\ \hline
\multicolumn{1}{|c|}{n} & \multicolumn{1}{c|}{${\widehat P}^2_D$} & \multicolumn{1}{c|}{CI} & \multicolumn{1}{c|}{${\widehat  P}^2_D$} & \multicolumn{1}{c|}{CI} & \multicolumn{1}{l|}{${\widehat  P}^2_D$} & \multicolumn{1}{l|}{CI} & \multicolumn{1}{l|}{${\widehat  P}^2_D$} & \multicolumn{1}{l|}{CI} \\ \hline
\multicolumn{1}{|r|}{500} & 0.1795 & 0.0038 & 0.1397 & 0.0087 & 0.1896 & 0.0029 & 0.1499 & 0.0074 \\ \hline
\multicolumn{1}{|r|}{1000} & 0.1898 & 0.0021 & 0.1517 & 0.0069 & 0.1973 & 0.0018 & 0.1612 & 0.0056 \\ \hline
%\multicolumn{1}{|r|}{1500} & 0.1958 & 0.0016 & 0.1600 & 0.0057 & 0.1984 & 0.0015 & 0.1706 & 0.0043 \\ \hline
\multicolumn{1}{|r|}{2000} & 0.1986 & 0.0013 & 0.1642 & 0.0051 & 0.1990 & 0.0014 & 0.1752 & 0.0036 \\ \hline
%\multicolumn{1}{|r|}{2500} & 0.1995 & 0.0010 & 0.1700 & 0.0043 & 0.2001 & 0.0011 & 0.1813 & 0.0028 \\ \hline
\multicolumn{1}{|r|}{5000} & 0.2001 & 0.0008 & 0.1850 & 0.0022 & 0.2000 & 0.0008 & 0.1937 & 0.0012 \\ \hline
\end{tabular}
\caption{Default-probability estimates over $200$ sample paths.}
\label{Table_withPD_p=0.03_p=0.05} 
 
 \vspace{4mm}
\centering
\begin{tabular}{r|r|r|r|r|r|r|r|r|}
\cline{2-9}
\multicolumn{1}{c|}{} & \multicolumn{2}{c|}{\textbf{Regular}} & \multicolumn{2}{c|}{\textbf{ER}} & \multicolumn{2}{c|}{\textbf{Regular}} & \multicolumn{2}{c|}{\textbf{ER}} \\ \cline{2-9} 
\multicolumn{1}{l|}{} & \multicolumn{4}{c|}{$p_2=0.03$} & \multicolumn{4}{c|}{$p_2=0.05$} \\ \hline
\multicolumn{1}{|c|}{n} & \multicolumn{1}{c|}{$ \widehat{E}[S]$} & \multicolumn{1}{c|}{CI} & \multicolumn{1}{c|}{$ \widehat{E}[S]$} & \multicolumn{1}{c|}{CI} & \multicolumn{1}{l|}{$ \widehat{E}[S]$} & \multicolumn{1}{l|}{CI} & \multicolumn{1}{l|}{$ \widehat{E}[S]$} & \multicolumn{1}{l|}{CI} \\ \hline
\multicolumn{1}{|r|}{500} & 5.9680 & 0.0254 & 5.9715 & 0.0298 & 5.9713 & 0.0249 & 5.9520 & 0.0297 \\ \hline
\multicolumn{1}{|r|}{1000} & 5.9878 & 0.0162 & 5.9564 & 0.0182 & 5.9656 & 0.0179 & 5.9668 & 0.0168 \\ \hline
%\multicolumn{1}{|r|}{1500} & 5.9714 & 0.0144 & 5.9497 & 0.0157 & 5.9714 & 0.0149 & 5.9585 & 0.0154 \\ \hline
\multicolumn{1}{|r|}{2000} & 5.9606 & 0.0125 & 5.9709 & 0.0128 & 5.9724 & 0.0137 & 5.9784 & 0.0131 \\ \hline
%\multicolumn{1}{|r|}{2500} & 5.9606 & 0.0099 & 5.9645 & 0.0113 & 5.9642 & 0.0107 & 5.9658 & 0.0098 \\ \hline
\multicolumn{1}{|r|}{5000} & 5.9642 & 0.0084 & 5.9637 & 0.0070 & 5.9657 & 0.0076 & 5.9655 & 0.0081 \\ \hline
\end{tabular}
\caption{ Expected-surplus estimates  over $200$ sample paths.}
\label{Table_withES_p=0.03_p=0.05}
\end{table}

We consider more sample paths in the next case-study in Tables \ref{Table_withPD_p=0.03_p=0.05}-\ref{Table_withES_p=0.03_p=0.05}.
We estimated the  default probability (${\widehat P}^2_D$) and the  expected surplus  ($ \widehat{E}[S]$) 
 for both   ER   and   regular graphs by averaging over $200$ sample paths. When the number of banks is small for ER graphs, the error between the estimated  default probability and the theoretical default probability ($P_D^2 = w =0.2$) is significantly high. However  the error reduces with the increase in the number of banks.  {\it Thus for  ER graphs, the rate of convergence of the performance  is slow.} On the other hand, the {\it same error is significantly small} for regular graphs. 
 Further the expected surplus and the aggregate clearing vector MC-estimates are close to the theoretical ones (provided in Table \ref{Table_1 withbarx}),  even for small values of $n$ and even for ER graphs (see Table \ref{Table_withES_p=0.03_p=0.05}). We tabulated only  $\widehat{E}[S]$   estimates, as the error/difference  is significantly small for one sample path itself as in  Table \ref{Table_1 withbarx}.

 %For ER graphs, the sample paths generated are irregular, i.e., the variance of deviations from the mean number of connections is high, which results in the {\it inaccuracy of the performance measure} for the lower value of $n$ (see Table \ref{Table_withPD_p=0.03_p=0.05}).

%We also generated the (correlated) regular graphs by discarding the samples if the mean number of connections deviated significantly from the true average links. Such a controlled generation of the random graphs leads to more regular graphs, i.e., the total number of borrowers and lenders is approximately equal in each sample path-wise. We observed that with a few hundred banks in the network, the estimated mean default probability and the expected surplus match the theoretical performance measure. 

\noindent {\bf Confidence intervals:} It is a standard practice to   compute the  $95 \%$ confidence interval\footnote{Confidence interval, CI = [estimated value - HW,  estimated value+HW].} by estimated-mean $\pm$ Half-width ($H_W$) of the estimated performance, where, $H_W$:=1.96 $\times$ $\sqrt{variance}/\sqrt{\#samples}$. For simpler representation, in all the tables  (Table \ref{Table_withPD_p=0.03_p=0.05}-\ref{Table_with_Regular_grpahp=0.05}), $H_W$  is shown as the confidence interval (CI).  
Once again regular graphs have  good CIs even for small values of $n$, while ER graphs are good only for  larger $n$. However the CIs related to expected surplus (in Table \ref{Table_withES_p=0.03_p=0.05}) are very good for both the types of graphs. Further CIs are better with bigger probability of connection $p_2$; and so are the estimated means. 

%In the same table, we observe that when the probability of connection increases, the estimated default probability performance improves both graphs. % while, the expected surplus performance improvement is not significant.

%\noindent  \textbf{Remark:}  When we generated the regular graphs with controlled experiments by discarding the sample if the paths are irregular. The performance is near to the true value. It is important to note that the assumption {\bf B.2} is satisfied with the regular graphs.

\begin{figure}[hbt!]
 \vspace{-2mm}
 \begin{subfigure}{\textwidth}
\begin{subfigure}{.5\textwidth}
  % include first image
  %\hspace{-20mm}
  \includegraphics[width=5.cm, height= 4.3cm]{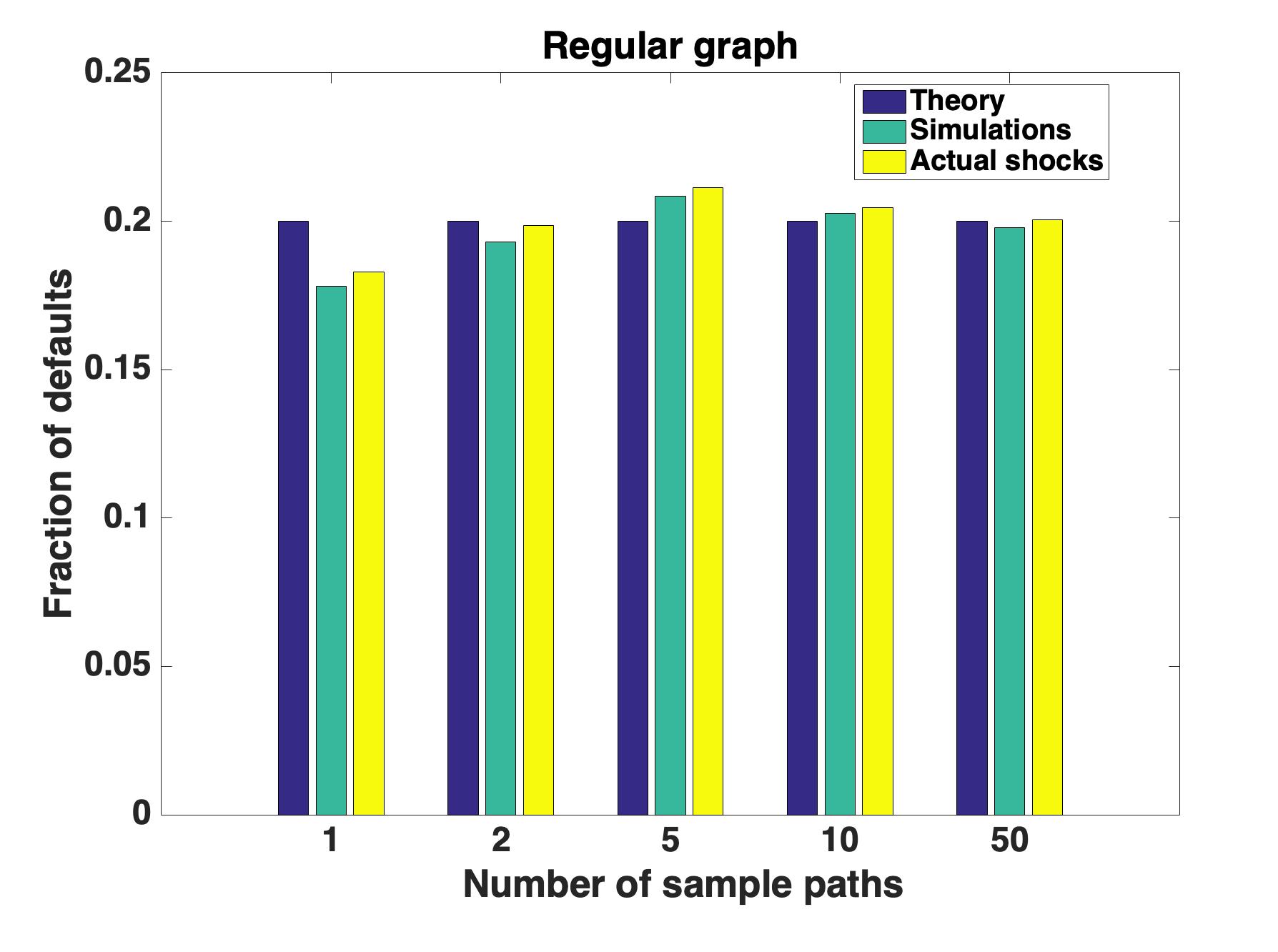} \end{subfigure}
  \hspace{.5cm}
\begin{subfigure}{.5\textwidth}
  % include second image
 % \hspace{-17mm}
  \includegraphics[width=5.cm, height= 4.35cm]{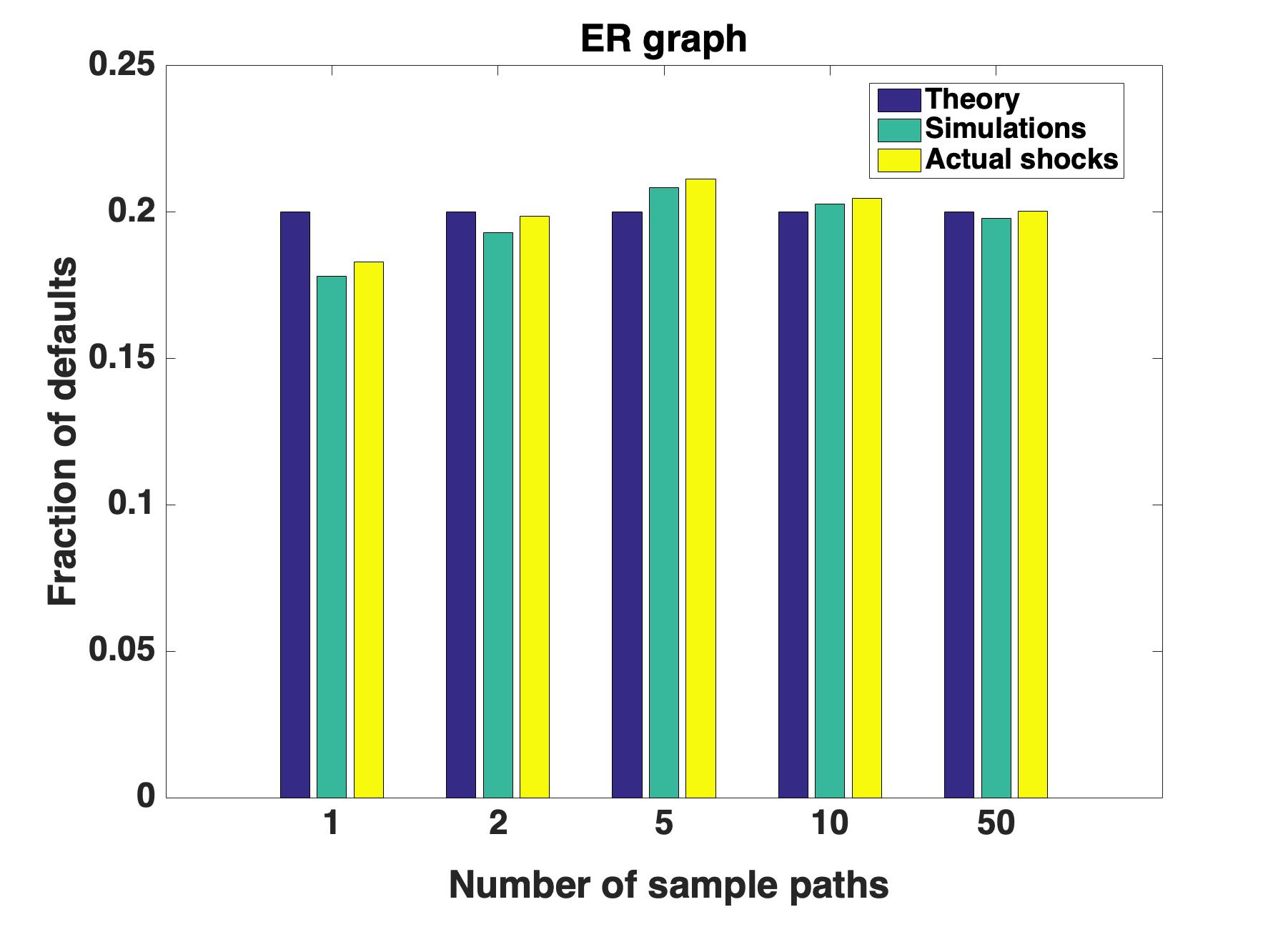} 
  \vspace{-1mm}
 \end{subfigure}
 \end{subfigure}
 \vspace{-4mm}
 \caption{ Few sample-path  estimates with $n=1000$: fraction of  defaults   and  fraction of banks  with shocks. }
\label{fig:bar_digram}
\end{figure}

In the bar charts of Figure \ref{fig:bar_digram}, we consider another case-study with small number of sample paths for default probability estimates. We consider average over  $1,2,5,10$ and $50$ sample paths.    The figure represents the theoretical default probability, simulated default probability, and the fraction of banks that receive the shocks for the regular graph (left sub-figure) and the ER graph (right sub-figure).   The correlation between the simulated fraction of defaults and the simulated fraction of banks that received the shocks  is higher in regular graphs, than that in ER graph. Nevertheless, in both the graphs there is good correlation between these two estimates, and both of them converge closer to theoretical estimate $w=0.2$ when the average is
over large number of sample paths (more than or equal to $10$). The theory indicates that, for the conditions of these case-studies, the number of defaults should equal the number of the banks that received the shocks and the same is well correlated in the simulations, even when the actual fraction of defaults is away from $0.2$.
 
 \begin{table}[hbt!]
\centering
\begin{tabular}{|l|l|l|l|l|l|}
\hline
$n$   & $\widehat{\bar{x}}$ & ${\widehat P}^2_D$     & CI      & $ \widehat{E}[S]$    & CI      \\ \hline
200 & 34.3821 & 0.1635 & 0.0062 & 5.9816 & 0.0375 \\ \hline
300 & 34.4105 & 0.1644 & 0.0059 & 5.9871 & 0.0351 \\ \hline
400 & 34.4327 & 0.1831 & 0.0038 & 5.9455 & 0.0283 \\ \hline
%500 & 34.4356 & 0.1778 & 0.0041 & 5.9748 & 0.0261 \\ \hline
\end{tabular}
\caption{ \textbf{Regular graphs} performance estimates over $200$ sample paths.}
\label{Table_with_Regular_grpahp=0.05}
\end{table}

Finally we consider further regular graphs, whose number of borrowers is also controlled and  with even smaller number of agents in Table \ref{Table_with_Regular_grpahp=0.05}. With $n=200$, $300$ and $400$ we respectively allowed number of borrowers  to be distributed between $10\pm 6$,  $15\pm 9$ and $20 \pm 12$, while the number of lenders is distributed respectively between 
$10\pm 1$, $15\pm2$ and $20\pm 2$. We observe a sufficiently good match between theoretical values  and the corresponding MC estimates, even for these small values of $n$.  We thus conclude that our theoretical estimates are sufficiently good matches even for few hundreds of banks when the number of connections across agents is not too diverse. Thus our results can provide a good method for estimating clearing vectors and further systemic-risk measures  for many practical scenarios.

 \section{Conclusions}
\label{sec_conclusion}
We consider a large dimensional fixed point equation, where the function corresponding to any component  depends upon a weighted aggregate of the values of its random neighbours. The underlying fixed point equations are random, and we provide a methodology to solve such equations under suitable graph structure(s);
our main contribution is to  solve these equations almost surely. We consider two different types of random graph models: the resources are shared equally across the (connected components of the) entire network in the first model. In contrast, in the second model, the resources are shared equally only within the group after allocating dedicated fractions to each group. In both the models,  the solution of the  random fixed point equation  converges   almost surely  to that of a limit system, and these solutions are asymptotically independent. This asymptotic simplification reduced the dimensionality of the  problem significantly; the almost-sure asymptotic solution is derived by solving a deterministic three-dimensional fixed point equation.

We apply the above results to a large financial network to study systemic risk related aspects. 
In such networks, the first object to be studied is clearing vector, a vector of maximum possible repayments (towards clearing their liabilities) by all the nodes of the network; this vector depends upon the random economic shocks received by a fraction of  agents, and the percolation of the influence of these shocks on the clearing capacity of their neighbours, neighbours of neighbours and so on.   
One of our  primary results is a procedure to compute  clearing vector of a variety of large-dimensional financial networks; in the example  considered in this paper, solution of a two-dimensional deterministic  fixed point equation is sufficient to study the clearing vector.

%accuracy. Simulations suggest that the approximation is well-matched with the Monte Carlo estimates.

 Considering a heterogeneous financial network with  two-period  framework,  binomial shocks and with two diverse groups (one aggressive to consider more risky avenues while the other is recessive), we derived  systemic-risk based performance measures, after deriving the clearing vectors.
 For majority of the cases,  closed-form expressions  are obtained for aggregate clearing vector, systemic risk performance measures viz. probability of default and expected surplus; for other cases we solved  the two-dimensional equations numerically. In the limiting network, we observed some phase transitions with respect to inter group-connectivity parameters: a) the existence of a regime of connectivity parameters,  wherein the surplus of both the groups improves with the   parameters; and b) in large shock regimes, the inter-connectivity has   adverse effect, at least on one of the groups.

 We performed exhaustive   Monte-Carlo simulations using  practical number of financial entities and Erd\H{o}s-R\'enyi graphs with an aim to study the accuracy of  approximation.  We observed good match for the surplus based performance measures and clearing vectors, even with few hundreds of banks. With more regular graphs (smaller variations in the number of connected components) there is a good match even between  default probability estimates and their theoretical counterparts.  
 
 \ignore{ 
    To conclude, our random fixed points asymptotic approximation leads to small-dimensional deterministic fixed point equations, which compute the financial entities' clearing vector by performing any iterative algorithm.
 
 %%%%%%%%%%%%%%%%%%%%%%%%% We could ignore the below part%%%%%%%%%%%%%%

  We observed  a certain phase transition   with respect to   inter lending parameters; below a certain threshold  defined by inter-lending  parameters, the   network is stable/resilient (no defaults in spite of  shocks);  once  the threshold is crossed  the  network is in default/collapse regime (percolation of shocks and eventual collapse). 
 
  We  illustrated the existence of a regime of  connectivity parameters,  wherein  the surplus of both the groups improves with inter-lending parameters; however,   in  large shock regimes, the inter-connectivity has an adverse effect (surplus of  at least one of the groups degrades).  

    %\item proved that if group 2 banks are resilient, then group 1 banks are also resilient.
}

%{\color{blue}
%We considered large dimensional fixed point equations which depend upon the performance of the neighbouring nodes in the network. The underlying fixed point equations are  random, and we solved such equations for a large number of nodes under suitable graph structure assumption. We showed that the aggregated fixed point solution converges to an almost sure constant. We apply the result to study the systemic risk in a financial network. We provided an analysis of how inter-connectivity plays a role in the stability of the system. We demonstrated critical phase change behaviour of the network and also identified the regime where both group banks surplus increases with respect to the inter-lending  parameter.
%We performed the Monte Carlo simulation, which suggests the approximation has an error within $10\%$.
%}

% \bibliographystyle{plain}
% \bibliography{biblography.bib}

\bibliographystyle{plain}

\ignore{
}

\ignore{
%\newpage
{\small
\begin{table}[htbp!]
%\label{Table_Notations}
	\begin{center}
		\begin{tabular}{|l|}
			\hline
			\\
			\centerline{\bf  Notations}
			\\	\hline \hline
		 $n$ represent the number of small nodes (or banks in the financial example).  
			\\ \hline
		
		$b$ represent the big node. 
			\\ \hline
		$n_1$, $n_2$ represent the  total number of nodes in each group.
		 \\ \hline
		 $N$ represent the set of the banks.
			\\ \hline
		$G^{m}_i$ = Represent the shock of the group $m\in\lbrace 1 ,2\rbrace $ for the $i$-th individual.  
		\\ \hline
		$\gamma$ be the positive fraction  of nodes in $\mathcal{G}_1$ and $(1-\gamma)$ be the nodes in $\mathcal{G}_2$.
			\\ \hline
			
			$\mathcal{G}_{m}$ represent the groups.
			\\ \hline
				$p_m$ are the group  and $p_{c_m}$ inter group  connection probability.
			\\ \hline
			$\{ I_{j,i} \}$ represent the indicator.
			
			\\ \hline
			
			$W_{j,i}$ weighted fraction from small node  $j$ to $i$.
			
		  \\ \hline 
		  $\{\eta_j^{sb} \}_j$,  $\{\eta_j^{bs} \}_j$ %$\{\zeta_j^{sb} \}_j$ 
are i.i.d.  random variable.
         \\ \hline
         
          $E(\eta_j^{sb} )= p^{sb}_m$ if $j\in \mathcal{G}_m$ and $E(\eta_j^{bs} )= p^{bs}_m$. 
          
         \\ \hline
          
          $\bar{X_i}^m$ represent the aggregate  for the small node $i$.
          \\ \hline
          $\bar{X_{b}}$ for the big node aggregate.
          
          \\ \hline 
          
          $\bf{x}$  is an appropriate  dimension  of the vector .
          
          \\ \hline
           $\bf{f}$ represent the vector valued function.
           
           \\  \hline
             
         $(f^{1},f^{2}, f^b)$ represent component wise functions.
         
         \\ \hline
           $(	X^{1},X^{2}, X^b)$ corresponding nodes performance vector.
           
           \\ \hline
           
             Fraction of  nodes in each group  through  $\gamma_1 := \gamma,  \  \gamma_2 := (1-\gamma)$. 
             
             \\ \hline
             $\gp{1} := \gamma p_1+ (1-\gamma)p_{c_1}  \mbox{  and  } \gp{2} :=  \gamma p_{c_2}+ (1-\gamma) p_2.$
             
			\\ \hline
			
			$K^{m}_i$ be the investments in outside network with $m =1,2$.
		    \\ \hline
		    $\bar{y}_m$ be the liability for group $m=1,2$.
		   
		   \\ \hline
		   
		   $\Omega_i$ be the investment in the outside security with $i =1 ,2 $.
		   
		   \\ \hline

		   $\kappa$: be the scaled tax factor or deposit.
		   \\ \hline
		
		\end{tabular}
		\caption{Summary of Notations}  
		\label{Table_Notations}
	\end{center}
\end{table}
}}
\ignore{
\begin{figure}
\tikzset{every picture/.style={line width=0.75pt}} %set default line width to 0.75pt        

\begin{tikzpicture}[x=0.65pt,y=0.65pt,yscale=-1,xscale=1]

%Straight Lines [id:da6311265433151858] 
\draw    (125,178.5) -- (122.05,48) ;
\draw [shift={(122,46)}, rotate = 448.7] [color={rgb, 255:red, 0; green, 0; blue, 0 }  ][line width=0.75]    (10.93,-3.29) .. controls (6.95,-1.4) and (3.31,-0.3) .. (0,0) .. controls (3.31,0.3) and (6.95,1.4) .. (10.93,3.29)   ;
%Straight Lines [id:da7254810126673821] 
\draw    (248.5,194) -- (171.5,193.03) ;
\draw [shift={(169.5,193)}, rotate = 360.73] [color={rgb, 255:red, 0; green, 0; blue, 0 }  ][line width=0.75]    (10.93,-3.29) .. controls (6.95,-1.4) and (3.31,-0.3) .. (0,0) .. controls (3.31,0.3) and (6.95,1.4) .. (10.93,3.29)   ;
%Shape: Rectangle [id:dp5499111928681848] 
\draw   (249.5,180) -- (346,180) -- (346,209.5) -- (249.5,209.5) -- cycle ;
%Straight Lines [id:da1642998262661275] 
\draw    (346.5,193) -- (434,194.47) ;
\draw [shift={(436,194.5)}, rotate = 180.96] [color={rgb, 255:red, 0; green, 0; blue, 0 }  ][line width=0.75]    (10.93,-3.29) .. controls (6.95,-1.4) and (3.31,-0.3) .. (0,0) .. controls (3.31,0.3) and (6.95,1.4) .. (10.93,3.29)   ;
%Shape: Rectangle [id:dp1251186125408894] 
\draw   (435,176) -- (530,176) -- (530,204.5) -- (435,204.5) -- cycle ;
%Shape: Rectangle [id:dp46143994775520336] 
\draw   (79.5,178.5) -- (171,178.5) -- (171,204.5) -- (79.5,204.5) -- cycle ;
%Straight Lines [id:da39124615307971644] 
\draw    (487,176.5) -- (486.51,52) ;
\draw [shift={(486.5,50)}, rotate = 449.77] [color={rgb, 255:red, 0; green, 0; blue, 0 }  ][line width=0.75]    (10.93,-3.29) .. controls (6.95,-1.4) and (3.31,-0.3) .. (0,0) .. controls (3.31,0.3) and (6.95,1.4) .. (10.93,3.29)   ;
%Shape: Rectangle [id:dp0850152941590252] 
\draw   (222,38.5) -- (361,38.5) -- (361,61.5) -- (222,61.5) -- cycle ;
%Straight Lines [id:da592184708782449] 
\draw    (361,50) -- (486.5,50) ;
%Straight Lines [id:da5809337676919744] 
\draw    (122,46) -- (221,46.5) ;

% Text Node
\draw (446.5,181.5) node [anchor=north west][inner sep=0.75pt]   [align=left] {Lemma \ref{lem: LLN_fixed point}};
% Text Node
\draw (262.5,187.5) node [anchor=north west][inner sep=0.75pt]   [align=left] {Lemma \ref{Master lemma}};
% Text Node
\draw (125.5,196.9) node [anchor=north west][inner sep=0.75pt]    {$$};
% Text Node
\draw (86.5,185.5) node [anchor=north west][inner sep=0.75pt]   [align=left] {Corollary \ref{Corrolary_supportinglemma}};
% Text Node
\draw (242.5,40.5) node [anchor=north west][inner sep=0.75pt]   [align=left] {Theorem \ref{Thm_MainGen 1}, \ref{Thm_MainGen 2}};
\end{tikzpicture}
\caption{Schematic view of the Proof of Theorem \ref{Thm_MainGen 1} and  \ref{Thm_MainGen 2}.}
\end{figure}
}
\section*{Appendix A: Proofs related to  Section \ref{Graphicalmodel}, random fixed points}

\ignore{
{\color{red}
\textbf{An auxiliary Lemma}
\begin{lemma}
\label{Lemma_Cauchy}
Let $\{z_n \}$ be a bounded infinite sequence i.e., $z_n \in (0,1)$, then $\bar{z}_n := \frac{1}{n}\sum_{i=1}^{n} z_i$ converges. Moreover the sequence $\{z_n \}$ have at most one limit point.
\end{lemma}
\textbf{Proof:} We consider the  following:
\begin{eqnarray*}
&& \bigg|\bar{z}_{n+1} -  \bar{z}_n \bigg| \\
& =& \bigg|\frac{1}{n+1}\sum_{i=1}^{n
+1} z_i -   \frac{1}{n}\sum_{i=1}^{n} z_i \bigg| \\
&=& \bigg|\frac{1}{n+1}
\displaystyle \sum_{i=1}^ {n}  z_i  + \frac{z_{n+1}}{n+1} -\frac{1}{n+1}
\displaystyle \sum_{i=1}^{n} z_i \frac{n+1}{n}\bigg| \\
&\le& \bigg|\frac{1}{n+1}
\displaystyle \sum_{n=1}^{n} z_i \bigg(1-\frac{n+1}{n} \bigg)\bigg| + \bigg|\frac{z_{n+1}}{n+1} \bigg|\\
&\le& \frac{y}{n+1} +\frac{y}{n+1} \le  \epsilon
\end{eqnarray*}
The last inequality is  due to the fact that for a fix $\epsilon > 0$   pick $n$  large enough such that $\frac{2y}{n+1}$ can be made arbitrary small $\forall n\ge n_k(\epsilon)$. Thus the sequence  $\{z_n \}$  is a Cauchy sequence and hence converges. The last part of the Lemma is immediate because $\{z_n \}$ is  Cauchy sequence and it has at most one limit point.
\eop

\noindent{\textbf{Proof of Lemma} \ref{lem: LLN_fixed point} :}
i) First observe that ${\bar f}^n_b ({\bar x}_b,{\x})$ converges  by applying the Lemma \ref{Lemma_Cauchy} individual term  with $z_i=  \xi^{m}_{j} ({\bar x}^{m}_j, x_b) W_{j, b}$  for $m=1,2$.
 Next we prove that ${\bar f}^{n, m}_i ( {\bar x}_b, \x )$ converges. It is enough to show  the individual term of ${\bar f}^{n, m}_i ( {\bar x}_b, \x )$ is a Cauchy sequence. Consider the following:

\vspace{-2mm} 
 {\scriptsize
 \begin{eqnarray*}
   &&\bigg|\sum_{j \in \mathcal{G}_m}   \xi^m_{j} ({\bar x}^m_j, x_b)  W_{j, i}(n+1)  - \sum_{j \in \mathcal{G}_m}   \xi^m_{j} ({\bar x}^m_j, x_b)  W_{j, i}(n)\bigg| \\
   & =&  \bigg|\sum_{j \in \mathcal{G}_m}   \xi^m_{j} ({\bar x}^m_j, x_b)  W_{j, i}(n+1)-   \sum_{j \in \mathcal{G}_m}    \xi^m_{j} ({\bar x}^m_j, x_b) \frac{p_m(1-\eta_j^{sb})}{(n+1) \gamma_{p_m}} + \sum_{j \in \mathcal{G}_m}    \xi^m_{j} ({\bar x}^m_j, x_b) \frac{p_m(1-\eta_j^{sb})}{(n+1) \gamma_{p_m}}\\
   && - \sum_{j \in \mathcal{G}_m}   \xi^m_{j} ({\bar x}^m_j, x_b)  W_{j, i}(n) +  \sum_{j \in \mathcal{G}_m}    \xi^m_{j} ({\bar x}^m_j, x_b) \frac{p_m(1-\eta_j^{sb})}{n\gamma_{p_m}} - \sum_{j \in \mathcal{G}_m}    \xi^m_{j} ({\bar x}^m_j, x_b) \frac{p_m(1-\eta_j^{sb})}{n \gamma_{p_m}}\bigg|  \\
   &\le& \bigg|\sum_{j \in \mathcal{G}_m}   \xi^m_{j} ({\bar x}^m_j, x_b)  W_{j, i}(n+1)-   \sum_{j \in \mathcal{G}_m}    \xi^m_{j} ({\bar x}^m_j, x_b) \frac{p_m(1-\eta_j^{sb})}{(n+1) \gamma_{p_m}}\bigg| +  \bigg|\sum_{j \in \mathcal{G}_m}    \xi^m_{j} ({\bar x}^m_j, x_b) \frac{p_m(1-\eta_j^{sb})}{(n+1) \gamma_{p_m}}\bigg| \\
   && + \bigg|\sum_{j \in \mathcal{G}_m}   \xi^m_{j} ({\bar x}^m_j, x_b)  W_{j, i}(n)-   \sum_{j \in \mathcal{G}_m}    \xi^m_{j} ({\bar x}^m_j, x_b) \frac{p_m(1-\eta_j^{sb})}{n \gamma_{p_m}}\bigg| +  \bigg|\sum_{j \in \mathcal{G}_m}    \xi^m_{j} ({\bar x}^m_j, x_b) \frac{p_m(1-\eta_j^{sb})}{n \gamma_{p_m}}\bigg| \\
   & \le& \bigg|\sum_{j \in \mathcal{G}_m}   \xi^m_{j} ({\bar x}^m_j, x_b) (1-\eta_j^{sb}) \bigg( \frac{1}{{\sum_{i \in   \mathcal{G}_1 \cup  \mathcal{G}_2 } I_{j, i}}} -\frac{1}{(n+1)\gamma_{_{p_m}}}\bigg) \bigg| +  \bigg|\sum_{j \in \mathcal{G}_m}    \xi^m_{j} ({\bar x}^m_j, x_b) \frac{p_m(1-\eta_j^{sb})}{(n+1) \gamma_{p_m}}\bigg| \\
   && +\bigg|\sum_{j \in \mathcal{G}_m}   \xi^m_{j} ({\bar x}^m_j, x_b) (1-\eta_j^{sb}) \bigg( \frac{1}{{\sum_{i \in   \mathcal{G}_1 \cup  \mathcal{G}_2 } I_{j, i}}} -\frac{1}{n\gamma_{_{p_m}}}\bigg) \bigg| +  \bigg|\sum_{j \in \mathcal{G}_m}    \xi^m_{j} ({\bar x}^m_j, x_b) \frac{p_m(1-\eta_j^{sb})}{n \gamma_{p_m}}\bigg| \\
   &\le&  \underbrace{\sum_{j \in \mathcal{G}_m} c_0 \bigg |    \bigg( \frac{1}{{\sum_{i \in   \mathcal{G}_1 \cup  \mathcal{G}_2 } I_{j, i}}} -\frac{1}{(n+1)\gamma_{_{p_m}}}\bigg) \bigg|}_{B.2}  + \underbrace{\bigg|\sum_{j \in \mathcal{G}_m}    \xi^m_{j} ({\bar x}^m_j, x_b) \frac{p_m(1-\eta_j^{sb})}{(n+1) \gamma_{p_m}}\bigg|}_{convergent~ as~ n\to \infty}  \\
   && \underbrace{\sum_{j \in \mathcal{G}_m} c_0 \bigg |    \bigg( \frac{1}{{\sum_{i \in   \mathcal{G}_1 \cup  \mathcal{G}_2 } I_{j, i}}} -\frac{1}{n\gamma_{_{p_m}}}\bigg) \bigg|}_{B.2}  + \underbrace{\bigg|\sum_{j \in \mathcal{G}_m}    \xi^m_{j} ({\bar x}^m_j, x_b) \frac{p_m(1-\eta_j^{sb})}{n \gamma_{p_m}}\bigg|}_{convergent~ as~ n\to \infty}\le \epsilon. 
   \end{eqnarray*}}
Therefore   ${\bar   f}^{n, m}_i ( {\bar x}_b, \x )$  is a    Cauchy sequence and hence converges.\\}}

\noindent{\textbf{Proof of Lemma} \ref{lem: LLN_fixed point}:}  First consider the term $\frac{1}{n} \sum_{j \in \mathcal{G}_1}   \xi^{1}_{j} ({\bar x}^{1}_j, x_b) W_{j, b}$ (recall $x_b = f^b({\bar x}_b)$ is deterministic). Now consider the constant sequences $({\x}^1, {\x}^2)$,  i.e., $\bar{x}_j^{m} = \bar{x}^{m}  $ (same for all  $j \in \mathcal{G}_m$) for each $m\in \{1,2\}$. Then $\{ \xi^{1}_{j} ({\bar x}^{1}, x_b) W_{j, b} \}_j $ is an i.i.d. sequence of random variables (see \eqref{Eqn_Weights 1} and \eqref{Eqn_xi 1}) and by law of large numbers (note $  E[W_{j, b}] = p_1^{sb}$ for $j \in {\cal G}_1$)
$$ \frac{1}{n} \sum_{j \in \mathcal{G}_1}   \xi^{1}_{j} ({\bar x}^{1}_j, x_b) W_{j, b} \to {\gamma} E_{ G^{1}_i, \eta_i^{bs}}  [\xi^{1}_i ({\bar x}^1, x_b)  ] p_1^{sb}  \mbox{ a.s., as  } n \to \infty \ \  (\mbox{for any } i). $$ 
Similarly, $
 \frac{1}{n} \sum_{j \in \mathcal{G}_2}   \xi^{2}_{j} ({\bar x}^{2}_j, x_b) W_{j, b} \to   (1-\gamma) E_{ G^{2}_i, \eta_i^{bs}}  [\xi^{2}_i ({\bar x}^2, x_b)  ]  p_2^{sb}  \mbox{ a.s. }
$ 
  
From \eqref{Eqn_bar_fixedpoint_randGen_1_more},    the aggregate (fixed point) function, for constant sequences is:
\begin{eqnarray}
{\bar f}^{n, m}_i ( {\bar x}_b, {\bar x}^1, {\bar x}^2) &=&  
\left \{  \begin{array}{lll}
\displaystyle \sum_{j \in \mathcal{G}_1}   \xi^1_{j} ({\bar x}^1, x_b)  W_{j, i} + \displaystyle \sum_{j \in \mathcal{G}_2}   \xi_{j}^2  ({\bar  x}^2, x_b)   {W}_{j, i} &\mbox{ if }  i \in \mathcal{G}_m   \\

 0  &\mbox{ else}   
\end{array} \right  .
\label{Eqn_fnm_repeat}  
\end{eqnarray}

Fix any $i \in \mathcal{G}_{m'}$. Define  $M^{m}_j = \xi^{m}_{j} ({\bar x}^{m}, x_b) I_{j,i} $ and observe these are i.i.d.  random variables for each $m$,   uniformly bounded by  constant $y$ of {\bf B.1}.  
By Lemma \ref{Master lemma}, each term of   (\ref{Eqn_fnm_repeat}) converges in almost sure sense on set  ${\cal E}$ of {\bf B.2}  to the following:
 $$ 
  \sum_{j \in \mathcal{G}_m}   \xi^m_{j} ({\bar x}^m x_b)  W_{j, i} \to  E_{G^{m}_i, \eta_i^{bs}} \left [ \xi^{m}_i ({\bar x}^m, x_b) \right] p_{m m'}  \gamma_m  \frac{(1- p^{sb}_{m})}{\gamma_{p_m}}.$$  Adding over all possible   $m, m'  \in \lbrace 1,2\rbrace$  we have the desired  convergence.
  %(see equation (\ref{Eqn_bar_fixedpoint_randGen_1}) for constant sequence). 
  \eop

\begin{lemma}
\label{Master lemma}
Assume {\bf B.2(C)}. Consider an i.i.d. sequence  $\{ M_j\}$,which is uniformly bounded by some constant  $c_0 < \infty$. Define the following random variables,  one for each $n$ and  for any $i \in \mathcal{G}_{m'}$ $m, m'  \in \{1, 2\}$ as below (see equation (\ref{Eqn_Weights 1})):
\begin{equation*}
{\zeta}^n  := \sum_{j \in \mathcal{G}_m}  M_j W_j, \  \  {\tilde \zeta}^n  := \sum_{j \in \mathcal{G}_m} I_{j,i} M_j W_j  \ \mbox{  with } W_j := \frac{  (1 - \eta_{j}^{sb})}{\displaystyle \sum_{i'
\le n} I_{j, i'}}.
\end{equation*}
%Then we have the following (see Table \ref{Table_Notations} for notations):\\ 
 i) we have the following convergence (${\cal E}$ defined in  equation (\ref{Eqn_calE})) 
\begin{equation*}
{\zeta}^n  \to  E[M_j]    \gamma_m \frac{1-p^{sb}_{m}}{ \gamma_{_{p_m}}    }    \mbox{, }  \  \  {\tilde \zeta}^n  \to  E[M_j]  p_{m'm}  \gamma_m \frac{1-p^{sb}_{m}}{ \gamma_{_{p_m}} }  \mbox { a.s. on  set, } {\cal E}   ; \mbox{ and, }
\end{equation*}
ii)  Under  {\bf B.2},  the above convergence is w.p.1.   
\end{lemma}
\ignore{
\AT{\begin{lemma}
\label{Master lemma}
Consider an i.i.d. sequence  $\{ M_j\}$,which is uniformly bounded by some constant  $c_0 < \infty$. Define the following random variables,  one for each $n$ and for any $m  \in \{1, 2\}$ as below (see equation (\ref{Eqn_Weights 1})):
\begin{equation*}
{\zeta}^n  := \sum_{j \in \mathcal{G}_m}  M_j W_j \mbox{  with } W_j := \frac{  (1 - \eta_{j}^{sb})c_m}{\displaystyle \sum_{i
\in \mathcal{G}_m} I_{j, i}}.
\end{equation*}
Then we have the following (see Table \ref{Table_Notations} for notations):\\ 
 i) we have the following convergence (${\cal E}$ defined in  equation (\ref{Eqn_calE})) 
\begin{equation*}
\begin{split}
{\zeta}^n  \to  E[M_j]    \gamma_m c_m \frac{1-p^m_{sb}}{p_m    }   \mbox { a.s. on  set, } {\cal E}; \mbox{ and, }
\end{split}
\end{equation*}
ii) \begin{equation*}
\tilde{\zeta}^n  := \sum_{j \in \mathcal{G}_m}  M_j \tilde{W}_j \mbox{  with } \tilde{W}_j := \frac{  (1-c_m)}{\displaystyle \sum_{i
\notin \mathcal{G}_m} I_{j, i}}.
\end{equation*}
 we have the following convergence (${\cal E}$ defined in  equation (\ref{Eqn_calE})) 
\begin{equation*}
\begin{split}
\tilde{\zeta}^n  \to  E[M_j]    \gamma_m (1-c_m) \frac{1-p^m_{sb}}{p_{c_m}    }   \mbox { a.s. on  set, } {\cal E}; \mbox{ and, }
\end{split}
\end{equation*}
 iii)  Under  {\bf B}.2,  the above convergence is w.p.1.   
\end{lemma}
}
}
\noindent {\textbf{Proof:}}  
By adding and subtracting appropriate terms and from  the equation (\ref{Eqn_Weights 1}):
\vspace{-4mm}
 {\tiny \begin{eqnarray}
 \label{Eqn_Lemma13}
 \bigg|{\zeta}^n  - E[M_j]  \gamma_m \frac{1-p^{sb}_{m}}{ \gamma_{_{p_m}}   } \bigg|  
 &=&   \bigg |\sum_{j \in \mathcal{G}_m}  M_j   W_{j} -\sum_{j \in \mathcal{G}_m} M_j \frac{(1- \eta_j^{sb})}{n\gamma_{_{p_m}} }+ \sum_{j \in \mathcal{G}_m} M_j \frac{(1- \eta_j^{sb})}{n\gamma_{_{p_m}}  } - E[M_j] \gamma_m  \frac{1-p^{sb}_{m}}{ \gamma_{_{p_m}}   } \bigg| \nonumber \\
 && \hspace{-20mm} = \  \bigg |\sum_{j \in \mathcal{G}_m}  M_j   \frac{  (1 - \eta_{j}^{sb})}{\displaystyle \sum_{i'
\in \mathcal{G}_1 \cup \mathcal{G}_2} I_{j, i'}} -\sum_{j \in \mathcal{G}_m} M_j \frac{(1- \eta_j^{sb})}{n\gamma_{_{p_m}} }+ \sum_{j \in \mathcal{G}_m} M_j \frac{(1- \eta_j^{sb})}{n\gamma_{_{p_m}}  } -E[M_j] \gamma_m  \frac{1-p^{sb}_{m}}{ \gamma_{_{p_m}}   } \bigg| \nonumber \\
 && \hspace{-20mm} \le  \ \bigg |\sum_{j \in \mathcal{G}_m}  M_j   (1- \eta_j^{sb})\bigg( \frac{1}{{\sum_{i \in   \mathcal{G}_1 \cup  \mathcal{G}_2 } I_{j, i}}} -\frac{1}{n\gamma_{_{p_m}}}\bigg) \bigg| +\bigg| \sum_{j \in \mathcal{G}_m} M_j \frac{(1- \eta_j^{sb})}{n\gamma_{_{p_m}}} -E[M_j] \gamma_m  \frac{1-p^{sb}_{m}}{ \gamma_{_{p_m}}   }\bigg| \nonumber\\
 && \hspace{-20mm} \le  \  \sum_{j \in \mathcal{G}_m} c_0 \bigg |    \bigg( \frac{1}{{\sum_{i \in   \mathcal{G}_1 \cup  \mathcal{G}_2 } I_{j, i}}} -\frac{1}{n\gamma_{_{p_m}}}\bigg) \bigg| +\bigg| \sum_{j \in \mathcal{G}_m} M_j \frac{(1- \eta_j^{sb})}{n\gamma_{_{p_m}}}- E[M_j]  \gamma_m \frac{1-p^{sb}_{m}}{ \gamma_{_{p_m}}   }\bigg|
\\
&& \xrightarrow [\text{B.2 and LLN}]{n\to \infty} 0 \mbox{ a.s. on set  ${\cal E}$ of {\bf B.2}.} \nonumber  \hspace{10mm}
\end{eqnarray}}
 The proof for ${\tilde \zeta}^n$ goes through in a similar way when $\{I_{j,i}\}$ are i.i.d. Otherwise it still goes through because of the following modified steps (in addition to the previous steps).
\vspace{-4mm}

{\tiny\begin{eqnarray*}
\bigg |
\sum_{j \in \mathcal{G}_m} M_j I_{j,i} \frac{(1- \eta_j^{sb})}{n\gamma_{_{p_m}} } -E[M_j]  p_{m',m} \gamma_m  \frac{1-p^{sb}_{m}}{ \gamma_{_{p_m}}   } \bigg| 
&\le& \bigg |
\sum_{j \in \mathcal{G}_m} M_j I_{j,i} \frac{(1- \eta_j^{sb})}{n\gamma_{_{p_m}} } -
\sum_{j \in \mathcal{G}_m} M_j  p_{m',m} \frac{(1- \eta_j^{sb})}{n\gamma_{_{p_m}} }  \bigg |  \\
&& + \bigg | \sum_{j \in \mathcal{G}_m} M_j  p_{m',m} \frac{(1- \eta_j^{sb})}{n\gamma_{_{p_m}} }  -E[M_j]  p_{m',m} \gamma_m  \frac{1-p^{sb}_{m}}{ \gamma_{_{p_m}}   } \bigg|  \\
&\le& \frac{c_0}{\gamma_{_{p_m}} }\bigg |
\sum_{j \in \mathcal{G}_m}  \frac{ I_{j,i} }{n } -
\sum_{j \in \mathcal{G}_m}   \frac{ p_{m',m} }{n }  \bigg |  \\
&& + \bigg | \sum_{j \in \mathcal{G}_m} M_j  p_{m',m} \frac{(1- \eta_j^{sb})}{n\gamma_{_{p_m}} }  -E[M_j]  p_{m',m} \gamma_m  \frac{1-p^{sb}_{m}}{ \gamma_{_{p_m}}   } \bigg|  \\
&\to&  0 \mbox{ a.s.}
\end{eqnarray*}}
because of assumption {\bf B.2(C)} and law of large numbers.   
\eop

\vspace{-6mm}

\section*{Appendix B: Proof of Theorem \ref{Thm_MainGen 1}}
    \noindent {\textbf{Proof of Theorem} \ref{Thm_MainGen 1}:}
% We showed in  Lemma  \ref{lem: LLN_fixed point} that the aggregate fixed point converges to that of its expected value for constant sequences.
 We consider the following norm for this proof on  the space of infinite sequences,    $[0,y]\times s^\infty \times  s^\infty $: 
 \begin{eqnarray} 
 \label{Eqn_infinynorm}
     || ({\bar x}_b ,\x)- ({\bar u}_b, \u)  ||_\infty %& =&\max \left \{ | {\bar x}_b -   {\bar u}_b |, \ \max_{m \in \lbrace 1, 2\rbrace} \sup_{i\in \mathcal{G}_m } |x_i - u_i| \right \}.\nonumber   \\
      &=&  \max_{m \in \lbrace 1, 2\rbrace} \sup_{i\in \mathcal{G}_m } \left (|\bar{x}^m_i - \bar{u}^m_i| +  | {\bar x}_b -   {\bar u}_b |\right)
 %    || ({\bar x}_b ,\x)- ({\bar u}_b, \u)  ||_2     &=& \max \left \{ | {\bar x}_b -   {\bar u}_b |, \ \max_{m \in \lbrace 1, 2\rbrace} \left (  \sum_{i\in \mathcal{G}_m } 2^{-i} |x_i - u_i|^2 \right )^{1/2} \right \}.
\end{eqnarray} 
%$l_2 (\mu)$ spaces,  little $l_2$ spaces ..

 From the equations \eqref{Eqn_xi 1}-(\ref{Eqn_bar_fixedpoint_randGen_1}) 
 and  assumption {\bf B.1},  the    finite $n$-systems
 have  aggregate fixed points  (in almost sure sense), 
 using Brouwer's fixed point Theorem.   
 Further the limit system has %unique 
 (aggregate) fixed  point as given by assumption {\bf B.5}.  
 Thus we define the following function whose zeros are the fixed points of the mappings ${\bar {\bf f}}^{n,m}$ given by (\ref{Eqn_bar_fixedpoint_randGen_1_more})   for any $n \le \infty$; 
 %{\color{blue} And, also by  Lemma \ref{Lemma_Contraction 1},  Assumption {\bf{B}}.2, and Corollary  \ref{Corrolary_supportinglemma} we have  almost sure existence of aggregate fixed points for all systems and uniqueness of the same for the limit system.} 
% In the next, we complete the remaining proof of Theorem.
define ${\cal N} = \{1, 2, \cdots, \infty\}$  and define\footnote{By boundedness assumption {\bf B.1},
  ${\bar f}(.)$'s are also bounded and hence the limit exists and so the  function $h$ is well defined for any $( {\bar x}_b, \x, \infty)$. } the   real valued function $h(\cdot)$ from $\left ( [0,   y] \times  s^\infty \times s^\infty \right ) \times {\cal N} \to  {\cal R}^+ $ as below    (recall $|{\cal G}_1 | + |{\cal G}_2 |= n$)
\begin{eqnarray*}
h({\bar x}_b ,\x, n) &=& | {\bar f}^{n}_b ({\bar x}_b, \x) - {\bar x}_b | + \sum_{m=1}^2 \left ( \sum_{i \in \mathcal{G}_m} 2^{-i}  | {\bar f}^{n,m}_i ({\bar x}_b,\x)  - {\bar x}^m_i  |   + \sum_{i > n} 2^{-i}   {\bar x}_i^m\right )
\\
h({\bar x}_b, \x, \infty) &=&  %\lim_{n \to \infty} 
\left ( |{\bar f}^{\infty}_b ({\bar x}_b,\x) - {\bar x}_b |   +  \sum_{m=1}^2 \sum_{i \in \mathcal{G}_m} 2^{-i}  | {\bar f}^{\infty,m}_i ({\bar x}_b, \x)  - {\bar x}^m_i  | \right ) .
\end{eqnarray*}
%It is clear that the zeros of the above functions are the fixed points of the mappings ${\bar {\bf f}}^{n,m}$ given by (\ref{Eqn_fnm})   for any $n \le \infty.$  And since the systems have fixed points, 
It is clear that the aggregate fixed points form the minimizers of the above functions. Our idea is to obtain the  convergence proof using continuity of optimizers as given by Maximum Theorem (e.g., \cite{Berge,Equilibria,Feinberg}). Towards this we begin with joint continuity of the objective function. \\
\noindent \underline{Joint continuity:} 
We require some  (sample path-wise)    inequalities  for showing the joint continuity of the  objective function and we begin with the same.
Towards this, we consider the set ${\cal C}$ of Corollary \ref{Corrolary_supportinglemma} (provided later in the same Appendix).
%here for convenience:
%\begin{eqnarray*}
%{\cal C}   =  \Bigg  \{ w :   \frac{1}{n}\sum_{j \in \mathcal{G}_1}  W_{j,b} (w) \stackrel{ n \to \infty}{\to}   \gamma p_1^{sb}  , \  \frac{1}{n}\sum_{j \in \mathcal{G}_2}  W_{j,b} (w) \stackrel{ n \to \infty}{\to}  (1-\gamma)  p_2^{sb}, \hspace{25mm}\\  %\\ & \hspace{-26mm}
%\sum_{j \in \mathcal{G}_1}  \frac{ 1 - \eta_{j}^{sb} } {\sum_{i' \in \mathcal{G}_1 \cup \mathcal{G}_2} I_{j, i'} }    \to  
%\frac{\gamma(1- p_1^{sb})}{\gp{1}}    \mbox{ and } 
%\sum_{j \in \mathcal{G}_2}  \frac{ 1 - \eta_{j}^{sb} } {\sum_{i' \in \mathcal{G}_1 \cup \mathcal{G}_2} I_{j, i'} } \to \frac{(1-\gamma) (1- p_2^{sb})}{\gp{2}}
%  \mbox{ and }\sum_{j \le n}   W_{j,i}  (w) \stackrel{ n \to \infty}{\to} 1-  p_{sb}  \mbox{ for all } i.  
 % \Bigg   \}.
%\end{eqnarray*}
By Corollary  \ref{Corrolary_supportinglemma} (under 
assumption {\bf B.2}),  $P({\cal C})=1$ and
%{\color{blue}Note by Lemma \ref{lem: LLN_fixed point} , Assumption {\bf B}.2,  and by using Corollary  \ref{Corrolary_supportinglemma}  $P({\cal C})=1$.} 
for any  $w \in {\cal C}$, there exists   an $ N_w < \infty$  
%(because of the convergences defined in  the proof of Lemma \ref{Master lemma}),
 such that for all $n \ge N_w$:

\vspace{-4mm}
{\small
\begin{eqnarray}
\label{eqn_inequality_joint1}
 \frac{1}{n} \sum_{j\in \mathcal{G}_1 } W_{j,b} + \frac{1}{n}  \sum_{j \in \mathcal{G}_2}  W_{j,b} \ \ \  &\le&  \  \   \  2  \gamma p_1^{sb} + 2 (1-\gamma) p_2^{sb}  \mbox{\normalsize,  and, }   \\
 \sum_{j \in \mathcal{G}_1}  \frac{1 - \eta_j^{sb}}{ \sum_{i' \in \mathcal{G}_1 \cup \mathcal{G}_2} I_{j, i'}} +  \sum_{j \in \mathcal{G}_2}  \frac{1 - \eta_j^{sb}}{ \sum_{i' \in \mathcal{G}_1 \cup \mathcal{G}_2} I_{j, i'}} 
  &\le&  \  \   \  %2 \lim_{n' \to \infty }   \sum_{j \le n'}  \frac{1 - \eta_j^{sb}}{ \sum_{i' \le n'} I_{j, i'}} =
  2  \bigg(   \frac{\gamma(1- p_1^{sb}) }{\gp{1}  } +  \frac{(1-\gamma)(1- p_2^{sb}) }{  \gp{2}}\bigg). \hspace{4mm}
  \label{eqn_inequality_joint2}
\end{eqnarray}
}
Thus for all  such  $n$ we have from equation (\ref{Eqn_Weights 1}), as $I_{j,i} \le 1$:

\vspace{-4mm}
{\small
\begin{eqnarray}
\label{eqn_inequality_forcontinuity}
\sum_{i \in \mathcal{G}_1 \cup {\mathcal G}_2} 2^{-i} \sum_{m=1}^2  \sum_{j \in \mathcal{G}_{m}} W_{j,i}    
& \le &    \sum_{i \in \mathcal{G}_1 \cup {\mathcal G}_2} 2^{-i} \bigg(\sum_{j \in \mathcal{G}_1 }  \frac{1 - \eta_j^{sb}}{ \sum_{i' \in \mathcal{G}_1\cup \mathcal{G}_2} I_{j, i'}}+ \sum_{j \in \mathcal{G}_2 }  \frac{1 - \eta_j^{sb}}{ \sum_{i' \in \mathcal{G}_1\cup \mathcal{G}_2} I_{j, i'}}\bigg) \nonumber    \\
& \le &  \hspace{-1 mm} \sum_{i \le \infty} 2^{-i}    
 2\bigg ( \frac{ \gamma(1-  p_1^{sb}) }{\gp{1} } + \frac{(1-\gamma) (1-  p_2^{sb} )}{\gp{2}}\bigg )\nonumber \\
 &=:&  2c < \infty,
 % \left (  2 \lim_{n'\to \infty }   \sum_{j \le n'}  \frac{1 - \eta_j^{sb}}{ \sum_{i' \le n'} I_{j, i'}}  + c'(w) \right ) := c(w)  <  \infty.
\end{eqnarray} 
}
where, $c$ is an appropriate constant.

In the second step we will show that $h(.)$ is a Lipschitz continuous function in $({\bar x}_b ,\x)$  for all such $w \in {\cal C}$    and that the co-efficient of Lipschitz continuity can be the same for all $n  \ge N_w$ and for $\infty$. 
For any $({\bar x}_b ,\x)$ and $({\bar u}_b,\u)$, using Lipschitz continuity assumption {\bf B.3} (note $\varsigma \le 1$),  we have\footnote{use triangular inequality, and then use inequalities like $|a|-|b| \le |a-b|$ etc.} 
 for all $n  > N_w$  (see equation (\ref{Eqn_bar_fixedpoint_randGen_1}) and \eqref{Eqn_infinynorm}): 

\vspace{-3mm}
{\small \begin{eqnarray}
\label{Eqn_Lip_cont}
| h({\bar x}_b ,\x, n)  - h({\bar u}_b, \u,  n) |  \nonumber \hspace{10mm}  \\ &  & \nonumber \hspace{-40mm}  \le \ \sum_{m=1}^2\sum_{i \in  \mathcal{G}_m } 2^{-i}   | {\bar f}^{n,m}_i ({\bar x}_b ,\x)  -  {\bar f}^{n,m}_i (  {\bar u}_b ,\u) | \\ \nonumber
&& \hspace{-35mm}
  + |  {\bar f}^{n}_b ({\bar x}_b, \x)-  {\bar f}^{n}_b ({\bar u}_b,\u) |  + | {\bar x}_b - {\bar u}_b | 
  %\\
%
%&& \hspace{-20mm}
+ \  \sum_{m=1}^2\sum_{i \in \mathcal{G}_m } 2^{-i}  | {\bar x}^m_i -  {\bar u}^m_i  | +
 \sum_{i >n}  2^{-i}  | {\bar x}^m_i -  {\bar u}^m_i  |
\\
&& \hspace{-40mm}
\le \  2 \sigma || ({\bar x}_b ,\x)  - ({\bar u}_b,\u  )||_\infty \hspace{-2mm}  \sum_{i \in  \mathcal{G}_1 \cup \mathcal{G}_2} 2^{-i} \hspace{-2mm} \sum_{j\in \mathcal{G}_1 \cup \mathcal{G}_2}
\hspace{-2mm} W_{j,i} 
%\\
%
%&& \hspace{-35mm}
+ 2   || ({\bar x}_b, \x)  - ( {\bar u}_b, \u)||_\infty    \left  (1+ \hspace{-2mm} \sum_{i \in \mathcal{G}_1 \cup \mathcal{G}_2 } \hspace{-2mm} 2^{-i} %+ 2 \sum_{i >n}  2^{-i} 
\right  ) \nonumber \\
&& \hspace{-35mm}
  + 2 \sigma || ({\bar x}_b,\x )  - ({\bar u}_b,\u  )||_\infty  \left (  \frac{1}{n} \sum_{j\in \mathcal{G}_1 } W_{j,b} + \frac{1}{n} \sum_{j\in \mathcal{G}_2 } W_{j,b}  
  \nonumber
  %\\
%
%&& \hspace{-32mm}
 \right )  +   || ({\bar x}_b,\x )  - ({\bar u}_b,\u  )||_\infty c'
 \\
&& \hspace{-40mm} \le  \ 4   || ({\bar x}_b ,\x)  - ({\bar u}_b ,\u )||_\infty   \left  (  c +   \gamma p_1^{sb} + (1-\gamma) p_2^{sb}+ 1  + \frac{c'}{4}\right ),
\end{eqnarray}}where $c'$ is an appropriate constant.
The above inequality is due to the bounds \eqref{eqn_inequality_joint1}-\eqref{eqn_inequality_forcontinuity} and observe that the upper bound is uniform in $n$.
Using exactly similar logic, we derive Lipschitz continuity for $n=\infty.$
By definition of limit system \eqref{Eqn_bar_fixedpoint_limtGen 1},  when the limit superior becomes limit, we  have:
$$
| h({\bar x}_{b}, \x,   n)  -  h({\bar x}_b,\x,\infty) | \to 0 \mbox{  as $n \to \infty$}.
$$Thus and using \eqref{Eqn_Lip_cont},
   if $({\bar x}_{n,b},\x_n,  n) \to  ({\bar x}_b, \x,  \infty) $, i.e., if  
$|| ( {\bar x}_{n,b}, \x_n) -  ({\bar x}_b,\x) ||_\infty \to 0 $  as $n \to \infty$,  with the limit $({\bar x}_b,\x) \in {\cal \tilde{D}} $, where,  
$$
{\cal \tilde{D}} := \bigg  \{ ({\bar x}_b , {\x} ) \in [0,   y] \times  s^\infty \times s^\infty  :  \mbox{ limit superior in \eqref{Eqn_bar_fixedpoint_limtGen 1} equals the limit}  \bigg \},
$$then, \vspace{-6mm} 
\begin{eqnarray}
\label{Eqn_h_conv} \hspace{1mm}
| h({\bar x}_{n,b}, \x_n,   n)  -  h({\bar x}_b,\x,\infty) | \hspace{-35mm} \nonumber\\
&\le& | h({\bar x}_{n,b}, \x_n,   n)  -  h({\bar x}_b,\x, n) | + | h({\bar x}_{b}, \x,   n)  -  h({\bar x}_b,\x,\infty) |
\to 0. 
\end{eqnarray} This implies joint (norm) continuity of the function $h(.)$    on set ${\cal \tilde{D}} \times {\cal N} $, with\footnote{\label{footnote_top_N}Equip ${\cal N}$ with topology obtained using Euclidean metric on ${\cal R}$ as well as include the following additional open sets around `$\infty$', $\{N, N+1, \cdots, \infty\}$  for any given $N < \infty$. Observe that any compact set in this topology is a set with finite number of elements or complement of a finite set, when it contains $\infty$. It is almost like  discrete topology, except that $\{\infty\}$ is not an open set.} ${\cal N} := \{1, 2, \cdots, \infty\}$. 
Observe that the fixed points of the finite systems (sequence has 0's after finite $n$) and constant sequence are easily in ${\cal \tilde{D}}$ by Lemma \ref{lem: LLN_fixed point}, and hence by {\bf B.5}  it suffices to consider minimizing $h (\cdot)$ over ${\cal \tilde{D}}$  to study the   required fixed points.

\ignore{
{\bf New Changes at one place !}
{\bf B.5} There exists at least one fixed point for the limit system \eqref{Eqn_bar_fixedpoint_limtGen 1} among constant sequences. 

Consider the following   domain  of optimization:
$$
{\cal D} := \bigg  \{ ({\bar x}_b , {\x} ) \in [0,   y] \times  s^\infty \times s^\infty  :  \mbox{ limit superior in \eqref{Eqn_bar_fixedpoint_limtGen 1} equals the limit}  \bigg \}.
$$
Observe that the fixed points of the finite systems (sequence has 0's after finite $n$) and constant sequence are easily in ${\cal D}$ by Lemma \ref{lem: LLN_fixed point}, and hence by {\cal B}.5  it suffices to consider this as the domain.   }

We now apply Maximum Theorem for non-compact sets given by \cite[Theorem 1.2]{Feinberg} to complete the proof: 

\noindent$\bullet$
We consider weak topology on $[0,   y] \times  s^\infty \times s^\infty $ generated\footnote{In such  a topology we say $({\bar x}_{b,n}, \x_n) \to ({\bar x}_b, \x) $ if and only if
$$
{\bar x}_{b,n} \to {\bar x}_{b} \mbox{ and }   x^{m}_{n,i}  \to  x^m_i
\mbox{ for all } (m,i).$$
} by the set of projections, i.e., the smallest topology that ensures the ($m,i$)-th projection mapping $({\bar x}_b , \x) \mapsto x^m_i$ from $[0,   y] \times  s^\infty \times s^\infty$ to ${\cal R}$  a continuous mapping for any  $(m, i)$. This is basically the product topology  and by the well known Tychonoff's theorem,   $[0,   y] \times  s^\infty \times s^\infty$ is a compact set under it.

\noindent$\bullet$ The set of the parameters ${\cal N} $ with topology as in footnote \ref{footnote_top_N} is clearly compactly generated topological space;

\noindent$\bullet$
We have a constant correspondence $\Phi$: for any $n \in {\cal N}$ the domain is $\Phi(n) = {\cal \tilde{D}}$;

\noindent$\bullet$
By definition of  ${\cal \tilde{D}}$ and by \eqref{Eqn_h_conv} the function $h (\cdot, \cdot)$ is strong continuous on ${\cal \tilde{D}}$, and hence is weak continuous (under product topology);
 and 

\noindent$\bullet$
Observe that for any compact set $K \subset {\cal N}$  the graph $Gr_K (\Phi) =    K \times {\cal \tilde{D}}$ (definitions as in \cite{Feinberg}). 
 Consider the following level set, 
$$
\bigg\{(n, {\bar x}_b, \x) \in K \times  {\cal \tilde{D}} : h( {\bar x}_b, \x, n) \le l \bigg \}, \mbox{ for any compact set $K$ and level $l$.}
$$The above set 
is   weak closed  (by   weak continuity of $h$) and hence is  weak compact, as it is a subset of $K\times [0, y] \times  s^\infty \times s^\infty$. Recall 
the domain ${\cal \tilde{D}} $ is a subset of weak (product) compact set $[0,   y] \times  s^\infty \times s^\infty$.  Thus 
the function $h$ is $K$-inf compact on $Gr_X (\Phi ) = X \times {\cal \tilde{D}}$ (definitions in \cite{Feinberg}).

\ignore{
{\color{blue}
Define the following neighbourhood of $(( {\bar q}_b, \bar{\bf q} ) ) \in [0,y]\cap {\cal Q} \times s^\infty \times s^\infty$:
$$
{O}_{i_1, \cdots, i_n, \epsilon} (( {\bar q}_b, \bar{\bf q} ) ) 
:= \left \{ 
( {\bar x}_b, \x) : 
|x^m_{i_k} - q^m_{i,k} | < \epsilon  \mbox{ for all } i_k, m \in \{1, 2\} 
\right \}
$$}

\begin{lemma}
\label{Lemma_w_compact}
The set  $[0,   y] \times  s^\infty \times s^\infty $ is compact in the w-topology.
\end{lemma}
{\bf Proof:} 
It is sufficient to prove  sequential  compactness (as $w$-topology is second countable).
Consider any sequence ${\x}_{n} $  in the given set. We omit here the first component for obvious reasons. 
Consider a sub-sequence of this sequence call it ${\x}_{n,1} $ such that the first component $(m= 1 $ and $i=1$) converges. Call the first element of this sub-sequence as $\bar{ \bf y}_1$. Now pick a further sub-sequence such that the second component $(m= 2 $ and $i=1$)  also converges. Pick the first element of this sub-sequence other than $\bar{ \bf y}_1$ and call it $\bar{ \bf y}_2$. Continue this and generate $\{\bar {\bf y}_k \}_k$. It is easy to verify that this sub-sequence converges component wise. Further it is easy to check that it is bounded by $y$. Thus the proof. \eop}

\ignore{
\noindent \underline{Weak topology:} d\footnote{ We say vector $( {\bar x}_{nb} ,\x_n) $ converges weakly to $({\bar x}_b,\x)$ (\cite{Berge,Equilibria}) if 
$$
 {\cal L } ( {\bar x}_{nb}, \x_n ) \to {\cal L} ({\bar x}_b, \x )   
$$ for any linear functional ${\cal L}$.}
The strong continuity implies continuity in weak sense.  
Further, the bounded set $[0, y] \times s^{\infty}\times s^{\infty}$ is weak compact.  
Now applying Maximum Theorem  to the function $h(.)$ using the topology generated by weak continuity and when the domain of optimization,  $[0, y]$ $\times s^{\infty}\times s^{\infty}$ is same for all $n$,}
Thus by applying non-compact Maximum Theorem\footnote{We reproduce  here  \cite[Theorem 1.2]{Feinberg}  in  our notations. 
Let $ \big ({\mathbb Z} , \Gamma_{prod} \big )$ be a topological space with  (say product) topology $\Gamma_{prod}$, $S({\mathbb Z})$ be the power set of ${\mathbb Z}$ and $K({\mathbb Z})$ be the class of compact subsets of ${\mathbb Z}$. Assume that \begin{enumerate}[(i)]
    \item ${\mathbb X}$ is a compactly generated topological space;
    \item $\Phi:{\mathbb X} \rightarrow S({\mathbb Z})$ is lower semi-continuous;
    \item  $h: {\mathbb X} \times {\mathbb Z} \rightarrow  \mathcal{R}$ is K-inf-compact and upper semi-continuous on $Gr_X(\Phi)$.
\end{enumerate}
Then the value function $v: {\mathbb X} \rightarrow \mathcal{R}$, defined by $v(x) := \max_{ {\bf z} \in \Phi(x)} h(x, {\bf z})$, is continuous and the solution   $\Phi^{*}:{\mathbb X} \rightarrow
K({\mathbb Z})$, where $\Phi^{*}(x) := \arg \max_{{\bf z} \in \Phi(x)} h(x, {\bf z})$,   is upper semi-continuous and compact-valued.}  \cite[Theorem 1.2]{Feinberg} to the function $h(.)$ 
we obtain that the set of minimizers of $h(.)$ form a upper hemi-continuous (with respect to $n$) compact correspondence,  under the  product-topology (\cite{Berge,Equilibria}).
%Let  the  unique fixed point of the limit system be represented by  $({\bar x}^{\infty*}_{b} ,\x^* )$.
By properties of upper hemi-continuous (with respect to $n$) compact correspondence (\cite{Berge,Equilibria}) we have the following result:
a) consider any sequence of  numbers $n_k \to \infty$ (it can also be $n \to \infty$);  b) consider (any) one fixed point for each $n_k$, call it 
$\big ({\bar x}^*_{b} (n_k) ,\x^*(n_k) \big )$; c) then we have the following  convergence along a sub-sequence $n_{k_l} \to \infty$ (under product topology): 
$$
\Big( {\bar x}^*_{b} (n_{k_l})  ,\x^* (n_{k_l}) \Big )  \xlongrightarrow{\text{component-wise}} \Big ( {\bar x}^{\infty*}_{b}, \x^{\infty*} \Big ),
$$where  $({\bar x}^{\infty*}_{b} ,\x^* )$ is a fixed point of the limit system. 
%Since the projection is a linear functional we  have the required result.  
The last statement of the Theorem is immediate once we have the convergence of the aggregate fixed points, component-wise. 
\eop
%\noindent \textbf{Note:} $\mathcal{D}$  is closes under strong topology but not closed under weak topology.
 %
    \begin{cor}
    \label{Corrolary_supportinglemma}
 The following set has measure  1 under the assumption $\textbf{B.2}$
\begin{eqnarray*}
{\cal C}   =  \Big  \{ w :   \frac{1}{n}\sum_{j \in \mathcal{G}_1}  W_{j,b} (w) \stackrel{ n \to \infty}{\to}   \gamma p_1^{sb}  , \frac{1}{n}\sum_{j \in \mathcal{G}_2}  W_{j,b} (w) \stackrel{ n \to \infty}{\to}  (1-\gamma)  p_2^{sb} \\  %\\ & \hspace{-26mm}
\sum_{j \in \mathcal{G}_1}  \frac{ 1 - \eta_{j}^{sb} } {\sum_{i' \in \mathcal{G}_1 \cup \mathcal{G}_2} I_{j, i'} }    \to  
\frac{\gamma(1- p_1^{sb})}{\gp{1}}    \mbox{ and }\\ 
 \sum_{j \in \mathcal{G}_2}  \frac{ 1 - \eta_{j}^{sb} } {\sum_{i' \in \mathcal{G}_1 \cup \mathcal{G}_2} I_{j, i'} } \to \frac{(1-\gamma) (1- p_2^{sb})}{\gp{2}}
%  \mbox{ and }\sum_{j \le n}   W_{j,i}  (w) \stackrel{ n \to \infty}{\to} 1-  p_{sb}  \mbox{ for all } i.  
  \Big   \}.
\end{eqnarray*}
\end{cor}
 \noindent {\textbf{Proof:}} 
 It follows directly by law of large numbers (for the first two quantities of the set $\mathcal{C}$) and by 
   Lemma \ref{Master lemma} with $M_j \equiv 1$ (for the last two quantities). 
 \eop 
 
\noindent \textbf{Proof of Theorem \ref{Lemma_limit_system_uniquness}:} We  consider  finite $n$-system and consider the following norm for this proof:
 \begin{eqnarray}
\label{Eqn_L_1norm_finitensystem}
|| (\x, {\bar x}_b) ||_1 := \frac{1}{n}
\sum_m \sum_{j \in {\cal G}_m}\left ( |x^m_j|+ \varsigma  |{\bar x}_b) | \right ) 
\end{eqnarray}
Observe that for any $ ( {\bar x}_b, \x ) $ (as $|{\cal G}_1|+|{\cal G}_2| = n$),

\vspace{-2mm}
{\small
\begin{eqnarray*}
||{\bar {\bf f}}^n ( {\bar x}_b, \x )  ||_1  \ =  \ \frac{1}{n}\sum_m 
\sum_{i \in {\cal G}_m}  \left ( |{\bar f}^{n, m}_i ( {\bar x}_b, \x )|+ \varsigma  |  {\bar f}^{n}_b ( {\bar x}_b, \x ) | \right ) \hspace{-75mm} \\ \\
&=&  \frac{1}{n}\sum_m 
\sum_{i \in {\cal G}_m}  \Bigg (  \sum_{j \in \mathcal{G}_1}   \xi^1_{j} ({\bar x}^1_j, x_b)  W_{j, i} + \displaystyle \sum_{j \in \mathcal{G}_2}   \xi_{j}^2  ({\bar  x}^2_j, x_b)   {W}_{j, i}   
\\
&& \hspace{30mm}
+ \varsigma \frac{1}{n}  \sum_{j \in \mathcal{G}_1}   \xi^{1}_{j} ({\bar x}^{1}_j, x_b) W_{j, b} + \varsigma  \frac{1}{n} 
\displaystyle \sum_{j \in \mathcal{G}_2}   \xi^{2}_{j} ({\bar x}^{2}_j, x_b) W_{j, b} \Bigg ) \\
&=&   \frac{1}{n}     \sum_{j \in \mathcal{G}_1}   \xi^1_{j} ({\bar x}^1_j, x_b)  \left ( \sum_m 
\sum_{i \in {\cal G}_m}  W_{j, i} +\varsigma W_{j,b} \right ) +   \frac{1}{n}  \sum_{j \in \mathcal{G}_2}   \xi_{j}^2  ({\bar  x}^2_j, x_b) \left (  \sum_m 
\sum_{i \in {\cal G}_m}   {W}_{j, i}   +\varsigma W_{j,b} \right ) 
 \\
&=&   \frac{1}{n}     \sum_{j \in \mathcal{G}_1}   \xi^1_{j} ({\bar x}^1_j, x_b) \left ( 1- \eta_j^{sb} + \varsigma \eta^{sb}_j  \right )
+ \frac{1}{n}  \sum_{j \in \mathcal{G}_2}   \xi_{j}^2  ({\bar  x}^2_j, x_b)  \left ( 1- \eta_j^{sb} + \varsigma \eta^{sb}_j  \right ) .
\end{eqnarray*}}
In the similar lines we have, with $
{\eta}^{\varsigma}_j := \left ( 1- \eta_j^{sb} + \varsigma   \eta^{sb}_j  \right )$:

\vspace{-4mm}
{\small 
\begin{eqnarray}
\label{Eqn_contraction_final_stepnsystem}
||{\bar {\bf f}}^n ( {\bar x}_b, \x ) - {\bar {\bf f}}^n ( {\bar u}_b, \u )  ||_1  \nonumber\hspace{-22mm} \\
& \le & \frac{1}{n}     \sum_{j \in \mathcal{G}_1}  \left |  \xi_{j}^1  ({\bar  x}^1_j, x_b) - 
 \xi_{j}^1  ({\bar  u}^1_j, u_b) \right | {\eta}^{\varsigma}_j + \frac{1}{n}  \sum_{j \in \mathcal{G}_2} |  \xi_{j}^2  ({\bar  x}^2_j, x_b) - 
 \xi_{j}^2  ({\bar  u}^2_j, u_b) | {\eta}^{\varsigma}_j \nonumber \\
 &\le & \sigma \frac{1}{n}  \sum_m \sum_{j \in {\cal G}_m} \left ( | {\bar x}_j^m - {\bar u}_j^m| +\varsigma | {\bar x}_b -  {\bar u}_b| \right )   {\eta}^{\varsigma}_j\nonumber   \\ 
 & \le&  \sigma \left ( 1- \underline{\eta} + \varsigma   \underline{\eta} \right )  || ( {\bar x}_b, \x )   - ( {\bar u}_b, \u )  ||_1 .
\end{eqnarray}}
Thus the finite $n$ -system is a strict contraction mapping for any $n$ and hence has a unique fixed point,  under the given hypothesis.

To prove the  second part,  we consider the following norm for the limit system:

\vspace{-4mm}
{\small
\begin{eqnarray}
\label{Eqn_limit_systemwithvarsigma_norm}
|| ({\bar x}_b ,\x)- ({\bar u}_b, \u)  ||_\infty 
      &=&  \max_{m \in \lbrace 1, 2\rbrace} \sup_{i\in \mathcal{G}_m } \left (|\bar{x}^m_i - \bar{u}^m_i| +  \varsigma| {\bar x}_b -   {\bar u}_b |\right) . %\mbox{with} \ 0 < \varsigma < 1.
\end{eqnarray} }
For any    any $i\in \mathcal{G}_m$ and $m$ we have: 

\vspace{-4mm}
{\small
\begin{eqnarray}
\label{Eqn_limitsystem_contraction_finalstep}
 |{\bar {\bf f}}_i^{\infty,m} ( {\bar x}_b, \x ) - {\bar {\bf f}}_i^{\infty,m} ( {\bar u}_b, \u )  | + \varsigma| {\bar f}^\infty_b ({\bar x}_b,{\x}) - {\bar f}^\infty_b ({\bar u}_b,{\u}) | \hspace{-65mm} &&  \nonumber\\
& =& |\limsup_n({\bar {\bf f}_i}^{n,m} ( {\bar x}_b, \x ) - {\bar {\bf f}}_i^{n,m} ( {\bar u}_b, \u ) ) | +\varsigma |\limsup_n ( {\bar f}^n_b ({\bar x}_b,{\x}) - {\bar f}^n_b ({\bar u}_b,{\u})| \nonumber \\
& \le& \ \limsup_n \displaystyle \sum_{j  \in \mathcal{G}_1  } |  \xi^{1}_j ({\bar x}^{1}_j, {\bar x}_b ) - \xi^{1}_j ({\bar u}^{1}_j, {\bar u}_b ) |   W_{j,i}+ \ \limsup_n  \displaystyle \sum_{j  \in \mathcal{G}_2  } |  \xi^{2}_j ({\bar x}^{2}_j, {\bar x}_b ) - \xi^{2}_j ({\bar u}^{2}_j, {\bar u}_b ) |   W_{j,i} \nonumber \\ 
 & &  + \varsigma \limsup_n \frac{1}{n} \left ( \displaystyle \sum_{j  \in \mathcal{G}_1  } |  \xi^{1}_j ({\bar x}^{1}_j, {\bar x}_b ) - \xi^{1}_j ({\bar u}^{1}_j, {\bar u}_b ) |   W_{j,b}+  \displaystyle \sum_{j  \in \mathcal{G}_2  } |  \xi^{2}_j ({\bar x}^{2}_j, {\bar x}_b ) - \xi^{2}_j ({\bar u}^{2}_j, {\bar u}_b ) |   W_{j,b} \right) \nonumber \\
 & \stackrel{a}{ \leq }&  \ \sigma  \max_{m \in \lbrace 1, 2\rbrace} \sup_{i\in \mathcal{G}_m } \left (|\bar{x}^m_i - \bar{u}^m_i| +  \varsigma| {\bar x}_b -   {\bar u}_b |\right) \lim_n \sum_{j \in \mathcal{G}_1 \cup \mathcal{G}_2} W_{ji}\nonumber \\
 & &+\sigma \varsigma\max_{m \in \{1,2\}}\sup_{i \in {\cal G}_m } \left ( |{\bar x}^{m}_i - {\bar u}^{m}_i | +  \varsigma|\bar{x}_b -\bar{u}_b |\right )  
\lim_n\frac{1}{n}\left ( \sum_{j \in\mathcal{G}_1}    W_{j, b} + \sum_{j \in\mathcal{G}_2}    W_{j, b} \right )\nonumber \\
 & \stackrel{b}{\leq} &  \  \sigma \varrho_\varsigma  || ({\bar x}_b,\x) - ({\bar u}_b, \u)  ||_\infty  
  \end{eqnarray}}
where,
  
\vspace{-4mm} 
 {\small
 \begin{eqnarray*} \hspace{2mm}
\varrho_\varsigma &:=&  \max \bigg \lbrace \frac{\gamma p_{c_1} (1-p_1^{sb})  }{\gp{1}} + \frac{(1-\gamma)p_2  (1-p_2^{sb}) }{\gp{2} } ,  \ \ \  \ \  
   \frac{\gamma p_1   (1-p_1^{sb}) }{\gp{1}} + \frac{(1-\gamma)p_{c_2}  (1-p_2^{sb})  }{\gp{2} } \bigg \rbrace  \\ && + \varsigma\left ( \gamma p_1^{sb} + (1-\gamma) p_2^{sb} \right ).
 \end{eqnarray*}}In the above,   the  inequality  $a$ is due to assumption  \textbf{B.3} and   the   inequality
$b$ is by Lemma  \ref{Master lemma} (with $M_j = I_{j,i} $). Therefore it is a contraction mapping under the assumption of $\sigma \varrho_\varsigma  <1$ (from the hypothesis either $\sigma < 1$ or $\varrho <1$). Hence the proof follows.
\eop   
   \ignore{
\begin{cor}
Assume {\bf B.1} to {\bf B.4} and also assume,  either  $\varrho <1$ or  $\sigma < 1$ in {\bf B.3}. Then we have unique fixed point   of the limit system \eqref{Eqn_bar_fixedpoint_limtGen 1}, which is a constant sequence, i.e., we have 
${\bar x}_i^{m\infty*} = {\bar x}^{m\infty*} $ for all $i \in {\cal G}_m$ in equation \eqref{Eqn_conv_aggr}. This limit is the fixed point of the three dimensional system given by \eqref{Eqn_barf_limit 1}
\end{cor}
\noindent \textbf{Proof:} By  assumption {\bf B.1}  and using Brouwer's fixed point Theorem,
  the finite $n$-system (\ref{Eqn_bar_fixedpoint_randGen_1}) has a  fixed point.  Consider  any $w \in {\cal C}$ and consider the norm as  in \eqref{Eqn_infinynorm}.
 From equation (\ref{Eqn_bar_fixedpoint_limtGen 1}) and using assumption {\bf B.3},  we observe the following for any  $i \in \mathcal{G}_1$ (see equation \eqref{Eqn_bar_fixedpoint_randGen_1}): 
 
 \vspace{-2mm}
 {\small
  \begin{eqnarray*}
 | {\bar f}_i^{\infty,1} ({\bar x}_b, \x)  - {\bar f}_i^{\infty,1} ( {\bar u}_b, \u )  |  \hspace{6mm} \\
  & & \hspace{-40mm} \le \ \lim_n \displaystyle \sum_{j  \in \mathcal{G}_1  } |  \xi^{1}_j ({\bar x}_j, {\bar x}_b ) - \xi^{1}_j ({\bar u}_j, {\bar u}_b ) |   W_{j,i}+ \ \lim_n  \displaystyle \sum_{j  \in \mathcal{G}_2  } |  \xi^{2}_j ({\bar x}_j, {\bar x}_b ) - \xi^{2}_j ({\bar u}_j, {\bar u}_b ) |   W_{j,i}   \nonumber \\
  & &  \hspace{-40mm} \le  \sigma|| ({\bar x}_b,\x) - ({\bar u}_b ,\u )  ||_\infty  \lim _n  \sum_{ j \in \mathcal{G}_1  \cup \mathcal{G}_2}W_{j,i}\\
& & \hspace{-40mm} \leq  \   || ({\bar x}_b,\x) - ({\bar u}_b, \u)  ||_\infty \bigg(  (1- p_1^{sb})  \frac{\gamma p_1 }{\gp{1}} + (1- p_2^{sb}) \frac{(1-\gamma)p_{c_2} }{ \gp{2}} \bigg)  . \end{eqnarray*}}The first equality holds due to assumption  \textbf{B.3} and   the last inequality
%the $\limsup$ becomes limit by the assumption  \textbf{B.2} and 
is by Lemma  \ref{Master lemma} (with $M_j = I_{j,i} $).
%{\color{blue}
%\vspace{-2mm}
%{\small
 % \begin{eqnarray*}
 %| {\bar f}_i^{\infty,2} ( \x, {\bar x}_b )  - {\bar %f}_i^{\infty,2} ( \u, {\bar u}_b )  | \hspace{9mm}  \\
  %& & \hspace{-40mm} \le \ \lim \sup_n  \sum_{j  \in \mathcal{G}_1  } |  \xi^{1}_j ({\bar x}_j, {\bar x}_b ) - \xi^{1}_j ({\bar u}_j, {\bar u}_b ) |   W_{j,i}+ \ \lim \sup_n  \sum_{j  \in \mathcal{G}_2  } |  \xi^{2}_j ({\bar x}_j, {\bar x}_b ) - \xi^{2}_j ({\bar u}_j, {\bar u}_b ) |   W_{j,i}   \nonumber \\
 % & &  \hspace{-40mm} \le  \sigma|| (\x, {\bar x}_b) - ( \u, {\bar u}_b )  ||_\infty  \lim \sup _n  \sum_{ j \in \mathcal{G}_1  \cup \mathcal{G}_2}W_{ji}\\
%& & \hspace{-40mm} \leq  \   || (\x, {\bar x}_b) - ( \u, {\bar u}_b )  ||_\infty  \bigg((1- p_1^{sb})  \frac{\gamma p_{c_1} }{\gp{1}} + (1- p_2^{sb})\frac{(1-\gamma)p_2 }{ \gp{2}} \bigg)    . 
%\end{eqnarray*} }
%}
 In a similar way,
 
 {\small \begin{eqnarray*}
   | {\bar f}_i^{\infty,2} ( {\bar x}_b, \x)  - {\bar f}_i^{\infty,2} ({\bar u}_b, \u )  |  
&   \leq  &   || ({\bar x}_b,\x) - ({\bar u}_b,\u)  ||_\infty  \bigg((1- p_1^{sb})  \frac{\gamma p_{c_1} }{\gp{1}} + (1- p_2^{sb})\frac{(1-\gamma)p_2 }{ \gp{2}} \bigg) \\
&& \hspace{40mm}\mbox{\normalsize for any } i \in {\cal G}_2, \mbox{\normalsize and, } \\
 | {\bar f}^\infty_b ({\bar x}_b,\x)  - {\bar f}^\infty_b ({\bar u}_b,\u)  |   &\le &  
    || ({\bar x}_b ,\x) - ({\bar u}_b,\u)  ||_\infty \bigg (\gamma  p_1^{sb}+ (1-\gamma)  p_2^{sb}\bigg)  .
 \end{eqnarray*} }
  Thus from assumption \textbf{B.4}
   \begin{eqnarray*}
% \label{Eqn_f_contract}
|| {\bar {\bf f}}^\infty ({\bar x}_b, \x) - {\bar {\bf f}}^\infty ({\bar u}_b, \u )  ||_\infty  \ForTR{\hspace{-40mm} \\ &= &
\max \left \{  | {\bar f}^\infty_b ( \x, {\bar u}_b )  - {\bar f}^\infty_b ( \u, {\bar u}_b )  | , \ 
\sup_{i }  | {\bar f}^\infty_i ( \x, {\bar u}_b )  - {\bar f}^\infty_i ( \u, {\bar u}_b )  | \right \}  }{   }
&\le &   \varrho|| (\x, {\bar x}_b) - ( \u, {\bar u}_b )  ||_\infty \ .   
 \end{eqnarray*}  
 Therefore, ${\bf{\bar f}}^\infty$ is  a strict contraction mapping with contraction co-efficient if $\varrho < 1$.
 
 \eop
}

\noindent\textbf{Proof of Corollary \ref{Cor_Three dimensional approximation}:}  
It is easy to observe that single valued upper semi-continuous correspondence is a continuous mapping (see \cite[Theorem 9.12]{Sundaram}).
Thus, because of unique fixed points for all $n$ as well as limit system, the set of minimizers of $h(\cdot)$ defined in the proof of Theorem \ref{Thm_MainGen 1} is a continuous mapping (as $n \to \infty$)   and hence the corollary.  \eop
 
%\begin{lemma}
%\label{Lemma_Contraction 1}
% Assume either  i) $0 < p_1^{sb}   < 1$ and $0 < p_2^{sb}   < 1$;     or  ii) $\sigma < 1$ in {\bf B.3}. Also assume   {\bf B}.4.
% The $n$-system (\ref{Eqn_bar_fixedpoint_randGen_1})   has (aggregate) fixed points, almost  surely.  
%Fix any $w \in {\cal C}$ defined in Corollary \ref{Corrolary_supportinglemma}. The limit  (aggregate) system  ${\bar f}^\infty$ given by (\ref{Eqn_bar_fixedpoint_limtGen 1}) for this sample  path  is a strict contraction  and  hence has unique fixed point.
 %\end{lemma}
 \ignore{
\noindent \textbf{Proof of Theorem \ref{Lemma_limit_system_uniquness}:} By  assumption {\bf B.1}  and using Brouwer's fixed point Theorem,
  the finite $n$-system (\ref{Eqn_bar_fixedpoint_randGen_1}) has a  fixed point.  Consider  any $w \in {\cal C}$ and consider the norm as  in \eqref{Eqn_infinynorm}.
 From equation (\ref{Eqn_bar_fixedpoint_limtGen 1}) and using assumption {\bf B.3},  we observe the following for any  $i \in \mathcal{G}_1$ (see equation \eqref{Eqn_bar_fixedpoint_randGen_1}) 
 
 \vspace{-2mm}
 {\small
  \begin{eqnarray*}
 | {\bar f}_i^{\infty,1} ({\bar x}_b, \x)  - {\bar f}_i^{\infty,1} ( {\bar u}_b, \u )  |  \hspace{6mm} \\
  & & \hspace{-40mm} \le \ \lim_n \displaystyle \sum_{j  \in \mathcal{G}_1  } |  \xi^{1}_j ({\bar x}_j, {\bar x}_b ) - \xi^{1}_j ({\bar u}_j, {\bar u}_b ) |   W_{j,i}+ \ \lim_n  \displaystyle \sum_{j  \in \mathcal{G}_2  } |  \xi^{2}_j ({\bar x}_j, {\bar x}_b ) - \xi^{2}_j ({\bar u}_j, {\bar u}_b ) |   W_{j,i}   \nonumber \\
  & &  \hspace{-40mm} \le  \sigma|| ({\bar x}_b,\x) - ({\bar u}_b ,\u )  ||_\infty  \lim _n  \sum_{ j \in \mathcal{G}_1  \cup \mathcal{G}_2}W_{j,i}\\
& & \hspace{-40mm} \leq  \   || ({\bar x}_b,\x) - ({\bar u}_b, \u)  ||_\infty \bigg(  (1- p_1^{sb})  \frac{\gamma p_1 }{\gp{1}} + (1- p_2^{sb}) \frac{(1-\gamma)p_{c_2} }{ \gp{2}} \bigg)  . \end{eqnarray*}}The first equality holds due to assumption  \textbf{B.3} and   the last inequality
%the $\limsup$ becomes limit by the assumption  \textbf{B.2} and 
is by Lemma  \ref{Master lemma} (with $M_j = I_{j,i} $).
%{\color{blue}
%\vspace{-2mm}
%{\small
 % \begin{eqnarray*}
 %| {\bar f}_i^{\infty,2} ( \x, {\bar x}_b )  - {\bar %f}_i^{\infty,2} ( \u, {\bar u}_b )  | \hspace{9mm}  \\
  %& & \hspace{-40mm} \le \ \lim \sup_n  \sum_{j  \in \mathcal{G}_1  } |  \xi^{1}_j ({\bar x}_j, {\bar x}_b ) - \xi^{1}_j ({\bar u}_j, {\bar u}_b ) |   W_{j,i}+ \ \lim \sup_n  \sum_{j  \in \mathcal{G}_2  } |  \xi^{2}_j ({\bar x}_j, {\bar x}_b ) - \xi^{2}_j ({\bar u}_j, {\bar u}_b ) |   W_{j,i}   \nonumber \\
 % & &  \hspace{-40mm} \le  \sigma|| (\x, {\bar x}_b) - ( \u, {\bar u}_b )  ||_\infty  \lim \sup _n  \sum_{ j \in \mathcal{G}_1  \cup \mathcal{G}_2}W_{ji}\\
%& & \hspace{-40mm} \leq  \   || (\x, {\bar x}_b) - ( \u, {\bar u}_b )  ||_\infty  \bigg((1- p_1^{sb})  \frac{\gamma p_{c_1} }{\gp{1}} + (1- p_2^{sb})\frac{(1-\gamma)p_2 }{ \gp{2}} \bigg)    . 
%\end{eqnarray*} }
%}
 In a similar way,
 
 {\small \begin{eqnarray*}
   | {\bar f}_i^{\infty,2} ( {\bar x}_b, \x)  - {\bar f}_i^{\infty,2} ({\bar u}_b, \u )  |  
&   \leq  &   || ({\bar x}_b,\x) - ({\bar u}_b,\u)  ||_\infty  \bigg((1- p_1^{sb})  \frac{\gamma p_{c_1} }{\gp{1}} + (1- p_2^{sb})\frac{(1-\gamma)p_2 }{ \gp{2}} \bigg) \\
&& \hspace{40mm}\mbox{\normalsize for any } i \in {\cal G}_2, \mbox{\normalsize and, } \\
 | {\bar f}^\infty_b ({\bar x}_b,\x)  - {\bar f}^\infty_b ({\bar u}_b,\u)  |   &\le &  
    || ({\bar x}_b ,\x) - ({\bar u}_b,\u)  ||_\infty \bigg (\gamma  p_1^{sb}+ (1-\gamma)  p_2^{sb}\bigg)  .
 \end{eqnarray*} }
 {\color{red}
 \begin{eqnarray*}
&&\varrho_\varsigma := \max \left \{ (1- p_1^{sb})  \frac{\gamma p_{c_1} }{\gp{1}} + (1- p_2^{sb})\frac{(1-\gamma)p_2 }{ \gp{2}}   ,    \frac{\gamma p_1   (1-p_1^{sb}) }{\gp{1}} + \frac{(1-\gamma)p_{c_2}  (1-p_2^{sb})  }{\gp{2} } \right \} \\
&& + \varsigma \bigg (\gamma  p_1^{sb}+ (1-\gamma)  p_2^{sb}\bigg) 
\end{eqnarray*}

 }
 
 Thus from assumption \textbf{B.4}
   \begin{eqnarray*}
% \label{Eqn_f_contract}
|| {\bar {\bf f}}^\infty ({\bar x}_b, \x) - {\bar {\bf f}}^\infty ({\bar u}_b, \u )  ||_\infty  \ForTR{\hspace{-40mm} \\ &= &
\max \left \{  | {\bar f}^\infty_b ( \x, {\bar u}_b )  - {\bar f}^\infty_b ( \u, {\bar u}_b )  | , \ 
\sup_{i }  | {\bar f}^\infty_i ( \x, {\bar u}_b )  - {\bar f}^\infty_i ( \u, {\bar u}_b )  | \right \}  }{   }
&\le &  || (\x, {\bar x}_b) - ( \u, {\bar u}_b )  ||_\infty \ \varrho .   
 \end{eqnarray*}  
 %
 %\begin{eqnarray*}
 %\varrho: =\max \bigg \{  \frac{\gamma p_{c_1} (1- p_1^{sb})  }{\gp{1}} + \frac{(1-\gamma)p_2  (1- p_2^{sb}) }{\gp{2}} , \hspace{-40mm}&\\
 %&
  % \frac{\gamma p_1   (1- p_1^{sb}) }{\gp{1}} +
 %\frac{(1-\gamma)p_{c_2}  (1- p_2^{sb})  }{ %\gp{2}}, \gamma p_1^{sb} + (1-\gamma)  p_2^{sb}   \bigg \}   \le  1 .
 %\end{eqnarray*}
Therefore, ${\bf{\bar f}}^\infty$ is  a strict contraction mapping with contraction co-efficient if $\varrho < 1$.
As above,  for any given $n$, we have (for example):

{\small 
 \begin{eqnarray*}
 | {\bar f}_i^{n,1} ( {\bar x}_b,\x)  - {\bar f}_i^{n,1} ( {\bar u}_b ,\u)  |   
  & \le  &   \sigma|| ({\bar x}_b ,\x) - ({\bar u}_b, \u  )  ||_\infty     \sum_{ j \in \mathcal{G}_1  \cup \mathcal{G}_2}W_{j,i}. \end{eqnarray*}}
Observe that by Lemma  \ref{Master lemma} (with $M_j = I_{j,i} $) and under the given assumptions:
\begin{eqnarray*}
  \sum_{ j \in \mathcal{G}_1  \cup \mathcal{G}_2}W_{j,i}  \to \bigg(  (1- p_1^{sb})  \frac{\gamma p_1 }{\gp{1}} + (1- p_2^{sb}) \frac{(1-\gamma)p_{c_2} }{ \gp{2}} \bigg)
  \mbox{ a.s., as  } n \to \infty. \end{eqnarray*} 
{\color{red}
Thus one can find (under the extra hypothesis of the theorem) an ${\bar N} $ such that: 
$$
\sigma \sum_{ j \in \mathcal{G}_1  \cup \mathcal{G}_2}W_{j,i} <  1 \mbox{ for all }  n \ge {\bar N},
$$and  hence 
$  {\bf {\bar f}}^n $ is a strict contraction  for all such $n$.
Thus we have unique fixed points for all $ n \ge {\bar N}$. 
}
 If $\sigma < 1$ using similar logic, it is clear that not only the limit system even the system for all but finite $n$ is a strict contraction.  
{\color{red}
It is easy to observe that single valued upper semi-continuous correspondence is a continuous mapping (see \cite[Theorem 9.12]{Sundaram}) and hence  by  (eventual) uniqueness of the fixed points,  convergence given by equations  \eqref{Eqn_conv_aggr} and \eqref{Eqn_Act_Fixed_with_m}  is  along the original sequence, i.e., as $n \to \infty$. }

Further from \eqref{Eqn_barf_limit 1} and  {\bf B.3} and with  $\varrho < 1$  or with  $\sigma <1$ we have a unique fixed point even when restricted to class of constant sequences; thus the unique fixed point of the limit system is a constant sequence.
\eop
}
\ignore{
{\color{red}
Consider the $L_1$ norm as follows:
$$||{\bar x}_b, \bar{x} )||_1 := \frac{1}{n} \sum_m \sum_{i\in \mathcal{G}_m }
\left (
|x^{m}_j|+\varsigma |x_b|  \right ) 
$$
We consider the $L_1$ norm to prove the contraction of the finite $n$-system. Consider the following:
\begin{eqnarray*}
&&\sum_m \sum_{i \in \mathcal{G}_m} |{\bar f}^{n, m}_i ( {\bar x}_b, \bar{x} ) - {\bar f}^{n, m}_i ( {\bar u}_b, \bar{u})|_1\\
&=&\sum_m  \sum_{i \in \mathcal{G}_m} \bigg| \displaystyle \sum_{j \in \mathcal{G}_1}   \xi^1_{j} ({\bar x}^1_j, x_b)  W_{j, i} + \displaystyle \sum_{j \in \mathcal{G}_2}   \xi_{j}^2  ({\bar  x}^2_j, x_b)   {W}_{j, i} -  \displaystyle \sum_{j \in \mathcal{G}_1}   \xi^1_{j} ({\bar u}^1_j, u_b)  W_{j, i} - \displaystyle \sum_{j \in \mathcal{G}_2}   \xi_{j}^2  ({\bar  u}^2_j, u_b)   {W}_{j, i} \bigg|_1 \\
&\le&   \displaystyle \sum_m  \sum_{i \in \mathcal{G}_m}\bigg|  \displaystyle \sum_{j \in \mathcal{G}_1}   \xi^1_{j} ({\bar x}^1_j, x_b)   - \displaystyle \sum_{j \in \mathcal{G}_1}   \xi^1_{j} ({\bar u}^1_j, u_b) \bigg|  W_{j, i} + \displaystyle \sum_{i \in \mathcal{G}_m}\bigg|  \displaystyle \sum_{j \in \mathcal{G}_2}   \xi^2_{j} ({\bar x}^2_j, x_b)   - \displaystyle \sum_{j \in \mathcal{G}_2}   \xi^2_{j} ({\bar u}^2_j, u_b) \bigg|  W_{j, i} \\
&\le & \sigma  \left ( \sum_{j \in \mathcal{G}_1} \left (  |x^1_j - u^1_j| +  \varsigma |x_b - u_b|  \right ) \sum_m  \sum_{i \in \mathcal{G}_m} W_{j,i} +  \sum_{j \in \mathcal{G}_2}  \left ( |x^2_j - u^2_j| + \varsigma |x_b - u_b|  \right )\sum_m  \sum_{i \in \mathcal{G}_m} W_{j,i}  \right ) 
\\
&= & \sigma  \left ( \sum_{j \in \mathcal{G}_1} \left (  |x^1_j - u^1_j| + \varsigma |x_b - u_b|  \right ) (1-\eta_j^{sb}) +  \sum_{j \in \mathcal{G}_2}  \left ( |x^2_j - u^2_j| + \varsigma |x_b - u_b|  \right )  (1-\eta^{sb}_j)  \right ) 
\end{eqnarray*}
And we have:
\begin{eqnarray*}
  |{\bar f}^{n}_b ( {\bar x}_b, \bar{x} ) - {\bar f}^{n}_b ( {\bar u}_b, \bar{u})|_1 
  &\le & \sigma \frac{1}{n}  \sum_{j \in {\cal G}_1 \cup {\cal G}_2}   
   \left ( |x^1_j - u^1_j| + |x_b - u_b|  \right )  \eta^{sb}_j
\end{eqnarray*}
Thus we have:
\begin{eqnarray*}
\sum_m \sum_{i \in \mathcal{G}_m} \left ( |{\bar f}^{n, m}_i ( {\bar x}_b, \bar{x} ) - {\bar f}^{n, m}_i ( {\bar u}_b, \bar{u})|_1 +  \varsigma |{\bar f}^{n}_b ( {\bar x}_b, \bar{x} ) - {\bar f}^{n}_b ( {\bar u}_b, \bar{u})|_1  \right)  \hspace{-80mm} \\
&=& \sum_m \sum_{i \in \mathcal{G}_m} \left ( |{\bar f}^{n, m}_i ( {\bar x}_b, \bar{x} ) - {\bar f}^{n, m}_i ( {\bar u}_b, \bar{u})|_1 \right)  + \varsigma n    |{\bar f}^{n}_b ( {\bar x}_b, \bar{x} ) - {\bar f}^{n}_b ( {\bar u}_b, \bar{u})|_1  \\ 
&=& 
\sigma \varsigma  \sum_m  \sum_{j \in \mathcal{G}_m} \left (  |x^m_j - u^m_j| + |x_b - u_b|  \right )   (1-\eta_j^{sb})   
\\
&& \hspace{5mm}
+ \sigma \sum_m \sum_{j \in {\cal G}_m}   
   \left ( |x^m_j - u^m_j| + |x_b - u_b|  \right )  \eta^{sb}_j \\
   &=& \sigma   \sum_m \sum_{j \in {\cal G}_m}   
   \left ( |x^m_j - u^m_j| + \varsigma |x_b - u_b|  \right ) 
   \left ( (1- \eta^{sb}_j)
    + \varsigma  \eta^{sb}_j \right )
\end{eqnarray*}
The claim is the following using LLN:
\begin{eqnarray*}
\frac{1}{n}
\left | 
\sigma   \sum_m \sum_{j \in {\cal G}_m}   
   \left ( |x^m_j - u^m_j| + \varsigma  |x_b - u_b|  \right ) 
   \left ( (1- \eta^{sb}_j)
    + \varsigma  \eta^{sb}_j \right ) - \sigma   \sum_m \sum_{j \in {\cal G}_m}   
   \left ( |x^m_j - u^m_j| + \varsigma |x_b - u_b|  \right ) 
   \left ( (1-  p^{sb}_m)
    + \varsigma   p^{sb}_m \right ) \right |  \to 0
\end{eqnarray*}
Thus we have that for large enough $n$
\begin{eqnarray*}
|| f(x..) - f(u...) ||_1 &  \le  &
\frac{1}{n}
 \sigma   \sum_m \sum_{j \in {\cal G}_m}   
   \left ( |x^m_j - u^m_j| + \varsigma |x_b - u_b|  \right ) 
   \left ( (1- \eta^{sb}_j)
    + \varsigma  \eta^{sb}_j \right )   \\ 
    &\le & \frac{1}{n}
     \sigma   \sum_m \sum_{j \in {\cal G}_m}   
   \left ( |x^m_j - u^m_j| + \varsigma |x_b - u_b|  \right ) 
   \left ( (1-  p^{sb}_m)
    + \varsigma   p^{sb}_m \right ) + \varepsilon \\
    &=& || (x..)  - (u..)||_1 \sigma \max_m \left ( (1-  p^{sb}_m)
    + \varsigma   p^{sb}_m \right ) + \varepsilon \\
     &=& || (x..)  - (u..)||_1 \sigma \max_m \left ( (1-  p^{sb}_m)
    + \varsigma   p^{sb}_m \right ) + \frac{\varepsilon}{|| (x..)  - (u..)||_1} || (x..)  - (u..)||_1
\end{eqnarray*}

If $\sigma < 1$ we have strict contraction and we have uniqueness.}

{\color{red} For alternate system --
 \begin{eqnarray*} 
     || ({\bar x}_b ,\x)- ({\bar u}_b, \u)  ||_\infty 
      &=&  \max_{m,m' \in \lbrace 1, 2\rbrace} \sup_{i\in \mathcal{G}_m } \left (|\bar{x}^{m,m'}_i - \bar{u}^{m,m'}_i| + | {\bar x}_b -   {\bar u}_b |\right) \\
       &=&  \max_{m \in \lbrace 1, 2\rbrace} \sup_{i\in \mathcal{G}_m } \left (|\bar{x}^{m1}_i - \bar{u}^{m1}_i|+ |\bar{x}^{m2}_i - \bar{u}^{m2}_i| + | {\bar x}_b -   {\bar u}_b |\right)
\end{eqnarray*} 
}
}

%\begin{eqnarray*}
%| {\bar f}^n_b ({\bar x}_b,{\x}) - {\bar f}^n_b ({\bar u}_b,{\u}) |
%&\le & 
%\max_{m \in \{1, 2\} }
%\sup_{j \in {\cal G}_m } \left ( |{\bar x}^{m1}_j - {\bar u}^{m1} | + |{\bar x}^{m2}_j - {\bar u}^{m2} |  +|\bar{x}_b -\bar{u}_b |\right )  
%\frac{1}{n}\left ( \sum_{j \in\mathcal{G}_1}    W_{j, b} + \sum_{j \in\mathcal{G}_2}    W_{j, b} \right )
%\end{eqnarray*}

\section*{Appendix C: Proofs related to Sub-section \ref{sec_assumptions}}
\noindent \textbf{Proof of Lemma \ref{Lemma_uniformconvergence}:}
By uniform convergence
%\footnote{$A^n_j := {\displaystyle\sum_{i \in \mathcal{G}_1 \cup \mathcal{G}_2}} I_{ji}$} 
in $j \in {\cal G}_1$  of $A_j^n/n$, for any given $\epsilon >0$,  $\exists$ an $N_\epsilon$ such that:
 \begin{eqnarray*}
\bigg |\frac{A^{n}_j}{n}-\gp{1} \bigg| <   \epsilon,  \forall n > N_{\epsilon}, \forall j\in \mathcal{G}_1.
\end{eqnarray*}
This implies   the  following: 
 \begin{equation*}
 -\epsilon + \gp{1} <  \frac{A^{n}_j}{n}  < \gp{1} + \epsilon 
 ~~ \forall n > N_{\epsilon}, \forall j \in \mathcal{G}_1.
 \end{equation*}
Now consider the term:

\vspace{-2mm}
  {\small
  \begin{multline*}
\displaystyle \sum_{j \in \mathcal{G}_1}    \left|  
    \frac{1 }{ \displaystyle  \sum_{i \in \mathcal{G}_1 \cup \mathcal{G}_2} I_{j, i}}  - \frac{1}{ n \gp{1}}    \right |
   = \frac{1}{n} \displaystyle \sum_{j \in \mathcal{G}_1}    \left|  
    \frac{n\gp{1} -A^n_j }{\gp{1} A^n_j}      \right |  \leq \frac{1}{n} \displaystyle \sum_{j \in \mathcal{G}_1}     
    \frac{\epsilon }{ \gp{1} \frac{A^n_j}{n}}    \  \leq \frac{1}{n} \displaystyle \sum_{j \in \mathcal{G}_1}     
    \frac{\epsilon }{ \gp{1} ( \gp{1}-\epsilon)}    \\
       =  \frac{|{\cal G}_1 | }{n}   \frac{\epsilon }{ \gp{1} ( \gp{1}-\epsilon)}  \to \frac{\epsilon \gamma }{\gp{1} (\gp{1} -\epsilon)}.    \\
   \end{multline*}}From above, by letting  $\epsilon \to 0$  one can prove the following  (a.s.) convergence:
\begin{eqnarray*}
 \lim_{n \to \infty } \displaystyle \sum_{j \in \mathcal{G}_1}    \left|  
    \frac{1 }{  \displaystyle \sum_{i \in   \mathcal{G}_1 \cup  \mathcal{G}_2 } I_{j, i}}  - \frac{1}{ n\gp{1} }    \right |  =0 .
\end{eqnarray*}
    In the similar line  one can show that  $\frac{A^{n}_k}{n} \stackrel{ n \to \infty} {\to} \gp{2} ~ \mbox{uniformly in k }$ implies the following:
    \begin{eqnarray*}
      \lim_{n \to \infty } \displaystyle \sum_{j \in \mathcal{G}_2}    \left|  \frac{1 }{  \displaystyle \sum_{i \in   \mathcal{G}_1 \cup  \mathcal{G}_2 } I_{j, i}}  - \frac{1}{ n\gp{2} }    \right |  =0 .
    \end{eqnarray*}
Using these above two we have the almost sure convergence on the set $\mathcal{E}$.
\eop

\noindent {\bf Proof of Theorem \ref{Thm_onsetD}:}  The proof follows by observing that the conclusions of all the required lemmas hold almost surely on set ${\cal E}$, and we have $P(\mathcal{C}\cap \mathcal{E}) = P(\mathcal{E})$. Finally the Theorem \ref{Thm_MainGen 1} holds true on the set $\mathcal{C}$ which by Lemma \ref{Lemma_uniformconvergence} hold on ${\cal D}$.
\eop

\section*{Appendix D: Proofs related to the Section \ref{Alternate_graphical_model}}    
\noindent \textbf{Proof of Theorem \ref{Thm_MainGen 2}:} 
The steps of the proof are exactly as in that of Theorem \ref{Thm_MainGen 1}, we would only mention the differences here.
We use the following norm for joint continuity:
\begin{eqnarray*} 
     || ({\bar x}_b ,\x)- ({\bar u}_b, \u)  ||_\infty 
       &=&  \max_{m \in \lbrace 1, 2\rbrace} \sup_{i\in \mathcal{G}_m } \left (|\bar{x}^{m1}_i - \bar{u}^{m1}_i|+ |\bar{x}^{m2}_i - \bar{u}^{m2}_i| + | {\bar x}_b -   {\bar u}_b |\right),
\end{eqnarray*}
The functions whose zeros provide the required fixed points are now given by: 

\vspace{-4mm}
{\small
\begin{eqnarray*}
h({\bar x}_b ,\x, n) &=& | {\bar f}^{n}_b ({\bar x}_b, \x) - {\bar x}_b | + \sum_{m=1}^2 \bigg( \sum_{i \in \mathcal{G}_m} 2^{-i}  | {\bar f}^{n,m1}_i ({\bar x}_b,\x)  - {\bar x}^{m1}_i  | \\ 
& & \hspace{20mm} +  \sum_{i \in \mathcal{G}_m} 2^{-i} | {\bar f}^{n,m2}_i ({\bar x}_b,\x)  - {\bar x}^{m2}_i |  + \sum_{i > n} 2^{-i}   {\bar x}^{m1}_i + \sum_{i > n} 2^{-i}   {\bar x}^{m2}_i\bigg)  \\
h({\bar x}_b, \x, \infty) & = & 
|{\bar f}^{\infty}_b ({\bar x}_b,\x) - {\bar x}_b |   +  \sum_{m=1}^2 \bigg( \sum_{i \in \mathcal{G}_m} 2^{-i}  | {\bar f}^{\infty,m1}_i ({\bar x}_b, \x)  - {\bar x}^{m1}_i  | 
  \\
& & \hspace{40mm}+ \sum_{i \in \mathcal{G}_m} 2^{-i}  | {\bar f}^{\infty,m2}_i ({\bar x}_b, \x)  - {\bar x}^{m2}_i  |\bigg ), 
\end{eqnarray*}}and these are defined over subsets of  $[0,y] \times  s^\infty \times s^\infty \times  s^\infty \times s^\infty$.
The rest of the  details  are as  in Theorem \ref{Thm_MainGen 1}.
\eop\\

\noindent \textbf{Proof of Theorem \ref{Thm_modified_forThm_2_Alt_system}:} 
Here again, we mention only the differences with respect to the proof of Theorem \ref{Lemma_limit_system_uniquness}. We would require  the following norm to show  that  the finite $n$- system is a  contraction mapping: 
\begin{eqnarray}
|| (\x, {\bar x}_b) ||_1 := \frac{1}{n}
\sum_m \sum_{j \in {\cal G}_m}\left ( |x^{m1}_j|+ |x^{m2}_j|+ \varsigma  |{\bar x}_b) | \right ) .
\end{eqnarray}
After going through the computations as 
 in the equation \eqref{Eqn_contraction_final_stepnsystem} the contraction coefficient now equals:
 $$\sigma\max_m \bigg(      (1-\underline{\eta}+\underline{\eta}\varsigma)\lambda_m +(1-\lambda_m) 1_{p_{c_m}>0}\bigg) .$$
 Clearly when $\sigma (1-\underline{\eta}+\underline{\eta}\varsigma) < 1$, the above is also less than 1 and this proves unique fixed point for any finite $n$ system.

To prove the  second part,  we consider the following norm for the limit system:

\vspace{-4mm}
{\small
\begin{eqnarray}
|| ({\bar x}_b ,\x)- ({\bar u}_b, \u)  ||_\infty 
      &=&  \max_{m \in \lbrace 1, 2\rbrace} \sup_{i\in \mathcal{G}_m } \left (|\bar{x}^{m1}_i - \bar{u}^{m1}_i| + (|\bar{x}^{m2}_i - \bar{u}^{m2}_i| + \varsigma| {\bar x}_b -   {\bar u}_b |\right) . %\mbox{with} \ 0 < \varsigma < 1.
\end{eqnarray} }
By similar steps as  in  equation \eqref{Eqn_limitsystem_contraction_finalstep} of the Theorem \ref{Lemma_limit_system_uniquness}, the contraction coefficient now modifies to $\sigma \rho_\varsigma$ where,

\vspace{-2mm}
{\small
\begin{eqnarray*}
\varrho_\varsigma : &=&\max \bigg\lbrace  \lambda_1(1- p_1^{sb})+ \frac{1-\gamma}{\gamma}(1-\lambda_2) 1_{p_{c_2}> 0}, \frac{\gamma}{1-\gamma}(1-\lambda_1)1_{p_{c_1}> 0}  + \lambda_2 (1-p_2^{sb})\bigg \rbrace \nonumber\\
 && + \varsigma\left ( \gamma \lambda_1 p_1^{sb} + (1-\gamma) \lambda_2 p_2^{sb} \right ).
\end{eqnarray*}}
Thus we have $\sigma \varrho_\varsigma  <1$   (as from the hypothesis either $\sigma < 1$ or $\varsigma < 1$ and  $\varrho \le 1$) and hence the proof follows.
\eop

\noindent \textbf{Proof of Corollary \ref{Cor_Three dimensional approximation_Alt_system}:} The first part of the proof (that the limit is given by five dimensional system) follows exactly as in   Corollary \ref{Cor_Three dimensional approximation}. By Theorem \ref{Thm_MainGen 2} the aggregate fixed point of the limit system is obtained under `constant sequences', and thus we need to solve a five dimensional fixed point equation given by  the  set of equations \eqref{Alt_Eqn_barf_limit 1}-\eqref{Alt_Eqn_barf_limit e}.
Further, from the same set of equations, it is easy to observe that:
\begin{eqnarray*}
{\bar f}_i^{\infty,12}  ({\bar x}_b,\x)  &=& {\bar f}_i^{\infty,22}  ({\bar x}_b,\x) \frac{1- \gamma}{\gamma} \frac{1-\lambda_2}{\lambda_2} \frac{1}{(1- p^{sb}_2)}  1_{p_{c_2}> 0},
\\
{\bar f}_i^{\infty,21}  ( {\bar x}_b,\x) &=&  {\bar f}_i^{\infty,11}  ({\bar x}_b,\x) \frac{\gamma}{1-\gamma} \frac{1-\lambda_1}{\lambda_1} \frac{1}{(1- p^{sb}_1)}1_{p_{c_1}> 0}.
\end{eqnarray*}
Thus we have a reduction  in the dimension, i.e., the effective fixed point equation reduces to three dimensional fixed point equation \eqref{Alt_eqn_limit_system} in terms of $({\bar x}^b, {\bar x}^{\infty, 11}, {\bar x}^{ \infty, 22 } $), while, the remaining   aggregate components  are given by \eqref{Eqn_aggrragate_al_system}.
\eop

\section*{Appendix E: Proofs related to the Section \ref{sec_asymptotic}}    
 \noindent \textbf{Proof of Lemma \ref{Lemma_G2 default}:}  We consider the following scenario's for the $\mathcal{G}_2$ banks with $v_2 < k_{d2}$. The  aggregate clearing vector for the $\mathcal{G}_2$ banks satisfies (see   \eqref{eqn_avgclearingvectorG2}):
 \begin{eqnarray}
{\bar x}_2^{ \infty} &=& \bigg( \min \left \{\bar{y}_2, \bigg (k_{d2}- v_2+ {\bar x}_2^{ \infty} \bigg)^+\right \}w  \nonumber \\ 
&& +\min \left \{\bar{y}_2,\bigg(k_{u1}- v_1+{\bar x}_2^{ \infty}  \bigg)^+\right \}(1-w) \bigg) (1-p^{sb}_2)\lambda_2.  \label{Eqn_FP_G2}
\end{eqnarray}
 \noindent \textbf{Case 1:} First consider the case when downward shock can be absorbed i.e., default probability is  $ P^2_D = 0$. If we have $k_{d2}-v_2 + {\bar x}_2^{ \infty} \geq \bar{y}_2$ 
 then  the aggregate clearing vector  ${\bar x}_2^{ \infty} =\bar{y}_2 (1-p^{sb}_2)\lambda_2$ and the above condition simplifies to the bound: $$\lambda_2  \ge  \frac{\bar{y}_2 +v_2-k_{d2} }{\bar{y}_2 (1-p_2^{sb})}.$$ \\
\noindent \textbf{Case 2:} Consider the case in which only the banks that receive  the shock will default, i.e., $P^2_D=w$ and the corresponding aggregate clearing vector equals:
\begin{eqnarray*}
{\bar x}_2^{ \infty} = \bigg(\bar{y}_2(1-w)+ ({\bar x}_2^{ \infty}  +k_{d2}-v_2)w\bigg) (1-p_2^{sb})\lambda_2, \mbox{ and satisfies, }  \\
(k_{d2}- v_2 + {\bar x}_2^{ \infty}) <\bar{y}_2 \mbox{ and } (k_{u2}- v_2 + {\bar x}_2^{ \infty}) > \bar{y}_2.
\end{eqnarray*}
In this case, the aggregate clearing vector  reduces to (equals that in hypothesis):
\begin{eqnarray*}
{\bar x}_2^{ \infty} =  \bigg(\frac{ \bar{y}_2 (1-w) + w( k_{d2} -v_2)   }{ 1 -w (1-p_2^{sb})\lambda_2  } \bigg)(1-p_2^{sb})\lambda_2,
\end{eqnarray*}   and  using the same  in the bounds,   we have:
\begin{equation*}
\beta_0  < \lambda_2
	< \frac{\bar{y}_2-k_{d2}+v_2}{\bar{y}_2 (1-p_2^{sb})}.
\end{equation*}
Observe in the above if ${\bar y}_2  - w(k_{u2}-k_{d2}) < 0$, then we will not require the lower bound, and hence the indicator in $\beta_0$ definition.  
 
\noindent\textbf{Case 3:} Consider the case with all default i.e., $P^2_D=1$.
 We compute ${\bar x}^{ \infty}_2$ which is obtained by solving following fixed point equation:
  
If we have $k_{u2} -v_2 + {\bar x}_2^{ \infty}   < \bar{y}_2$  then from equation \eqref{Eqn_FP_G2} the aggregate clearing vector  reduces to :
\begin{eqnarray*}
{\bar x}_2^{ \infty}  =  \bigg (\frac{ ( \bar{l}_2 -v_2)^+}{ 1- (1-p_2^{sb})\lambda_2 }\bigg)(1-p_2^{sb})\lambda_2.
\end{eqnarray*}
Substituting ${\bar x}_2^{ \infty}$  in the above condition we have the  following bound:
\begin{eqnarray*}
 \lambda_2 < \beta_0 .
  \hspace{5mm}  \mbox{ \eop }
\end{eqnarray*}

\noindent \textbf{Proof of Lemma \ref{Lemma_G2 default with v2> kd}:}
First observe that with $v_2 > k_{d2}$, one can never have resilience regime, i.e., $P^2_{D} \ge w$, as in this case,  we always have
$(k_{d2}-v_2 + {\bar x}_2^{ \infty}) < {\bar y_2}$ (recall ${\bar x}_2^{ \infty} \le {\bar y}_2$).

We have the following   sub cases under $k_{d2}< v_ 2$. Towards the end of the proof we provide two flow diagrams and some relevant explanations to show that the following four cases exhaust all possibles parameters that satisfy the given hypothesis. \\
\textbf{Case 1:} First consider the scenario when only the banks with shock  default,  i.e., when $P^{2}_D =w$. In this case the aggregate clearing vector  is obtained as (see   \eqref{eqn_avgclearingvectorG2}):
\begin{eqnarray*}
{\bar x}_2^{ \infty}& =&  \bigg(\bar{y}_2(1-w) +w(k_{d2}-v_2 + {\bar x}_2^{ \infty})^+ \bigg) (1-p_2^{sb})\lambda_2. \end{eqnarray*}
If further $k_{d2} -v_2 + {\bar x}_2^{ \infty} \leq 0$  and  $k_{u2} -v_2 + {\bar x}_2^{ \infty} > \bar{y}_2$, then we will have:
\begin{eqnarray*}
   {\bar x}_2^{ \infty}& =&  \bar{y}_2(1-w)(1-p_2^{sb})\lambda_2, 
\end{eqnarray*}
and the above conditions 
  simplify to 
$$
\beta_4 < \lambda_2 \leq \beta_1.
$$
\textbf{Case 2:} Consider the regime when all default i.e., if  $P^{2}_D =1$ and the aggregate clearing vector is obtained as:
\begin{eqnarray*}
{\bar x}_2^{ \infty}& =& \bigg (\bigg (k_{u2}-v_2 + {\bar x}_2^{ \infty}\bigg)^+(1-w) +\bigg (k_{d2}-v_2 + {\bar x}_2^{ \infty}\bigg)^+w \bigg) (1-p^{sb}_2) \lambda_2 \end{eqnarray*}
If $k_{d2}-v_2 + {\bar x}_2^{ \infty} < 0$ and $k_{u2} -v_2 + {\bar x}_2^{ \infty} < \bar{y}_2$, then as before, we will have: 
\begin{eqnarray*}
  {\bar x}_2^{ \infty} &=& \frac{(k_{u2}-v_2)^+(1-w)(1-p^{sb}_2)\lambda_2}{1-(1-p^{sb}_2)\lambda_2(1-w)}.
\end{eqnarray*}
The above conditions are satisfied  when:
\begin{eqnarray*}
 \lambda_2 < \min  \bigg \lbrace \beta_4, \beta_3  \bigg\rbrace .
\end{eqnarray*}
 
\noindent \textbf{Case 3:} Again consider the  case with $P^{2}_D=w$.  In this case the aggregate clearing vector  is obtained as:
\begin{eqnarray*}
{\bar x}_2^{ \infty}& =& \bigg (\bar{y}_2(1-w) +\bigg (k_{d2}-v_2 + {\bar x}_2^{ \infty}\bigg)w \bigg) (1-p^{sb)}_2 \lambda_2
\end{eqnarray*}

If we have, $k_{u2}-v_2 +{\bar x}_2^{ \infty} > \bar{y}_2$ 
and $0< k_{d2}-v_2 + {\bar x}_2^{ \infty}$. In this case,
\begin{eqnarray*}
 {\bar x}_2^{ \infty}& =& 
\bigg(\frac{(k_{d2}-v_2)w +\bar{y}_2(1-w) }{1-(1-p^{sb}_2)w\lambda_2} \bigg)(1-p^{sb}_2)\lambda_2.
\end{eqnarray*}
The above two conditions are equivalent to the  following bound:
\begin{eqnarray*}
\lambda_2 > \max \bigg \lbrace \beta_2,\beta_1 \bigg \rbrace.  
\end{eqnarray*}
\textbf{Case 4:} Now consider the scenario with all default i.e., $P^2_D =1$. In this sub- case the bank can survive with the downward shock and the corresponding aggregate clearing vector in this regime is given by the following:
\begin{eqnarray*}
{\bar x}_2^{ \infty} = \bigg (\bigg (k_{u2}-v_2 + {\bar x}_2^{ \infty}\bigg)(1-w) +\bigg (k_{d2}-v_2 + {\bar x}_2^{ \infty}\bigg)w \bigg) (1-p^{sb)}_2 \lambda_2
\end{eqnarray*}
If we have, $k_{u2}-v_2+ {\bar x}_2^{ \infty} < \bar{y}_2$ and $k_{d2} -v_2 + {\bar x}_2^{ \infty} > 0$. Then the aggregate  clearing vector obtained as follows:
\begin{eqnarray*}
{\bar x}_2^{ \infty} = \frac{(\bar{l}_2-v_2)(1-p^{sb}_2)\lambda_2}{1-(1-p^{sb}_2)\lambda_2}
\end{eqnarray*}
The above conditions reduces to the following bound:
\begin{eqnarray*}
 \beta_3 < \lambda_2 < \beta_2. \hspace{5mm}   
\end{eqnarray*}
% Case 1  picture
{\bf Sufficiency of the four cases:}
  First consider   $\beta_4 > \beta_1$. By definitions, $\beta_3 > \beta_1$ and $\beta_4 > \max\{\beta_1 ,\beta_2\}$. If further,  $\min \{\beta_4,\beta_3 \} < \lambda_2 < \max \{\beta_2, \beta_1 \} $, then  
\begin{eqnarray}
\label{Eqn_beta_regime}
 \beta_1< \min \{\beta_4,\beta_3 \} < \lambda_2 < \max \{\beta_2, \beta_1 \} < \beta_4  \implies \beta_3 < \lambda_2 < \beta_2.
\end{eqnarray}
All the sub-cases that fall into this category are covered by the four sub-cases as explained in the flow chart of Figure \ref{Fig_regime1}.

\begin{figure}
\begin{tikzpicture}[x=0.6pt,y=0.6pt,yscale=-0.9,xscale=0.9]

%uncomment if require: \path (0,845); %set diagram left start at 0, and has height of 845

%Straight Lines [id:da7029366993018923]
\draw    (124,111) -- (228,111) ;
%Straight Lines [id:da5403419566964796]
\draw    (381,110) -- (474,112) ;
%Straight Lines [id:da7809696799930331]
\draw    (124,111) -- (124.98,202) ;
\draw [shift={(125,204)}, rotate = 269.38] [color={rgb, 255:red, 0; green, 0; blue, 0 }  ][line width=0.75]    (10.93,-3.29) .. controls (6.95,-1.4) and (3.31,-0.3) .. (0,0) .. controls (3.31,0.3) and (6.95,1.4) .. (10.93,3.29)   ;
%Shape: Rectangle [id:dp08179322156210245]
\draw   (89,207) -- (165,207) -- (165,247) -- (89,247) -- cycle ;
%Straight Lines [id:da41207941446750285]
\draw    (474,112) -- (473.51,194) ;
\draw [shift={(473.5,196)}, rotate = 270.34000000000003] [color={rgb, 255:red, 0; green, 0; blue, 0 }  ][line width=0.75]    (10.93,-3.29) .. controls (6.95,-1.4) and (3.31,-0.3) .. (0,0) .. controls (3.31,0.3) and (6.95,1.4) .. (10.93,3.29)   ;
%Straight Lines [id:da5281895169237881]
\draw    (298,226.75) -- (381,225.75) ;
%Straight Lines [id:da06810646320066527]
\draw    (566,225.75) -- (659,224.75) ;
%Straight Lines [id:da37989373510214297]
\draw    (298,226.75) -- (298.97,300.75) ;
\draw [shift={(299,302.75)}, rotate = 269.25] [color={rgb, 255:red, 0; green, 0; blue, 0 }  ][line width=0.75]    (10.93,-3.29) .. controls (6.95,-1.4) and (3.31,-0.3) .. (0,0) .. controls (3.31,0.3) and (6.95,1.4) .. (10.93,3.29)   ;
%Straight Lines [id:da3502101909651032]
\draw    (659,223) -- (659.97,296) ;
\draw [shift={(660,298)}, rotate = 269.24] [color={rgb, 255:red, 0; green, 0; blue, 0 }  ][line width=0.75]    (10.93,-3.29) .. controls (6.95,-1.4) and (3.31,-0.3) .. (0,0) .. controls (3.31,0.3) and (6.95,1.4) .. (10.93,3.29)   ;
%Shape: Rectangle [id:dp03337768559815324]
\draw   (620,418) -- (690,418) -- (690,453) -- (620,453) -- cycle ;
%Shape: Rectangle [id:dp8695873329134886]
\draw   (265,302) -- (335,302) -- (335,342) -- (265,342) -- cycle ;
%Shape: Diamond [id:dp9605295269845976]
\draw   (304.5,62) -- (381,110) -- (304.5,158) -- (228,110) -- cycle ;
%Shape: Diamond [id:dp633293378632396]
\draw   (473.5,196) -- (566,225.75) -- (473.5,255.5) -- (381,225.75) -- cycle ;
%Shape: Rectangle [id:dp9382364755331015]
\draw   (594,300) -- (700,300) -- (700,360) -- (594,360) -- cycle ;
%Straight Lines [id:da6124506205590861]
\draw    (659,361) -- (660.95,415) ;
\draw [shift={(661,417)}, rotate = 268.49] [color={rgb, 255:red, 0; green, 0; blue, 0 }  ][line width=0.75]    (10.93,-3.29) .. controls (6.95,-1.4) and (3.31,-0.3) .. (0,0) .. controls (3.31,0.3) and (6.95,1.4) .. (10.93,3.29)   ;

% Text Node
\draw (235,101) node [anchor=north west][inner sep=0.75pt]   [align=left] {$\lambda_2 > \min \{\beta_4,\beta_1\}$};
% Text Node
\draw (101,218) node [anchor=north west][inner sep=0.75pt]   [align=left] {Case 2};
% Text Node
\draw (628,428) node [anchor=north west][inner sep=0.75pt]   [align=left] {Case 4};
% Text Node
\draw (275,311) node [anchor=north west][inner sep=0.75pt]   [align=left] {Case 3};
% Text Node
\draw (418,115) node [anchor=north west][inner sep=0.75pt]   [align=left] {Yes};
% Text Node
\draw (167,113) node [anchor=north west][inner sep=0.75pt]   [align=left] {No};
% Text Node
\draw (325,229) node [anchor=north west][inner sep=0.75pt]   [align=left] {No};
% Text Node
\draw (605,227) node [anchor=north west][inner sep=0.75pt]   [align=left] {Yes};
% Text Node
\draw (401,217) node [anchor=north west][inner sep=0.75pt]   [align=left] {$\lambda_2 < \max\{\beta_2,\beta_1\}$};
% Text Node
\draw   (601,311)  node [anchor=north west][inner sep=0.9pt]   [align=left] {By  \eqref{Eqn_beta_regime} \\ $\beta_3<\lambda_2<\beta_2$};
\end{tikzpicture}
\vspace{-4mm}
\caption{With $\beta_4 > \beta_1$ \label{Fig_regime1}}
\end{figure}

Now we consider the left-over regime, i.e., with   $\beta_4 < \beta_1$. This implies $\beta_1 > \max \{\beta_4 , \beta_3 \}$  and $\beta_4 < \min \{\beta_1 ,\beta_2 \}$. The details of this regime are in the flow-chart of Figure \ref{Fig_regime2}.  \eop

\begin{figure}
\begin{tikzpicture}[x=0.55pt,y=0.55pt,yscale=-0.96,xscale=0.96]
%uncomment if require: \path (0,845); %set diagram left start at 0, and has height of 845

%Straight Lines [id:da7029366993018923]
\draw    (148,113) -- (252,113) ;
%Straight Lines [id:da7809696799930331]
\draw    (148,113) -- (148.97,182) ;
\draw [shift={(149,184)}, rotate = 269.19] [color={rgb, 255:red, 0; green, 0; blue, 0 }  ][line width=0.75]    (10.93,-3.29) .. controls (6.95,-1.4) and (3.31,-0.3) .. (0,0) .. controls (3.31,0.3) and (6.95,1.4) .. (10.93,3.29)   ;
%Shape: Rectangle [id:dp08179322156210245]
\draw   (111,185) -- (187,185) -- (187,225) -- (111,225) -- cycle ;
%Straight Lines [id:da41207941446750285]
\draw    (328.5,160) -- (329.48,256) ;
\draw [shift={(329.5,258)}, rotate = 269.42] [color={rgb, 255:red, 0; green, 0; blue, 0 }  ][line width=0.75]    (10.93,-3.29) .. controls (6.95,-1.4) and (3.31,-0.3) .. (0,0) .. controls (3.31,0.3) and (6.95,1.4) .. (10.93,3.29)   ;
%Straight Lines [id:da5281895169237881]
\draw    (175,288.75) -- (250,287.75) ;
%Straight Lines [id:da06810646320066527]
\draw    (409,287.75) -- (504,288) ;
%Straight Lines [id:da37989373510214297]
\draw    (175,288.75) -- (175.97,362.75) ;
\draw [shift={(176,364.75)}, rotate = 269.25] [color={rgb, 255:red, 0; green, 0; blue, 0 }  ][line width=0.75]    (10.93,-3.29) .. controls (6.95,-1.4) and (3.31,-0.3) .. (0,0) .. controls (3.31,0.3) and (6.95,1.4) .. (10.93,3.29)   ;
%Straight Lines [id:da3502101909651032]
\draw    (504,288) -- (504.97,361) ;
\draw [shift={(505,363)}, rotate = 269.24] [color={rgb, 255:red, 0; green, 0; blue, 0 }  ][line width=0.75]    (10.93,-3.29) .. controls (6.95,-1.4) and (3.31,-0.3) .. (0,0) .. controls (3.31,0.3) and (6.95,1.4) .. (10.93,3.29)   ;
%Shape: Rectangle [id:dp03337768559815324]
\draw   (563,468) -- (633,468) -- (633,502) -- (563,502) -- cycle ;
%Shape: Rectangle [id:dp8695873329134886]
\draw   (247,463) -- (317,463) -- (317,503) -- (247,503) -- cycle ;
%Shape: Diamond [id:dp9605295269845976]
\draw   (328.5,64) -- (405,112) -- (328.5,160) -- (252,112) -- cycle ;
%Shape: Diamond [id:dp633293378632396]
\draw   (329.5,258) -- (409,287.75) -- (329.5,317.5) -- (250,287.75) -- cycle ;
%Shape: Diamond [id:dp9350442376728605]
\draw   (176,364.75) -- (254,397.38) -- (176,430) -- (98,397.38) -- cycle ;
%Straight Lines [id:da3981428829103487]
\draw    (282,398) -- (283.94,460) ;
\draw [shift={(284,462)}, rotate = 268.21] [color={rgb, 255:red, 0; green, 0; blue, 0 }  ][line width=0.75]    (10.93,-3.29) .. controls (6.95,-1.4) and (3.31,-0.3) .. (0,0) .. controls (3.31,0.3) and (6.95,1.4) .. (10.93,3.29)   ;
%Straight Lines [id:da22878076985807683]
\draw    (61,398) -- (61,465) ;
\draw [shift={(61,467)}, rotate = 270] [color={rgb, 255:red, 0; green, 0; blue, 0 }  ][line width=0.75]    (10.93,-3.29) .. controls (6.95,-1.4) and (3.31,-0.3) .. (0,0) .. controls (3.31,0.3) and (6.95,1.4) .. (10.93,3.29)   ;
%Shape: Rectangle [id:dp16630001006467576]
\draw   (26,466) -- (96,466) -- (96,506) -- (26,506) -- cycle ;
%Straight Lines [id:da8233894495515137]
\draw    (61,398) -- (98,397.38) ;
%Straight Lines [id:da07694955545630622]
\draw    (254,397.38) -- (282,398) ;
%Shape: Diamond [id:dp6035085314408115]
\draw   (506,363) -- (564,398) -- (506,433) -- (448,398) -- cycle ;
%Straight Lines [id:da482099224039773]
\draw    (606,400) -- (606.97,458) ;
\draw [shift={(607,460)}, rotate = 269.05] [color={rgb, 255:red, 0; green, 0; blue, 0 }  ][line width=0.75]    (10.93,-3.29) .. controls (6.95,-1.4) and (3.31,-0.3) .. (0,0) .. controls (3.31,0.3) and (6.95,1.4) .. (10.93,3.29)   ;
%Straight Lines [id:da7785417536635518]
\draw    (564,398) -- (606,400) ;
%Straight Lines [id:da6769598050925204]
\draw    (395,399) -- (448,398) ;
%Straight Lines [id:da022533831264768822]
\draw    (395,399) -- (395.97,462) ;
\draw [shift={(396,464)}, rotate = 269.12] [color={rgb, 255:red, 0; green, 0; blue, 0 }  ][line width=0.75]    (10.93,-3.29) .. controls (6.95,-1.4) and (3.31,-0.3) .. (0,0) .. controls (3.31,0.3) and (6.95,1.4) .. (10.93,3.29)   ;
%Shape: Rectangle [id:dp2747241031101324]
\draw   (361,466) -- (431,466) -- (431,506) -- (361,506) -- cycle ;
%Shape: Rectangle [id:dp2654632965768702]
\draw   (122,250) -- (258,250) -- (258,280) -- (122,280) -- cycle ;
%Shape: Rectangle [id:dp8083995714348291]
\draw   (425,255) -- (575,255) -- (575,281) -- (425,281) -- cycle ;

% Text Node
\draw (272,102) node [anchor=north west][inner sep=0.75pt]   [align=left] {$\beta_4<\lambda_2 < \beta_1$};
% Text Node
\draw (123,196) node [anchor=north west][inner sep=0.75pt]   [align=left] {Case 1};
% Text Node
\draw (571,475) node [anchor=north west][inner sep=0.75pt]   [align=left] {Case 3};
% Text Node
\draw (257,472) node [anchor=north west][inner sep=0.75pt]   [align=left] {Case 4};
% Text Node
\draw (208,117) node [anchor=north west][inner sep=0.75pt]   [align=left] {Yes};
% Text Node
\draw (332,199) node [anchor=north west][inner sep=0.75pt]   [align=left] {No};
% Text Node
\draw (209.5,289.25) node [anchor=north west][inner sep=0.75pt]   [align=left] {No};
% Text Node
\draw (461,290) node [anchor=north west][inner sep=0.75pt]   [align=left] {Yes};
% Text Node
\draw (283,278) node [anchor=north west][inner sep=0.75pt]   [align=left] {$\lambda_2 >\beta_1$};
% Text Node
\draw (425,256) node [anchor=north west][inner sep=0.75pt]   [align=left] {$\lambda_2> \beta_1> \beta_4,\beta_3$};
% Text Node
\draw (132,389) node [anchor=north west][inner sep=0.75pt]   [align=left] {$\lambda_2>\beta_3$};
% Text Node
\draw (245,365) node [anchor=north west][inner sep=0.75pt]   [align=left] {Yes};
% Text Node
\draw (80,370) node [anchor=north west][inner sep=0.75pt]   [align=left] {No};
% Text Node
\draw (38,480) node [anchor=north west][inner sep=0.75pt]   [align=left] {Case 2};
% Text Node
\draw (465,388) node [anchor=north west][inner sep=0.75pt]   [align=left] {$\lambda_2>\beta_2$};
% Text Node
\draw (584,377) node [anchor=north west][inner sep=0.75pt]   [align=left] {Yes};
% Text Node
\draw (414,376) node [anchor=north west][inner sep=0.75pt]   [align=left] {No};
% Text Node
\draw (368,476) node [anchor=north west][inner sep=0.75pt]   [align=left] {Case 4};
% Text Node
\draw (121,257) node [anchor=north west][inner sep=0.75pt]   [align=left] {$\lambda_2< \beta_4< \beta_1,\beta_2$};
\end{tikzpicture}
\caption{With $\beta_4 \le \beta_1$ \label{Fig_regime2}}
\end{figure}

\ignore{

\tikzset{every picture/.style={line width=0.75pt}} %set default line width to 0.75pt
{\scriptsize
\begin{tikzpicture}[x=0.52pt,y=0.52pt,yscale=-1,xscale=1]
%uncomment if require: \path (0,699); %set diagram left start at 0, and has height of 699

%Shape: Diamond [id:dp4795468814791437]
\draw (351.5,46) -- (437,101) -- (351.5,156) -- (266,101) -- cycle ;
%Straight Lines [id:da4601379438225859]
\draw (117,100) -- (265,101) ;
%Straight Lines [id:da16356284727526438]
\draw (437,101) -- (570,103) ;
%Straight Lines [id:da6740743344171112]
\draw (117,100) -- (118.95,184) ;
\draw [shift={(119,186)}, rotate = 268.67] [color={rgb, 255:red, 0; green, 0; blue, 0 } ][line width=0.75] (10.93,-3.29) .. controls (6.95,-1.4) and (3.31,-0.3) .. (0,0) .. controls (3.31,0.3) and (6.95,1.4) .. (10.93,3.29) ;
%Shape: Rectangle [id:dp8389447253655358]
\draw (50,185) -- (193,185) -- (193,265) -- (50,265) -- cycle ;
%Straight Lines [id:da13943032925641585]
\draw (122,268) -- (37.21,379.41) ;
\draw [shift={(36,381)}, rotate = 307.27] [color={rgb, 255:red, 0; green, 0; blue, 0 } ][line width=0.75] (10.93,-3.29) .. controls (6.95,-1.4) and (3.31,-0.3) .. (0,0) .. controls (3.31,0.3) and (6.95,1.4) .. (10.93,3.29) ;
%Shape: Rectangle [id:dp23378991034449625]
\draw (4,383) -- (175,383) -- (175,479) -- (4,479) -- cycle ;
%Straight Lines [id:da9064103275127751]
\draw (122,268) -- (519.08,384.44) ;
\draw [shift={(521,385)}, rotate = 196.34] [color={rgb, 255:red, 0; green, 0; blue, 0 } ][line width=0.75] (10.93,-3.29) .. controls (6.95,-1.4) and (3.31,-0.3) .. (0,0) .. controls (3.31,0.3) and (6.95,1.4) .. (10.93,3.29) ;
%Shape: Rectangle [id:dp3148420253513802]
\draw (223,385) -- (365,385) -- (365,480) -- (223,480) -- cycle ;
%Straight Lines [id:da7666863451108787]
\draw (122,268) -- (291.35,384.86) ;
\draw [shift={(293,386)}, rotate = 214.61] [color={rgb, 255:red, 0; green, 0; blue, 0 } ][line width=0.75] (10.93,-3.29) .. controls (6.95,-1.4) and (3.31,-0.3) .. (0,0) .. controls (3.31,0.3) and (6.95,1.4) .. (10.93,3.29) ;
%Shape: Rectangle [id:dp9312569262606369]
\draw (430,385) -- (618,385) -- (618,474) -- (430,474) -- cycle ;
%Straight Lines [id:da22713429149773368]
\draw (570,103) -- (569.02,182) ;
\draw [shift={(569,184)}, rotate = 270.71] [color={rgb, 255:red, 0; green, 0; blue, 0 } ][line width=0.75] (10.93,-3.29) .. controls (6.95,-1.4) and (3.31,-0.3) .. (0,0) .. controls (3.31,0.3) and (6.95,1.4) .. (10.93,3.29) ;
%Shape: Rectangle [id:dp8953178948993405]
\draw (486,185) -- (638,185) -- (638,259) -- (486,259) -- cycle ;

% Text Node
\draw (314,91) node [anchor=north west][inner sep=0.75pt] [align=left] {$\bar{y}_2 > k_{u2}-v_2$};
% Text Node
\draw (55,205) node [anchor=north west][inner sep=0.75pt] [align=left] {when $\beta_4 > \beta_2$, \\
$\beta_3 > \beta_1$\\ $P_D^2 = w$ or $1$};
% Text Node
\draw (16,400) node [anchor=north west][inner sep=0.75pt] [align=left] {if $\beta_4 < \lambda_2 < \beta_1$\\$P_D^2 =w$\\else\\
$P_D^2 =1$};
% Text Node
\draw (242,399) node [anchor=north west][inner sep=0.75pt] [align=left] {if $\lambda_2 < \beta_2$\\ $P_D^2 =w$\\else\\
$P_D^2 =1$};
% Text Node
\draw (444,399) node [anchor=north west][inner sep=0.75pt] [align=left] {when $\beta_1 < \lambda_2 < \beta_3$\\if $\lambda_2 < \beta_2$ \\
$P_D^2 =w$\\else
$P_D^2 =1$};
% Text Node
\draw (122,143) node [anchor=north west][inner sep=0.75pt] [align=left] {Yes};
% Text Node
\draw (502,74) node [anchor=north west][inner sep=0.75pt] [align=left] {No};
% Text Node
\draw (496,202) node [anchor=north west][inner sep=0.75pt] [align=left] {$P_D^2 = w$ };

\end{tikzpicture}
 }
 }
%\textbf{Remarks:} In the above lemma, the assumption  $\bar{y}_2  > w(k_{u2}-k_{d2})$  is reasonable; This is because the probability of having downward shock $w$ is a typically small value, i.e., close to zero and the total liability is non-zero.
   
\noindent \textbf{Proof of Lemma \ref{Lemma_G1 default}:} 
Let $\beta:= \mu_1 {\bar x}_2^\infty.$
Recall that in general ${\bar x}_1^{ \infty}$ has to satisfy (see \eqref{eqn_avgclearingvectorG1}):

\vspace{-7mm}
{\small \begin{eqnarray}
\frac{{\bar x}_1^{ \infty}}{(1-p_1^{sb})} =  \min \left \{\bar{y}_1, \bigg (k_{d1}- v_1+ {\bar x}_1^{ \infty}+ \beta \bigg)^+ \right \} w \nonumber \\ + \min \left \{\bar{y}_1,\bigg(k_{u1}- v_1+{\bar x}_1^{ \infty}+ \beta \bigg)^+\right \}(1-w).  \label{Eqn_FP_G1}
\end{eqnarray}} 
We prove this lemma in the following sub-cases: \\
\noindent \textbf{Case 1:} First consider the case when downward shock can be absorbed i.e., default probability, $ P^1_D = 0$. In this case,   $k_{d1}-v_1 + {\bar x}_1^{ \infty} + \beta \ge \bar{y}_1$, and  then the aggregate clearing vector ${\bar x}_1^{ \infty} =\bar{y}_1 (1-p_1^{sb})$. Substituting  ${\bar x}_1^{ \infty} $ in the above we have the following bound:
\vspace{-4mm}
\begin{eqnarray*}
\beta  \ge  (v_1-k_{d1} + \bar{y}_1 p_1^{sb}).
\end{eqnarray*}
\noindent \textbf{Case 2:} Consider the case, in which,  only the banks that receive shock will default, i.e., $P^1_D=w$. In this case,
 $k_{d1}- v_1 + {\bar x}_1^{ \infty}+ \beta <\bar{y}_1$ and  $k_{u1}- v_1 +  {\bar x}_1^{ \infty}+\beta \ge \bar{y}_1.$ Then the aggregate clearing vector reduces to: 
\begin{equation}
{\bar x}_1^{ \infty} = \Bigg(\frac{\bar{y}_1(1-w) +w\bigg (k_{d1}-v_1 + \beta  \bigg) }{1-w(1-p_1^{sb})} \Bigg) (1-p_1^{sb}). \nonumber
\end{equation}
Under the given hypothesis, we always have, 
$k_{d1}- v_1 + {\bar x}_1^{ \infty}+ \beta \ge 0$.  
\noindent Substituting  ${\bar x}_1^{ \infty}$:
\begin{equation*}
v_1-\bar{l}_1 +p_1^{sb} (\bar{y}_1+ w(k_{d1} -k_{u1})) \le \beta
	< (v_1-k_{d1} + \bar{y}_1 p_1^{sb}).
\end{equation*}
\ignore{
 we need that $k_{d1} -v_1 +\bar{x}_1^\infty +\beta \ge 0 $ which simplifies to the following after substituting $\bar{x}_1^\infty$:
 $$\beta + k_{d1}- v_1 +\bar{y}_1(1-w)(1-p_1^{sb}) \geq 0$$
 Or we need
 $$
 \beta
	- ( (v_1-k_{d1} + \bar{y}_1 p_1^{sb}) ) + \bar{y}_1 (1-w)+ w \bar{y}_1 p_1^{sb}  > 0
$$ 

Or using the lower bound of $\beta$ in the above we  have the following:
$$
(1-w)(k_{d1}-k_{u1}) + \bar{y}_1 (1-w)+ w \bar{y}_1 p_1^{sb}  + p_1^{sb}  w(k_{d1} -k_{u1}) > 0
$$

$$
(1-w + w p_1^{sb})( \bar{y}_1 + k_{d1}-k_{u1})   > 0
$$
 Thus it is reduces to the condition:  $$\bar{y}_1 > (k_{u1} -k_{d1})$$
 }
 \ignore{
{\color{red}\noindent\textbf{Case 3:}
We are now left with the case
when $\beta < v_1 -  \bar{l}_1 + p_1^{sb} (\bar{y}_1+ w(k_{d1} -k_{u1}))$. Then, $k_{d1}- v_1+   \beta < 0$. If further 
$ k_{d1}- v_1+{\bar x}_1^{ \infty}+ \beta \le 0$, then clearly\footnote{If further $k_{u1}- v_1+ \beta \le 0$, no positive ${\bar x}_1^{ \infty}$ can satisfy the required equation} from \eqref{Eqn_FP_G1}:
$${\bar x}_1^{ \infty}  =  \bigg(k_{u1}- v_1+ \beta \bigg)^+ \frac{(1-w)(1-p_1^{sb)}}{1-(1-w)(1-p_1^{sb})}  .$$
This obviously translates to the following condition:
$$
 k_{d1}- v_1+\bigg(k_{u1}- v_1+ \beta \bigg)^+ \frac{(1-w)(1-p_1^{sb})}{1-(1-w)(1-p_1^{sb})} + \beta \le 0
$$which is trivially true with $\beta \le v_1 -  \bar{l}_1 +p_1^{sb}(1-w)(k_{u1} -k_{d1})$. 

Also we need that $k_{u1} -v_1 +\bar{x}_1^{\infty} +\beta \ge 0$ and this reduces to the  following:
$$
\beta \le \bar{y}_1(1-(1-w)(1-p_1^{sb}))+ v_1- k_{u1}.
$$

On the other hand, if $k_{d1}- v_1+{\bar x}_1^{ \infty}+ \beta > 0$, then the fixed point equation \eqref{Eqn_FP_G1} is satisfied only when $v_1 - {\bar l}_1  = \beta$, but in this sub-case we have $\beta < v_1 -  \bar{l}_1$.\\}
}
\textbf{Case 3:} We are now left with the case
when
%$\beta < v_1 -  \bar{l}_1 + p_1^{sb} (\bar{y}_1+ w(k_{d1} -k_{u1}))$ and 
the default probability is, $P_D^1=1$. In this case   we will have  $k_{u1}  -v_1 +\bar{x}_1^{\infty} +\beta < \bar{y}_1$. Then the aggregate clearing vector reduces to:
\begin{eqnarray*}
\bar{x}_1^{\infty} = \frac{(\bar{l}_1 -v_1 +\beta)(1-p_1^{sb})}{p_1^{sb}} \ \mbox{where},\  p_1^{sb} > 0.
\end{eqnarray*}
Substituting $\bar{x}_1^{\infty}$  in the required condition,  $k_{u1}  -v_1 +\bar{x}_1^{\infty} +\beta < \bar{y}_1$,  we have  $\beta< v_1 - \bar{l}_1 + p_1^{sb} (\bar{y}_1+ w(k_{d1} -k_{u1}))$. 
\eop

\ignore{
 In this we first calculate ${\bar x}_1^{ \infty}$ which is obtained by solving following fixed point equation:
  $${\bar x}_1^{ \infty} =  \bigg (k_{d1}- v_1+ {\bar x}_1^{ \infty}+ \beta \bigg)^+w + \bigg(k_{u1}- v_1+{\bar x}_1^{ \infty}+ \beta \bigg)^+(1-w) $$ 
\noindent $ \implies {\bar x}_1^{ \infty}  =  \bigg(k_{u1}- v_1+ \beta \bigg)^+ \frac{1-w}{w}  .$\\
This is true when the following two conditions are satisfied:
$$k_{u1}- v_1+ \beta > 0 \mbox{ and }
 \bar{l}_1 - v_1 +\beta    \le 0
$$
   In this case the default probability is $P^1_D = 1$.
   The regime satisfies if the following hold
    \begin{eqnarray*}
(v_1-k_{u1}) <  \beta_1 \le   (v_1-\bar{l}_1). \hspace{5mm}  \mbox{ \eop }
    \end{eqnarray*} }
   
 %{\color{blue}
%\noindent\underline{Case 3:}
 %In this we first calculate ${\bar x}_1^{ \infty}$ which is obtained by solving following fixed point equation:
  %$${\bar x}_1^{ \infty} =  \bigg  (k_{d1}- v_1+ {\bar x}_1^{ \infty}+ \beta \bigg )^+ w +
  %\bigg (k_{u1}- v_1+{\bar x}_1^{ \infty}+ \beta \bigg)^+ (1-w) $$ 
%\noindent $ \implies {\bar x}_1^{ \infty} =  0.$

  % In this case the default probability is $P^1_D = 1$.
 %  The regime satisfies if the following hold
   %\vspace{-2mm}
   % \begin{eqnarray*}
  % k_{d1} -v_1 + {\bar x}_1^\infty + \beta < 0 \\ 
   %\implies \beta <  (v_1-k_{d1}). %\hspace{5mm}  \mbox{ \eop }
   % \end{eqnarray*}}
%\noindent \textbf{Note:} In the last regime i.e.,  when $P_D^{1}=1$, a even  more tighter bound can be   obtained by using  ${\bar x}_1^\infty < \frac{(v_1- \bar{l}_1)\gamma}{1-\gamma}$. Then  the last regime  can be modified by  follows:
%\begin{equation*}
   %(1-p_2^{sb})\lambda_2{\bar x}_1^\infty < \frac{\gamma}{(1-\gamma)}\bigg( \bar{y}_1 - (k_{u1}-v_1)-\frac{\gamma}{(1-\gamma)}(v_1-\bar{l}_1)\bigg).
%end{equation*}
%{\color{red}For notational  convenience we define another constant, $\kappa := \kappa' / (1+u-d_c)$, i.e., 
%$v_m = \kappa (1+u-d_c)  \Omega_m$ for all $m$.\\}

\noindent \textbf{Proof of Lemma \ref{Lemma_default}:} 
The  resilient condition  for group $\mathcal{G}_2$  is that there is no default even when its banks  receive shock and hence  with ${\bar x}_2^\infty = {\bar y}_2 \lambda_2 (1-p_2^{sb})$   the following should be satisfied (see \eqref{eqn_avgclearingvectorG2}):
\begin{equation}
\label{eqn_resilient}
    k_{d2}-v_2 + \bar{y}_2 \lambda_2(1-p^{sb}_2) \geq  \bar{y}_2 . \\
%\mbox{, or equivalently, } (1-p^{sb}_2) \geq \frac{\bar{y}_2 +v_2-k_{d2}}{y_2(1+r_2)}.
\end{equation}
%In the above 
Note  
from Lemma \ref{Lemma_G2 default with v2> kd} that one can't have $P_D^2 = 0$ with $v_2 > k_{d_2}$, thus for resilience of   $\mathcal{G}_2$ it is required that $v_2 < k_{d_2}$.
Also from \eqref{eqn_omega_2},  
$$v_2- k_{d2}=(k_0 +y_c+ y_2 p^{sb}_2 ) \bigg( \kappa  - (1+d-d_c) \bigg).$$
Thus   condition (\ref{eqn_resilient})   equivalently modifies to the following requirement (recall $\bar{y}_2:= (y_2+y_c)(1+r_2)$, $\lambda_2 = y_2 / (y_2 + y_c)$):
\begin{eqnarray*}
    y_c\left ((1+r_2) +\kappa  - (1+d-d_c) \right )   
    & \leq &-    {\bigg( \kappa  - (1+d-d_c)\bigg)(k_0 +y_2 p^{sb}_2) - y_2 p^{sb}_2(1+r_2)  }  ,
    \end{eqnarray*}
    Note that $r_2 > d$, hence the left hand side (LHS) is non-negative, and thus
observe  that the right hand side (RHS) should be a non-negative quantity, and this is possible only when the first term in the RHS is positive, and,
\begin{eqnarray}
  (1+d-d_c) - \kappa   \ge \frac{(y_c+ y_2 p^{sb}_2)(1+r_2)}{k_0 +y_c+y_2 p^{sb}_2} > 0, \label{Eqn_cond_resilience}
\end{eqnarray}
as the second term is definitely negative and $y_c$ has to be a non-negative quantity.  
Under the above condition, consider   a   $\mathcal{G}_1$ bank  that 
receives shock. Then its clearing value is governed by the following:
 \begin{equation}
     \bigg(k_{d1}-v_1 + {\bar x}_1^\infty  + {\bar y}_2 \lambda_2 (1-p_2^{sb}) \frac{1- \gamma}{\gamma} \frac{y_c}{y_2} \frac{1}{(1-p^{sb}_2)}\bigg).
     \label{eqn_resilient cond for g1}
 \end{equation}
 Recall $v_1- k_{d1}= \bigg( \kappa  - (1+d-d_c) \bigg)\bigg(k_0 +y_1p_1^{sb} -\frac{1-\gamma}{\gamma}y_c\bigg)$, and hence the $\mathcal{G}_1$ bank (with shock) also does   not default because:

\vspace*{-2mm}
 {\small
\begin{eqnarray*}
     k_{d1}-v_1 + \bar{y}_1(1-p_1^{sb})  +  (1+r_2) y_c \frac{1- \gamma}{\gamma} \hspace{-55mm} \\
     &=& \bigg( (1+d-d_c) -\kappa \bigg) \left (k_0 +y_1p_1^{sb} -\frac{1-\gamma}{\gamma}y_c \right  ) + \bar{y}_1(1-p_1^{sb})  +  (1+r_2) y_c \frac{1- \gamma}{\gamma}  \\
     &=& y_c \frac{1- \gamma}{\gamma} \left ((1+r_2) + \bigg( \kappa  - (1+d-d_c) \bigg) \right ) + {\bar y}_1(1-p_1^{sb}) + \bigg( 1+d-d_c -\kappa \bigg) (k_0+y_1p_1^{sb})  \\
          &=& y_c \frac{1- \gamma}{\gamma} \left ( r_2 +   \kappa  -  d+ d_c   \right ) + {\bar y}_1(1-p_1^{sb}) + \bigg(  (1+d-d_c)  -\kappa \bigg) (k_0+y_1p_1^{sb})   .
\end{eqnarray*}}
From  equation \eqref{Eqn_cond_resilience} we have

\vspace{-2mm}
{\small
\begin{eqnarray*}
\bigg(  1+d-d_c  -\kappa \bigg) k_0  >  (y_c + y_2 p_2^{sb} ) \bigg (r_2-d+d_c +\kappa ) \bigg)
> y_1p_1^{sb} \bigg (r_2-d+d_c +\kappa ) \bigg)
\end{eqnarray*}}
and we have following

\vspace{-2mm}
{\small
 \begin{eqnarray*}
  \bigg( 1+d-d_c  -\kappa \bigg) (k_0+y_1p_1^{sb})
      > \bigg (r_2-d+d_c +\kappa  \bigg) y_1p_1^{sb}+ (1+d-d_c -\kappa)y_1p_1^{sb} > \bar{y}_1 p_1^{sb}.
 \end{eqnarray*} } 
Using this lower bound of the above, the aggregate clearing vector for the $\mathcal{G}_1$ banks simplifies to  the following:
\begin{eqnarray*}
k_{d1}-v_1 + \bar{y}_1(1-p_1^{sb})  +  (1+r_2) y_c \frac{1- \gamma}{\gamma} \hspace{-55mm} \\
&>& y_c \frac{1- \gamma}{\gamma} (r_2 -d+d_c+\kappa) + {\bar y}_1(1-p_1^{sb}) + \bar{y}_1p_1^{sb}\\
&>&  \bar{y}_1.
\end{eqnarray*}
Basically if the taxes  ($\kappa$) are not   high the condition \eqref{Eqn_cond_resilience} is satisfied and one can have resilience. This completes  the proof.  \eop

\noindent \textbf{Proof of Lemma \ref{lemma_surplus1}:} 
With resilient regime   we  will have $v_2 < k_{d2}$, as in previous proof. And for this case, using \eqref{eqn_limitsurplussb}-\eqref{Eqn_surplus_atu_G_2}, the expected surplus as well as SaU simplify to the following:

\vspace{-2mm}
{\small
\begin{eqnarray*} 
%\label{Eqn_surplus 1}
E[S_1]&=& E[K_i^1] -v_1+y_c(1+r_2)\frac{1-\gamma}{\gamma} -\bar{y}_1 p_1^{sb} \hspace{-75mm}\\
&=& \bigg(k_0 +y_1p_1^{sb}- \frac{1-\gamma}{\gamma}y_c \bigg) \bigg(1-d_c +{\bar r}_r- \kappa \bigg) + y_c(1+r_2)\frac{1-\gamma}{\gamma} -\bar{y}_1 p_1^{sb} \\
&=& k_0( 1+ \bar{r}_r  -(d_c + \kappa) ) - \frac{1-\gamma}{\gamma} y_c(\Delta_r - d_c -\kappa)  +y_1p_1^{sb}(\bar{r}_r- d_c-\kappa-r_1),
\\
\label{Eqn_surplus 2}
E[S_2]&=&  E[K_i^2]- v_2 +\bar{x}_2^{\infty} -\bar{y}_2 \\
&=& (k_0+ y_c + y_2p^{sb}_2 ) (1-d_c +\bar{r}_r-\kappa)  -(1+r_2)(y_c +y_2 p^{sb}_2 ) \\
&=& k_0 (1+\bar{r}_r-(d_c +\kappa))  + \bigg (y_2 p_2^{sb} + y_c \bigg  )  (\Delta_r - d_c -\kappa)  ,
\\
E[S_2]-E[S_1] &=& \left (\Delta_r - (d_c +\kappa) \right  ) \bigg (y_2 p_2^{sb} + \frac{y_c}{\gamma} \bigg  )    - y_1p_1^{sb}(\bar{r}_r- d_c-\kappa-r_1), \\
\hat{S}_{2,u}& =& k_{u2}+ \bar{x}_2^{\infty} -v_2- \bar{y}_2 \\
%&=& \Omega_2\bigg ((1+u-d_c) -\kappa\bigg) - (1+r_2) \bigg(y_c+ y_2 p^{sb}_2 \bigg) \\
&=& \bigg(k_0 + y_2 p_2^{sb} + y_c \bigg ) \bigg ((1+u-d_c) -\kappa\bigg) - (1+r_2) \bigg(y_c+ y_2 p^{sb}_2 \bigg) \\
&=& k_0 (1+u-d_c -\kappa )
+ \bigg (y_2 p_2^{sb} + y_c \bigg  ) \left (\Delta_u - (d_c +\kappa) \right  ).
\end{eqnarray*}}
All the above quantities are linear in $y_c$ and hence the lemma. \eop

%%%%%%%%%%%%%%%%%%%%%%%%%%

%Game Over

\ignore{
\begin{table}[hbt!]
\centering
\begin{tabular}{l|l|l|l|l|}
\cline{2-5}
\multicolumn{1}{c|}{} & \multicolumn{2}{c|}{\textbf{Regular Graph}} & \multicolumn{2}{c|}{\textbf{ER Graph}} \\ \hline
\multicolumn{1}{|c|}{$n$} & \multicolumn{1}{c|}{${\bar P}^2_D$} & \multicolumn{1}{c|}{Confidence Interval} & \multicolumn{1}{c|}{${\bar P}^2_D$} & \multicolumn{1}{c|}{Confidence Interval}  \\ \hline
\multicolumn{1}{|l|}{500}  &  0.1893& (0.1892, 0.1906) & 0.1485  & (0.1451, 0.1518) \\ \hline
\multicolumn{1}{|l|}{1000} & 0.1981 & (0.1981, 0.1989) & 0.1619 & (0.1595, 0.1644)  \\ \hline
\multicolumn{1}{|l|}{1500} & 0.1995  &(0.1995, 0.2001)  & 0.1700 & (0.1680, 0.1720)   \\ \hline
\multicolumn{1}{|l|}{2000} & 0.1998  & (0.1998, 0.2004) &  0.1772&  (0.1756, 0.1787)  \\ \hline
\multicolumn{1}{|l|}{2500} & 0.1998 & (0.1998, 0.2003) & 0.1812 &(0.1799, 0.1825)   \\ \hline
\multicolumn{1}{|l|}{5000} & 0.2008 & (0.2004, 0.2011) & 0.1936 &(0.1931, 0.1941)   \\ \hline
\end{tabular}
\caption{Average over $1000$ sample paths :  $u=0.2$, $d=-0.6$, $r_2= 0.12$,
   $\kappa= 0.56$, $\Omega_2= 12.5$, $\bar{y}_2 =35$, $v_2=7$, $w= 0.2$, $p=0.05$, $p_2^{sb}=0.001$, $y_c=0$, $d_c=0$.}
\label{Table_with PD}   
\end{table}
\begin{table}[hbt!]
\centering
\begin{tabular}{l|l|l|l|l|}
\cline{2-5}
\multicolumn{1}{c|}{} & \multicolumn{2}{c|}{\textbf{Regular Graph}} & \multicolumn{2}{c|}{\textbf{ER Graph}} \\ \hline
\multicolumn{1}{|c|}{$n$} & \multicolumn{1}{c|}{$\bar{E}[S]_{sim}$} & \multicolumn{1}{c|}{Confidence Interval} & \multicolumn{1}{c|}{$\bar{E}[S]_{sim}$} & \multicolumn{1}{c|}{Confidence Interval}  \\ \hline
\multicolumn{1}{|l|}{500}  & 5.9724 & (5.9611, 5.9836) & 5.9632  & (5.9507, 5.9657) \\ \hline
\multicolumn{1}{|l|}{1000} & 5.9581 & (5.9504, 5.9658) & 5.9599 & (5.9516, 5.9610)  \\ \hline
\multicolumn{1}{|l|}{1500} &  5.9611& (5.9547, 5.9675)  &5.9673  & (5.9609, 5.9679)   \\ \hline
\multicolumn{1}{|l|}{2000} & 5.9633  & (5.9577, 5.9689) &5.9595  &  (5.9536, 5.9601)  \\ \hline
\multicolumn{1}{|l|}{2500} & 5.9670 & (5.9619, 5.9720) & 5.9678 &(5.9626, 5.9682)   \\ \hline
\multicolumn{1}{|l|}{5000} & 5.9654 & (5.9618, 5.9655) & 5.9653 &(5.9626, 5.9682)   \\ \hline
\end{tabular}
\caption{Average over $1000$ sample paths:  $u=0.2$, $d=-0.6$, $r_2= 0.12$,
   $\kappa= 0.56$, $\Omega_2= 12.5$, $\bar{y}_2 =35$, $v_2=7$, $w= 0.2$, $p=0.05$, $p_2^{sb}=0.001$, $y_c=0$, $d_c=0$.}
   \label{Table_withsurplus}
\end{table}

}

\ignore{
  
  %% Table for 
 
   \begin{table}[hbt!]
\centering
\begin{center}
\begin{tabular}{|l|l|l|l|l|l|l|l|}
\hline
 $n$ & ${\bar P}^2_D$ &$var({\bar P}^2_D)$  & CI & 
 $ \bar{E}[S]_{sim}$& $ var(E[S]_{sim})$  &  CI  \\  \hline

  500  &  0.148476&  0.002951  &   &5.963185   &  0.040578  &  \\ \hline
   1000&  0.161940  &  0.001603  &  & 5.959884    & 0.017797 &  \\ \hline
   1500&    0.169995& 0.001002  &  & 5.967263   &0.010641   &  \\ \hline
   2000&  0.177182 &   0.000611&  & 5.959510  &   0.009048& \\ \hline
   2500& 0.181178  & 0.000424 &  &  5.967784  & 0.006950  &  \\ \hline
\end{tabular}
\end{center}
\caption{Average over $1000$ sample paths with Erdos Renyi graph:  $u=0.2$, $d=-0.6$, $r_2= 0.12$,
   $\kappa= 0.56$, $\Omega_2= 12.5$, $\bar{y}_2 =35$, $v_2=7$, $w= 0.2$, $p=0.05$, $p_2^{sb}=0.001$, $y_c=0$, $d_c=0$.}
\label{Table_3}
\end{table}

 \begin{table}[hbt!]
\centering
\begin{center}
\begin{tabular}{|l|l|l|l|l|l|l|l|}
\hline
 $n$ & ${\bar P}^2_D$ &$var({\bar P}^2_D)$  & CI & 
 $ \bar{E}[S]_{sim}$& $ var(E[S]_{sim})$  &  CI  \\  \hline

  500  & 0.189274 &  0.000449 &   &  5.972367 & 0.033049  &  \\ \hline
   1000& 0.198122  &  0.000167  &  &   5.958070 & 0.015455  &  \\ \hline
   1500& 0.199471 & 0.000110 &  &   5.961123  &   0.010724 &  \\ \hline
   2000& 0.199781  &  0.000085  &  & 5.963309 & 0.008262 & \\ \hline
   2500&  0.199796 &0.000066  &  & 5.966977  &  0.006604 &  \\ \hline
\end{tabular}
\end{center}
\caption{Average over $1000$ sample paths with regular graph:  $u=0.2$, $d=-0.6$, $r_2= 0.12$,
   $\kappa= 0.56$, $\Omega_2= 12.5$, $\bar{y}_2 =35$, $v_2=7$, $w= 0.2$, $p=0.05$, $p_2^{sb}=0.001$, $y_c=0$, $d_c=0$.}
\label{Table_4}
\end{table}
}
%%%%%%%%%%%%%%%%%%%

% Alternate case 2 picture without boxes 
\ignore{

\begin{tikzpicture}[x=0.54pt,y=0.54pt,yscale=-0.95,xscale=0.95]
%uncomment if require: \path (0,845); %set diagram left start at 0, and has height of 845

%Straight Lines [id:da7029366993018923]
\draw    (167,113) -- (271,113) ;
%Straight Lines [id:da7809696799930331]
\draw    (167,113) -- (167.97,182) ;
\draw [shift={(168,184)}, rotate = 269.19] [color={rgb, 255:red, 0; green, 0; blue, 0 }  ][line width=0.75]    (10.93,-3.29) .. controls (6.95,-1.4) and (3.31,-0.3) .. (0,0) .. controls (3.31,0.3) and (6.95,1.4) .. (10.93,3.29)   ;
%Shape: Rectangle [id:dp08179322156210245]
\draw   (130,185) -- (206,185) -- (206,225) -- (130,225) -- cycle ;
%Straight Lines [id:da41207941446750285]
\draw    (347.5,160) -- (348.48,256) ;
\draw [shift={(348.5,258)}, rotate = 269.42] [color={rgb, 255:red, 0; green, 0; blue, 0 }  ][line width=0.75]    (10.93,-3.29) .. controls (6.95,-1.4) and (3.31,-0.3) .. (0,0) .. controls (3.31,0.3) and (6.95,1.4) .. (10.93,3.29)   ;
%Straight Lines [id:da5281895169237881]
\draw    (173,288.75) -- (269,287.75) ;
%Straight Lines [id:da06810646320066527]
\draw    (428,287.75) -- (534,286.75) ;
%Straight Lines [id:da37989373510214297]
\draw    (173,288.75) -- (173.97,362.75) ;
\draw [shift={(174,364.75)}, rotate = 269.25] [color={rgb, 255:red, 0; green, 0; blue, 0 }  ][line width=0.75]    (10.93,-3.29) .. controls (6.95,-1.4) and (3.31,-0.3) .. (0,0) .. controls (3.31,0.3) and (6.95,1.4) .. (10.93,3.29)   ;
%Straight Lines [id:da3502101909651032]
\draw    (534,285) -- (534.97,358) ;
\draw [shift={(535,360)}, rotate = 269.24] [color={rgb, 255:red, 0; green, 0; blue, 0 }  ][line width=0.75]    (10.93,-3.29) .. controls (6.95,-1.4) and (3.31,-0.3) .. (0,0) .. controls (3.31,0.3) and (6.95,1.4) .. (10.93,3.29)   ;
%Shape: Rectangle [id:dp03337768559815324]
\draw   (626,461) -- (696,461) -- (696,500) -- (626,500) -- cycle ;
%Shape: Rectangle [id:dp8695873329134886]
\draw   (247,463) -- (317,463) -- (317,503) -- (247,503) -- cycle ;
%Shape: Diamond [id:dp9605295269845976]
\draw   (347.5,64) -- (424,112) -- (347.5,160) -- (271,112) -- cycle ;
%Shape: Diamond [id:dp633293378632396]
\draw   (348.5,258) -- (428,287.75) -- (348.5,317.5) -- (269,287.75) -- cycle ;
%Shape: Diamond [id:dp9350442376728605]
\draw   (176,364.75) -- (254,397.38) -- (176,430) -- (98,397.38) -- cycle ;
%Straight Lines [id:da3981428829103487]
\draw    (282,398) -- (283.94,460) ;
\draw [shift={(284,462)}, rotate = 268.21] [color={rgb, 255:red, 0; green, 0; blue, 0 }  ][line width=0.75]    (10.93,-3.29) .. controls (6.95,-1.4) and (3.31,-0.3) .. (0,0) .. controls (3.31,0.3) and (6.95,1.4) .. (10.93,3.29)   ;
%Straight Lines [id:da22878076985807683]
\draw    (61,398) -- (61,465) ;
\draw [shift={(61,467)}, rotate = 270] [color={rgb, 255:red, 0; green, 0; blue, 0 }  ][line width=0.75]    (10.93,-3.29) .. controls (6.95,-1.4) and (3.31,-0.3) .. (0,0) .. controls (3.31,0.3) and (6.95,1.4) .. (10.93,3.29)   ;
%Shape: Rectangle [id:dp16630001006467576]
\draw   (26,466) -- (96,466) -- (96,506) -- (26,506) -- cycle ;
%Straight Lines [id:da8233894495515137]
\draw    (61,398) -- (98,397.38) ;
%Straight Lines [id:da07694955545630622]
\draw    (254,397.38) -- (282,398) ;
%Shape: Diamond [id:dp6035085314408115]
\draw   (535,360) -- (593,395) -- (535,430) -- (477,395) -- cycle ;
%Straight Lines [id:da482099224039773]
\draw    (660,397) -- (660.97,455) ;
\draw [shift={(661,457)}, rotate = 269.05] [color={rgb, 255:red, 0; green, 0; blue, 0 }  ][line width=0.75]    (10.93,-3.29) .. controls (6.95,-1.4) and (3.31,-0.3) .. (0,0) .. controls (3.31,0.3) and (6.95,1.4) .. (10.93,3.29)   ;
%Straight Lines [id:da7785417536635518]
\draw    (593,395) -- (660,396) ;
%Straight Lines [id:da6769598050925204]
\draw    (424,396) -- (477,395) ;
%Straight Lines [id:da022533831264768822]
\draw    (424,396) -- (424.97,459) ;
\draw [shift={(425,461)}, rotate = 269.12] [color={rgb, 255:red, 0; green, 0; blue, 0 }  ][line width=0.75]    (10.93,-3.29) .. controls (6.95,-1.4) and (3.31,-0.3) .. (0,0) .. controls (3.31,0.3) and (6.95,1.4) .. (10.93,3.29)   ;
%Shape: Rectangle [id:dp2747241031101324]
\draw   (390,463) -- (460,463) -- (460,503) -- (390,503) -- cycle ;

% Text Node
\draw (291,102) node [anchor=north west][inner sep=0.75pt]   [align=left] {$\beta_4<\lambda_2 < \beta_1$ };
% Text Node
\draw (142,196) node [anchor=north west][inner sep=0.75pt]   [align=left] {Case 1};
% Text Node
\draw (634,471) node [anchor=north west][inner sep=0.75pt]   [align=left] {Case 3};
% Text Node
\draw (257,472) node [anchor=north west][inner sep=0.75pt]   [align=left] {Case 4};
% Text Node
\draw (227,117) node [anchor=north west][inner sep=0.75pt]   [align=left] {Yes};
% Text Node
\draw (351,199) node [anchor=north west][inner sep=0.75pt]   [align=left] {No};
% Text Node
\draw (228.5,289.25) node [anchor=north west][inner sep=0.75pt]   [align=left] {No};
% Text Node
\draw (461,290) node [anchor=north west][inner sep=0.75pt]   [align=left] {Yes};
% Text Node
\draw (302,278) node [anchor=north west][inner sep=0.75pt]   [align=left] {$\lambda_2 >\beta_1$};
% Text Node
\draw (432,264) node [anchor=north west][inner sep=0.75pt]   [align=left] {$\lambda_2>\beta_1>\beta_4,\beta_3$};
% Text Node
\draw (132,389) node [anchor=north west][inner sep=0.75pt]   [align=left] {$\lambda_2>\beta_3$};
% Text Node
\draw (245,365) node [anchor=north west][inner sep=0.75pt]   [align=left] {Yes};
% Text Node
\draw (80,370) node [anchor=north west][inner sep=0.75pt]   [align=left] {No};
% Text Node
\draw (38,475) node [anchor=north west][inner sep=0.75pt]   [align=left] {Case 2} ;
% Text Node
\draw (494,385) node [anchor=north west][inner sep=0.75pt]   [align=left] {$\lambda_2>\beta_2$};
% Text Node
\draw (613,374) node [anchor=north west][inner sep=0.75pt]   [align=left] {Yes};
% Text Node
\draw (443,373) node [anchor=north west][inner sep=0.75pt]   [align=left] {No};
% Text Node
\draw (397,473) node [anchor=north west][inner sep=0.75pt]   [align=left] {Case 4};
% Text Node
\draw (154,266) node [anchor=north west][inner sep=0.75pt]   [align=left] {$\lambda_2<\beta_4<\beta_1,\beta_2$};
\end{tikzpicture}
}

\ignore{
We now extend this proof for the limit system. First define the following norm for the infinite dimensional space as below:
\begin{eqnarray}
\label{Eqn_infinite_norm_finitensystem}
|| (\x, {\bar x}_b) ||_{\infty ,1} := \sup _n \frac{1}{n}
\sum_m \sum_{j \in {\cal G}_m}\left ( |x^m_j|+ \varsigma  |{\bar x}_b) | \right ) 
\end{eqnarray}
Observe that if the infinite sequence $ (\x, {\bar x}_b) $ is non-zero only for the first $k$-components then:
\begin{eqnarray}
\label{Eqn_infinite_norm_finitensystem_upto_k}
|| (\x, {\bar x}_b) ||_{\infty ,1} = \sup _{n \le k} \frac{1}{n}
\sum_m \sum_{j \in {\cal G}_m}\left ( |x^m_j|+ \varsigma  |{\bar x}_b) | \right ) .
\end{eqnarray}
With this, and following similar steps as before (after observing that $\sum_m 
\sum_{i \in {\cal G}_m, i \le k}  W_{j, i} \le (1-\eta_j)$ for all $j$ by non-negativity) we have that (for any $n$:
{\small 
\begin{eqnarray*}
||{\bar {\bf f}}^n ( {\bar x}_b, \x ) - {\bar {\bf f}}^n ( {\bar u}_b, \u )  ||_{\infty, 1} 
& \le &     || ( {\bar x}_b, \x )   - ( {\bar u}_b, \u )  ||_{\infty, 1} \left ( 1- \underline{\eta} + \varsigma   \underline{\eta} \right ).
\end{eqnarray*}}

Now consider the following:
\begin{eqnarray*}
 ||{\bar {\bf f}}^\infty ( {\bar x}_b, \x ) - {\bar {\bf f}}^\infty ( {\bar u}_b, \u )  ||_{\infty ,1} \hspace{-20mm} \\
& = &||\limsup({\bar {\bf f}}^n ( {\bar x}_b, \x ) - {\bar {\bf f}}^n ( {\bar u}_b, \u ) ) ||_{\infty,1} \\
&=& \lim_n ||\sup_{k \ge n}({\bar {\bf f}}^k ( {\bar x}_b, \x ) - {\bar {\bf f}}^k ( {\bar u}_b, \u ) ) ||_{\infty, 1} \\
&\le&  \lim_n \sup_{k\ge n} ||{\bar {\bf f}}^k ( {\bar x}_b, \x ) - {\bar {\bf f}}^k ( {\bar u}_b, \u ) ) ||_{\infty,1} \\
& \le &  \lim_n \sup_{k \ge n} \frac{1}{n} \sum_{j \in \mathcal{G}_1}   |  \xi_{j}^1  ({\bar  x}^1_j, x_b) - 
 \xi_{j}^1  ({\bar  u}^1_j, u_b) | \\
 &&+ \frac{1}{n}  \sum_{j \in \mathcal{G}_2} |  \xi_{j}^2  ({\bar  x}^2_j, x_b) - 
 \xi_{j}^2  ({\bar  u}^2_j, u_b) | \left ( 1- \eta_j^{sb} + \varsigma   \eta^{sb}_j  \right ) \\
  &\le & \lim_n \sup_{k \ge n} \sigma \frac{1}{n}  \sum_m \sum_{j \in {\cal G}_m} \left ( | {\bar x}_j^m - {\bar u}_j^m| +\varsigma | {\bar x}_b -  {\bar u}_b| \right )    \left ( 1- \eta_j^{sb} + \varsigma   \eta^{sb}_j  \right ) \\
   & \le&    || ( {\bar x}_b, \x )   - ( {\bar u}_b, \u )  ||_1 \left ( 1- \underline{\eta} + \varsigma   \underline{\eta} \right )
\end{eqnarray*}
Therefore the limit system  is strict contraction have unique fixed point with contraction coefficient$ (1- \underline{\eta} + \varsigma   \underline{\eta} )  < 1$.
}

\ignore{
a) With resilient regime   we  will have $v_2 < k_{d2}$, as in previous proof. Then the expected surplus simplifies and we have: 
\begin{eqnarray*}
E[S_1]- E[S_2] &=&\bigg (E[K_1] - v_1 +y_c(1+r_2)\frac{1-\gamma}{\gamma} \bigg)-  \\
& &\bigg (E[K_2] - v_2 -(1+r_2)(y_c+ y_2 p^{sb}_2) \bigg) \\
&=& (E[K_1]- E[K_2]) +(v_2-v_1) + (1+r_2) \bigg(y_2 p^{sb}_2 + \frac{y_c}{\gamma} \bigg)  \\
&=& \bigg(y_2 p^{sb}_2 + \frac{y_c}{\gamma} \bigg) \bigg (-({\bar r}_r+1-d_c) + \kappa +(1+r_2) \bigg).
\end{eqnarray*}
Therefore  $E[S_1]- E[S_2]  \geq  0$  if and only if  ${\bar r}_r   \leq  r_2+d_c +\kappa$. Also we used the following for the above computations:
\begin{eqnarray*}
v_2-v_1 &=& \kappa(\Omega_{2}-  \Omega_{1})
= \kappa \bigg(y_2 p^{sb}_2 + \frac{y_c}{\gamma} \bigg),\\
 E[K_1] - E[K_2] &=& w(k_{d1}- k_{d2})+(1-w) (k_{u1}- k_{u2}) \\
 &=& -({\bar r}_r +1-d_c)  \bigg(y_2 p^{sb}_2 + \frac{y_c}{\gamma} \bigg).
\end{eqnarray*}
b)   The surplus of the $\mathcal{G}_1$ banks in the resilient regime is  given by :
\begin{eqnarray*} 
\label{Eqn_surplus 1}
E[S_1]&=& E[K_1] -v_1+y_c(1+r_2)\frac{1-\gamma}{\gamma}, \mbox{ with, } \\ E[K_1]-v_1 &=& \bigg(k_0- \frac{1-\gamma}{\gamma}y_c \bigg) \bigg(1-d_c +{\bar r}_r- \kappa \bigg).
\end{eqnarray*}
\ignore{
Now consider  the following function  and its derivative with respect to $y_c$:
\begin{eqnarray}
g_1(y_c) &:=& \bigg(k_0- \frac{1-\gamma}{\gamma}y_c \bigg) \bigg(1-d_c +{\bar r}_r- \kappa \bigg) +y_c(1+r_2)\frac{1-\gamma}{\gamma}, \mbox{ and, } \nonumber \\ 
   g'_1(y_c)& =& \frac{1-\gamma}{\gamma} \bigg(1+r_2 - \big( 1-d_c +{\bar r}_r-\kappa\big) \bigg)\label{eqn_diffrention_with_yc}.
\end{eqnarray}
}
The surplus of $\mathcal{G}_1$ bank is again linear in $y_c$ and this  increases with $y_c$ if and only if   $\Delta_r < d_c + \kappa$, and remain constant  $\Delta_r = d_c + \kappa$.
%It increases with inter-lending amount $y_c$ if   ${\bar r}_r <\re$ and decreases if   ${\bar r}_r > \re$.\\
%\noindent The  next part of the Lemma  is immediate from  the equation \eqref{eqn_diffrention_with_yc}.\\

\noindent c)   
Now consider the surplus of $\mathcal{G}_2$ at upward movement:
\begin{eqnarray*}
\hat{S}_{2,u}& =& k_{u2}+ \bar{x}_2^{\infty} -v_2- \bar{y}_2 \\
&=& \Omega_2\bigg ((1+u-d_c) -\kappa\bigg) - (1+r_2) \bigg(y_c+ y_2 p^{sb}_2 \bigg) \\
&=& \bigg(k_0 + y_2 p_2^{sb} + y_c \bigg ) \bigg ((1+u-d_c) -\kappa\bigg) - (1+r_2) \bigg(y_c+ y_2 p^{sb}_2 \bigg).
\end{eqnarray*}
Once again the above function is linear in $y_c$ and 
thus the proof. 
%\noindent d) The expected surplus of the $\mathcal{G}_1$ banks in the resilient regime is  increases with $y_c$  from part(b) if   
%${\bar r}_r < \re$. 
\eop
}
%This is increasing with  inter lending amount $y_c$ if  $u  > \re$. Thus combining the two regime we have the proof.
%
% \noindent b) The proof of part(b) is similar to part(a) and hence omitted.
%\eop
%
%
\ignore{
Recall that we have,
\begin{eqnarray}
\bar {\bf f}^\infty ({\bar x}_b,\x)  &=&  ({\bar f}^\infty_b, \ \ \  {\bar f}^{\infty,1}_1,  \cdots , \ \ \ {\bar f}^{\infty,2}_1, \cdots) \mbox{ with} \\
{\bar f}_b^{\infty} &: =& \displaystyle \lim \sup_{n}{\bar f}^n_b ({\bar x}_b,\x) \  \mbox{and}, \nonumber \\
 {\bar f}_i^{\infty,m} ( {\bar x}_b,\x ) &:=
&\displaystyle \lim \sup_{n} {\bar f}^{n,m}_i ({\bar x}_b,\x) \mbox{ for all } i \mbox{ and } m.\nonumber
\end{eqnarray}

\textbf{Assumption:} Assume that
$\xi^{1}_{j} ({\bar x}^{1}_j, x_b)$ and  $\xi^{2}_{j} ({\bar x}^{2}_j, x_b)$ random variables are uniformly bounded by a constant $c_0$.

Our aim is to extend the Lemma \ref{lem: LLN_fixed point} for any bounded sequence. Show that for a bounded sequence $[0,y]$ $\bar {\bf f}^\infty ({\bar x}_b,\x)$  and  ${\bar f}_i^{\infty,m} ( {\bar x}_b,\x )$ converges.\\

 Recall we have 
 \begin{eqnarray*}
 {\bar f}^n_b ({\bar x}_b,{\x}) &:=&  
\frac{1}{n}
\displaystyle \sum_{j \in \mathcal{G}_1}   \xi^{1}_{j} ({\bar x}^{1}_j, x_b) W_{j, b} +  \frac{1}{n}
\displaystyle \sum_{j \in \mathcal{G}_2}   \xi^{2}_{j} ({\bar x}^{2}_j, x_b) W_{j, b}\\
&=& \frac{1}{n}
\displaystyle \sum_{j \in \mathcal{G}_1}  A_j + \frac{1}{n}
\displaystyle \sum_{j \in \mathcal{G}_2}  B_j
\end{eqnarray*}

${\bar f}^n_b ({\bar x}_b,{\x})$ converges because of the following:
\begin{eqnarray*}
&& \bigg|{\bar f}^{n+1}_b ({\bar x}_b,{\x}) -  {\bar f}^n_b ({\bar x}_b,{\x}) \bigg| \\
& =& \bigg|\frac{1}{n+1}
\displaystyle \sum_{j \in \mathcal{G}_1}  A_j + \frac{1}{n+1}
\displaystyle \sum_{j \in \mathcal{G}_2}  B_j - \frac{1}{n}
\displaystyle \sum_{j \in \mathcal{G}_1}  A_j - \frac{1}{n}
\displaystyle \sum_{j \in \mathcal{G}_2}  B_j \bigg| \\
&\le& \bigg|\frac{1}{n+1}
\displaystyle \sum_{j \in \mathcal{G}_1}  A_j - \frac{1}{n}
\displaystyle \sum_{j \in \mathcal{G}_1}  A_j \bigg| + \bigg|\frac{1}{n+1}
\displaystyle \sum_{j \in \mathcal{G}_2}  B_j - \frac{1}{n}
\displaystyle \sum_{j \in \mathcal{G}_2}  B_j \bigg|\\
&\le& \bigg|\frac{1}{n+1}
\displaystyle \sum_{j \in \mathcal{G}_1}  A_j  + \frac{A_{n\gamma+1}}{n+1} -\frac{1}{n+1}
\displaystyle \sum_{j \in \mathcal{G}_1}  A_j \frac{n+1}{n}\bigg| \\
&&+ \bigg|\frac{1}{n+1}
\displaystyle \sum_{j \in \mathcal{G}_2}  B_j  + \frac{A_{n(1-\gamma)+1}}{n+1} -\frac{1}{n+1}
\displaystyle \sum_{j \in \mathcal{G}_2}  B_j \frac{n+1}{n}\bigg| \\
&\le& \bigg|\frac{1}{n+1}
\displaystyle \sum_{j \in \mathcal{G}_1}  A_j (1-\frac{n+1}{n})\bigg| + \bigg|\frac{A_{n\gamma+1}}{n+1} \bigg|\\
&& + \bigg|\frac{1}{n+1}
\displaystyle \sum_{j \in \mathcal{G}_2}  B_j (1-\frac{n+1}{n})\bigg| + \bigg|\frac{B_{n(1-\gamma)+1}}{n+1} \bigg| \\
&\le& \frac{(1+\gamma)c_0}{n+1} +\frac{(2-\gamma)c_0}{n+1} \le \frac{\epsilon}{2} +\frac{\epsilon}{2} =\epsilon
\end{eqnarray*}
The last inequality due to the fact choose $n$ large  enough such that $\frac{(1+\gamma)c_0}{n+1}$ can be made arbitrarily small  $\forall n \ge N_1(\epsilon)$ and similarly the other term can be made arbitrarily small for  $\forall n \ge N_2(\epsilon)$.
Thus we have $\lbrace{\bar f}^n_b ({\bar x}_b,{\x})\rbrace$ is a cauchy sequence and hence converges.

Also  recall for the small nodes,
\begin{eqnarray}
{\bar f}^{n, m}_i ( {\bar x}_b, \x ) &:=&  
\left \{  \begin{array}{lll}
 \displaystyle \sum_{j \in \mathcal{G}_1}   \xi^1_{j} ({\bar x}^1_j, x_b)  W_{j, i} + \displaystyle \sum_{j \in \mathcal{G}_2}   \xi_{j}^2  ({\bar  x}^2_j, x_b)   {W}_{j, i} &\mbox{ if }  i \in \mathcal{G}_m ,    \\ \\
 0  &\mbox{\normalsize else,}  \mbox{ \normalsize and, }  \\
\end{array} \right.   
\end{eqnarray}
 Using the similar logic as above $\lbrace {\bar f}^{n, m}_i ( {\bar
 x}_b, \x ) \rbrace$ converges.}

 \ignore{
 % First figure start from here
 
\begin{tikzpicture}[x=0.75pt,y=0.75pt,yscale=-0.9,xscale=0.9]
%uncomment if require: \path (0,502); %set diagram left start at 0, and has height of 502

%Shape: Diamond [id:dp042926846690521936] 
\draw   (325.5,22) -- (411.5,57) -- (325.5,92) -- (239.5,57) -- cycle ;
%Straight Lines [id:da7617962014015986] 
\draw    (239.5,57) -- (102.5,58) ;
%Straight Lines [id:da845473854487296] 
\draw    (411.5,57) -- (539.5,58) ;
%Straight Lines [id:da31869413337731056] 
\draw    (102.5,58) -- (102.02,116.5) ;
\draw [shift={(102,118.5)}, rotate = 270.47] [color={rgb, 255:red, 0; green, 0; blue, 0 }  ][line width=0.75]    (10.93,-3.29) .. controls (6.95,-1.4) and (3.31,-0.3) .. (0,0) .. controls (3.31,0.3) and (6.95,1.4) .. (10.93,3.29)   ;
%Straight Lines [id:da7208770158084672] 
\draw    (539.5,58) -- (540.47,122) ;
\draw [shift={(540.5,124)}, rotate = 269.13] [color={rgb, 255:red, 0; green, 0; blue, 0 }  ][line width=0.75]    (10.93,-3.29) .. controls (6.95,-1.4) and (3.31,-0.3) .. (0,0) .. controls (3.31,0.3) and (6.95,1.4) .. (10.93,3.29)   ;
%Straight Lines [id:da7389590056510802] 
\draw    (111.5,197) -- (190.27,298.42) ;
\draw [shift={(191.5,300)}, rotate = 232.16] [color={rgb, 255:red, 0; green, 0; blue, 0 }  ][line width=0.75]    (10.93,-3.29) .. controls (6.95,-1.4) and (3.31,-0.3) .. (0,0) .. controls (3.31,0.3) and (6.95,1.4) .. (10.93,3.29)   ;
%Straight Lines [id:da9829575407191531] 
\draw    (111.5,197) -- (353.65,298.23) ;
\draw [shift={(355.5,299)}, rotate = 202.69] [color={rgb, 255:red, 0; green, 0; blue, 0 }  ][line width=0.75]    (10.93,-3.29) .. controls (6.95,-1.4) and (3.31,-0.3) .. (0,0) .. controls (3.31,0.3) and (6.95,1.4) .. (10.93,3.29)   ;
%Straight Lines [id:da7584753483073938] 
\draw    (111.5,197) -- (70.22,303.14) ;
\draw [shift={(69.5,305)}, rotate = 291.25] [color={rgb, 255:red, 0; green, 0; blue, 0 }  ][line width=0.75]    (10.93,-3.29) .. controls (6.95,-1.4) and (3.31,-0.3) .. (0,0) .. controls (3.31,0.3) and (6.95,1.4) .. (10.93,3.29)   ;
%Shape: Rectangle [id:dp09611174071817041] 
\draw   (16.5,311) -- (131.5,311) -- (131.5,423) -- (16.5,423) -- cycle ;
%Shape: Rectangle [id:dp47739619437048453] 
\draw   (161.5,301.5) -- (268.5,301.5) -- (268.5,413) -- (161.5,413) -- cycle ;
%Shape: Rectangle [id:dp7697152580348556] 
\draw   (313,301) -- (428.5,301) -- (428.5,409) -- (313,409) -- cycle ;
%Shape: Rectangle [id:dp40978223136595093] 
\draw   (34,116.5) -- (173,116.5) -- (173,196) -- (34,196) -- cycle ;
%Shape: Rectangle [id:dp31450665767357133] 
\draw   (476,122) -- (608,122) -- (608,180.5) -- (476,180.5) -- cycle ;

% Text Node
\draw (204,47) node  [xscale=0.75,yscale=0.75] [align=left] {\begin{minipage}[lt]{68pt}\setlength\topsep{0pt}
Yes
\end{minipage}};
% Text Node
\draw (509,36.75) node  [xscale=0.75,yscale=0.75] [align=left] {\begin{minipage}[lt]{68pt}\setlength\topsep{0pt}
No
\end{minipage}};
% Text Node
\draw (122.5,125.9) node [anchor=north west][inner sep=0.75pt]  [xscale=0.75,yscale=0.75]  {$$};
% Text Node
\draw (277.5,46.5) node [anchor=north west][inner sep=0.75pt]  [font=\normalsize,xscale=0.75,yscale=0.75] [align=left] {{\fontsize{0.8em}{0.96em}\selectfont $\bar{y}_2 > k_{u2} -v_2$}};
% Text Node
\draw (561.5,156.75) node  [font=\tiny,xscale=0.75,yscale=0.75] [align=left] {\begin{minipage}[lt]{63.24pt}\setlength\topsep{0pt}
{\fontsize{1.6em}{1.92em}\selectfont $P^2_D =w$}
\end{minipage}};
% Text Node
\draw (84.5,146.25) node  [font=\normalsize,xscale=0.75,yscale=0.75] [align=left] {\begin{minipage}[lt]{64.60000000000001pt}\setlength\topsep{0pt}
{\fontsize{0.8em}{0.96em}\selectfont $P^2_D= w$ or $1$ }\\{\fontsize{0.8em}{0.96em}\selectfont $\beta_4 > \beta_2$}
{\fontsize{0.8em}{0.96em}\selectfont $\beta_3 > \beta_1$}

\end{minipage}};
% Text Node
\draw (70.5,336) node  [xscale=0.75,yscale=0.75] [align=left] {\begin{minipage}[lt]{68pt}\setlength\topsep{0pt}
{\fontsize{0.8em}{0.96em}\selectfont   $\lambda_2 < \beta_1$ if $\lambda_2 >\beta_4$}\\{\fontsize{0.8em}{0.96em}\selectfont $P^2_D = w$ else $P^2_D =1$ }
\end{minipage}};
% Text Node
\draw (221,334.5) node  [xscale=0.75,yscale=0.75] [align=left] {\begin{minipage}[lt]{74.80000000000001pt}\setlength\topsep{0pt}
{\fontsize{0.8em}{0.96em}\selectfont if $\lambda_2 < \beta_2$}\\{\fontsize{0.8em}{0.96em}\selectfont $P^2_D= w$ else $P^2_D=1$}\\\\
\end{minipage}};
% Text Node
\draw (367.25,339) node  [font=\tiny,xscale=0.75,yscale=0.75] [align=left] {\begin{minipage}[lt]{61.540000000000006pt}\setlength\topsep{0pt}
{\fontsize{0.8em}{0.96em}\selectfont  $\beta_1< \lambda_2< \beta_3$}\\{\fontsize{0.8em}{0.96em}\selectfont if $\lambda_2 > \beta_2$}\\{\fontsize{0.8em}{0.96em}\selectfont $P^2_D =w$ else $P^2_D =1$}
\end{minipage}};

\end{tikzpicture}
% First figure ends here
}

 \end{document}

\ignore{

% The following material did not work and the analysis for default regime, currently we are not using the materials

{\color{red}

\newpage
\begin{lemma}
\label{lemma_surplus1_defaultregime}
Assume   proportional taxes, i.e., $v =  \kappa \Omega $, and   let ${\bar r}_r = u(1-w) +dw$,  
$\re := r_2+d_c +\kappa$ and $\beta  : ={\bar x}_2^\infty\frac{1- \gamma}{\gamma} \frac{y_c}{y_2(1-p^{sb}_2)}$. Under the default regime i..e, if $P^1_D \le w$ and $P^2_D \in \{w,1\}$
%The expected surplus of  $\mathcal{G}_1$  under the resilient regime is :
\begin{enumerate}[a)]
    \item  the expected surplus $E[S_1] \geq E[S_2]$  if  and only if  ${\bar r}_r  < \re $ 
    \item the expected surplus ($E[S_1]$) of  $\mathcal{G}_1$ is increasing with inter lending amount $y_c$ if  $ \frac{d\beta}{dy_c} > \frac{1-\gamma}{\gamma} \bigg(1-d_c -\kappa +wd +(1-w)\bar{r}_r \bigg)   $;
    \item $E[S_1]$ is decreasing with $y_c$ if  $ \frac{d\beta}{dy_c} < \frac{1-\gamma}{\gamma} \bigg(1-d_c -\kappa +wd +(1-w)\bar{r}_r \bigg)   $;and remains constant when $  \frac{d\beta}{dy_c} = \frac{1-\gamma}{\gamma} \bigg(1-d_c -\kappa +wd +(1-w)\bar{r}_r \bigg) $  ; and; 
    \item if further $  \frac{d\beta}{dy_c}   > \max \bigg\{\bigg(\frac{\beta+ \frac{1-\gamma}{\gamma} \frac{y^2_c(\re -u)}{y_2(1-p^{sb}_2)}}{y_c} \bigg), \frac{1-\gamma}{w\gamma} \bigg(1-d_c -\kappa +wd +(1-w)\bar{r}_r \bigg)  \bigg\} $
    then,  SaU
$\hat{S}_{2,u}$   of $\mathcal{G}_2$  and the expected surplus $E[S_1]$ of $\mathcal{G}_1$ are both  increasing with inter lending amount $y_c$ .
\end{enumerate}
\end{lemma}
\textbf{Proof:} a) For this we divide the proof in the following cases.\\
\textbf{Case 1:} First consider the case  with $P^1_D=0$ and $P^2_D\in \{w,1\}$  and $\beta  : ={\bar x}_2^\infty\frac{1- \gamma}{\gamma} \frac{y_c}{y_2(1-p^{sb}_2)}$.
\begin{eqnarray*}
E[S_1] -E[S_2] &=& E\bigg ( K_{1,i}  +   {\bar x}_1^\infty +  {\bar x}_2^\infty \frac{1-\gamma}{\gamma}\frac{y_c}{y_2(1-p^{sb}_2)}-v_1- \bar{y}_1 \bigg )^+ \\
& &- E \bigg ( K_{2,i}  +  {\bar x}_2^\infty   -v_2- \bar{y}_2 \bigg )^+\\ 
&=& E[K_1] -E[K_2] +(v_2-v_1)+\beta + (\bar{y}_2 - {\bar x}_2^\infty) \\
&=& -({\bar r}_r +1-d_c)  \bigg(y_2 p^{sb}_2 + \frac{y_c}{\gamma} \bigg) + \kappa \bigg(y_2 p^{sb}_2 + \frac{y_c}{\gamma} \bigg) + \beta  + \bigg(\bar{y}_2 - {\bar x}_2^\infty \bigg) \\
&=& \underbrace{\bigg(y_2 p^{sb}_2 + \frac{y_c}{\gamma} \bigg)}_{>0}(\kappa - {\bar r}_r -1+d_c) + \underbrace{  \bigg(\bar{y}_2 - {\bar x}_2^\infty \bigg) +\beta}_{> 0 }  
\end{eqnarray*}
Thus $E[S_1] \ge E[S_2]$ if and only if ${\bar r}_r  \leq  \re -(r_2+ 1)< \re $. \\
\textbf{Case 2:} Now consider the case with $P^1_D=w$ and $P^2_D\in \{w,1\}$ 
\begin{eqnarray*}
E[S_1] -E[S_2] &=&  E[K_1] -E[K_2] +(v_2-v_1)+ {\bar x}_1^\infty  +\beta -\bar{y}_1+ (\bar{y}_2 - {\bar x}_2^\infty) \\
&=& E[K_1] -E[K_2] +(v_2-v_1)+ (\bar{y}_2 - {\bar x}_2^\infty) + {\bar x}_1^\infty  +\beta -\bar{y}_1\\
&=& E[K_1] -E[K_2] +(v_2-v_1)+ (\bar{y}_2 - {\bar x}_2^\infty)   +\frac{\beta+w(k_{d1}-v_1)}{1-w} \\
&=& \underbrace{\bigg(y_2 p^{sb}_2 + \frac{y_c}{\gamma} \bigg)}_{>0}(\kappa - {\bar r}_r -1+d_c) + \underbrace{(\bar{y}_2 - {\bar x}_2^\infty)}_{> 0}+ \underbrace{\frac{\beta+w(k_{d1}-v_1)}{1-w}}_{> 0,\ \mbox{as w is small}}
\end{eqnarray*}
Thus combining two cases $E[S_1] \ge E[S_2]$ if and only if ${\bar r}_r  \leq  \re -(r_2+ 1)< \re$. \\
b) Similar to the above we divide the proof in the following  cases:\\
\textbf{Case 1:} When $P^1_D=0$ and $P^2_D \in\{w,1\}$.
\begin{eqnarray*}
E[S_1] &=& E[K_1]-v_1 +{\bar x}_1^\infty  +\beta -\bar{y}_1\\
&=& E[K_1]-v_1 +\beta \\
&=& \bigg(k_0 -\frac{1-\gamma}{\gamma}y_c\bigg)\bigg(\bar{r}_r+1-d_c -\kappa\bigg)+\beta
\end{eqnarray*}
Therefore $E[S_1]$ is increasing with the inter lending  parameter $y_c$ if  and only  if 
$$
\frac{d\beta}{dy_c} > \frac{1-\gamma}{\gamma}\bigg(\bar{r}_r+1-d_c -\kappa\bigg).
$$
\textbf{Case 2:} We consider the expected surplus of the $\mathcal{G}_1$ when $P^1_D= w$ and $P^2_D \in \{w,1\}$.
{\small
\begin{eqnarray}
E[S_1]&=&E[K_1] -v_1 + {\bar x}_1^\infty  +\beta -\bar{y}_1  \nonumber \\
&=& E[K_1]-v_1+ \frac{w}{1-w}(k_{d1}-v_1+\beta)+\beta \nonumber \\
&=& \bigg(k_0 -\frac{1-\gamma}{\gamma}y_c\bigg) \bigg \{\bigg(1-d_c+\bar{r}_r -\kappa\bigg) + \frac{w}{1-w}\bigg(1+d-d_c-\kappa \bigg)\bigg\} +\frac{\beta}{1-w} 
\label{Eqn_defultRegimeforG1}
\end{eqnarray}}
Thus $E[S_1]$ is increasing  with the inter lending parameter $y_c$ if 
$$ 
\frac{d\beta}{dy_c} > \frac{1-\gamma}{\gamma} \bigg(1-d_c -\kappa +wd +(1-w)\bar{r}_r \bigg) . 
$$
c) The proof is immediate from the equation \eqref{Eqn_defultRegimeforG1}.\\
d) Recall surplus of $\mathcal{G}_2$ at upward movement:
\begin{eqnarray*}
\hat{S}_{2,u}& =& k_{u2}+ \bar{x}_2^{\infty} -v_2- \bar{y}_2 \\
&=& \Omega_2\bigg ((1+u-d_c) -\kappa\bigg) + \bar{x}_2^{\infty} - \bar{y}_2  \\
&=& \bigg(k_0 + y_2 p^{sb}_2 + y_c \bigg ) \bigg ((1+u-d_c) -\kappa\bigg) +\frac{\gamma y_2(1-p^{sb}_2)\beta}{y_c(1-\gamma)} - (y_2+y_c)(1+r_2)
\end{eqnarray*}
$\hat{S}_{2,u}$ is increasing with  the inter-lending parameter $y_c$ if 
$$ 
\frac{d\beta}{dy_c}  > \bigg(\frac{\beta+ \frac{1-\gamma}{\gamma} \frac{y^2_c(\re -u)}{y_2(1-p^{sb}_2)}}{y_c} \bigg).
$$
Thus combining part(b) and the above we have the proof.
\eop
\ignore{
 Define the performance of the $\mathcal{G}_1$ with the upward movement, i.e., banks are interested in the shock-free return from the outside investment. Let $\hat{S}_{1,u}$ be the surplus  of the $\mathcal{G}_1$ banks, and it is defined as:
\begin{eqnarray*}
\hat{S}_{1,u} = \bigg(k_{u1}+ \bar{x}_1^{\infty} +{\bar x}_2^\infty\frac{1- \gamma}{\gamma} \frac{y_c}{y_2(1-p^{sb}_2)} -v_1- \bar{y}_1 \bigg)
\end{eqnarray*}
}
\begin{lemma}
Consider the regime where $P^{1}_D =0$ and $P^{2}_D =w$, the expected surplus $E[S_1]$  and $\hat{S}_{2,u}$ is increasing with the  inter lending parameter $y_c$ if and only if $y_{*} > (1-d_c+\bar{r}_r -\kappa)\bigg((y_2+y_c) -wy_2(1-p^{sb}_2)\bigg)^2$ and 

$ \frac{d {\bar x}_2^\infty}{d\beta} > \re -u$
respectively, where 
\begin{eqnarray*}
y_{*}&=& (1-w) \bigg( y_c(y_2+y_c)(1+r_2) + y_2 \bar{y}_2 -w(1-p^{sb}_2) y_2(\bar{y}_2 +y_2(1+r_2)) \bigg)\\ 
&&+ w \bigg((k_{d2} -v_2)y_2 (1-w(1-p^{sb}_2)) +y_c(1+d-d_c -\kappa)\bigg((y_2+y_c) -wy_2(1-p^{sb}_2) \bigg)\bigg).
\end{eqnarray*}
\end{lemma}
\textbf{Proof:} Consider the expected surplus of the $\mathcal{G}_1$  banks  when  $P^{1}_D =0$ and $P^{2}_D =w$.
\begin{eqnarray*}
E[S_1]&=&E\bigg (  K_{1,i}  +   {\bar x}_1^\infty +  {\bar x}_2^\infty \frac{1-\gamma}{\gamma}\frac{y_c}{y_2(1-p^{sb}_2)}-v_1- \bar{y}_1 \bigg )^+ \\
&=&E[K_1] -v_1 + \beta\\
&=& \bigg(k_0 -\frac{1-\gamma}{\gamma} y_c\bigg)(1-d_c+\bar{r}_r -\kappa)+ \frac{1-\gamma}{\gamma}y_c \bigg( \frac{\bar{y}_2(1-w) +(k_{d2}-v_2)w}{(y_2+y_c) -w(1-p^{sb}_2) y_2}\bigg) 
\end{eqnarray*}
Thus $E[S_1]$ increases with the inter lending parameter $y_c$ if and only if
\begin{eqnarray*}
&&\bigg((1-w)(\bar{y}_2+y_c(1+r_2))+ w(y_c(1+d-d_c-\kappa) +k_{d2}-v_2) \bigg)\\
&&\bigg((y_2+y_c) -w(1-p^{sb}_2)y_2)\bigg) -y_c\bigg(\bar{y}_2(1-w) +w(k_{d2}-v_2)\bigg)\\
&&>  (1-d_c+\bar{r}_r -\kappa)\bigg((y_2+y_c) -wy_2(1-p^{sb}_2)\bigg)^2
\end{eqnarray*}
which simplifies to 
$$y_{*} > (1-d_c+\bar{r}_r -\kappa)\bigg((y_2+y_c) -wy_2(1-p^{sb}_2)\bigg)^2.$$

Now consider the surplus at upward movement of the $\mathcal{G}_1$ banks
\begin{eqnarray*}
\hat{S}_{2,u}& =& k_{u2}+ \bar{x}_2^{\infty} -v_2- \bar{y}_2 \\
&=& \Omega_2\bigg ((1+u-d_c) -\kappa\bigg) + \bar{x}_2^{\infty} - \bar{y}_2   \\
&=&( k_0+y_c+ y_2p^{sb}_2)\bigg ((1+u-d_c) -\kappa\bigg) + \bar{x}_2^{\infty} -(y_2+y_c)(1+r_2)  
\end{eqnarray*}
Using the similar steps as above  the surplus $\hat{S}_{2,u}$ increasing with the  inter-lending parameter $y_c$ if and only if 
 $$
 \frac{d {\bar x}_2^\infty}{d\beta} > \re -u .
 $$
\eop

\noindent \textbf{Note:} A simpler representation of the above Lemma  could be in terms of the variable $\beta$ i.e., we have   $E[S_1]$ increases with  the inter lending parameter $y_c$ if and only if  $\frac{d\beta}{dy_c} > \frac{1-\gamma}{\gamma} (1-d_c +\bar{r}_r -\kappa)$. Also it is easy to compute that: 
\begin{eqnarray*}
\beta &=& \frac{1-\gamma}{\gamma}y_c \bigg( \frac{\bar{y}_2(1-w) +(k_{d2}-v_2)w}{(y_2+y_c) -w(1-p^{sb}_2) y_2}\bigg)  \\
\beta &=&\frac{1-\gamma}{\gamma} \frac{A_0 y^2_c+A_1 y_c}{(y_c+ A_2)} \\
\frac{d\beta}{dy_c}& =&  \frac{1-\gamma}{\gamma} \frac{y^2_cA_0 +2y_cA_0A_2+A_2 A_1}{(y_c+A_2)^2} > 0
\end{eqnarray*}
\begin{eqnarray*}\frac{dE[S_1]}{dy_c} &=& \frac{1-\gamma}{\gamma}
\left ( \frac{y^2_cA_0 +2y_cA_0A_2+A_2 A_1 
- (y_c^2 + 2y_c A_2 + A_2^2) (1+d-d_c-\kappa) }{(y_c+A_2)^2} \right )  \\
&= &
 \frac{1-\gamma}{\gamma} (1-w) (\re-d)  \frac{y_c^2  
 + 2y_c A_2  + 
 A_2 ( y_2   +  \frac{(1+d-d_c-\kappa)}{(1-w) (\re-d)} (k_0 w - y_2 p^{sb}_2 (1-w) )   )
 }{(y_c+A_2)^2}
\end{eqnarray*}
where the constants are 
\begin{eqnarray*}
A_0 : &=&(1+r_2)(1-w)+ w(1+d-d_c-\kappa) \\
A_1: &=& y_2(1+r_2)(1-w) +(k_0+y_2p^{sb}_2)(1+d-d_c-\kappa)w \\
A_2 :&=&y_2(1-w(1-p^{sb}_2)).
\end{eqnarray*}
Thus the above condition simplifies to  $(1-d_c +\bar{r}_r -\kappa) < 0$.\\

Recall the average clearing vector of the $\mathcal{G}_2$ bank  in the default regime is:
\begin{eqnarray*}
\bar{x}^{\infty}_2 &=&  \bigg( \frac{\bar{y}_2(1-w) +(k_{d2}-v_2)w}{(y_2+y_c) -w(1-p^{sb}_2) y_2}\bigg)(1-p^{sb}_2) y_2 \\
\frac{d\bar{x}^{\infty}_2}{dy_c}&=& = \frac{(A_0 y_c +A_1)}{(y_c+ A_2)^2}(1-p^{sb}_2) y_2.
\end{eqnarray*}
}
}